\documentclass[11pt,a4]{article}
\usepackage[T2A]{fontenc}
\usepackage[cp866]{inputenc}
\usepackage[mathscr]{eucal}
\usepackage{amssymb}
\usepackage{amsmath}
\renewcommand{\leq}{\leqslant}
\renewcommand{\geq}{\geqslant}

\renewcommand{\C}{{\mathbb C}}

\newcommand{\R}{{\mathbb R}}

\newcommand{\K}{{\mathbb K}}

\renewcommand{\k}{\rule{0.7em}{0.7em}}
\usepackage{fancyhdr}
\begin{document}
\sloppy

{\normalsize

\thispagestyle{empty}

\mbox{}
\\[-2.25ex]
\centerline{
{%\large
\bf
FIRST INTEGRALS  OF  ORDINARY LINEAR  DIFFERENTIAL SYSTEMS
}
}
\\[2.5ex]
\centerline{
\bf 
V.N. Gorbuzov$\!{}^{*},$   A.F. Pranevich$\!{}^{**}$
}
\\[2ex]
\centerline{
\it 
Yanka Kupala Grodno State University
{\rm(}Ozeshko 22, Grodno, 230023, Belarus{\rm)}
}
\\[1.5ex]
\centerline{
E-mail: gorbuzov@grsu.by$\!{}^{*},$
pronevich@tut.by$\!{}^{**}$
}
\\[5.5ex]
\centerline{{\large\bf Abstract}}
\\[1ex]
\indent
The spectral method for building first integrals of ordinary linear differential systems is elaborated. 
Using this method, we obtain bases of first integrals for linear differential systems with constant coefficients,
for linear nonautonomous differential systems integrable in closed form (algebraic reducible systems, triangular systems,
the Lappo-Danilevskii systems), and for reducible ordinary differential systems with respect to various transformation groups.
\\[1.5ex]
\\[1.25ex]
\\[5ex]
\centerline{{\large\bf Contents}}
\\[1.5ex]
{\bf  Introduction}                   \dotfill\ 2
\\[0.5ex]
{\bf 1. 
Integrals of ordinary linear differential system with constant coefficients}\;\!\!\! \footnote{
The main results of this Section were originally published in
{\it Vestnik of the Yanka Kupala Grodno State Univ.}, 
2002, Ser. 2,  No. 2(11), 23-29; 
2003, Ser. 2,  No. 2(22), 50-60; 
2008, Ser. 2,  No. 2(68), 5-10.}
                                                 \dotfill \ 4
\\[0.5ex]
\mbox{}\hspace{1em}
1.1. Linear homogeneous differential system
                                                 \dotfill \ 4
\\
\mbox{}\hspace{2.75em}
1.1.1. Linear partial integral
                                                 \dotfill \ 4\\
\mbox{}\hspace{2.75em}
1.1.2. Autonomous first integrals      
                                                 \dotfill \ 4\\
\mbox{}\hspace{2.75em}
1.1.3. Nonautonomous first integrals
                                                 \dotfill \ 14
\\[0.5ex]
\mbox{}\hspace{1em}
1.2. Linear nonhomogeneous differential system
                                                 \dotfill \ 17\\
\mbox{}\hspace{2.75em}
1.2.1. Case of simple elementary divisors
                                                 \dotfill \ 17\\
\mbox{}\hspace{2.75em}
1.2.2. Case of multiple elementary divisors
                                                 \dotfill \ 18
\\[0.75ex]
\noindent
{\bf 2. 
Integrals of ordinary linear nonautonomous differential systems 
\\
\mbox{}\hspace{0.9em}
integrable in closed form}
                                                 \dotfill \ 21
\\[0.5ex]
\mbox{}\hspace{1em}
2.1. Algebraic reducible systems\;\!\!\! \footnote{
The results of this Subsection on the multidimensional case has been published in the journal
{\it Vestnik of the Belarusian State Univ.}, 2008, Ser. 1, No. 2, 75-79.
}
                                                 \dotfill \ 21
\\
\mbox{}\hspace{2.75em}
2.1.1. Partial integrals
                                                 \dotfill \ 21
\\
\mbox{}\hspace{2.75em}
2.1.2. First integrals
                                                 \dotfill \ 22
\\[0.5ex]
\mbox{}\hspace{1em}
2.2. Triangular systems
                                                 \dotfill \ 26
\\[0.5ex]
\mbox{}\hspace{1em}
2.3. The Lappo-Danilevskii systems\;\!\!\! \footnote[3]{
The definitive version of this Subsection has been published in the journals
{\it Differential equations and control processes}, 2001, No. 3, 17-45 (http://www.neva.ru/journal) and
{\it Vestnik of the Yanka Kupala Grodno State Univ.}, 2008, Ser. 2,  No. 3(73), 79-83. 
}
                                                 \dotfill \ 30
\\
\mbox{}\hspace{2.75em}
2.3.1. Linear homogeneous differential system
                                                 \dotfill \ 30
\\
\mbox{}\hspace{2.75em}
2.3.2. Linear nonhomogeneous differential system
                                                 \dotfill \ 51
\\[0.5ex]
\noindent
{\bf 3. 
Integrals of reducible ordinary differential systems}
                                                 \dotfill \ 57
\\[0.35ex]
\mbox{}\hspace{1em}
3.1. Linear homogeneous differential system
                                                 \dotfill \ 57
\\[0.25ex]
\mbox{}\hspace{1em}
3.2. Linear nonhomogeneous differential system
                                                 \dotfill \ 64
\\
\mbox{}\hspace{2.75em}
3.2.1. Case of simple elementary divisors
                                                 \dotfill \ 64
\\
\mbox{}\hspace{2.75em}
3.2.2. Case of multiple elementary divisors
                                                 \dotfill \ 66
\\[0.35ex]
{\bf References}
                                              \dotfill \ 72

\newpage

\mbox{}
\\[-1.75ex]
\centerline{\large\bf  Introduction}
\\[1.5ex]
\indent
One of the most important problems of the general theory of differential systems 
is the problem of finding first integrals. 
In 1878 the French mathematician J.G. Darboux showed how general integrals of 
the first-order ordinary differential equations possessing sufficient in\-va\-ri\-ant curves
are constructed [1]. 
His investigation gave the classical problem {\sl {\rm (}the Dar\-bo\-ux problem}) 
about building of first integrals by known partial integrals. 
The review of the litera\-tu\-re and the current situation of the theory of integrals 
are given in the monographies [2\! -- 5]. 
\\[0.25ex]
\indent
Linear differential systems are of interest for mathematicians both per se and as a tool 
for studying nonlinear differential equations by means of the method of linearization. 
Note also that systems of linear differential equations 
play a broad and fundamental role in electrical, mechanical, 
chemical and aerospace engineering, communications, and signal processing.
There exist many good works on linear differential systems (see, for example, [6 -- 24]).

In this paper we study the Dar\-bo\-ux problem of the existence of first integrals for 
main classes of linear ordinary differential systems. 
Using the method of partial integrals for polynomial differential systems [3, pp. 187 -- 226; 28 -- 31],
we obtain {\sl the spectral method} for building first integrals of linear differential systems [25 -- 27].

At a later time, the spectral method has been applied to  
linear multidimensional differential systems [32 -- 38] and 
to nonlinear Jacobi's differential systems [39 -- 44].
\vspace{0.25ex}

The material of this paper is made on the base of our articles  [25 -- 27; 36; 37; 45; 46].
\vspace{0.25ex}

To avoid ambiguity, we stipulate the following notions (using the article [47]).

Consider an ordinary differential system of $n\!$-th order
\\[2ex]
\mbox{}\hfill                                   %(0.1)
$
\dfrac{dx}{dt} = X(t,x),
$
\hfill (0.1)
\\[2.25ex]
where $x\in {\mathbb R}^n,\ t\in {\mathbb R},$
\vspace{1ex}
the column vector
$\dfrac{dx}{dt}={\rm colon}\Bigl(\dfrac{dx_1^{}}{dt}\,,\ldots, \dfrac{dx_n^{}}{dt}\Bigr),$
the coordinates of the vector function
\vspace{0.75ex}
$X(t,x)={\rm colon}\bigl(X_{1}^{}(t,x),\ldots, X_{n}^{}(t,x)\bigr)$ 
are the continuously differentiable on a domain $G\subset {\mathbb R}^{n+1}$ scalar functions
\vspace{0.35ex}
$
X_{i}^{}\colon G \to {\mathbb R},
\ i=1,\ldots, n.
$

We recall that by domain we mean open arcwise connected set.
\vspace{0.25ex}

A continuously differentiable scalar function $F \colon G^{\,\prime} \to {\mathbb R}$
is called a \textit{first integral} on a domain $G^{\,\prime}\subset G$ of the 
ordinary differential system {\rm (0.1)} if 
\\[1.5ex]
\mbox{}\hfill                               
$
{\frak X}\;\! F(t,x) = 0
$
\ for all 
$
(t,x) \in G^{\,\prime},
\hfill
$
\\[1.75ex]
where the linear differential operator 
$
{\frak X}(t,x)=
\partial_{t}^{}+ 
\sum\limits_{i=1}^{n}
X_{i}^{}(t,x)\,\partial^{}_{x^{}_i}$
for all 
$
(t,x)\in G
$
is the {\it operator of differentiation by virtue of system} (0.1).
\vspace{0.35ex}

A continuously differentiable scalar function $F \colon G^{\,\prime} \to {\mathbb R}$
\vspace{0.25ex}
is a first integral on a domain $G^{\;\!\prime}\subset G$ of the 
ordinary differential system {\rm (0.1)} if and only if
\vspace{0.25ex}
the function $F \colon G^{\,\prime} \to {\mathbb R}$ is constant along any solution
\vspace{0.5ex}
$x\colon t \to x(t)$ for all $t\in J$ of system {\rm(0.1)}, 
where $x\colon t \to x(t)$ is such that 
$(t,x(t))\in G^{\;\!\prime}$ for all $t\in J\subset {\mathbb R},$ i.e.,
\vspace{0.5ex}
$F(t,x(t))=C$ for all $t\in J,\ C={\rm const}.$

A continuously differentiable scalar function $w \colon G^{\,\prime} \to {\mathbb R}$ 
\vspace{0.25ex}
is said to be a {\it partial integral} on a domain $G^{\;\!\prime}\subset G$ of the 
ordinary differential system {\rm (0.1)} if
\\[1.75ex]
\mbox{}\hfill                                     
$
{\frak X}\;\! w(t,x) = \Phi(t,x)
$
\ for all 
$
(t,x)\in G^{\;\!\prime},
\hfill
$
\\[1.75ex]
where $\Phi\colon G^{\;\!\prime}\to {\mathbb R}$ is a scalar function such that
$
\Phi(t,x)_{\displaystyle |_{w(t,x)=0}} = 0
$
for all 
$
(t,x)\in G^{\;\!\prime}.
$
\vspace{1ex}

A continuously differentiable scalar function $w \colon G^{\,\prime} \to {\mathbb R}$
\vspace{0.25ex}
is a partial integral on a domain $G^{\;\!\prime}\subset G$ of the 
\vspace{0.25ex}
ordinary differential system {\rm (0.1)} if and only if
the function $w \colon G^{\,\prime} \to {\mathbb R}$ vanishes identically along any solution
\vspace{0.5ex}
$x\colon t \to x(t)$ for all $t\in J\subset\! {\mathbb R}$ 
of system {\rm(0.1)}, where $x\colon t \to x(t)$ is
such that $(t,x(t))\in G^{\;\!\prime}$ for all $t\in J,$ i.e.,
$
w(t,x(t))=0
$ 
for all 
$t\in J.$

\newpage

A set of the functionally independent on a domain $G^{\;\!\prime}\subset G$ 
\vspace{0.25ex}
first integrals 
$
F_l^{} \colon G^{\;\!\prime}\to {\mathbb R}, 
\linebreak
l= 1,\ldots, k,
$
of system {\rm (0.1)} is called a {\it basis of first integrals} (or {\it integral basis}) 
\vspace{0.25ex}
on the domain $G^{\;\!\prime}$ of system {\rm (0.1)} if
for any first integral $\Psi \colon G^{\;\!\prime} \to {\mathbb R}$ of system {\rm (0.1)},
we have 
\\[1.5ex]
\mbox{}\hfill
$
\Psi(t,x) = \Phi\bigl(F_{1}^{}(t,x), \ldots, F_{k}^{}(t,x)\bigr)
$
\ for all 
$
(t,x) \in G^{\;\!\prime},
\hfill
$
\\[1.5ex]
where 
\vspace{0.35ex}
$\Phi$ is some continuously differentiable function on the range 
of the vector function 
$F\colon (t,x)\to \bigl(F_{1}^{}(t,x), \ldots, F_{k}^{}(t,x)\bigr)$ for all $(t,x) \in G^{\;\!\prime}.$
\vspace{0.35ex}
The number $k$ is said to be the \textit{di\-men\-si\-on} of 
ba\-sis of first integrals on the domain $G^{\;\!\prime}$ for system {\rm (0.1)}.
\vspace{0.35ex}

The ordinary differential system {\rm(0.1)} on a neighbourhood of any point of the domain $G$ 
has a basis of first integrals of dimension $n$ [48, pp.  175 -- 177; 49, pp. 256 -- 263].
\vspace{0.35ex}

If the ordinary differential 
system {\rm(0.1)} is autonomous 
(the vector function $X\colon x\to X(x)$ for all $x\in {\mathscr X}\subset {\mathbb R}^n$ 
not de\-pends on the independent variable $t),$ then the ordinary differential system (0.1) 
in an $n\!$-di\-men\-si\-o\-nal integral basis has 
$n-1$ functionally independent on the domain ${\mathscr X}$ 
autonomous first integrals [50, pp.  161 -- 169].

The paper is organized as follows.

{\it In Section} 1 the spectral method for building first integrals 
for linear homogeneous and nonhomogeneous
differential systems with real constant coefficients is developed.\! 
Here we pay spe\-ci\-al attention to the construction 
of bases of {\it real} first integrals for these systems.

The Euler method [51, pp. 350 -- 360; 52, pp. 93 -- 101] and 
the matrix method [51, pp. 329 -- 349; 53, pp. 166 -- 172] are the main methods for 
finding solutions of linear homogeneous differential systems with constant coefficients.
For building first integrals we know the method of integrable combinations [48, pp. 171 -- 173]
and the N.P. Erugin --- N.A. Zboichik method [54, pp. 464 -- 469; 55]. 
These approaches show only ways of constructing first integrals 
for linear differential systems, but they do not build first integrals in explicit form.

Using the method of partial integrals [3, pp. 187 -- 226; 28 -- 31]
for polynomial differential systems, we obtain the 
spectral method of building first integrals in explicit form for linear differential systems 
with constant coefficients [26; 27; 45]. 

Note also that the results of Subsection 1.1.2 are co-ordinated with the articles [20 -- 24]. 
%for linear homogeneous autonomous differential systems.
The approach of finding autonomous first integrals with using the method of Jordan ca\-no\-ni\-cal 
forms was developed by M.~Falconi and J. Llibre [20].
The existence of common first integrals for two coupled third-order linear differential systems 
(which satisfy the Frobenius compatibility condition) is discussed in [21].
The problem of building rational first integrals for linear systems was considered 
by the French mathematician J.A. Weil in the paper [22].
The class of linear autonomous differential systems with no rational first integrals 
is determined by A. Nowicki [23].
First integrals are constructed using integrating factors in [24].
\vspace{0.35ex}

{\it In Section} 2 we investigate the problem of the exis\-tence of first in\-teg\-rals for 
linear nonautonomous differential systems integrable in closed form. 

In the general theory of differential systems we know some results about 
reducibility of 
%linear nonautonomous 
differential systems to some special forms. 
For example, any linear nonautonomous differential system is 
reducible to an algebraic reducible system or to a triangular system [56; 8, p. 33].

Three classes of linear nonautonomous 
differential systems integrable in closed form (algebraic reducible systems, triangular systems,
the Lappo-Danilevsky systems) are considered.
The regular method of building first integrals for these systems is elaborated [36; 46].
Also, the N.P. Erugin problem of the existence of autonomous first integrals [54, p. 469]
for an nonautonomous homogeneous Lappo-Danilevskii differential system is solved [25].
\vspace{0.35ex}

{\it In Section} 3 we consider reducible linear nonautonomous differential systems 
with respect to various tran\-s\-for\-ma\-ti\-on groups 
(periodic, po\-ly\-no\-mi\-al,  or\-tho\-go\-nal, Liapunov, exponential, and other groups)
to the linear differential systems with constant coefficients.

The spectral method for building first integrals of reducible systems is elaborated. 

The notion of reducible systems has been entered by A.M. Lyapunov [57].
The development of the theory of reducible systems is associated with the name of N.P.Erugin [58].
%The general theory of reducible systems was developed by N.P.Erugin in 1946 [58].

In addition, in this article some examples are given to illustrate the obtained results.

\newpage

\mbox{}
\\[-2.5ex]
\centerline{
\large\bf  
1. Integrals of ordinary linear differential system 
}
\\[0.5ex]
\centerline{
\large\bf  
with constant coefficients
}
\\[1.25ex]
\centerline{
{\bf  1.1. 
Linear homogeneous differential system}
}
\\[1.25ex]
\indent
Consider an ordinary linear homogeneous differential system 
with constant coefficients
\\[1.75ex]
\mbox{}\hfill                                             % (1.1)
$
\dfrac{dx}{dt} =  A\;\! x,
$
\hfill (1.1)
\\[2ex]
where $x = \mbox{colon}(x_{1}^{},\ldots,x_{n}^{})\in {\mathbb R}^n$ and 
\vspace{0.75ex}
$A = \bigl\|a_{\xi j}^{}\bigr\|$ is an square constant matrix of order $n$ 
with entries $a_{\xi j}^{}\in {\mathbb R},\ \xi=1,\ldots,n,\ j=1,\ldots,n.$
\vspace{0.5ex}
The differential system (1.1) is induced the autonomous linear differential ope\-ra\-tor of first order
\\[1.25ex]
\mbox{}\hfill
$
\displaystyle
{\frak A} (x) = \sum\limits_{\xi=1}^{n} A_{\xi}^{}x\,\partial_{x_\xi^{}}^{}
$
\ for all 
$
x\in {\mathbb R}^{n},
\hfill
$
\\[1.25ex]
where the vectors $A_\xi^{} = (a_{\xi 1}^{},\ldots,a_{\xi n}^{}),\ \xi=1,\ldots, n.$ 
\vspace{0.35ex}

The differential system (1.1) on a domain $G\subset {\mathbb R}^{n+1}$ 
has a basis of first integrals of  dimension  $n.$ 
Moreover, the differential system (1.1) has also $n-1$ autonomous functionally independent first integrals 
(autonomous integral basis) on a domain ${\mathscr X}\subset {\mathbb R}^{n}.$ 

The aim of Section~1 is to build these bases of first integrals for system (1.1).
\\[1.75ex]
\indent
{\bf 1.1.1. Linear partial integral}
\\[0.75ex]
\indent
A complex-valued linear homogeneous function
\\[1.25ex]
\mbox{}\hfill                                % (1.2)
$
\displaystyle
p\colon x\to
\sum\limits_{\xi=1}^{n} \nu_{\xi}^{}\;\! x_{\xi}^{}
\ 
$
for all $x\in {\mathbb R}^{n}
\quad 
(\nu_{\xi}^{}\in {\mathbb C},\  \xi=1,\ldots, n)
$
\hfill (1.2)
\\[1.25ex]
is a partial integral of system (1.1) if and only if
\\[1.5ex]
\mbox{}\hfill                                  % (1.3)
$
{\frak A} p(x) =  \lambda p(x)
\ 
$
for all 
$
x\in \R^{n},
\quad 
\lambda\in {\mathbb C}.
$
\hfill (1.3)
\\[1.75ex]
\indent
The identity (1.3) is equivalent to the linear homogeneous system 
\\[1.25ex]
\mbox{}\hfill                         % (1.4)
$
(B -  \lambda E)\nu = 0,
$
\hfill (1.4)
\\[1.25ex]
where the column vector
\vspace{0.25ex}
$\nu=\mbox{colon}(\nu_{1}^{},\ldots,\nu_{n}^{})\in {\mathbb C}^n, \ E$ is 
the $n\times n$ identity matrix,
and the matrix $B$ is the transpose of the matrix $A.$

The linear homogeneous system (1.4) has a nontrivial solution if and only if
\\[1.25ex]
\mbox{}\hfill                                             % (1.5)
$
\det (B-\lambda E) =  0.
$
\hfill (1.5)
\\[1.25ex]
\indent
We shall say that the equation (1.5) is the {\it integral characteristic equation} of system (1.1), 
and its roots are {\it integral characteristic roots} of system (1.1). 
Besides, a solution $\nu$ of the li\-ne\-ar system (1.4) is an eigenvector of the matrix $B$ 
corresponding to the eigenvalue $\lambda.$

Thus we have proved the following statement.
\vspace{0.35ex}

{\bf Lemma 1.1.}                                  
{\it
The linear homogeneous function {\rm (1.2)} is 
a partial integral of the differential system {\rm (1.1)} if and only if
the vector $\nu\in {\mathbb C}^n$ is an eigenvector of the matrix $B.$ 
\vspace{0.5ex}
}

Therefore a linear partial integral of the linear homogeneous differential system {\rm (1.1)} is generated by  
an eigenvector of the matrix $B.$ 

Lemma 1.1 is base for ours spectral method of building first integrals of  system (1.1).
In addition, first integrals of  the differential system (1.1) are building by eigenvectors and 
eigenvalues of the matrix $B$ with using orders of elementary divisors.
\\[1.75ex]
\indent
{\bf 1.1.2. Autonomous first integrals}
\\[0.5ex]
\indent
{\sl Case of real eigenvalues}.
If the matrix $B$ has two linearly independent real eigenvectors, then we can find an 
autonomous first integral of the differential system (1.1) by using following assertions
(Theorem 1.1, Corollaries 1.1 and 1.2).
%\vspace{0.75ex}

{\bf Theorem 1.1.}                                   % 1
{\it
Suppose $\nu^{1}$ and $\nu^{2}$ are real eigenvectors of  the matrix $B$ 
corresponding to the distinct eigenvalues $\lambda_{1}^{}$ and $\lambda_{2}^{} \ (\lambda_{1}^{}\ne\lambda_{2}^{}),$ respectively. 
\vspace{0.25ex}
Then the ordinary linear homogeneous differential system {\rm (1.1)} has the autonomous first integral
\\[2ex]
\mbox{}\hfill                                      % (1.6)
$
F\colon x\to \bigl|\nu^{1}x\bigr|^{h_1^{}}
\bigl|\nu^{2}x{\bigr|}^{h_2^{}} 
$ 
\ for all 
$
x\in {\mathscr  X},
$
\hfill {\rm (1.6)}
\\[2ex]
where ${\mathscr  X}$ is a domain of the domain of function ${\rm D}F\subset {\mathbb R}^n,$ 
\vspace{0.25ex}
the numbers $h_1^{}$ and $h_2^{}$ are a real solution to the equation 
$\lambda_1^{}h_1^{}+\lambda_2^{}h_2^{}=0$
with $|h_1^{}|+|h_2^{}| \not=0.$ 
}
\vspace{0.75ex}

{\sl Proof}. By Lemma 1.1, the linear functions
\\[1.5ex]
\mbox{}\hfill
$
p_k^{}\colon x\to \nu^k x
$
\ for all 
$
x\in {\mathbb R}^n,
\quad k=1,\ k=2,
\hfill
$
\\[1.5ex]
are partial integrals of system (1.1). Hence, 
\\[1.5ex]
\mbox{}\hfill
$
{\frak A}\,\nu^k x =\lambda_k^{}\,\nu^k x
$
\ for all 
$
x\in {\mathbb R}^n,
\quad k=1,\ k=2.
\hfill
$
\\[1.5ex]
\indent
Using these identities, we have
\\[2ex]
\mbox{}\hfill
$
{\frak A} F(x)
=\bigl|\nu^{1}x\bigr|^{h_{1}^{}-1}\bigl|\nu^{2}x\bigr|^{h_{2}^{}-1}\,
\bigl(h_{1}^{}\,\mbox{sgn}(\nu^{1}x)\,
\bigl|\nu^{2}x\bigr|\ {\frak A}\,\nu^{1}x +
h_{2}^{}\,\mbox{sgn}(\nu^{2}x)\,\bigl|\nu^{1}x\bigr|\
{\frak A}\, \nu^{2}x\;\!\bigr)=
\hfill
$
\\[2.75ex]
\mbox{}\hfill
$
=(\lambda_{1}^{}h_{1}^{}+\lambda_{2}^{}h_{2}^{})\;\! F(x)
$
\ for all 
$
x\in {\mathscr  X}\subset {\rm D}F.
\hfill
$
\\[2.25ex]
\indent
If $h_1^{}$ and $h_2^{}$ are real numbers such that 
\vspace{0.5ex}
$\lambda_{1}^{}h_{1}^{}+\lambda_{2}^{}h_{2}^{}=0$ with $|h_1^{}|+|h_2^{}| \not=0,$ 
then the scalar function (1.6) is an autonomous first integral of system (1.1). \k
\vspace{0.75ex}

{\bf Corollary 1.1.}                         
{\it
If $\nu$ is a real eigenvector of the matrix $B$ corresponding to the eigenvalue $\lambda=0,$
then the linear homogeneous function
\\[1.5ex]
\mbox{}\hfill                           % (1.7)
$
F\colon x\to \nu x
$
\ for all 
$
x\in {\mathbb R}^n
$ 
\hfill {\rm (1.7)}
\\[1.5ex]
is an autonomous first integral of the ordinary linear homogeneous differential system {\rm (1.1).}
}
\vspace{0.5ex}

{\sl Indeed}, 
suppose $\nu$ is a real eigenvector of the matrix $B$ corresponding to the eigenvalue $\lambda=0.$
Then, from the identity (1.3) it follows that 
${\frak A}\;\!\nu x = 0$ for all $x\in {\mathbb R}^{n}.$
\vspace{0.25ex}
This yields that the function (1.7) is an autonomous first integral on the space ${\mathbb R}^n$ of system {\rm (1.1)}. \k
\vspace{0.75ex}

{\bf Corollary 1.2.} 
{\it
Let $\lambda\ne 0$ be an eigenvalue of the matrix $B$ corresponding to two 
real linearly independent eigenvectors $\nu^{1}$ and $\nu^{2}.$ Then 
the ordinary linear homogeneous differential system {\rm (1.1)} has the autonomous first integral
\\[2ex]
\mbox{}\hfill                                            % (1.8)
$
F\colon x\to \dfrac{\nu^1 x}{\nu^2 x}
$
\ for all 
$
x\in {\mathscr X},
$
\hfill {\rm (1.8)}
\\[2ex]
where ${\mathscr X}$ is a domain from the set $\{x\colon \nu^2 x\ne 0\}\subset {\mathbb R}^n.$
}
\vspace{0.75ex}

Let us remark that the first integrals (1.7) and (1.8) are algebraic.
Note also that algebraicity of the first integral (1.6) depends on the numbers $h_1^{}$ and $h_2^{}.$ 
For example, if the numbers $h_1^{}$ and $h_2^{}$ are rational, then the first integral (1.6) is algebraic.
But this condition isn't necessary for algebraicity of the first integral (1.6). 
At the same time we have
\vspace{0.5ex}

{\bf Property 1.1}
{\sl {\rm(}sufficient condition for algebraicity of basis of autonomous first integrals{\rm)}.}
{\it
If all eigenvalues of the matrix $B$ are simple and rational, then 
the differential system {\rm (1.1)} has a basis of autonomous algebraic first integrals.
}
\vspace{0.75ex}

{\bf Example 1.1.}
Consider the fourth-order ordinary linear differential system
\\[2.5ex]
\mbox{}\hfill                              % (1.9)
$
\begin{array}{ll}
\dfrac{dx_1^{}}{dt}=   x_{1}^{} - 2x_{2}^{} -  x_{4}^{}, &
\dfrac{dx_2^{}}{dt} = {}- x_{1}^{} + 4x_{2}^{} - x_{3}^{} + 2x_{4}^{},
\\[4.25ex]
\dfrac{dx_3^{}}{dt} =  2x_{2}^{} +  x_{3}^{} + x_{4}^{},\ \ \ \   &
\dfrac{dx_4^{}}{dt} =  2x_{1}^{} - 4x_{2}^{} + 2x_{3}^{} - 2x_{4}^{}.
\end{array}
$
\hfill (1.9)
\\[2.75ex]
\indent
We claim that the matrix
\\[1ex]
\mbox{}\hfill
$
B = \left\|\!\!
\begin{array}{rrrr}
   1& {}-1&    0&    2
\\[0.5ex]
{}- 2 &    4&    2& {}- 4
\\[0.5ex]
   0& {}-1&    1&    2
\\[0.5ex]
{}- 1&    2&    1& {}- 2
\end{array}
\!\! \right\|
\hfill
$
\\[2ex]
has the eigenvalues $\lambda_{1}^{}=0,\  \lambda_{2}^{}=\lambda_{3}^{}=1,$ and $\lambda_{4}^{}=2.$
Indeed, the characteristic equation 
\\[1.75ex]
\mbox{}\hfill
$
\left|\!\!
\begin{array}{rrrr}
1-\lambda & {}-1 & 0 & 2
\\[0.75ex]
{}-2 & 4-\lambda & 2 & -4
\\[0.75ex]
0 & {}-1 & 1-\lambda & 2
\\[0.75ex]
{}-1 & 2 & 1 & {}-2-\lambda
\end{array}
\!\!\right| =0
\iff
\lambda (\lambda-1)^{2}(\lambda-2)=0.
\hfill
$
\\[2ex]
\indent
The rank of the matrix $B-\lambda_{2}^{}E$ is equal 2.
\vspace{0.25ex}
Therefore the double eigenvalue $\lambda_{2}^{}=1$ of the matrix $B$ has 
$\kappa=4-2=2$ simple elementary divisors $\lambda-1$ and $\lambda-1.$
\vspace{0.25ex}

The matrix $B$ has four simple elementary divisors $\lambda,\ \lambda-1,\ \lambda-1,$ and $\lambda-2.$
\vspace{0.25ex}

The linear homogeneous system 
\\[1ex]
\mbox{}\hfill
$
(B-\lambda_{1}^{}E)\,\mbox{colon}(\nu_{1}^{},\ldots,\nu_{4}^{}) = 0
\!\iff\!
\left\{\!\!
\begin{array}{r}
\nu_{1}^{}-\nu_{2}^{}+2\nu_{4}^{}=0,
\\[1ex]
{}-2\nu_{1}^{}+4\nu_{2}^{}+2\nu_{3}^{}-4\nu_{4}^{}=0,
\\[1ex]
{}-\nu_{2}^{}+\nu_{3}^{}+2\nu_{4}^{}=0,
\\[1ex]
{}-\nu_{1}^{}+2\nu_{2}^{}+\nu_{3}^{}-2\nu_{4}^{}=0\,
\end{array}
\right.
\!\iff\!
\left\{\!\!
\begin{array}{l}
\nu_{1}^{}={}-\nu_{4}^{},
\\[1ex]
\nu_{2}^{}=\nu_{4}^{},
\\[1ex]
\nu_{3}^{}={}-\nu_{4}^{}.
\end{array}
\right.
\hfill
$
\\[2ex]
\indent
Hence $\nu^{1}=({}-1,1,{}-1,1)$  is a real eigenvector corresponding  to the eigenvalue $\lambda_{1}^{}=0$
of the matrix $B.$

The scalar function (by Corollary 1.1)
\\[1.5ex]
\mbox{}\hfill
$
F_{1}^{}\colon x\to {}-x_{1}^{}+x_{2}^{}-x_{3}^{}+x_{4}^{} 
$
\ for all 
$
x\in {\mathbb R}^4
\hfill
$
\\[1.5ex]
is an autonomous linear first integral of the ordinary differential system (1.9).
\vspace{0.35ex}

The linear homogeneous system 
\\[1ex]
\mbox{}\hfill
$
(B-\lambda_{2}^{}E)\,\mbox{colon}(\nu_{1}^{},\ldots,\nu_{4}^{}) = 0
\!\iff\!
\left\{\!\!
\begin{array}{r}
{}-\nu_{2}^{}+2\nu_{4}^{}=0,
\\[1ex]
{}-2\nu_{1}^{}+3\nu_{2}^{}+2\nu_{3}^{}-4\nu_{4}^{}=0,
\\[1ex]
{}-\nu_{2}^{}+2\nu_{4}^{}=0,
\\[1ex]
{}-\nu_{1}^{}+2\nu_{2}^{}+\nu_{3}^{}-3\nu_{4}^{}=0\,
\end{array}
\right.
\!\!\!\iff\!
\left\{\!\!
\begin{array}{l}
\nu_{1}^{} = \nu_{3}^{}+\nu_{4}^{},
\\[1ex]
\nu_{2}^{} = 2\nu_{4}^{}.
\end{array}
\right.
\hfill
$
\\[2ex]
\indent
Hence $\nu^{2}=(2,2,1,1)$ and $\nu^{3}=(1,0,1,0)$  are two 
\vspace{0.25ex}
linearly independent real eigenvectors corresponding  to the double eigenvalue 
$\lambda_{2}^{}=1$ of the matrix $B.$
\vspace{0.25ex}

The scalar function (by Corollary 1.2)
\\[2ex]
\mbox{}\hfill
$
F_{23}^{}\colon x\to
\dfrac{2x_{1}^{}+2x_{2}^{}+x_{3}^{}+x_{4}^{}}{x_{1}^{}+x_{3}^{}}
$
\ for all 
$
x\in {\mathscr X}_1^{},
\hfill
$
\\[2.25ex]
where a domain ${\mathscr X}_1^{}\subset \{x\colon x_{1}^{}+x_{3}^{}\ne 0\},$
\vspace{0.35ex}
is an autonomous first integral of system (1.9).

The linear homogeneous system 
\\[1ex]
\mbox{}\hfill
$
(B-\lambda_{4}^{}E)\,\mbox{colon}(\nu_{1}^{},\ldots,\nu_{4}^{}) = 0
\iff
\left\{\!\!
\begin{array}{r}
{}-\nu_{1}^{}-\nu_{2}^{}+2\nu_{4}^{}=0,
\\[1ex]
{}-2\nu_{1}^{}+2\nu_{2}^{}+2\nu_{3}^{}-4\nu_{4}^{}=0,
\\[1ex]
{}-\nu_{2}^{}-\nu_{3}^{}+2\nu_{4}^{}=0,
\\[1ex]
{}-\nu_{1}^{}+2\nu_{2}^{}+\nu_{3}^{}-4\nu_{4}^{}=0\,
\end{array}
\right.
\iff
\left\{\!\!
\begin{array}{l}
\nu_{1}^{} = \nu_{3}^{},
\\[1ex]
\nu_{2}^{} = 2\nu_{4}^{}.
\end{array}
\right.
\hfill
$
\\[2ex]
\indent
Hence $\nu^{4}=(0,2,0,1)$ is a real eigenvector corresponding  to the eigenvalue $\lambda_{4}^{}=2$
of the matrix $B.$

Using the linearly independent real eigenvectors $\nu^{2}$ and $\nu^{4}$ of the matrix $B,$ we can build (by Theorem 1.1) an
autonomous first integral of the differential system (1.9).
Since the equation $h_{2}^{}+2h_{4}^{} = 0,$ we have, for example, $h_{2}^{}=2,\  h_{4}^{}={}-1.$
The scalar function 
\\[2.25ex]
\mbox{}\hfill
$
F_{24}^{}\colon x\to
\dfrac{(2x_{1}^{}+2x_{2}^{}+x_{3}^{}+x_{4}^{})^{2}}{2x_{2}^{}+x_{4}^{}} 
$
\ for all 
$
x\in {\mathscr X}_2^{},
\hfill
$
\\[2.25ex]
where a domain ${\mathscr X}_2^{}\subset \{x\colon 2x_{2}^{}+x_{4}^{}\ne 0\},$ 
\vspace{0.35ex}
is an autonomous first integral of system (1.9).

The functionally independent first integrals $F_{1}^{},\ F_{23}^{},$ and $F_{24}^{}$ are 
\vspace{0.35ex}
an autonomous integral basis of system (1.9) on any domain ${\mathscr X}$ from the set
\vspace{1.25ex}
$\{x\colon x_{1}^{}+x_{3}^{}\ne 0,\ 2x_{2}^{}+x_{4}^{}\ne 0\}\subset {\mathbb R}^4.$

{\sl Case of complex eigenvalues}.
Let the function (1.2) be a complex-valued partial integral of system (1.1).
Then, from the identity (1.3)  it follows that
\\[1.75ex]
\mbox{}\hfill                                         % (1.10)
$
{\frak A}\,{\rm Re}\, p(x) =
{\stackrel{*}{\lambda}}\, {\rm Re}\,p(x) \, - \ 
\widetilde{\lambda}\; {\rm Im}\,p(x)
$
\ for all 
$
x\in {\mathbb R}^n,
\hfill
$
\\[-0.25ex]
\mbox{}\hfill (1.10)
\\
\mbox{}\hfill
$
{\frak A}\,{\rm Im}\,p(x) =
\widetilde{\lambda}\, {\rm Re}\,p(x) \, + \
{\stackrel{*}{\lambda}}\;
{\rm Im}\, p(x)
$
\ for all 
$
x\in {\mathbb R}^n,
%\ \ \lambda = {\stackrel{*}{\lambda}} +
%\widetilde{\lambda}\,i.
\hfill
$
\\[1.75ex]
where the real numbers ${\stackrel{*}{\lambda}}={\rm Re}\,\lambda,\ \widetilde{\lambda}={\rm Im}\,\lambda.$
Thus we have the following criteria for the existence of a complex-valued partial integral of system (1.1).
\vspace{0.5ex}

{\bf Lemma 1.2.}
{\it
The function {\rm (1.2)} is a complex-valued partial integral of system {\rm (1.1)} if and only if 
the system of identities {\rm (1.10)} holds.
}
\vspace{0.5ex}

Using Lemma 1.2, we may establish the following propositions.
\vspace{0.5ex}

{\bf Property 1.2.}
{\it
If the system {\rm(1.1)} has a complex-valued partial integral {\rm (1.2)}, then
the complex conjugate function $\overline{p}$ is a complex-valued partial integral of system {\rm (1.1)} and 
\\[1.5ex]
\mbox{}\hfill
$
{\frak A}\,\overline{p}\;\!(x) =
\overline{\lambda}\; \overline{p}(x)
$
\ for all 
$
x\in {\mathbb R}^n,
\hfill
$
\\[1.5ex]
where $\overline{\lambda}$ is 
the conjugate of the complex number $\lambda$ from the identity {\rm (1.3)}.
}
\vspace{0.5ex}

{\sl Proof}. Using the system of identities (1.10), we get
\\[1.5ex]
\mbox{}\hfill                        
$
{\frak A}\,\overline{p}\;\!(x) = 
{\frak A}\,{\rm Re}\;\!p(x) -\;\! 
i\,{\frak A}\,{\rm Im}\;\!p(x)  =
{\stackrel{*}{\lambda}}\, {\rm Re}\,p(x)  - \;\! 
\widetilde{\lambda}\; {\rm Im}\,p(x) -
i\;\!\bigl(\widetilde{\lambda}\, {\rm Re}\,p(x)  + \;\! 
{\stackrel{*}{\lambda}}\;{\rm Im}\, p(x)\bigr)=
\hfill                        
$
\\[1.75ex]
\mbox{}\hfill                        
$
=\bigl({\stackrel{*}{\lambda}}- i\,\widetilde{\lambda}\,\bigr)
\bigl({\rm Re}\,p(x)-i\,{\rm Im}\,p(x)\bigr)=
\overline{\lambda}\; \overline{p}(x)
\quad
$
\ for all 
$
x\in \R^n.
\, \k
\hfill                        
$
\\[1.35ex]
\indent
Let $\K$ be either the real number field $\R$ or the complex number field $\C.$
\vspace{0.75ex}

{\bf Property 1.3.}
{\it
\vspace{0.15ex}
The product $u_1^{}u_2^{}$ of the polynomials $u_1^{}\colon {\mathbb R}^n\to{\mathbb K}$ and 
$u_2^{}\colon{\mathbb R}^n\to{\mathbb K}$ is a polynomial partial integral of system {\rm(1.1)} 
if and only if the polynomials $u_1^{}$ and  $u_2^{}$ are 
polynomial partial integrals of system {\rm(1.1)}. 
}

{\sl Proof}. From definition of partial integral and
\\[1.5ex]
\mbox{}\hfill                        
$
{\frak A}\;\!\bigl(u_1^{}(x)\,u_2^{}(x)\bigr)=
u_2^{}(x)\,{\frak A}\,u_1^{}(x)+
u_1^{}(x)\,{\frak A}\,u_2^{}(x)
$
\ for all 
$
x\in {\mathbb R}^n,
\hfill                        
$
\\[1.5ex]
we get the assertion of Property 1.3. \k
\vspace{0.5ex}

{\bf Property 1.4.}
{\it
The system {\rm(1.1)} has the complex-valued polynomial partial integral {\rm (1.2)}
if and only if 
the real polynomial
\\[1.5ex]
\mbox{}\hfill                                    
$
P\colon x\to \ 
{\rm Re}^2\,p(x) +\,  {\rm Im}^2\,p(x)
$
\ for all 
$
x\in {\mathbb R}^n
\hfill 
$
\\[1.5ex]
is a partial integral of system {\rm (1.1)}. Moreover, the following identity holds
\\[1.25ex]
\mbox{}\hfill                                    
$
{\frak A}\,P(x) =2\;\! {\stackrel{*}{\lambda}}\;\!P(x)
$
\ for all 
$
x\in\R^n,
\hfill 
$
\\[1.25ex]
where the number $\lambda = {\stackrel{*}{\lambda}} + i\,\widetilde{\lambda}\ 
({\stackrel{*}{\lambda}}={\rm Re}\,\lambda,\ \widetilde{\lambda}={\rm Im}\,\lambda)$
is given by the identity {\rm (1.3)}.
}
\vspace{0.5ex}

{\sl Proof.} Using Properties 1.2 and 1.3, and the identity
\\[1.5ex]
\mbox{}\hfill                        
$
p(x)\,\overline{p}\;\!(x)=
{\rm Re}^2\,p(x)+{\rm Im}^2\,p(x)
$
\ for all $
x\in {\mathbb R}^n,
\hfill                        
$
\\[1.5ex]
we have the polynomial $P$ is a real partial integral of system {\rm (1.1)} if and only if 
the linear function (1.2) is a complex-valued partial integral of system {\rm (1.1)}. 

By the system of identities (1.10), it follows that
\\[1.5ex]
\mbox{}\hfill                        
$
{\frak A}\bigl({\rm Re}^2\,p(x)+{\rm Im}^2\,p(x)\bigr)=
2\;\!{\rm Re}\,p(x)\,{\frak A}\,{\rm Re}\,p(x)+
2\;\!{\rm Im}\,p(x)\,{\frak A}\,{\rm Im}\,p(x)=
\hfill                        
$
\\[1.75ex]
\mbox{}\hfill                        
$
=2\;\!{\rm Re}\,p(x)
\bigl({\stackrel{*}{\lambda}}\, {\rm Re}\,p(x) -  
\widetilde{\lambda}\; {\rm Im}\,p(x)\bigr)  + 
2\;\!{\rm Im}\,p(x)\bigl(\widetilde{\lambda}\, {\rm Re}\,p(x) + 
{\stackrel{*}{\lambda}}\;{\rm Im}\, p(x)\bigr)=
\hfill                        
$
\\[1.75ex]
\mbox{}\hfill                        
$
= 2\;\!{\stackrel{*}{\lambda}}\;\!\bigl({\rm Re}^2\,p(x)+{\rm Im}^2\,p(x)\bigr)=
2\;\!{\stackrel{*}{\lambda}}\,P(x)
$
\ for all
$
x\in \R^n.
\ \k
\hfill                        
$
\\[2ex]
\indent
{\bf Property 1.5.}
{\it
Let the function {\rm (1.2)} be a complex-valued partial integral of system {\rm (1.1)}. 
Then the Lie derivative of the function
\\[1.75ex]
\mbox{}\hfill
$
\varphi\colon x\to\ \arctan\dfrac{{\rm Im}\,p(x)}{{\rm Re}\, p(x)}
$
\ for all 
$
x\in {\mathscr X}\subset \{x\colon {\rm Re}\,p(x)\ne 0\}
\hfill
$
\\[1.75ex]
by virtue of system {\rm (1.1)} is equal to 
\\[1.25ex]
\mbox{}\hfill                                   
$
{\frak A}\,\varphi(x) =
\widetilde{\lambda}
$
\ for all 
$
x\in {\mathscr X},
\hfill
$
\\[1.75ex]
where the number $\lambda = {\stackrel{*}{\lambda}} + i\, \widetilde{\lambda}\ 
({\stackrel{*}{\lambda}}={\rm Re}\,\lambda,\ \widetilde{\lambda}={\rm Im}\,\lambda)$
is given by the identity {\rm (1.3)}.
}
\vspace{1ex}

{\sl Proof.} 
Using the identities (1.10), we obtain
\\[2.25ex]
\mbox{}\hfill                        
$
{\frak A}\,\varphi(x)\, =\, 
{\frak A}\,\arctan\dfrac{{\rm Im}\,p(x)}{{\rm Re}\,p(x)} \, =\,
\dfrac{1}{1+\dfrac{{\rm Im}^2\,p(x)}{{\rm Re}^2\,p(x)}}
\cdot
\dfrac{{\rm Re}\,p(x)\,{\frak A}\,{\rm Im}\,p(x)- {\rm Im}\,p(x)\,{\frak A}\,{\rm Re}\,p(x)}
{{\rm Re}^2\,p(x)}\ =
\hfill                        
$
\\[2.25ex]
\mbox{}\hfill                        
$
=
\dfrac{{\rm Re}\,p(x)
\bigl({\stackrel{*}{\lambda}}\, {\rm Re}\,p(x) -  \widetilde{\lambda}\; {\rm Im}\,p(x)\bigr)  -
{\rm Im}\,{\frak m}(x)
\bigl(\widetilde{\lambda}\, {\rm Re}\,p(x) + 
{\stackrel{*}{\lambda}}\;{\rm Im}\, p(x)\bigr)}
{{\rm Re}^2\,p(x)+{\rm Im}^2\,p(x)}=
\widetilde{\lambda}
$
\ for all 
$
x\in {\mathscr X}.
\ \k
\hfill                        
$
\\[2.25ex]
\indent
If the matrix $B$ has an imaginary eigenvalue, then we can build an 
autonomous first integral of system (1.1) by using following assertions
(Theorems 1.2, 1.3, and 1.4).
\vspace{0.5ex}

{\bf Theorem 1.2.}\!
\vspace{0.25ex}
{\it
Suppose $\nu={\stackrel{*}{\nu}}+\widetilde{\nu}\,i\
({\rm Re}\,\nu={\stackrel{*}{\nu}},\ {\rm Im}\,\nu=\widetilde{\nu}\,)$ is an 
eigenvector of the mat\-rix $B$ corresponding  to the imaginary eigenvalue
\vspace{0.5ex}
$\lambda={\stackrel{*}{\lambda}}+\widetilde{\lambda}\,i\ 
({\rm Re}\,\lambda={\stackrel{*}{\lambda}},\ {\rm Im}\, \lambda=\widetilde{\lambda}\not=0).\!$ 
Then the ordinary differential system {\rm (1.1)}  has the autonomous first integral}
\\[1.75ex]
\mbox{}\hfill                        % (1.11)
$
F\colon x\to
\bigl((\,{\stackrel{*}{\nu}}x)^2
+(\,\widetilde{\nu}x)^2\,\bigr)
\exp\biggl(
{}-2 \
\dfrac{{\stackrel{*}{\lambda}}}{\widetilde{\lambda}} \
\arctan\dfrac{\widetilde{\nu}x}{{\stackrel{*}{\nu}}x}
\,\biggr)
$
\ for all 
$
x\in {\mathscr X},
$
\hfill (1.11)
\\[1.75ex]
{\it where ${\mathscr X}$ is a domain from the set 
$\bigl\{x\colon {\stackrel {*}{\nu}}x\ne 0\bigr\}.$ 
}
\vspace{0.75ex}

{\sl Proof}. Taking into account Properties 1.4 and 1.5, we get
\\[1.5ex]
\mbox{}\hfill
$
{\frak A}F(x)=
\exp\biggl(\!\!
{}-2\,
\dfrac{{\stackrel{*}{\lambda}}}{\widetilde{\lambda}}\;\!
\arctan\dfrac{\widetilde{\nu}x}{{\stackrel{*}{\nu}}x}\biggr)\,
{\frak A}\bigl(({\stackrel{*}{\nu}}x)^2+(\widetilde{\nu}x)^2\bigr) + 
\bigl(({\stackrel{*}{\nu}}x)^2+(\widetilde{\nu}x)^2\bigr)\,
{\frak A}
\exp\biggl(\!\!
{}-2\,\dfrac{{\stackrel{*}{\lambda}}}{\widetilde{\lambda}}\;\!
\arctan\dfrac{\widetilde{\nu}x}{{\stackrel{*}{\nu}}x}\biggr)
=
\hfill
$
\\[2.25ex]
\mbox{}\hfill
$
=\biggl(2{\stackrel{*}{\lambda}}-
2\,\dfrac{{\stackrel{*}{\lambda}}}{\widetilde{\lambda}}\ \widetilde{\lambda}\biggr)
\bigl((\,{\stackrel{*}{\nu}}x)^2
+(\,\widetilde{\nu}x)^2\,\bigr)
\exp\biggl(
{}-2 \
\dfrac{{\stackrel{*}{\lambda}}}{\widetilde{\lambda}} \
\arctan\dfrac{\widetilde{\nu}x}{{\stackrel{*}{\nu}}x}
\,\biggr) =0
$
\ for all 
$
x\in {\mathscr X}.
\hfill
$
\\[2ex]
\indent
This implies that the function (1.11) is an autonomous first integral of system (1.1). \k
\vspace{1.25ex}

Transcendency of the first integral (1.11) of systems (1.1) depends on 
\vspace{0.25ex}
the imaginary eigenvalue $\lambda={\stackrel{*}{\lambda}}+\widetilde{\lambda}\,i$
(if ${\stackrel{*}{\lambda}}=0,$ then the first integral (1.11) of system (1.1) is algebraic).
\vspace{1ex}

The proof of Theorems 1.3 and 1.4 is similar to that one in Theorem 1.2. 
\vspace{0.5ex}

\newpage

{\bf Theorem 1.3.}\!
\vspace{0.25ex}
{\it
Let $\nu^{1}\!={\stackrel{*}{\nu}}{}^{1}+ \widetilde{\nu}{}^{\,1}\,i
\ ({\rm Re}\,\nu^1\!={\stackrel{*}{\nu}}{}^{\,1},\, {\rm Im}\,\nu^1\!=\widetilde{\nu}{}^{\,1})$
be an eigenvector of the matrix $B$ cor\-res\-pon\-ding to the imaginary eigenvalue
\vspace{0.75ex}
$\lambda_1^{}={\stackrel{*}{\lambda}}_1^{}+\widetilde{\lambda}_1^{}\,i
\ ({\rm Re}\,\lambda_1^{}\!=\!{\stackrel{*}{\lambda}}_1^{},\ 
{\rm Im}\,\lambda_1^{}=\widetilde{\lambda}_1^{}\not= 0),\ 
\nu^2$ be a real eigenvector of the matrix $B$ cor\-res\-pon\-ding to the eigenvalue 
\vspace{0.35ex}
$\lambda_2^{}\ne 0.$ Then the ordinary differential sys\-tem {\rm (1.1)} 
has the autonomous first integral}
\\[1.25ex]
\mbox{}\hfill
$
F\colon x\to
\nu^2x\,\exp\biggl(
{}-\dfrac{\lambda_2^{}}{\widetilde{\lambda}_1^{}} \
\arctan\dfrac{\widetilde{\nu}{}^{\,1}x}
{{\stackrel{*}{\nu}}{}^1x}\,\biggr)
$
\ for all 
$
x\in {\mathscr X},
\hfill
$
\\[1.25ex]
{\it
where ${\mathscr X}$ is a domain from the set
$\bigl\{x\colon {\stackrel{*}{\nu}}{}^1x\ne 0\bigr\}.$
}
\vspace{1ex}

{\bf Theorem 1.4.}\!
\vspace{0.25ex}
{\it
Let
$\nu^1={\stackrel{*}{\nu}}{}^1+\widetilde{\nu}{}^{\,1}\,i$
and $\nu^2={\stackrel{*}{\nu}}{}^2+\widetilde{\nu}{}^{\,2}\,i\ 
\ ({\rm Re}\,\nu^{\tau}={\stackrel{*}{\nu}}{}^{\,{\tau}},\, {\rm Im}\,\nu^{\tau}=\widetilde{\nu}{}^{\,{\tau}},\ \tau=1, 2)$
be two linearly independent eigenvectors of the matrix $B$ corresponding to the imaginary eigenvalues
\vspace{0.35ex}
$\lambda_1^{}\!=\!
{\stackrel{*}{\lambda}}_1^{} +\widetilde{\lambda}_1^{}\,i$  and
$\lambda_2^{}\!=\!{\stackrel{*}{\lambda}}_2^{} + \widetilde{\lambda}_2^{}\,i\,   
({\rm Re}\,\lambda_{\tau}\!=\!{\stackrel{*}{\lambda}}{}_{\,{\tau}},\ 
{\rm Im}\,\lambda_{\tau}\!=\!\widetilde{\lambda}{}_{\,{\tau}},\tau\!=\!1, 2)$
\vspace{0.5ex}
with $\lambda_1^{}\ne \overline{\lambda}_2^{}.$ 
Then the or\-di\-na\-ry dif\-fe\-ren\-ti\-al sys\-tem {\rm (1.1)} 
has the autonomous first integral
\\[2ex]
\mbox{}\hfill
$
F\colon x\to\
\widetilde{\lambda}_{1}\,
\arctan\dfrac{\widetilde{\nu}{}^{\,2}x}
{{\stackrel{*}{\nu}}{}^{2}x}
\ -\ \widetilde{\lambda}_{2}\,
\arctan\dfrac{\widetilde{\nu}{}^{\,1}x}
{{\stackrel{*}{\nu}}{}^{1}x}
$ 
\ for all 
$
x\in {\mathscr X},
\hfill
$
\\[2ex]
where ${\mathscr X}$ is a domain from the set
\vspace{1ex}
$\bigl\{x\colon  {\stackrel{*}{\nu}}{}^{1}x\ne 0,\ 
{\stackrel{*}{\nu}}{}^{2}x\ne 0\bigr\}.$
}

{\bf Example 1.2.}
%\marginpar{\No\ 02\,443}
The autonomous differential system
\\[2.25ex]
\mbox{}\hfill                              % (1.12)
$
\dfrac{dx_1^{}}{dt} =  2x_{1}^{} + x_{2}^{},
\qquad 
\dfrac{dx_2^{}}{dt} = x_{1}^{} + 3x_{2}^{} -  x_{3}^{},
\qquad 
\dfrac{dx_3^{}}{dt} =  {}-x_{1}^{} + 2x_{2}^{} + 3x_{3}^{}
$
\hfill (1.12)
\\[2.25ex]
has the eigenvalues $\lambda_{1}^{}=3+i$ and $\lambda_{2}^{}=2$ corresponding 
\vspace{0.25ex}
to the eigenvectors $\nu^{1}=(1,i,-1)$ and $\nu^{2}=(3,-1,-1),$ respectively. 
The first integrals (by Theorem 1.2)
\\[2ex]
\mbox{}\hfill                        
$
F_{1}^{}\colon x\to\,
\bigl( (x_{1}^{}-x_{3}^{})^2+x_{2}^2\,\bigr)
\exp\Bigl( {}-6\;\!\arctan\dfrac{x_2}{x_1^{}-x_3^{}}\,\Bigr)
$
\ for all $x\in {\mathscr X}
\hfill 
$
\\[1.75ex]
and (by Theorem 1.3)
\\[2ex]
\mbox{}\hfill                          
$
F_{2}^{}\colon x\to\
(3x_1^{}-x_2^{}-x_3^{})\exp\Bigl( {}-2\;\!
\arctan\dfrac{x_2^{}}{x_1^{}-x_3^{}}\Bigr)
$
\ for all 
$
x\in {\mathscr X}
\hfill
$
\\[2ex]
are an autonomous integral basis of system (1.12) 
on any domain ${\mathscr X}\subset \{x\colon x_1^{}-x_3^{}\!\ne 0\}.$ 
\vspace{1ex}

{\bf Example 1.3.}
%\marginpar{\No\ 02\,445}
The autonomous differential system
\\[2ex]
\mbox{}\hfill                              % (1.13)
$
\begin{array}{ll}
\dfrac{dx_1^{}}{dt} = {}-3x_1^{} + x_2^{} + 4x_3^{}+ 2x_4^{},
\quad
&
\dfrac{dx_2^{}}{dt} =  8x_1^{} - 3x_2^{} -2x_3^{}+ 6x_4^{},
\\[3.75ex]
\dfrac{dx_3^{}}{dt} = {}-9x_1^{} + 3x_2^{} + 4x_3^{}- 4x_4^{},
&
\dfrac{dx_4^{}}{dt} = 6x_1^{} - 3x_2^{} - 4x_3^{}+ 2x_4^{}
\end{array}
\hfill (1.13)
$
\\[2.25ex]
has the eigenvalues $\!\lambda_1^{}\!=\!i,\, \lambda_2^{}\!=\!2i$ corresponding 
\vspace{0.25ex}
to the eigenvectors
$\!\nu^{1}\!=\!(1-i, -1+2i,2i,2),$ $\nu^{2}=(i,{}-1,i,1+2i).$  
The functionally independent first integrals (by Theorem 1.2)
\\[1.25ex]
\mbox{}\hfill                           
$
F_{1}^{}\colon x\to
(x_{1}^{}-x_{2}^{}+2x_{4}^{})^2+({}-x_{1}^{}+2x_{2}^{}+2x_{3}^{})^2
$
\ for all 
$
x\in {\mathbb R}^4,
\hfill
$
\\[2.5ex]
\mbox{}\hfill                              
$
F_{2}^{}\colon x\to ({}-x_2^{}+x_4^{})^{2}+(x_1^{}+x_3^{}+2x_4^{})^2
$
\ for all 
$
x\in {\mathbb R}^4,
\hfill
$
\\[1.75ex]
and (by Theorem 1.4)
\\[2ex]
\mbox{}\hfill                          
$
F_{3}^{}\colon x\to\
\arctan\dfrac{x_{1}^{}+x_{3}^{}+2x_{4}^{}}{{}-x_{2}^{}+x_{4}^{}}
-2\arctan\dfrac{{}-x_{1}^{}+2x_{2}^{}+2x_{3}^{}}{x_{1}^{}-x_{2}^{}+2x_{4}^{}}
$
\ for all 
$
x\in {\mathscr X}
\hfill
$
\\[1.75ex]
are an autonomous integral basis of the ordinary differential system (1.13) 
\vspace{0.25ex}
on any domain ${\mathscr X}$ from the set
$\{x\colon x_1^{}-x_2^{}+2x_4^{}\ne 0,\ x_2^{}-x_4^{}\ne 0\}\subset {\mathbb R}^4.$
\vspace{1.25ex}

{\sl Case of multiple elementary divisors}.
\vspace{0.25ex}

{\bf Definition 1.1.}
{\it
Let $\nu^{0}$ be an eigenvector of the matrix $B$ corresponding to the eigen\-va\-lue $\lambda$ with 
the elementary divisor of multiplicity $m.$ 
A non-zero vector $\nu^{k}\in {\mathbb C}^n$ is called a
\textit{\textbf{ge\-ne\-ra\-li\-zed eigenvector of order}} {\boldmath $k$} 
if the vector $\nu^{k}$ satisfying}
\\[1.25ex]
\mbox{}\hfill                                            % (1.14)
$
(B-\lambda E)\,\nu^{k}=\, k \;\! \nu^{k-1},
\quad 
k=1,\ldots, m-1.
$
\hfill (1.14)
\\[1.5ex]
\indent
In this case we can build first integrals of system (1.1) by using following assertions.
\vspace{0.35ex}

{\bf Theorem 1.5.}
{\it
Let $\lambda$ be the eigenvalue of the matrix $B$ 
with elementary divisor of multiplicity $m\ (m\geq 2)$
corresponding to the real eigenvector $\nu^{0}$ and to the real {\rm 1}-st order 
generalized eigenvector $\nu^{1}.$ 
Then the sys\-tem {\rm (1.1)} has the autonomous first integral
\\[1.5ex]
\mbox{}\hfill
$
F\colon x\to\,
\nu^{0}x\,
\exp\biggl( {}-\lambda\, \dfrac{\nu^{1}x}{\nu^{0}x}\,\biggr)
$
\ for all 
$
x\in {\mathscr X},
\quad 
{\mathscr X}\subset \{x\colon  \nu^0x\ne 0\}. 
\hfill
$
\\[1.5ex]
}
\indent
{\sl Proof}. Using the equalities (1.14) under the condition $k=1,$ we get 
\\[1.5ex]
\mbox{}\hfill
$
{\frak A}\,\nu^1x=\lambda\nu^1x+\nu^0x
$
\ for all 
$
x\in\R^n.
\hfill
$
\\[1.5ex]
\indent
Then, the Lie derivative of the function $F$ by virtue of system (1.1) is
\\[1.5ex]
\mbox{}\hfill
$
{\frak A}\,F(x)=
\Bigl(\lambda-\lambda\dfrac{(\lambda\nu^1x+\nu^0x)\nu^0x-\lambda\,\nu^0x\,\nu^1x}
{(\nu^0x)^2}\,\Bigr)F(x)=0
$
\ for all 
$
x\in {\mathscr X}.
\ \ \k
\hfill
$
\\[1.75ex]
\indent
From Theorem 1.5, we have the following
\vspace{0.75ex}

{\bf Corollary 1.3.}
\vspace{0.25ex}
{\it
Let $\lambda ={\stackrel{*}{\lambda}}+\widetilde{\lambda}\,i
\ ({\rm Re}\,\lambda={\stackrel{*}{\lambda}},\ {\rm Im}\,\lambda=\widetilde{\lambda}\ne 0)$ be
the complex eigenvalue of the matrix $B$ with elementary divisor of multiplicity $m\ (m\geq 2)$
corresponding to the eigenvector 
\vspace{0.25ex}
$\nu^{0}={\stackrel{*}{\nu}}{}^{0}+\widetilde{\nu}{}^{\,0}\,i
\ ({\rm Re}\,\nu^0={\stackrel{*}{\nu}}{}^{\,0},\ {\rm Im}\,\nu^0\!=\!\widetilde{\nu}{}^{\,0})$
and to the {\rm 1}-st order generalized eigen\-vec\-tor 
\vspace{0.25ex}
$\nu^{1}\!={\stackrel{*}{\nu}}{}^{1}\!+\widetilde{\nu}{}^{\,1}\,i
\ ({\rm Re}\,\nu^1\!={\stackrel{*}{\nu}}{}^{1}\!,\, {\rm Im}\,\nu^1\!=\widetilde{\nu}{}^{\,1}\!).\!\!$
Then the ordinary differential sys\-tem {\rm (1.1)} on a domain
${\mathscr X}\subset \bigl\{x\colon {\stackrel{*}{\nu}}{}^{0}x\ne 0\bigr\}$  
has the autonomous first integrals
\\[1.5ex]
\mbox{}\hfill
$
F_1^{}\colon x\to\
\Bigr(\bigl({\stackrel{*}{\nu}}{}^{\,0}x\bigr)^2+
\bigl(\widetilde{\nu}{}^{\,0}x\bigr)^2\,\Bigr)
\exp\biggl( {}-2\,\dfrac{{\stackrel{*}{\lambda}}\,\alpha(x)-
\widetilde{\lambda}\,\beta(x)}
{\bigl({\stackrel{*}{\nu}}{}^{\,0}x\bigr)^2+
\bigl(\widetilde{\nu}{}^{\,0}x\bigr)^2}\,\biggr)
$
\ for all 
$
x\in {\mathscr X}
\hfill
$
\\[1.25ex]
and
\\[1ex]
\mbox{}\hfill
$
F_2^{}\colon x\to\
\arctan\dfrac{\widetilde{\nu}{}^{\,0}x}{{\stackrel{*}{\nu}}{}^{0}x}
\ - \
\dfrac{\widetilde{\lambda}\,\alpha(x)+
{\stackrel{*}{\lambda}}\,\beta(x)}
{\bigl({\stackrel{*}{\nu}}{}^{\,0}x\bigr)^2+
\bigl(\,\widetilde{\nu}{}^{\,0}x\bigr)^2}
$
\ for all 
$
x\in {\mathscr X},
\hfill
$
\\[1.5ex]
where the polynomials
\vspace{0.75ex}
$
\alpha\colon x\to
{\stackrel{*}{\nu}}{}^{0}x\,{\stackrel{*}{\nu}}{}^{1}x +
\widetilde{\nu}\,{}^{\,0}x\,\widetilde{\nu}{}^{\,1}x,
\ 
\beta\colon x\to
{\stackrel{*}{\nu}}{}^{0}x\,\widetilde{\nu}{}^{\,1}x -
\widetilde{\nu}{}^{\,0}x\,{\stackrel{*}{\nu}}{}^{1} x
$
for all 
$
x\in\R^n.
$
}

The proof of Theorem 1.6 is similar to that one in Theorem 1.5. 
\vspace{0.5ex}

{\bf Theorem 1.6.}
{\it
Suppose $\nu^{0}$ and $\nu^{1}$ are a real eigenvector and 
a real {\rm 1}-th order generalized eigenvector of the matrix $B$ corresponding to the
eigenvalue $\lambda_1^{}\!=0$ with elementary divisor of multiplicity $m\geq 2,$
and $\nu^{2}$ is a real eigenvector of the matrix $B$ corresponding to the 
eigenvalue $\lambda_2.$ Then the system {\rm (1.1)} on a domain ${\mathscr X}$ 
has the autonomous first integral
\\[1.5ex]
\mbox{}\hfill
$
F\colon x\to\
\nu^{2}x\;\!
\exp\biggl({}-\lambda_{2}^{}\ \dfrac{\nu^{1}x}{\nu^{0}x}\,\biggr)
$
\ for all 
$
x\in {\mathscr X},
\quad 
{\mathscr X}\subset \{x\colon \nu^0x\ne 0\}.
\hfill
$
\\[1.75ex]
}
\indent
If $\lambda_2^{}=0,$ then from Theorem 1.6, we have Corollary 1.1.
\vspace{0.5ex}

{\bf Example 1.4.}
%\marginpar{\No\ 02\,447}
The autonomous system of ordinary linear differential equations 
\\[2ex]
\mbox{}\hfill                              % (1.15)
$
\dfrac{dx_1^{}}{dt} =4x_1^{} - 5x_2^{}+ 2x_3^{},
\quad \ 
\dfrac{dx_2^{}}{dt} = 5x_1^{} - 7x_2^{}+ 3x_3^{},
\quad \
\dfrac{dx_3^{}}{dt} = 6x_1^{} - 9x_2^{} + 4x_3^{}
$
\hfill (1.15)
\\[2ex]
has the eigenvalue $\lambda_1^{}=0$ with elementary divisor $\lambda^2$ of multiplicity $2$ 
corresponding to the eigenvector $\nu^{1}=(1,{}-2,1)$ and 
to the {\rm 1}-st order generalized eigenvector  $\nu^{2}=(0,{}-1,1),$
and the simple eigenvalue $\lambda_{3}^{}=1$ with elementary divisor $\lambda-1$ 
corresponding to the eigenvector $\nu^{3}=(3,{}-3,1).$ 
The scalar functions (by Theorem 1.5)
\\[1.5ex]
\mbox{}\hfill                  
$
F_1^{}\colon x\to\ 
 x_1^{}-2x_2^{}+x_3^{}
$
\ for all 
$
x\in {\mathbb R}^3
\hfill
$
\\[1ex]
and (by Theorem 1.6)
\\[1ex]
\mbox{}\hfill                   
$
F_{2}^{}\colon x\to\
(3x_1^{}-3x_2^{}+x_3^{})\exp\dfrac{x_2^{}-x_3^{}}{x_1^{}-2x_2^{}+x_3^{}}
$ 
\ for all 
$
x\in {\mathscr X}
\hfill 
$
\\[2ex]
are an autonomous integral basis of system (1.15) on a 
\vspace{0.5ex}
domain ${\mathscr X}\subset \{x\colon x_1^{}-2x_2^{}+x_3^{}\ne 0\}.$

From Theorem 1.6, we obtain
\vspace{0.5ex}

{\bf Corollary 1.4.}
\vspace{0.25ex}
{\it
Suppose $\nu^{0}$ and $\nu^{1}$ are a real eigenvector and 
a real {\rm 1}-th order generalized eigenvector of the matrix $B$ corresponding to the
\vspace{0.25ex}
eigenvalue $\lambda_1^{}\!=0$ with elementary divisor of multiplicity $m\geq 2,$
and 
$\nu^{2}={\stackrel{*}{\nu}}{}^{2}+\widetilde{\nu}{}^{\,2}\,i
\ ({\rm Re}\;\!\nu^2={\stackrel{*}{\nu}}{}^{\,2},\, {\rm Im}\;\!\nu^2=\widetilde{\nu}{}^{\,2})$
is an eigenvector of the matrix $B$ corresponding to the complex eigenvalue 
\vspace{0.5ex}
$
\lambda_{2}^{}\!={\stackrel{*}{\lambda}}_2^{}+\widetilde{\lambda}_2^{}\,i
\ ({\rm Re}\,\lambda_2^{}\!={\stackrel{*}{\lambda}}_2^{},\, 
{\rm Im}\,\lambda_2^{}\!=\widetilde{\lambda}_2^{}\!\ne 0).$
Then the system {\rm (1.1)} on a domain ${\mathscr X}$ 
has the autonomous first integrals
\\[1.75ex]
\mbox{}\hfill
$
F_1^{}\colon x\to \
\Bigl( \bigl(\,{\stackrel{*}{\nu}}{}^{2}x\bigr)^{2} +
\bigl(\,\widetilde{\nu}{}^{\,2}x\bigr)^2\,\Bigr)
\exp\biggl({}-2\,
{\stackrel{*}{\lambda}}_{2}^{} \ \dfrac{\nu^{1}x}{\nu^{0}x}\,\biggr)
$
\ for all 
$
x\in {\mathscr X}
\hfill
$
\\[1.35ex]
and
\\[1ex]
\mbox{}\hfill
$
F_2^{}\colon x\to \
\arctan\dfrac{\widetilde{\nu}{}^{\,2}x}{{\stackrel{*}{\nu}}{}^2x}
\ - \
\widetilde{\lambda}_2^{} \ \dfrac{\nu^{1}x}{\nu^{0}x}
$
\ for all 
$
x\in {\mathscr X},
\quad 
{\mathscr X}\subset \{x\colon \nu^0x\ne 0,\ {\stackrel{*}{\nu}}{}^2x\ne 0\}.
\hfill
$
\\[2.25ex]
}
\indent
{\bf Theorem 1.7.}\!
\vspace{0.25ex}
{\it 
Suppose $\nu^{01}$ and $\nu^{11}$ are a real eigenvector and 
a real {\rm 1}-th order generalized eigenvector of the matrix $B$ corresponding to the
\vspace{0.25ex}
eigenvalue $\lambda_1^{}$ with elementary divisor of multiplicity $m_1^{}\geq 2,$
\vspace{0.25ex}
and 
$\nu^{02}$ and $\nu^{12}$ are a real eigenvector and 
a real {\rm 1}-th order generalized eigenvector of the matrix $B$ corresponding to the
\vspace{0.25ex}
eigenvalue $\lambda_2^{}$ with elementary divisor of multiplicity $m_2^{}\geq 2.$
Then the system {\rm (1.1)} has the autonomous first integral
\\[1.5ex]
\mbox{}\hfill
$
F\colon x\to  \
\dfrac{\nu^{11}x}{\nu^{01}x}\ -\ \dfrac{\nu^{12}x}{\nu^{02}x}
$
\ for all 
$
x\in {\mathscr X},
\quad 
{\mathscr X}\subset \{x\colon \nu^{01}x\ne 0,\ \nu^{02}x\ne 0\}.
\hfill
$
\\[2ex]
}
\indent
{\sl Indeed}, 
the Lie derivative of the function $F$ on a domain ${\mathscr X}$ by virtue of system (1.1) is
\\[1.75ex]
\mbox{}\hfill
$
{\frak A}\,F(x)=
\dfrac{(\lambda_1^{}\nu^{11}x+\nu^{01}x)\nu^{01}x-\lambda_1^{}\,\nu^{01}x\,\nu^{11}x}
{(\nu^{01}x)^2}  -
\dfrac{(\lambda_2^{}\nu^{12}x+\nu^{02}x)\nu^{02}x-\lambda_2^{}\,\nu^{02}x\,\nu^{12}x}
{(\nu^{02}x)^2} =0.
\ \k
\hfill
$
\\[2.25ex]
\indent
{\bf Theorem 1.8.}
{\it
Suppose $\lambda$ is the eigenvalue with elementary divisor of multiplicity $m\geq 2$ 
of the matrix $B$ corresponding to an eigenvector $\nu^{0}$ and to 
generalized eigenvectors $\nu^{k},\ k=1,\ldots, m-1.$ 
Then the system {\rm (1.1)} on a domain ${\mathscr X}\subset \{x\colon \nu^0x\ne 0\}$ 
has the fun\-c\-ti\-o\-na\-l\-ly independent autonomous first integrals
\\[2ex]
\mbox{}\hfill                                       % (1.16)
$
F_{\zeta}^{}\colon x\to\, \Psi_{\zeta}^{}(x)
$
\ for all 
$
x\in {\mathscr X},
\ \ \zeta=2,\ldots, m-1,
$
\hfill {\rm(1.16)}
\\[2.25ex]
where the functions $\Psi_{\zeta}^{}\colon {\mathscr X}\to {\mathbb R},\ \zeta=2,\ldots, m-1,$ 
are the solution to the system
\\[1.5ex]
\mbox{}\hfill   % (1.17)
$
\displaystyle
\nu^{k}x =
\sum\limits_{\tau=1}^{k}{\textstyle\binom{k-1}{\tau-1}}
\Psi_{\tau}^{}(x)\,\nu^{k-\tau}x
$
\ for all 
$
x\in {\mathscr X},
\ \ k=1,\ldots, m-1.
$
\hfill {\rm (1.17)}
\\[2ex]
}
{\sl Proof}. 
Using the identities (1.3) and the equalities (1.14), we obtain 
\\[1.75ex]
\mbox{}\hfill                                     % (1.18)
$
{\frak A}\,\nu^0x = \lambda\,\nu^0x
$
\ for all 
$
x\in \R^n,
\hfill
$
\\[-0.35ex]
\mbox{}\hfill (1.18)
\\[-0.25ex]
\mbox{}\hfill
$
{\frak A}\,\nu^{k}x =
\lambda\,\nu^{k}x+k\,\nu^{k-1}x
$
\ for all 
$
x\in \R^n,
\ \ k=1,\ldots, m-1.
\hfill
$
\\[1.75ex]
\indent
The system (1.17) has the determinant $\,(\nu^0x)^{m-1}$ such that 
\vspace{0.25ex}
$\,(\nu^0x)^{m-1}\ne 0$ for all $x$ from a domain ${\mathscr X}\subset \{x\colon \nu^0x\ne 0\}.$ 
\vspace{0.25ex}
Therefore there exists the solution $\Psi_{\tau}^{},\ \tau=1,\ldots, m-1,$ 
on the domain ${\mathscr X}$ of the functional system (1.17).
Let us show that
\\[2.5ex]
\mbox{}\hfill                                        % (1.19)
$
{\frak A}\Psi_{k}^{}(x)=
\left[\!\!
\begin{array}{lll}
1\! & \text{for all}\ \, x\in {\mathscr X}, &  k=1,
\\[1.25ex]
0\! & \text{for all}\ \, x\in {\mathscr X}, & k=2,\ldots, m-1,
\end{array}
\right.
$
\hfill (1.19)
\\[2.5ex]
\indent
The proof of identities (1.19) is by induction on $m.$

For $m=2$ and $m=3,$ the assertions (1.19) follows from the identities (1.18).
\vspace{0.35ex}

Assume that the identities (1.19) for $m=\varepsilon$ is true. 
Using the system of identities (1.18) and the system (1.17) for 
$m=\varepsilon+1,\ m=\varepsilon,$ we get
\\[2ex]
\mbox{}\hfill
$
{\frak A}\,\nu^{\varepsilon}x =
\lambda\ {\displaystyle \sum\limits_{\tau=1}^{\varepsilon}}
{\textstyle\binom{\varepsilon-1}{\tau-1}}\Psi_{\tau}^{}(x)\, \nu^{\varepsilon-\tau}x +
(\varepsilon-1) {\displaystyle \sum\limits_{\tau=1}^{\varepsilon-1}}
{\textstyle\binom{\varepsilon-2}{\tau-1}}\Psi_{\tau}^{}(x)\, \nu^{\varepsilon-\tau-1}x +
\nu^{\varepsilon-1}x + 
\nu^{0}x\, {\frak A}\Psi_{\varepsilon}^{}(x)
\hfill
$
\\[2.5ex]
\mbox{}\hfill
for all 
$
x\in {\mathscr X},
\quad
{\mathscr X}\subset \{x\colon \nu^0x\ne 0\}.
\hfill
$
\\[2ex]
\indent
Now taking into account the system (1.17) with $k=\varepsilon-1$ and  $k=\varepsilon,$
the identity (1.18) with $k=\varepsilon,$ and $\nu^{0}x\ne 0$ for all $x\in {\mathscr X},$ we have
\\[1.75ex]
\mbox{}\hfill
$
{\frak A}\Psi_{\varepsilon}^{}(x)=0
$
\ for all 
$
x\in {\mathscr X}.
\hfill
$
\\[1.5ex]
\indent
This implies that the identities (1.19) for $m=\varepsilon+1$ are true.
So by the principle of mathematical induction, the statement (1.19) is true for every 
natural number $m\geq 2.$ 

Thus the functions (1.16) 
\vspace{1ex}
are functionally independent first integrals of system (1.1). $\k$

Theorem 1.8 is true both for the case of the real eigenvalue $\lambda$ and 
for the case of the com\-p\-lex eigenvalue $\lambda\ ({\rm Im}\,\lambda\ne 0).$
In the complex case, from the complex-valued first integrals (1.16) of 
system (1.1), we obtain the real-valued first integrals of system (1.1)
\\[1.75ex]
\mbox{}\hfill
$
F_{\zeta}^{1}\colon x\to \mbox{Re}\,\Psi_\zeta^{}(x),
\quad 
F_{\zeta}^{2}\colon x\to \mbox{Im}\,\Psi_\zeta^{}(x)
$
\ for all 
$
x\in {\mathscr X}, 
\ \ \zeta=2,\ldots, m-1,
\hfill
$
\\[1.5ex]
where ${\mathscr X}$ is a domain from the set
$\bigl\{x\colon \bigl(\,{\stackrel{*}{\nu}}{}^0x\bigr)^2 +
\bigl(\,\widetilde{\nu}{}^{\,0}x\bigr)^2\ne 0\bigr\}\subset {\mathbb R}^n.$
\vspace{1.5ex}

{\bf Example 1.5.}
%\marginpar{\No\ 03\,173}
Consider the linear autonomous differential system
\\[2.25ex]
\mbox{}\hfill                                   % (1.20)
$
\dfrac{dx_1^{}}{dt} =   4x_1^{} - x_2^{},
\quad \
\dfrac{dx_2^{}}{dt} =  3x_1^{} + x_2^{} - x_3^{},
\quad \
\dfrac{dx_3^{}}{dt} =    x_1^{} + x_3^{}.
$
\hfill (1.20)
\\[2.75ex]
\indent
The matrix 
\vspace{1.5ex}
$
B=\left\|\!\!
\begin{array}{rrc}
4 & 3 & 1
\\[0.35ex]
{}-1 & 1 & 0
\\[0.35ex]
0 & {}-1 & 1
\end{array}
\!\right\|
$
has one triple eigenvalue $\lambda_{1}^{}=2.$ 

The rank of the matrix $B-\lambda_1^{} E$ is equal 2.
\vspace{0.5ex}
Therefore the eigenvalue $\lambda_{1}^{}=2$ has 
$\kappa_1^{}=3-2=1$ elementary divisor $(\lambda-2)^3$ of multiplicity 3.
\vspace{1ex}

The system
\vspace{1ex}
$
(B-\lambda_{1}^{}E)\,\mbox{colon}(\nu_{1}^{},\nu_2^{},\nu_{3}^{}) = 0
\Leftrightarrow
\left\{\!\!
\begin{array}{r}
2\nu_{1}^{}+3\nu_{2}^{}+\nu_{3}^{}=0,
\\[0.5ex]
{}-\nu_{1}^{}-\nu_{2}^{}=0,
\\[0.5ex]
{}-\nu_{2}^{}-\nu_{3}^{}=0\,
\end{array} 
\right.
\Leftrightarrow
\left\{\!\!
\begin{array}{l}
\nu_{1}^{}=\nu_3^{},
\\[0.75ex]
\nu_{2}^{}=-\;\!\nu_{3}^{}.
\end{array} 
\right.
$
Hence $\nu^{0}=(1,{}-1,1)$ is an eigenvector of $B$ corresponding to 
the eigenvalue $\lambda_{1}^{}=2.$

The 1-st order generalized eigenvector $\nu^{1}$ of the matrix $B$ corresponding to 
the eigenvalue $\lambda_1=2^{}$ is a solution of the system
\\[1.25ex]
\mbox{}\hfill
$
(B-\lambda_{1}^{}E)\,\mbox{colon}(\nu_{1}^{},\nu_2^{},\nu_{3}^{}) =
\mbox{colon}(1,-1,1)
\Leftrightarrow
\left\{\!\!
\begin{array}{r}
2\nu_{1}^{}+3\nu_{2}^{}+\nu_{3}^{}=1,
\\[0.5ex]
{}-\nu_{1}^{}-\nu_{2}^{}={}-1,
\\[0.5ex]
{}-\nu_{2}^{}-\nu_{3}^{}=1\,
\end{array}
\right.
\Leftrightarrow
\left\{\!\!
\begin{array}{l}
\nu_{1}^{}=-\,\nu_{2}^{}+1,
\\[0.5ex]
\nu_{3}^{}=-\,\nu_{2}^{}-1.
\end{array}
\right.
\hfill
$
\\[1.5ex]
Hence $\nu^{1}=(1,0,{}-1)$ is a generalized eigenvector of the 1-st order
of the matrix $B$ corresponding to the eigenvalue $\lambda_1^{}=2.$

The 2-nd order generalized eigenvector $\nu^{2}$ of the matrix $B$ corresponding to 
the eigenvalue $\lambda_1^{}=2$ is a solution of the system
\\[1ex]
\mbox{}\hfill
$
(B-\lambda_{1}^{}E)\,\mbox{colon}(\nu_{1}^{},\nu_2^{},\nu_{3}^{}) =
2\,\mbox{colon}(1,0,-1)
\Leftrightarrow
\left\{\!\!
\begin{array}{r}
2\nu_{1}^{}+3\nu_{2}^{}+\nu_{3}^{}=2,
\\[0.5ex]
{}-\nu_{1}^{}-\nu_{2}^{}=0,
\\[0.5ex]
{}-\nu_{2}^{}-\nu_{3}^{}={}-2\,
\end{array} 
\right.
\!\!
\Leftrightarrow
\left\{\!\!
\begin{array}{l}
\nu_{1}^{}=-\,\nu_{2}^{},
\\[0.5ex]
\nu_{3}^{}=-\,\nu_2^{}+2.
\end{array} 
\right.
\hfill
$
\\[1.25ex]
Hence $\nu^{2}=(0,0,2)$ is a generalized eigenvector of the 2-nd order
of the matrix $B$ corresponding to the eigenvalue $\lambda_1^{}=2.$
\vspace{0.35ex}

The functionally independent scalar functions (by Theorem 1.5)
\\[1.5ex]
\mbox{}\hfill                            
$
F_{1}^{}\colon (x_1^{},x_2^{},x_3^{})\to\,
(x_1^{}-x_2^{}+x_3^{})\exp\Bigl({}-2\, \dfrac{x_1^{}-x_3^{}}{x_1^{}-x_2^{}+x_3^{}}\,\Bigr)
$
\ for all 
$
(x_1^{},x_2^{},x_3^{})\in {\mathscr X},
\hfill 
$
\\[1.75ex]
\mbox{}\hfill                             
$
F_{2}^{}\colon  (x_1^{},x_2^{},x_3^{})\to
\dfrac{(x_1^{}-x_3^{})^2-2x_3^{}(x_1^{}-x_2^{}+x_3^{})}{(x_1^{}-x_2^{}+x_3^{})^2}
$
for all 
$
(x_1^{},x_2^{},x_3^{})\in {\mathscr X}
$
(by Theorem 1.8)\hfill\mbox{}
\\[1.5ex]
are a basis of autonomous first integrals on a domain 
${\mathscr X}\subset \{(x_1^{},x_2^{},x_3^{})\colon x_1^{}-x_2^{}+x_3^{}\ne 0\}$ of
the ordinary differential system (1.20).
\vspace{0.5ex}

{\bf Example 1.6.}
%\marginpar{\No\ 00\,000}
The sixth-order ordinary autonomous linear differential system 
\\[1.5ex]
\mbox{}\hfill                                   % (1.21)
$
\begin{array}{ll}
\dfrac{dx_1^{}}{dt} = x_1^{} - 2x_2^{}+x_3^{}-2x_6^{},
\quad &
\dfrac{dx_2^{}}{dt} = 3x_2^{}-x_3^{}-x_5^{}+2x_6^{},
\\[3ex]
\dfrac{dx_3^{}}{dt} = {}-x_1^{} +x_3^{}+2x_4^{}+2x_5^{},
\quad &
\dfrac{dx_4^{}}{dt} = {}-x_1^{}+x_4^{}+x_5^{}+x_6^{},
\\[3ex]
\dfrac{dx_5^{}}{dt} = x_1^{} +x_2^{}+x_5^{},
\quad &
\dfrac{dx_6^{}}{dt} = x_1^{}-x_2^{}+x_3^{}-x_4^{}-x_6^{}
\end{array}
$
\hfill (1.21)
\\[2ex]
has the triple complex eigenvalue $\lambda_{1}^{}=1+i$ 
\vspace{0.25ex}
with the elementary divisor $(\lambda-1-i)^3$ 
of mul\-ti\-p\-li\-ci\-ty $3$ corresponding to the eigenvector $\nu^{0}=(1,1,0,0,i,0)$ and to the 
\vspace{0.35ex}
generalized eigenvectors $\nu^{1}=(0,1,0,i,i,1), \ \nu^{2}=(0,1,i,0,i,0).$
The functions
\\[2ex]
\mbox{}\hfill                              
$
F_{1}^{}\colon x\to\ P(x)\,\exp({}-2\varphi (x))
$
\ for all 
$
x\in {\mathscr X}
$
\ \ \ (by Theorem 1.2),
\hfill\mbox{} 
\\[2.25ex]
\mbox{}\hfill                             
$
F_{2}^{}\colon x\to\ 
P(x)\exp\biggl({}-2\,\dfrac{\alpha(x)-\beta(x)}{P(x)}\,\biggr)
$
\ for all 
$
x\in {\mathscr X}
$
\ \ \ (by Corollary 1.3),
\hfill\mbox{}
\\[2ex]
\mbox{}\hfill                              
$
F_{3}^{}\colon x\to \
\varphi(x)- \dfrac{\alpha(x)+\beta(x)}{P(x)}
$
\ for all 
$
x\in {\mathscr X}
$
\ \ \ (by Corollary 1.3),
\hfill\mbox{} 
%(1.22)
\\[2ex]
\mbox{}\hfill                             
$
F_{4}^{}\colon x\to\
\dfrac{\gamma(x)P(x)+\beta^2(x)-\alpha^2(x)}{P^2(x)}
$
\ for all 
$
x\in {\mathscr X}
$
\ \ \ (by Theorem 1.8),
\hfill\mbox{}
\\[2ex]
\mbox{}\hfill                              
$
F_{5}^{}\colon x\to\
\dfrac{\delta(x)P(x)-2\alpha(x)\beta(x)}{P^2(x)}
$
\ for all 
$
x\in {\mathscr X}
$
\ \ \ (by Theorem 1.8),
\hfill\mbox{} 
\\[1.5ex]
where the polynomials
\\[1.5ex]
\mbox{}\hfill
$
P\colon x\to (x_1^{}+x_2^{})^2 + x_5^2, 
\quad
\alpha\colon x\to (x_1^{}+x_2^{})(x_2^{}+x_6^{}) + x_5^{}(x_4^{}+x_5^{}),
\hfill
$
\\[1.5ex]
\mbox{}\hfill
$
\beta\colon x\to (x_1^{}+x_2^{})(x_4^{}+x_5^{}) - x_5^{}(x_2^{}+x_6^{}), 
\ \ 
\gamma\colon x\to x_2^{}(x_1^{}+x_2^{}) + x_5^{}(x_3^{}+x_5^{}),
\hfill
$
\\[1.5ex]
\mbox{}\hfill
$
\delta\colon x\to (x_1^{}+x_2^{})(x_3^{}+x_5^{}) - x_2^{}x_5^{}
$
\ for all 
$
x\in {\mathbb R}^6,
\hfill
$
\\[1.5ex]
the scalar function
\vspace{0.75ex}
$
\varphi\colon x\to \arctan\dfrac{x_5^{}}{x_1^{}+x_2^{}}
$
for all 
$
x\in {\mathscr X},
$
are first integrals on a domain ${\mathscr X}$ from the set 
$\{x\colon x_1^{}+x_2^{}\ne 0\}\subset {\mathbb R}^6$ 
of the ordinary differential system (1.21). 
\vspace{0.5ex}

The functionally independent first integrals $F_1^{},\ldots, F_5^{}$ are an autonomous integral basis 
on a domain ${\mathscr X}$ of the ordinary differential system (1.21).
\\[2.25ex]
\indent
{\bf 1.1.3. Nonautonomous first integrals}
\\[0.75ex]
\indent
The ordinary differential system (1.1) 
\vspace{0.25ex}
is induced the nonautonomous linear differential operator of first order
$
{\frak B}(t,x)=\partial_{t}^{}+{\frak A}(x)
$
for all 
$
(t,x)\in\R^{n+1}.
$
\vspace{0.25ex}

Adding one nonautonomous first integral of system (1.1) 
to an autonomous integral basis of system (1.1),
we can construct a basis of first integrals for system (1.1).  
\vspace{0.25ex}

Such procedure every time can be carried out on the base of the following statements.
\vspace{0.5ex}

{\bf Theorem 1.9.}
\vspace{0.25ex}
{\it
Suppose $\nu$ is a real eigenvector of  the matrix $B$ 
corresponding to the eigenvalue $\lambda.$  
Then the ordinary differential system {\rm (1.1)} has the first integral}
\\[1.5ex]
\mbox{}\hfill
$
F\colon (t,x)\to\ \nu x\,\exp({}-\lambda\,t)
$
\ for all 
$
(t,x)\in {\mathbb R}^{n+1}.
\hfill
$
\\[1.5ex]
\indent
{\sl Indeed},  
the Lie derivative of the function $F$ by virtue of system (1.1) is
\\[1.75ex]
\mbox{}\hfill
$
{\frak B}F(t,x)=\partial_{t}F(t,x)+{\frak A}F(t,x)=
{}-\lambda F(t,x)+\lambda F(t,x)=0
$
\ for all 
$
(t,x)\in\R^{n+1}.
\ \k
\hfill
$
\\[2ex]
\indent
{\bf Example 1.7}
%\marginpar{\No\ 02\,442}
(continuation of Example 1.1). 
\vspace{0.35ex}
By Theorem 1.9,
using the eigenvalue $\lambda_2^{}=1$ corresponding to the eigenvector $\nu^2=(2,2,1,1),$ 
\vspace{0.5ex}
we can build the first integral
$
F_2^{}\colon (t,x)\to (2x_1+2x_2+x_3+x_4)\,e^{{}-t}
$
for all $(t,x)\in\R^5
$
of system (1.9).
\vspace{0.5ex}

The functions $F_1^{},\, F_{23}^{},\ F_{24}^{}$ (constructed in Example 1.1), 
and the function $F_2^{}$ 
\vspace{0.25ex}
are an integral basis of the autonomous differential system (1.9) 
on a domain ${\mathbb R}\times {\mathscr X},$
\vspace{0.25ex}
where ${\mathscr X}$ is a domain from the set 
$\{x\colon x_1^{}+x_3^{}\ne 0,\  2x_2^{}+x_4^{}\ne 0\}\subset {\mathbb R}^4.$
\vspace{0.75ex}

{\bf Example 1.8}
%\marginpar{\No\ 02\,443}
(continuation of Example 1.2).
The autonomous differential system (1.12) has the eigenvector $\nu^2=(3,{}-1,{}-1)$
corresponding to the eigenvalue $\lambda_2^{}=2.$ 
\vspace{0.35ex}

By Theorem 1.9, the differential system (1.12) has the first integral
\\[1.5ex]
\mbox{}\hfill
$
F_3^{}\colon (t,x_1^{},x_2^{},x_3^{})\to\ (3x_1^{}-x_2^{}-x_3^{})\,e^{{}-2t}
$
\ for all 
$
(t,x_1^{},x_2^{},x_3^{})\in {\mathbb R}^4.
\hfill
$
\\[1.5ex]
\indent
The functionally independent first integrals $F_1^{},\ F_2^{}$ (constructed in Example 1.1), 
\vspace{0.25ex}
and $F_3^{}$ are a basis of first integrals for the differential system (1.12) on a domain 
\vspace{0.35ex}
${\mathbb R}\times {\mathscr X},$ 
where ${\mathscr X}$ is a domain from the set
$\{(x_1^{},x_2^{},x_3^{})\colon x_1^{}-x_3^{}\ne 0\}.$
\vspace{0.75ex}

{\bf Corollary 1.5.}
\vspace{0.35ex}
{\it
Let $\nu={\stackrel{*}{\nu}}+\widetilde{\nu}\,i\
({\rm Re}\,\nu\!=\!{\stackrel{*}{\nu}},\, {\rm Im}\,\nu\!=\!\widetilde{\nu}\,)\!\!$ be
an eigenvector of the matrix $B$ corresponding to the complex eigenvalue
\vspace{0.35ex}
$\lambda ={\stackrel{*}{\lambda}}+ \widetilde{\lambda}\;\!i\
({\rm Re}\,\lambda={\stackrel{*}{\lambda}},\ {\rm Im}\,\lambda=\widetilde{\lambda}\ne 0).$
Then the or\-di\-na\-ry differential system {\rm (1.1)} has the first integrals
\\[1.75ex]
\mbox{}\hfill
$
F_1^{}\colon (t,x)\to\
\bigl( ({\stackrel{*}{\nu}}x)^2
+ (\widetilde{\nu}x)^2\bigr)
\exp\bigl({}-2{\stackrel{*}{\lambda}}\,t\bigr)
$
\ for all 
$
(t,x)\in {\mathbb R}^{n+1}
\hfill
$
\\[1.25ex]
and
\\[1ex]
\mbox{}\hfill
$
F_2^{}\colon (t,x)\to\
\arctan\dfrac{\widetilde{\nu}x}{{\stackrel{*}{\nu}}x}\ -\
\widetilde{\lambda}\,t
$
\ for all 
$
(t,x)\in {\mathbb R}\times {\mathscr X},
\hfill
$
\\[1.5ex]
where ${\mathscr X}$ is a domain from the set
\vspace{0.75ex}
$\bigl\{x\colon {\stackrel{*}{\nu}}x\ne 0\bigr\}\subset {\mathbb R}^n.$
}

{\sl Proof}. Using Properties 1.4 and 1.5, we get
\\[1.5ex]
\mbox{}\hfill
$
{\frak B}\;\!F_1^{}(t,x)=
\exp({}-2\;\! {\stackrel{*}{\lambda}}\,t)\, 
{\frak B}\;\!\bigl( ({\stackrel{*}{\nu}}x)^2+(\widetilde{\nu}x)^2\bigr) +\, 
\bigl( ({\stackrel{*}{\nu}}x)^2+(\widetilde{\nu}x)^2\bigr)
\,{\frak B}\exp({}-2\;\! {\stackrel{*}{\lambda}}\,t)=
\hfill
$
\\[2ex]
\mbox{}\hfill
$
=2\;\! {\stackrel{*}{\lambda}}\;\!
\bigl(\bigl({\stackrel{*}{\nu}}x\bigr)^2
+ \bigl(\widetilde{\nu}x\bigr)^2\bigr)
\exp\bigl({}-2{\stackrel{*}{\lambda}}\,t\bigr)
+ 
\bigl( ({\stackrel{*}{\nu}}x)^2+(\widetilde{\nu}x)^2\bigr)
\partial_{t}^{}\exp({}-2\;\! {\stackrel{*}{\lambda}}\,t)=0
$
for all 
$\!
(t,x)\!\in\! {\mathbb R}^{n+1},
\hfill
$
\\[2.25ex]
\mbox{}\hfill                                            
$
\displaystyle
{\frak B}\;\!F_2^{}(t,x)= 
{\frak B}\;\!\arctan\dfrac{\widetilde{\nu}x}{{\stackrel{*}{\nu}}x}\, - \,
{\frak B}\;\!\bigl(\;\!\widetilde{\lambda}\,t\bigr) \, =\, 
\widetilde{\lambda}\, -\,
\partial_{t}^{}\bigl(\;\!\widetilde{\lambda}\,t\bigr)=0
$
\ for all 
$
(t,x)\in {\mathbb R}\times {\mathscr X}.
\hfill
$
\\[2ex]
\indent
Therefore the scalar functions $F_1^{}\colon {\mathbb R}^{n+1}\to {\mathbb R}$ and 
\vspace{0.35ex}
$F_2^{}\colon {\mathbb R}\times {\mathscr X}\to {\mathbb R}$ 
are first integrals of the ordinary linear autonomous differential system (1.1). \k
\vspace{1ex}

\newpage

{\bf Example 1.9}
%\marginpar{\No\ 02\,445}
(continuation of Example 1.3).
\vspace{0.25ex}
Using the eigenvalue $\lambda_1^{}=i$ corresponding to the eigenvector 
$\nu^1=(1-i,{}-1+2i,2i,2),$
\vspace{0.5ex}
we can construct the first integral
\\[2ex]
\mbox{}\hfill
$
F_4^{}\colon (t,x)\to\
\arctan\dfrac{{}-x_1^{}+2x_2^{}+2x_3^{}}{x_1^{}-x_2^{}+2x_4^{}} \  -\ t
$
\ for all 
$
(t,x)\in {\mathbb R}\times {\mathscr X}_1^{}
$
\ \ (by Corollary 1.5)
\hfill\mbox{}
\\[2ex]
of the differential system (1.13) on a domain ${\mathscr X}_1^{}$ from the set 
$\{x\colon x_1^{}-x_2^{}+2x_4^{}\ne 0\}.$
\vspace{0.35ex}

The functionally independent first integrals
\vspace{0.25ex}
$F_1^{},\, F_{2}^{},\ F_{3}^{}$ (constructed in Example 1.3), 
and the function $F_4^{}$ 
are an integral basis of the differential system (1.13) 
\vspace{0.25ex}
on a domain ${\mathbb R}\times {\mathscr X},$
where ${\mathscr X}$ is a domain from the set 
$\{x\colon x_1^{}-x_2^{}+2x_4^{}\ne 0,\ x_4^{}-x_2^{}\ne 0\}.$
\vspace{1ex}

{\bf Theorem 1.10.}
{\it
Let $\lambda$ be the eigenvalue with elementary divisor of multiplicity $m\geq 2$ 
of the matrix $B$ corresponding to a real eigenvector $\nu^0$ and to 
a real generalized eigenvector $\nu^1$ of the {\rm 1}-st order.
Then the system {\rm (1.1)} has the first integral
\\[1.5ex]
\mbox{}\hfill                                % (1.22)
$
F\colon (t,x)\to\ \dfrac{\nu^1x}{\nu^0x} - t
$ 
\ for all 
$
(t,x)\in {\mathbb R}\times {\mathscr X},
\quad 
{\mathscr X}\subset \{x\colon \nu^0x\ne 0\}.
$
\hfill {\rm (1.22)}
\\[2ex]
}
\indent
{\sl Indeed},  
the Lie derivative of $F$ on a domain ${\mathbb R}\times {\mathscr X}$ 
by virtue of system (1.1) is
\\[1.5ex]
\mbox{}\hfill
$
{\frak B}F(t,x)=\partial_{t}F(t,x)+{\frak A}F(t,x)=
{}-1+\dfrac{(\lambda \nu^1x+\nu^0x)\nu^0x-\nu^1x\;\! \lambda\nu^0x}{(\nu^0x)^2}=0. 
\ \ \k 
\hfill
$
\\[2ex]
\indent
{\bf Example 1.10}
%\marginpar{\No\ 02\,447}
(continuation of Example 1.4).
The system (1.15) has the eigenvalue $\lambda_1^{}=0$ corresponding to the eigenvector
$\nu^1=(1,{}-2,1)$ and to the 1-st order generalized eigenvector $\nu^2=(0,{}-1,1).$
By Theorem 1.10, the scalar function
\\[1.5ex]
\mbox{}\hfill
$
F_3^{}\colon (t,x_1^{},x_2^{},x_3^{})\to\ 
\dfrac{x_3^{}-x_2^{}}{x_1^{}-2x_2^{}+x_3^{}}\ - \ t
$
\ for all 
$
(t,x_1^{},x_2^{},x_3^{})\in {\mathbb R}^4
\hfill
$
\\[1.5ex]
is a first integral of the linear differential system (1.15).

The functionally independent first integrals
\vspace{0.25ex}
$F_1^{},\, F_{2}^{}$ (constructed in Example 1.4), 
and $F_3^{}$ are an integral basis of the differential system (1.15) 
\vspace{0.25ex}
on a domain ${\mathbb R}\times {\mathscr X},$
where ${\mathscr X}$ is a domain from the set 
$\{(x_1^{},x_2^{},x_3^{})\colon x_1^{}-2x_2^{}+x_3^{}\ne 0\}.$
\vspace{0.75ex}

{\bf Example 1.11}
%\marginpar{\No\ 03\,173}
(continuation of Example 1.5).
The system (1.20) has the eigenvalue 
$\lambda_1^{}=2$ corresponding to the eigenvector
$\nu^0=(1,{}-1,1)$ and to the 1-st order generalized eigenvector $\nu^1=(1,0,{}-1).$
Using Theorem 1.10, we can build the first integral 
\\[1.5ex]
\mbox{}\hfill
$
F_3^{}\colon (t,x_1^{},x_2^{},x_3^{})\to\ 
\dfrac{x_1^{}-x_3^{}}{x_1^{}-x_2^{}+x_3^{}}\ - \ t
$
\ for all 
$
(t,x_1^{},x_2^{},x_3^{})\in {\mathbb R}^4.
\hfill
$
\\[1.75ex]
of the linear autonomous differential system (1.20)

The functionally independent first integrals
\vspace{0.25ex}
$F_1^{},\, F_{2}^{}$ (constructed in Example 1.5), 
and $F_3^{}$ are an integral basis of the differential system (1.20) 
\vspace{0.25ex}
on a domain ${\mathbb R}\times {\mathscr X},$
where ${\mathscr X}$ is a domain from the set 
$\{(x_1^{},x_2^{},x_3^{})\colon x_1^{}-x_2^{}+x_3^{}\ne 0\}.$
\vspace{1ex}

If $\lambda$ is a complex number, then the function (1.22) is a
complex-valued first integral of system (1.1).
In this case, from the complex-valued first integral (1.22) of 
system (1.1), we obtain two real-valued first integrals of system (1.1)
\\[1.75ex]
\mbox{}\hfill
$
F_1^{}\colon (t,x)\to\
\dfrac{{\stackrel{*}{\nu}}\,{}^{0}x\,{\stackrel{*}{\nu}}{}^{1}x +
\widetilde{\nu}{}^{\,0}x\,\widetilde{\nu}{}^{\,1}x}
{\bigl({\stackrel{*}{\nu}}\,{}^{0}x\bigr)^{2} +
\bigl(\widetilde{\nu}{}^{\,0}x\bigr)^{2}} \ - t
$
\ for all 
$
(t,x)\in {\mathbb R}\times{\mathscr X}
\hfill
$
\\[1ex]
and
\\[1ex]
\mbox{}\hfill
$
F_2^{}\colon (t,x)\to\
\dfrac{{\stackrel{*}{\nu}}\,{}^{0}x\,\widetilde{\nu}{}^{\,1}x
-\widetilde{\nu}{}^{\,0}x\,{\stackrel{*}{\nu}}{}^{1}x}
{\bigl({\stackrel{*}{\nu}}\,{}^{0}x\bigr)^{2} +
\bigl(\widetilde{\nu}{}^{\,0}x\bigr)^{2}}
$
\ for all 
$
(t,x)\in {\mathbb R}\times{\mathscr X},
\hfill
$
\\[1.75ex]
where ${\mathscr X}$ is a domain from the set
$\bigl\{x\colon \bigl({\stackrel{*}{\nu}}\,{}^0x\bigr)^2 +
\bigl(\widetilde{\nu}{}^{\,0}x\bigr)^2\ne 0\bigr\}.$
\vspace{1.5ex}

{\bf Example 1.12}
%\marginpar{\No\ 00\,000}
(continuation of Example 1.6).
\vspace{0.25ex}
Using the eigenvector 
$\nu^{0}=(1,1,0,0,i,0)$ and 1-st order generalized eigenvector
\vspace{0.25ex}
$\nu^{1}=(0,1,0,i,i,1)$ corresponding to the eigenvalue $\lambda_1^{}=1+i,$ 
we can construct the first integral of system (1.21)
\\[1.75ex]
\mbox{}\hfill
$
F_6^{}\colon (t,x)\to\
\dfrac{(x_1^{}+x_2^{})(x_2^{}+x_6^{})+x_5^{}(x_4^{}+x_5^{})}{(x_1^{}+x_2^{})^2+x_5^2}\ - \ t
$
\ for all 
$
(t,x)\in {\mathbb R}\times {\mathscr X}.
\hfill
$
\\[1.75ex]
\indent
The functionally independent first integrals
\vspace{0.25ex}
$F_1^{},\ldots, F_5^{}$ (constructed in Example 1.6), 
and $F_6^{}$ are a basis of first integrals for the differential system (1.21) 
\vspace{0.25ex}
on a domain ${\mathbb R}\times {\mathscr X},$
where ${\mathscr X}$ is a domain from the set 
$\{x\colon x_1^{}+x_2^{}\ne 0\}.$
\vspace{0.75ex}

{\sl Let us remark that} using Theorems 1.8, 1.9, 1.10, and Corollary 1.5, we can always 
build an integral basis of system (1.1) on the base of nonautonomous first integrals.
\vspace{0.5ex}

{\bf Example 1.13.}
The fifth-order ordinary linear autonomous differential system
\\[2ex]
\mbox{}\hfill                               % (1.23)
$
\dfrac{dx_1^{}}{dt} =   2(x_1^{} - 3x_2^{} +2x_3^{}+x_5^{}),
\qquad
\dfrac{dx_2^{}}{dt} =  2x_1^{} - x_2^{} +2x_3^{}+2x_4^{},
\qquad
\dfrac{dx_3^{}}{dt} =  x_2^{} - x_3^{},
\hfill
$
\\[0.75ex]
\mbox{}\hfill (1.23)
\\[0.5ex]
\mbox{}\hfill
$
\dfrac{dx_4^{}}{dt} =  {}-3x_1^{} + 5x_2^{} -4x_3^{}-x_4^{}-2x_5^{},
\qquad
\dfrac{dx_5^{}}{dt} =  4x_2^{} - 2x_3^{}+x_4^{}
\hfill
$
\\[2ex]
has the eigenvalues $\lambda_{1}^{}={}-1,\ \lambda_{2}^{}=1-i,$ and $\lambda_{3}^{}=1+i$
\vspace{0.25ex}
with elementary divisors $(\lambda+1)^3,$ $\lambda-1+i,$ and $\lambda-1-i,$ respectively,
\vspace{0.25ex}
corresponding to the eigenvectors
$\nu^{01}=(1, 0, 1, 1, 0),$ $\nu^{2}=(1, 1+i, 0, 1+i, 1-i),\ \nu^{3}=(1, 1-i, 0, 1-i, 1+i)$ 
\vspace{0.25ex}
and to the generalized eigenvectors 
$\nu^{11}=(1, 1/2, 1, 1, 0),\ \nu^{21}=(0, 1, 1, 0, 0).$
\vspace{0.35ex}

Using the eigenvector $\nu^{01}\!=(1, 0, 1, 1, 0),$ the generalized eigenvector
\vspace{0.25ex}
$\nu^{11}\!=(1, 1/2, 1, 1, 0)$ of the 1-st order, and the generalized eigenvector 
\vspace{0.25ex}
$\nu^{21}=(0, 1, 1, 0, 0)$ of the 2-nd order corresponding to the 
eigenvalue $\lambda_1^{}={}-1,$ we can build the first integrals of system (1.23):
\\[2ex]
\mbox{}\hfill
$
F_1^{}\colon (t,x)\to \,
(x_1^{}+x_3^{}+x_4^{})\;\!e^{t}
$
\ for all 
$
(t,x)\in {\mathbb R}^6
$
\quad (by Theorem 1.9),
\hfill\mbox{}
\\[2.75ex]
\mbox{}\hfill
$
F_2^{}\colon (t,x)\to \
\dfrac{2x_1^{}+x_2^{}+2x_3^{}+2x_4^{}}{2(x_1^{}+x_3^{}+x_4^{})}\,-\, t
$
\ for all 
$
(t,x)\in {\mathbb R}\times {\mathscr X}_1^{}
$
\quad (by Theorem 1.10),
\hfill\mbox{}
\\[1.75ex]
and (by Theorem 1.8)
\\[2ex]
\mbox{}\hfill
$
F_3^{}\colon (t,x)\to \
\dfrac{4(x_2^{}+x_3^{})(x_1^{}+x_3^{}+x_4^{})-(2x_1^{}+x_2^{}+2x_3^{}+2x_4^{})^2}
{4(x_1^{}+x_3^{}+x_4^{})^2}
$
\ for all 
$
(t,x)\in {\mathbb R}\times {\mathscr X}_1^{},
\hfill
$
\\[2.25ex]
where ${\mathscr X}_1^{}$ is a domain from the set 
$\{x\colon x_1^{}+x_3^{}+x_4^{}\ne 0\}\subset {\mathbb R}^5.$
\vspace{0.75ex}

Using the eigenvector $\nu^{2}=(1, 1+i, 0, 1+i, 1-i)$ corresponding to the 
\vspace{0.25ex}
complex eigenvalue $\lambda_{2}^{}=1-i,$ we can construct (by Corollary 1.5) 
the first integrals of system (1.23): 
\\[2ex]
\mbox{}\hfill
$
F_4^{}\colon (t,x)\to \,
\bigl((x_1^{}+x_2^{}+x_4^{}+x_5^{})^2+(x_2^{}+x_4^{}-x_5^{})^2\;\!\bigr)\;\!e^{{}-2t}
$
\ for all 
$
(t,x)\in {\mathbb R}^6
\hfill
$
\\[1.25ex]
and
\\[1.5ex]
\mbox{}\hfill
$
F_5^{}\colon (t,x)\to \
\arctan\dfrac{x_2^{}+x_4^{}-x_5^{}}{x_1^{}+x_2^{}+x_4^{}+x_5^{}}\,+\, t
$
\ for all 
$
(t,x)\in {\mathbb R}\times {\mathscr X}_2^{},
\hfill
$
\\[2.25ex]
where ${\mathscr X}_2^{}$ is a domain from the set  
\vspace{0.75ex}
$\{x\colon x_1^{}+x_2^{}+x_4^{}+x_5^{}\ne 0\}\subset {\mathbb R}^5.$

The functionally independent first integrals $F_1^{},\ldots, F_5^{}$ 
\vspace{0.25ex}
are an integral basis of the linear autonomous differential system (1.23) 
on a domain ${\mathbb R}\times {\mathscr X},$ 
\vspace{0.25ex}
where ${\mathscr X}$ is a domain from the set 
$\{x\colon x_1^{}+x_3^{}+x_4^{}\ne 0,\ x_1^{}+x_2^{}+x_4^{}+x_5^{}\ne 0\}\subset {\mathbb R}^5.$
\vspace{1ex}

\newpage

\mbox{}
\\[-1.75ex]
\centerline{
\bf  
1.2. Linear nonhomogeneous differential system
}
\\[1.5ex]
\indent
Consider an ordinary linear nonhomogeneous differential system 
with constant coefficients
\\[1.75ex]
\mbox{}\hfill                                          % (1.24)
$
\dfrac{dx}{dt}=Ax+f(t),
$
\hfill (1.24)
\\[2ex]
where 
\vspace{0.5ex}
$f\colon t\to {\rm colon}(f_1^{}(t),\ldots,f_n^{}(t))$ for all $t\in J$ is a continuous function 
on an in\-ter\-val $J\subset {\mathbb R}.$ The differential system (1.1) is the corresponding 
\vspace{0.35ex}
homogeneous system of (1.24).

The system (1.24) is induced the linear differential ope\-ra\-tor of first order
\\[1ex]
\mbox{}\hfill
$
\displaystyle
{\frak C}(t,x)=\partial_{t}+
\sum\limits_{\xi=1}^n (A_{\xi}^{} x+f_{\xi}^{}(t))\,\partial_{x_\xi} 
$
for all 
$
(t,x)\in J\times {\mathbb R}^{n},
\ \,
A_\xi^{}\!=\!(a_{\xi 1}^{},\ldots,a_{\xi n}^{}),
\
\xi\!=\!1,\ldots, n.
\hfill
$
\\[1.5ex]
\indent
{\bf 1.2.1. Case of simple elementary divisors}
\\[0.75ex]
\indent
If the matrix $B$ has simple structure, then we can build  
first integrals of the differential system (1.24) by using following assertions
(Theorem 1.11 and Corollary 1.6).
\vspace{0.75ex}

{\bf Theorem 1.11.}
{\it
Let $\nu$ be a real eigenvector  of the matrix $B$ corresponding to the eigenvalue $\lambda.$ 
Then a first integral of system {\rm (1.24)} is the scalar function 
\\[1.5ex]
\mbox{}\hfill                                            % (1.25)
$
\displaystyle
F\colon (t,x)\to\  \nu x\,\exp({}-\lambda\;\! t)
-\int\limits_{t_0^{}}^t \nu f(\zeta)\,\exp({}-\lambda\;\! \zeta)\,d\zeta
$
\ for all 
$
(t,x)\in J\times {\mathbb R}^n,
$
\hfill {\rm (1.25)}
\\[1.5ex]
where $t_0^{}$ is a fixed point from the in\-ter\-val $J.$
}
\vspace{0.75ex}

{\sl Proof}. Using Lemma 1.1, we get
\\[1.5ex]
\mbox{}\hfill                                           
$
\displaystyle
{\frak C}\;\!F(t,x)=
\partial_{t}\;\!F(t,x)+{\frak A}\;\!F(t,x)+
\sum\limits_{\xi=1}^{n}f_\xi(t)\;\!\partial_{x_\xi^{}}^{}\!F(t,x)=
{}-\lambda\;\!\nu x\exp({}-\lambda\;\! t)
-\;\! \nu f(t)\exp({}-\lambda\;\! t) \, +
\hfill 
$
\\[1.75ex]
\mbox{}\hfill                                           
$
\displaystyle
+\ \lambda\,\nu x\,\exp({}-\lambda\;\! t) +
\sum\limits_{\xi=1}^{n}\nu_\xi^{} f_\xi^{}(t)\,\exp({}-\lambda\;\! t)=0
$
\ for all 
$
(t,x)\in J\times {\mathbb R}^n.
\hfill 
$
\\[1.5ex]
\indent
Therefore the function (1.25) is a first integral of the differential system (1.24). $\k$
\vspace{0.5ex}

{\bf Corollary 1.6.}
\vspace{0.35ex}
{\it
Suppose $\nu={\stackrel{*}{\nu}}+\widetilde{\nu}\,i\
({\rm Re}\,\nu={\stackrel{*}{\nu}},\ {\rm Im}\,\nu=\widetilde{\nu}\,)$ is an eigenvector of the matrix  $B$
corresponding to the complex eigenvalue 
\vspace{0.5ex}
$\lambda={\stackrel{*}{\lambda}}+\widetilde{\lambda}\,i\ 
({\rm Re}\,\lambda={\stackrel{*}{\lambda}},\ {\rm Im}\,\lambda=\widetilde{\lambda}\ne 0).$
Then first integrals of the differential system {\rm (1.24)} are the scalar functions 
\\[1.5ex]
\mbox{}\hfill                                           % (1.26)
$
\displaystyle
F_{{}_{\scriptstyle \theta}}\colon (t,x)\to\ 
\alpha_{{}_{\scriptstyle \theta}}(t, x)-
\int\limits_{t_0^{}}^{t} \alpha_{{}_{\scriptstyle \theta}}(\zeta, f(\zeta))\,d\zeta
$
\ for all 
$
(t,x)\in J\times {\mathbb R}^n,
\ \ \theta=1,2,
$
\hfill {\rm(1.26)}
\\[1.25ex]
where $t_0^{}$ is a fixed point from the in\-ter\-val $J,$ the functions
\\[1ex]
\mbox{}\hfill                                     
$
\displaystyle
\alpha_1^{}\colon (t,x)\to 
\bigl(\,{\stackrel{*}{\nu}}x\,\cos\widetilde{\lambda}\,t +
\widetilde{\nu}x\,\sin\widetilde{\lambda}\,t \bigr)\,
\exp\bigl({}-{\stackrel{*}{\lambda}}\,t\bigr)
$
\ for all 
$
(t,x)\in {\mathbb R}^{n+1},
\hfill
$
\\[1.35ex]
\mbox{}\hfill                                     
$
\displaystyle
\alpha_{2}^{}\colon (t,x)\to
\bigl(\,\widetilde{\nu}x\,\cos\widetilde{\lambda}\,t -
{\stackrel{*}{\nu}}x\,\sin\widetilde{\lambda}\,t\bigr)\,
\exp\bigl({}-{\stackrel{*}{\lambda}}\,t\bigr)
$
\ for all 
$(t,x)\in {\mathbb R}^{n+1}.
\hfill
$
\\[1.75ex]
}
\indent
{\sl Proof}. From Theorem 1.11 it follows that the function (1.25) is
a complex-valued first integral of the differential system (1.24). Then, 
the real and imaginary parts of this function are the real-valued 
first integrals (1.26) of the differential system (1.24).$\k$
\vspace{0.75ex}

{\bf Example 1.13.} 
Let us consider the linear nonhomogeneous differential system
\\[1.75ex]
\mbox{}\hfill                                                              % (1.27)                     
$
\dfrac{dx_1^{}}{dt} =  2x_{1}^{} + x_{2}^{}+2e^{2t},
\quad
\dfrac{dx_2^{}}{dt} = x_{1}^{} + 3x_{2}^{} -  x_{3}^{}+10,
\quad
\dfrac{dx_3^{}}{dt} =  {}-x_{1}^{} + 2x_{2}^{} + 3x_{3}^{}+e^{3t}.
$
\hfill (1.27)
\\[1.75ex]
\indent
The differential system (1.12) is the corresponding homogeneous system of (1.27).
\vspace{0.35ex}

Using the eigenvectors $\nu^{1}=(1,i,{}-1),\ \nu^{2}=(3,{}-1,{}-1)$ corresponding 
\vspace{0.35ex}
to the eigenvalues $\lambda_{1}^{}=3+i,\ \lambda_{2}^{}=2,$ respectively, we can build 
\vspace{0.35ex}
(by Corollary 1.6 and Theorem 1.11)
the functionally independent first integrals (an integral basis) of system (1.27)
\\[1.75ex]
\mbox{}\hfill                         
$
F_{1}^{}\colon (t,x_1^{},x_2^{},x_3^{})\to
\bigl( (x_{1}^{}-x_{3}^{}+1)\cos t+(x_{2}^{}+3)\sin t\bigr) e^{{}-3t} +
(\cos t-\sin t)e^{{}-t}+\sin t,
\hfill
$
\\[2.25ex]
\mbox{}\hfill                         
$
F_{2}^{}\colon (t,x_1^{},x_2^{},x_3^{})\to
\bigl( (x_{2}^{}+3)\cos t+(x_{3}^{}-x_{1}^{}-1)\sin t\bigr) e^{{}-3t} -
(\cos t+\sin t)e^{{}-t}+\cos t,
\hfill
$
\\[2.25ex]
\mbox{}\hfill                          
$
F_{3}^{}\colon (t,x_1^{},x_2^{},x_3^{})\to
(3x_1^{}-x_2^{}-x_3^{}-5)e^{{}-2t}+e^{t}-6t
$
\ for all 
$
(t,x_1^{},x_2^{},x_3^{})\in {\mathbb R}^4.
\hfill
$
\\[2ex]
\indent
{\sl Remark}. Under thr conditions of Corollary 1.6, we have the scalar function
\\[1ex]
\mbox{}\hfill                                     
$
\displaystyle
F\colon (t,x)\to \
\Bigr(\bigl({\stackrel{*}{\nu}}x\bigr)^2 +
\bigl(\widetilde{\nu}x\bigr)^2\Bigr)
\exp\bigl({}-2\,{\stackrel{*}{\lambda}}\,t\bigr)\, -
\, 2\biggl(\alpha_1^{}(t,x)\int\limits_{t_0^{}}^{t} \alpha_1^{}(\zeta,f(\zeta))\,d\zeta\ + 
\hfill                                     
$
\\[1.25ex]
\mbox{}\hfill                                     
$
\displaystyle
+\, \alpha_2^{}(t,x)\int\limits_{t_0^{}}^{t}\! \alpha_2^{}(\zeta,f(\zeta))\;\!d\zeta\biggr)\, +
\, \biggl(\, \int\limits_{t_0^{}}^{t}\! \alpha_1^{}(\zeta,f(\zeta))\;\!d\zeta\biggr)^2 + \,
\biggl(\, \int\limits_{t_0^{}}^{t}\! \alpha_2^{}(\zeta,f(\zeta))\;\!d\zeta\biggr)^2
\hfill
$
\\[2ex]
is also a first integral on the domain $J\times {\mathbb R}^n$
of the linear differential system (1.24).
\\[2ex]
\indent
{\bf 1.2.2. Case of multiple elementary divisors}
\\[0.75ex]
\indent
If the matrix $B$ has multiple elementary divisors, then we can construct 
first integrals of system (1.24) by using following assertions
\vspace{0.75ex}
(Theorem 1.12 and Corollary 1.7).

{\bf Theorem 1.12.}
\vspace{0.25ex}
{\it
Let $\lambda$ be the eigenvalue with elementary divisor of multiplicity $m\geq 2$ 
of the matrix $B$ corresponding to a real eigenvector $\nu^{0}$ and to real
\vspace{0.35ex}
generalized eigenvectors $\nu^{k},\, k=1,\ldots, m-1.$ 
Then the system {\rm (1.24)} has the fun\-c\-ti\-o\-na\-l\-ly independent first integrals
\\[1.75ex]
\mbox{}\hfill                                   % (1.28)
$
\displaystyle
F_{k+1}^{}\colon (t,x)\to \ 
\nu^kx\;\!\exp({}-\lambda\;\! t) \, - \,  
\sum\limits_{\tau=0}^{k-1}\!
{\textstyle \binom{k}{\tau}}\;\! t^{k-\tau} F_{\tau+1}^{}(t,x) \, -\,  C_{k}^{}(t)
\hfill
$
\\
\mbox{}\hfill {\rm (1.28)}
\\
\mbox{}\hfill
for all 
$
(t,x)\in J\times {\mathbb R}^n,
\quad
k=1,\ldots, m-1,
\mbox{}\hfill
$
\\[2ex]
where the integral 
\vspace{0.35ex}
$
F_1^{}\colon\! (t,x)\!\to
\nu^0 x\exp(-\;\!\lambda\;\! t)-C_0^{}(t)
\!$
for all\! 
$
(t,x)\!\in\! J\!\times\! {\mathbb R}^n\!
$
{\rm(}by Theorem\! {\rm 1.11)},
the scalar functions  
\\[1.5ex]
\mbox{}\hfill
$
\displaystyle
C_k^{}\colon t \to \
\int\limits_{t_0^{}}^{t} \bigl(\nu^kf(\zeta)\exp({}-\lambda\;\!\zeta) + k\,C_{k-1}^{}(\zeta)\bigr)\, d\zeta
$
\ for all 
$
t\in J, 
\ \ k=0,\ldots, m-1,
\ \ t_0^{}\in J.
\hfill
$
\\[1.75ex]
}
\indent
{\sl Proof}. The proof of Theorem 1.12 is by induction on $m.$
\vspace{0.25ex}

By the equalities (1.14), it follows that
\\[2ex]
\mbox{}\hfill                       % (1.29)
$
{\frak C}\;\!\bigl(\nu^{\,\varepsilon}x\,\exp({}-\lambda t)\bigr)=
\bigl(\varepsilon\,\nu^{\,\varepsilon-1}x+\nu^{\varepsilon}f(t)\bigr)\exp({}-\lambda\;\! t)
\hfill
$
\\
\mbox{}\hfill (1.29)
\\
\mbox{}\hfill
for all $(t,x)\in J\times {\mathbb R}^n,
\ \ 
\varepsilon=1\ldots, m-1.
\hfill 
$
\\[2ex]
\indent
Let $m=2.$ Using the system of identities (1.29), we get
\\[2ex]
\mbox{}\hfill                                        
$
\displaystyle
{\frak C}\;\!F_2^{}(t,x)=
{\frak C}\;\!\bigl(
\nu^{1}x\,\exp({}-\lambda t)-t\, F_1^{}(t,x) - C_1^{}(t)\bigr)=
\hfill                                        
$
\\[2ex]
\mbox{}\hfill                                        
$
\displaystyle
=\bigl(\nu^{0}x+\nu^{1}f(t)\bigr)\exp({}-\lambda t)- F_1^{}(t,x) -
\bigl(\nu^{1}f(t)\,\exp({}-\lambda t) +C_0^{}(t)\bigr)=
\hfill                                        
$
\\[2ex]
\mbox{}\hfill                                        
$
\displaystyle
=\bigl(\nu^{0}x\,\exp({}-\lambda t)- C_0^{}(t)\bigr) - F_1^{}(t,x)=0
$
\ for all 
$
(t,x)\in J\times {\mathbb R}^n.
\hfill
$
\\[2.25ex]
\indent
Therefore the function
\vspace{0.35ex}
$F_2^{}\colon J\times{\mathbb R}^n\to {\mathbb R}$ 
is a first integral of system (1.24).

Assume that the functions (1.28) for $m=\mu$ are first integrals of system (1.24).
\vspace{0.25ex}
Then, from the identity (1.29) for the function $F_{\mu+1}^{}$ it follows that
\\[1.5ex]
\mbox{}\hfill                                  
$
\displaystyle
{\frak C}F_{\mu+1}^{}(t,x)\!=\! 
{\frak C}\Bigl(\!\nu^{\mu} x\exp(-\lambda t) -  
\sum\limits_{\tau=0}^{\mu-1}\!
{\textstyle \binom{\mu}{\tau}}\;\! t^{\mu-\tau}\! F_{\tau+1}(t,x)  -  C_{\mu}^{}(t)\!\Bigr)\!=\!
\bigl(\mu\nu^{\mu-1}x+ \nu^{\mu}f(t)\!\bigr)\!\exp(-\lambda t) -
\hfill                              
$
\\[2ex]
\mbox{}\hfill                                  
$
\displaystyle
- \
\mu\  \sum\limits_{\tau=0}^{\mu-2}
{\textstyle \binom{\mu-1}{\tau}}\, t^{\,{}^{\scriptstyle \mu-\tau-1}}\, F_{\tau+1}^{}(t,x)-
\mu\, F_{\mu}^{}(t,x)  -
\bigl(\nu^{\mu}f(t)\,\exp({}-\lambda t)+
\mu\;\! C_{\mu-1}^{}(t)\bigr)=
\hfill 
$
\\[2ex]
\mbox{}\hfill                                  
$
\displaystyle
=\mu\Bigl( \nu^{\mu-1}x-  
\sum\limits_{\tau=0}^{\mu-2}
{\textstyle \binom{\mu-1}{\tau}}\;\! t^{\,{}^{\scriptstyle \mu-\tau-1}}\;\! F_{\tau+1}^{}(t,x)  -
C_{\mu-1}^{}(t)\Bigr)
 - \mu\;\! F_{\mu}(t,x) =0
$ 
for all 
$
(t,x)\in J\times {\mathbb R}^n.
\hfill 
$
\\[2ex]
\indent
This implies that the scalar function $F_{\mu+1}^{}\colon J\times{\mathbb R}^n\to {\mathbb R}$ 
\vspace{0.35ex}
(for $m=\mu+1)$ is a first integral of the linear nonhomogeneous differential system (1.24).

So by the principle of mathematical induction, the scalar functions (1.28) are 
first integrals of the differential system (1.24) for every natural number $m\geq 2.\ \k$
\vspace{0.75ex}

{\bf Example 1.14.}
Consider the linear nonhomogeneous differential system
\\[1.75ex]
\mbox{}\hfill                    % (1.30)               
$
\dfrac{dx_1^{}}{dt} =   4x_1^{} - x_2^{}+e^{3t},
\quad
\dfrac{dx_2^{}}{dt} =  3x_1^{} + x_2^{} - x_3^{}+8t,
\quad
\dfrac{dx_3^{}}{dt} =    x_1^{} + x_3^{}+4.
$
\hfill (1.30)
\\[2.25ex]
\indent
%The system (1.20) 
%\vspace{0.35ex}
%is the corresponding homogeneous differential system of system (1.30).
The differential system (1.30) has the eigenvalue $\lambda_{1}=2$ 
\vspace{0.35ex}
corresponding to the eigenvector $\nu^{0}=(1,{}-1,1)$ and to the generalized eigenvectors 
$\nu^{1}=(1,0,{}-1),\ \nu^{2}=({}-2,2,0).$
\vspace{0.35ex}

An integral basis of system (1.30) is the functions (by Theorems 1.11 and 1.12)
\\[1.5ex]
\mbox{}\hfill                           
$
F_{1}^{}\colon (t, x)\to\
(x_1-x_2+x_3-4t)e^{{}-2t}-e^{t}
$
\ for all 
$
(t, x)\in {\mathbb R}^4,
\hfill
$
\\[2ex]
\mbox{}\hfill                           
$
F_{2}^{}\colon (t, x)\to\
(x_1-x_3+2t-1) e^{{}-2t}-t\, F_1^{}(t,x)-2e^{t}
$
\ for all 
$
(t, x)\in {\mathbb R}^4,
\hfill
$
\\[2ex]
\mbox{}\hfill                           
$
F_{3}^{}\colon (t, x)\to\
2(x_2-x_1+3t+2) e^{{}-2t}-t^2\,F_1^{}(t,x)-2t\,F_2^{}(t,x)-2e^{t}
$ 
\ for all 
$
(t, x)\in {\mathbb R}^4.
\hfill
$
\\[2ex]
\indent
{\bf Corollary 1.7.}
{\it
Let 
\vspace{0.25ex}
$\lambda={\stackrel{*}{\lambda}}+\widetilde{\lambda}\,i
\ ({\rm Re}\,\lambda={\stackrel{*}{\lambda}},\ {\rm Im}\,\lambda=\widetilde{\lambda}\ne 0)$
be the complex eigenvalue of the matrix $B$ with elementary divisor of multiplicity $m\geq 2$ 
\vspace{0.25ex}
corresponding to an eigenvector
$\nu^0={\stackrel{*}{\nu}}\;\!{}^{0}+\widetilde{\nu}\;\!{}^{0}\,i
\ ({\rm Re}\,\nu^0={\stackrel{*}{\nu}}\;\!{}^{0},\ {\rm Im}\,\nu^0=\widetilde{\nu}\;\!{}^{0})$ 
and to generalized eigenvectors  
$\nu^k={\stackrel{*}{\nu}}\;\!{}^{k}+\widetilde{\nu}\;\!{}^{k}\,i$ 
$({\rm Re}\,\nu^k={\stackrel{*}{\nu}}\;\!{}^{k},\ {\rm Im}\,\nu^k=\widetilde{\nu}\;\!{}^{k}),
\ k=1,\ldots, m-1.$ 
\vspace{0.35ex}
Then first integrals of the linear nonhomogeneous differential system {\rm(1.24)} 
are the functions
\\[1.5ex]
\mbox{}\hfill                                  % (1.31)
$
\displaystyle
F_{{}_{\scriptstyle \theta, k+1}}\colon (t,x)\to \
\alpha_{{}_{\scriptstyle \theta k}}(t,x)-\,
\sum\limits_{\tau=0}^{k-1}
{\textstyle \binom{k}{\tau}}\, t^{k-\tau}\, F_{{}_{\scriptstyle \theta, \tau+1}}(t,x) -
C_{{}_{\scriptstyle \theta k}}(t)
\hfill                              
$
\\[0.25ex]
\mbox{}\hfill      {\rm (1.31)}                            
\\[0.25ex]
\mbox{}\hfill                                  
for all 
$
(t,x)\in J\times {\mathbb R}^n,
\quad
k=1,\ldots, m-1,
\quad
\theta=1,2,
\hfill                              
$
\\[1.75ex]
where the first integrals {\rm(}by Corollary {\rm 1.6)}
\\[1.75ex]
\mbox{}\hfill                                     
$
F_{{}_{\scriptstyle \theta\;\!1}}^{}\colon (t,x)\to \
\alpha_{{}_{\scriptstyle \theta\;\! 0}}(t, x)-
C_{{}_{\scriptstyle \theta\;\! 0}}(t)
$
\ for all 
$
(t,x)\in J\times {\mathbb R}^n,
\ \ \theta=1,2,
\hfill
$
\\[1.75ex]
the scalar functions
\\[1.5ex]
\mbox{}\hfill                                     
$
\displaystyle
\alpha_{{}_{\scriptstyle 1k}}\colon (t, x)\to
\bigl(\,{\stackrel{*}{\nu}}\;\!{}^kx\,\cos\widetilde{\lambda}\,t +
\widetilde{\nu}\;\!{}^kx\,\sin\widetilde{\lambda}\,t \bigr)
\exp\bigl({}-{\stackrel{*}{\lambda}}\,t\bigr)
$
for all 
$
(t,x)\in {\mathbb R}^{n+1},
\ k=0,\ldots, m-1,
\hfill
$
\\[2ex]
\mbox{}\hfill                                     
$
\displaystyle
\alpha_{{}_{\scriptstyle 2k}}\colon (t,x)\to
\bigl(\,\widetilde{\nu}\;\!{}^k x\,\cos\widetilde{\lambda}\,t -
{\stackrel{*}{\nu}}\;\!{}^kx\,\sin\widetilde{\lambda}\,t\bigr)
\exp\bigl({}-{\stackrel{*}{\lambda}}\,t\bigr)
$
for all 
$(t,x)\in {\mathbb R}^{n+1},
\ k=0,\ldots, m-1,
\hfill                              
$
\\[1.75ex]
\mbox{}\hfill                                  
$
\displaystyle
C_{{}_{\scriptstyle \theta k}}\colon t\to\,
\int\limits_{t_0^{}}^{t}\bigl(\alpha_{{}_{\scriptstyle \theta k}}(\zeta, f(\zeta)) +
k\;\!C_{{}_{\scriptstyle \theta, k-1}}(\zeta)\bigr)\;\! d\zeta
$
for all 
$
t\in J,
\ \,
k\!=\!0,\ldots, m\!-\!1,
\ \,
\theta\!=\!1,2,
\ \, 
t_0^{}\in J.
\hfill                              
$
\\[1.5ex]
}
\indent
{\sl Proof}. 
Formally using Theorem 1.12, we get 
the complex-valued first integrals (1.28) of the differential system (1.24). Then, 
the real and imaginary parts of the functions (1.28) are the real-valued 
first integrals (1.31) of the differential system (1.24).$\k$
\vspace{0.75ex}

{\bf Example 1.15.}
The linear nonhomogeneous differential system
\\[2ex]
\mbox{}\hfill                               % (1.32)
$
\dfrac{dx_1^{}}{dt} =   {}- x_3^{} +x_4^{}+4a\cos t,
\qquad
\dfrac{dx_2^{}}{dt} =  x_1^{} + x_3^{} +4a\sin t,
\hfill
$
\\[0.25ex]
\mbox{}\hfill  (1.32)
\\[0.25ex]
\mbox{}\hfill
$
\dfrac{dx_3^{}}{dt} =  {}-x_2^{} + x_3^{}-x_4^{}+b\;\!t,
\quad \
\dfrac{dx_4^{}}{dt} =  {}-x_1^{} + x_2^{} +x_3^{}-x_4^{}+\dfrac{c}{\sin^2 t}\,,
\hfill
$
\\[2.25ex]
where $a,\, b,$ and $c$ are some real numbers,
\vspace{0.25ex}
has the complex eigenvalues $\lambda_{1}^{}=i,\ \lambda_{2}^{}={}-i$ with
elementary divisors $(\lambda-i)^2,\ (\lambda+i)^2,$ respectively.
\vspace{0.5ex}

The number $\lambda_{1}^{}=i$ 
\vspace{0.5ex}
is the eigenvalue corresponding to the eigenvector
$\nu^{0}=(1, i, 1, 0)$ and to the generalized eigenvector 
$\nu^{1}=({}-1+i, 0, 0, i)$ of the 1-st order.
\vspace{0.35ex}

Using the real numbers
\vspace{0.35ex}
${\stackrel{*}{\lambda}}_1^{}=0,\ \widetilde{\lambda}_1^{}=1,$  
the real vectors
${\stackrel{*}{\nu}}{}^{0}=(1, 0, 1, 0),\ \widetilde{\nu}{}^{\,0}=(0, 1, 0, 0),$ 
${\stackrel{*}{\nu}}{}^{1}=({}-1, 0, 0, 0),\ \widetilde{\nu}{}^{\,1}=(1, 0, 0, 1),$
and the scalar functions
\\[2ex]
\mbox{}\hfill
$
\alpha_{{}_{\scriptstyle 10}}\colon (t,x)\to\,
\cos t\, (x_1^{}+x_3^{})+\sin t\, x_2^{},
\quad
\alpha_{{}_{\scriptstyle 20}}\colon (t,x)\to\,
\cos t\, x_2^{}-\sin t\,(x_1^{}+x_3^{}),
\hfill
$
\\[2.25ex]
\mbox{}\hfill
$
\alpha_{{}_{\scriptstyle 11}}\!\colon (t,x)\to
{}-\cos t\, x_1^{} +\sin t\, (x_1^{}+x_4^{}),
\quad 
\alpha_{{}_{\scriptstyle 21}}\!\colon (t,x)\to
\cos t\, (x_1^{}+x_4^{})+\sin t\,x_1^{}
\hfill
$
\\[2.25ex]
\mbox{}\hfill
for all 
$
(t,x)\in J_l^{}\times {\mathbb R}^4,
\quad 
J_l^{}=(\pi l;\pi (l+1))
$
\ for all $l\in {\mathbb Z},
\hfill
$
\\[2.5ex]
\mbox{}\hfill
$
\displaystyle
C_{{}_{\scriptstyle 10}}(t)=
\int \!
(4a+ b\;\!t\cos t)\;\!dt=
4at+b(\cos t+t\sin t),
\ \ 
C_{{}_{\scriptstyle 20}}(t)=
{}-\int \!
b\;\!t\sin t\;\!dt=
b(t\cos t-\sin t),
\hfill
$
\\[2.75ex]
\mbox{}\hfill
$
\displaystyle
C_{{}_{\scriptstyle 11}}(t)=
\int
\Bigl({}-4a\cos^2 t+ 4a\cos t\sin t+\dfrac{c}{\sin t}+4at+b(\cos t+t\sin t)\Bigr)\;\!dt=
\hfill
$
\\[2ex]
\mbox{}\hfill
$
\displaystyle
=
2a(t^2-t+\sin^2 t-\sin t\cos t) +b(2\sin t -t\cos t)+c\ln\Bigl|\tan\dfrac{t}{2}\Bigr|
$ 
\ for all 
$
t\in J_l^{},
\hfill
$
\\[2.75ex]
\mbox{}\hfill
$
\displaystyle
C_{{}_{\scriptstyle 21}}(t)=
\int
\Bigl(4a\cos^2 t+ \dfrac{c\cos t}{\sin^2 t}+4a\cos t\sin t+b(t\cos t-\sin t)\Bigr)\;\!dt=
\hfill
$
\\[2ex]
\mbox{}\hfill
$
\displaystyle
=
2a(t+\sin^2 t+\sin t\cos t) +b(2\cos t +t\sin t)-\dfrac{c}{\sin t}
$ 
\ for all 
$
t\in J_l^{},
\hfill
$
\\[2.25ex]
we can build (by Corollaries 1.6 and 1.7) the first integrals of system (1.32)
\\[2ex]
\mbox{}\hfill
$
F_{{}_{\scriptstyle 11}}\colon (t,x)\to\
\cos t\,(x_1^{}+x_3^{})+\sin t\, x_2^{} -4at -b(\cos t+t\sin t)
$ 
\ for all 
$
(t,x)\in J_l^{}\times {\mathbb R}^4,
\hfill
$
\\[2.5ex]
\mbox{}\hfill
$
F_{{}_{\scriptstyle 21}}\colon (t,x)\to\
\cos t\, x_2^{}-\sin t\, (x_1^{}+x_3^{})+b(\sin t-t\cos t)
$
\ for all 
$
(t,x)\in J_l^{}\times {\mathbb R}^4,
\hfill
$
\\[2.5ex]
\mbox{}\hfill
$
F_{{}_{\scriptstyle 12}}\colon (t,x)\to\
{}-\cos t\,x_1^{}+\sin t\, (x_1^{}+x_4^{}) -tF_{{}_{\scriptstyle 11}}(t,x)-
2a(t^2-t+\sin^2 t-\sin t\cos t)\ -
\hfill
$
\\[2ex]
\mbox{}\hfill
$
-\ b(2\sin t -t\cos t)-c\ln\Bigl|\tan\dfrac{t}{2}\Bigr|
$ 
\ for all 
$
(t,x)\in J_l^{}\times {\mathbb R}^4,
\hfill
$
\\[2.5ex]
\mbox{}\hfill
$
F_{{}_{\scriptstyle 22}}\colon (t,x)\to\
\cos t\,(x_1^{}+x_4^{})+\sin t\,x_1^{}  -tF_{{}_{\scriptstyle 21}}(t,x)-
2a(t+\sin^2 t+\sin t\cos t) \ -
\hfill
$
\\[2ex]
\mbox{}\hfill
$
-\ b(2\cos t +t\sin t)+\dfrac{c}{\sin t}
$
\ for all 
$(t,x)\in J_l^{}\times {\mathbb R}^4.
\hfill
$
\\[2ex]
\indent
The functionally independent first integrals
\vspace{0.75ex}
$F_{11}^{},\, F_{21}^{},\, F_{12}^{},$ and $F_{22}^{}$ are
a basis of first integrals for the differential system (1.32) on a domain $J_l^{}\times {\mathbb R}^4.$

\newpage

\mbox{}
\\[-1.75ex]
\centerline{                               %{Параграф 2}
\large\bf
2. Integrals of ordinary linear nonautonomous differential systems 
}
\\[0.75ex]
\centerline{
\large\bf
integrable in closed form
}
\\[2.5ex]
\centerline{                               %{П. 1, Параграф 1}
\bf 2.1.
Algebraic reducible systems
}
\\[1.25ex]
\indent
Consider an ordinary real linear nonhomogeneous differential system of the $n\!$-th order
\\[1.75ex]
\mbox{}\hfill                                          % (2.1)
$
\dfrac{dx}{dt}=A(t)\,x+f(t),
$
\hfill (2.1)
\\[2ex]
where 
\vspace{0.5ex}
$x={\rm colon}(x_1^{},\ldots,x_n^{})\in {\mathbb R}^n,$ 
the continuous on an interval $J\subset {\mathbb R}$ coefficient matrix
$A\colon t\to A(t)$ for all $t\in J$ is diagonalizable by 
\vspace{0.25ex}
a constant semilarity matrix [59, p. 61],
the vector function
$
f\colon t\to {\rm colon}(f_1^{}(t),\ldots,f_n^{}(t)) 
$
for all 
$t\in J$ 
is continuous.
\vspace{0.35ex}

The corresponding homogeneous system of the 
nonhomogeneous system (2.1) is
\\[1.75ex]
\mbox{}\hfill                                          % (2.2)
$
\dfrac{dx}{dt}=A(t)\,x.
$
\hfill (2.2)
\\[1.75ex]
\indent
By [48, p. 186], 
the nonhomogeneous differential system (2.1) and 
the corresponding homogeneous differential system (2.2) are called {\it algebraic reducible}.
\\[1.25ex]
\indent
{\bf 2.1.1. Partial integrals}
\\[0.5ex]
\indent
A complex-valued linear homogeneous function
\\[1.25ex]
\mbox{}\hfill                                % (2.3)
$
\displaystyle
p\colon x\to \nu x
$ 
\ for all 
$
x\in {\mathbb R}^{n},
\quad
\nu\in {\mathbb C}^{\;\!n},
$
\hfill (2.3)
\\[1.25ex]
is a partial integral of the algebraic reducible system (2.2) 
if and only if 
the identity holds
\\[1.5ex]
\mbox{}\hfill                                  % (2.4)
$
{\frak A} p(x) =  \lambda(t)\;\! p(x)
$
\ for all 
$
(t,x)\in J\times {\mathbb R}^{n},
$
\hfill (2.4)
\\[1.75ex]
where the linear differential operator
\vspace{0.35ex}
$
{\frak A}(t,x)=\partial_{t}^{}+
A(t)\,x\,\partial_{x}^{} 
$
for all 
$
(t,x)\in J\times {\mathbb R}^{n},
$
the scalar function 
$\lambda\colon J\to {\mathbb R}.$
The identity (2.4) is equivalent to the linear system
\\[1.25ex]
\mbox{}\hfill                         % (2.5)
$
\bigl(B(t) -  \lambda(t)\;\! E\bigr)\;\!\nu = 0,
$
\hfill (2.5)
\\[1.25ex]
where $E$ is the $n\times n$ identity matrix,
and the matrix $B\colon t\to B(t)$ for all $t\in J$
is the transpose of the matrix $A\colon t\to A(t)$ for all $t\in J.$
\vspace{0.15ex}

The following basic propositions (Lemmas 2.1 and 2.2)
are base for the method of building first integrals of 
the algebraic reducible systems (2.1) and (2.2).
\vspace{0.35ex}

{\bf Lemma 2.1.}
{\it 
Suppose $\nu$ is a real eigenvector of the matrix 
$B\colon t\to B(t)$ for all $t\in J$ 
corresponding to the eigenfunction $\lambda\colon t\to \lambda (t)$ for all $t\in J.$ 
Then the linear function {\rm (2.3)} is a partial integral 
of the algebraic reducible differential system {\rm (2.2)}. 
}
\vspace{0.35ex}

{\sl Indeed}, if $\nu\in {\mathbb R}^{n}$ is an eigenvector of the matrix 
$B\colon t\to B(t)$ for all $t\in J$ corresponding to the eigenfunction 
$\lambda\colon t\to \lambda (t)$ for all $t\in J,$ then  
$\nu$ is a solution to the system (2.5).
This implies that the identity (2.4) holds. 
Thus the linear function (2.3) is a partial integral 
of the algebraic reducible differential system {\rm (2.2)}. \k
\vspace{0.5ex}

{\bf Lemma 2.2.}
{\it
Let 
$\nu={\stackrel{*}{\nu}}+\widetilde{\nu}\,i\ 
({\stackrel{*}{\nu}}={\rm Re}\,\nu,\ \widetilde{\nu}={\rm Im}\,\nu\ne 0)$ 
be a complex eigenvector of the matrix $B\colon t\to B(t)$ for all $t\in J$ 
\vspace{0.35ex}
corresponding to the eigenfunction
$\lambda\colon t\to {\stackrel{*}{\lambda}}(t)+\widetilde{\lambda}(t)\;\!i$ for all 
$t\in J\ 
({\stackrel{*}{\lambda}}\colon t\to {\rm Re}\,\lambda(t), \ 
\widetilde{\lambda}\colon t\to {\rm Im}\,\lambda(t)$ for all $t\in J),$
\vspace{0.75ex}
let the functions $P$ and $\psi$ be defined by
\vspace{0.35ex}
$
P\colon x\to 
({\stackrel{*}{\nu}}x)^2+(\widetilde{\nu}x)^2
$
for all 
$
x\in {\mathbb R}^n
$
and
$
\displaystyle
\psi\colon x\to 
\arctan\dfrac{\widetilde{\nu}x}{{\stackrel{*}{\nu}}x}
$
for all 
$
x\in {\mathscr X},
$
where ${\mathscr X}$ is a domain  
from the set 
$\bigl\{x\colon {\stackrel{*}{\nu}}x\ne 0\bigr\}\subset {\mathbb R}^n.$ 
Then, we have
\\[2ex]
\mbox{}\hfill                                            
$
{\frak A} P(x)=2\;\!{\stackrel{*}{\lambda}}(t)\;\!P(x)
$
for all 
$
(t,x)\in J\times {\mathbb R}^n
$
\
and
\ 
$
\displaystyle
{\frak A}\psi(x)=
\widetilde{\lambda}(t)
$ 
for all 
$
(t,x)\in J\times {\mathscr X}.
\hfill 
$
}
\\[2ex]
\indent
{\sl Proof}. 
Formally using Lemma 2.1, we get 
the complex-valued function (2.3) is a partial integral 
of the differential system (2.2) and the following identity holds
\\[1.25ex]
\mbox{}\hfill
$
{\frak A}\;\!
\bigl({\stackrel{*}{\nu}}x+i\,\widetilde{\nu} x\bigr)=
\bigl({\stackrel{*}{\lambda}}(t)+i\,\widetilde{\lambda}(t)\bigr)
\bigl({\stackrel{*}{\nu}}x+i\,\widetilde{\nu} x\bigr)
$
\ for all 
$
(t,x)\in J\times {\mathbb R}^{n}.
\hfill
$
\\[1.35ex]
\indent
This complex identity is equivalent to the real system of identities
\\[1.5ex]
\mbox{}\hfill
$
{\frak A}\;\!{\stackrel{*}{\nu}} x =
{\stackrel{*}{\lambda}}(t)\,{\stackrel{*}{\nu}} x -
\widetilde{\lambda}(t)\,\widetilde{\nu} x,
\quad \
{\frak A}\;\!\widetilde{\nu} x=
{\stackrel{*}{\lambda}}(t)\,\widetilde{\nu} x+
\widetilde{\lambda}(t)\,{\stackrel{*}{\nu}} x
$
\ for all 
$
(t,x)\in J\times {\mathbb R}^{n}.
\hfill
$
\\[1.5ex]
\indent
Using this system of identities, we obtain
\\[1.5ex]
\mbox{}\hfill                        
$
{\frak A}\;\!P(x)=
{\frak A}\bigl(
({\stackrel{*}{\nu}} x)^2+(\widetilde{\nu} x)^2\bigr)=
2\;\!{\stackrel{*}{\nu}} x\,{\frak A}\;\!{\stackrel{*}{\nu}} x +
2\;\!\widetilde{\nu} x\,{\frak A}\;\!\widetilde{\nu} x=
2\;\!{\stackrel{*}{\nu}} x\,
\bigl(\;\!{\stackrel{*}{\lambda}}(t)\,{\stackrel{*}{\nu}} x -
\widetilde{\lambda}(t)\,\widetilde{\nu} x\bigr)  \ + 
\hfill                        
$
\\[1.75ex]
\mbox{}\hfill                        
$
+\ 2\;\!\widetilde{\nu} x\,
\bigl(\;\!{\stackrel{*}{\lambda}}(t)\,\widetilde{\nu} x+
\widetilde{\lambda}(t)\,{\stackrel{*}{\nu}} x\bigr)=
2\;\!{\stackrel{*}{\lambda}}(t)\;\!
\bigl(({\stackrel{*}{\nu}} x)^2+(\widetilde{\nu} x)^2\bigr)=
2\;\!{\stackrel{*}{\lambda}}(t)\,P(x)
$
\ for all 
$
(t,x)\in J\times {\mathbb R}^{n}
\hfill                        
$
\\[1.5ex]
and                        
\\[1.5ex]
\mbox{}\hfill                        
$
{\frak A}\,\psi(x)\, =\, 
%{\frak A}\,\arctg\dfrac{\widetilde{\nu} x}{{\stackrel{*}{\nu}} x} \, =\,
%\dfrac{1}{1+\dfrac{(\widetilde{\nu} x)^2}{({\stackrel{*}{\nu}} x)^2}}
%\cdot
%\dfrac{{\stackrel{*}{\nu}} x\, {\frak A}\;\!\widetilde{\nu} x - 
%\widetilde{\nu} x\, {\frak A}\;\!{\stackrel{*}{\nu}} x}
%{({\stackrel{*}{\nu}} x)^2}\ =
%\hfill                        
%$
%\\[2.25ex]
%\mbox{}\hfill                        
%$
%=
\dfrac{
{\stackrel{*}{\nu}} x\;\! 
\bigl(\;\!{\stackrel{*}{\lambda}}(t)\,\widetilde{\nu} x+
\widetilde{\lambda}(t)\,{\stackrel{*}{\nu}} x\bigr)  -\,
\widetilde{\nu} x\;\! 
\bigl(\;\!{\stackrel{*}{\lambda}}(t)\,{\stackrel{*}{\nu}} x -
\widetilde{\lambda}(t)\,\widetilde{\nu} x\bigr)}
{({\stackrel{*}{\nu}} x)^2+(\widetilde{\nu} x)^2}=
\widetilde{\lambda}(t)
$
\ for all 
$
(t,x)\in J\times {\mathscr X}.
\ \k
\hfill                        
$
\\[2.5ex]
\indent
{\bf 2.1.2. First integrals}
\vspace{0.5ex}

Using Theorems 2.1 and 2.2, we can obtain first integrals of the 
algebraic reducible nonhomogeneous differential system (2.1). 
\vspace{0.75ex}

{\bf Theorem 2.1.}\!
{\it
Suppose $\nu$ is a real eigenvector of the matrix 
$B\colon t\to B(t)$ for all $t\in J$ 
corresponding to the eigenfunction $\lambda\colon t\to \lambda (t)$ for all $t\in J.$ 
Then the algebraic reducible differential system {\rm (2.1)} has the first integral  
\\[1.5ex]
\mbox{}\hfill                                            % (2.6)
$
\displaystyle
F\colon (t,x)\to\, 
\nu x\;\! \varphi(t) \, -\, 
\int\limits_{t_0^{}}^{t} \nu f(\zeta)\;\! \varphi(\zeta)\;\!d\zeta
$
\ for all 
$
(t,x)\in J\times\R^n,
$
\hfill {\rm (2.6)}
\\[1.5ex]
where the exponential function
\\[1.5ex]
\mbox{}\hfill
$
\varphi\colon t\to \exp\biggl(\!{}-{\displaystyle \int\limits_{t_0^{}}^{t}} 
\lambda(\zeta)\,d\zeta\!\biggr)
$ 
\ for all 
$t\in J,
\hfill
$ 
\\[1.5ex]
and $t_0^{}$ is a fixed point from the interval $J.$
}
\vspace{0.75ex}

{\sl Proof}. 
From Lemma 2.1, we get the following
\\[1.75ex]
\mbox{}\hfill
$
\displaystyle
{\frak B}F(t,x)=
{\frak A}F(t,x)+ f(t)\,\partial_{x}^{}F(t,x)= 
\nu x\, \partial_{t}^{}\;\! \varphi(t) - 
\partial_{t}^{}\int\limits_{t_0^{}}^{t} \nu f(\zeta)\;\! \varphi(\zeta)\;\!d\zeta
+\varphi(t)\,{\frak A}\;\!\nu x \ +
\hfill
$
\\[2.25ex]
\mbox{}\hfill
$
\displaystyle
+ \ f(t)\,\partial_{x}^{}\;\! \nu x\, \varphi(t)=
{}-\lambda(t)\,\nu x\;\!\varphi(t) -
\nu f(t)\;\!\varphi(t)+
\lambda(t)\,\nu x\;\!\varphi(t)+\nu f(t)\;\!\varphi(t)=0
\hfill
$
\\[2ex]
\mbox{}\hfill
for all 
$
(t,x)\in J\times\R^n,
\hfill 
$
\\[2ex]
where the linear differential operator
\vspace{0.35ex}
$
{\frak B}(t,x)=\partial_{t}^{}+ \bigl(A(t)\,x+f(t)\bigr)\,\partial_{x}^{} 
$
for all 
$
(t,x)\in J\times {\mathbb R}^{n}
$
is the operator of differentiation by virtue of system (2.1).
\vspace{0.25ex}

Therefore the function (2.6) is a first integral of the algebraic reducible system (2.1). \k
\vspace{1ex}

{\bf Example 2.1.}
Let us consider the linear differential system of the third-order
\\[2ex]
\mbox{}\hfill
$
\dfrac{dx_1^{}}{dt} =\alpha_3^{}(t)\,x_1^{}+(\alpha_1^{}(t)-\alpha_3^{}(t))(x_3^{}-x_2^{})+f_1^{}(t),
\qquad \, \mbox{}
\hfill
$
\\[2.25ex]
\mbox{}\hfill
$
\dfrac{dx_2^{}}{dt} = (\alpha_2^{}(t)-\alpha_3^{}(t))\,(x_3^{}-x_1^{}) +\alpha_3^{}(t)\,x_2^{}+f_2^{}(t),
$
\hfill {\rm (2.7)}
\\[2ex]
\mbox{}\hfill
$
\dfrac{dx_3^{}}{dt} = (\alpha_3^{}(t)-\alpha_2^{}(t))\,x_1^{} + (\alpha_3^{}(t)-\alpha_1^{}(t))\,x_2^{}  +
(\alpha_1^{}(t)+\alpha_2^{}(t)-\alpha_3^{}(t))\,x_3^{}+f_3^{}(t),
\hfill 
$
\\[2ex]
where the scalar functions
\vspace{0.5ex}
$\alpha_{\xi}^{}\colon J\to {\mathbb R}$ and $f_{\xi}^{}\colon J\to{\mathbb R},\ \xi=1,2, 3,$ 
are continuous on an interval $J\subset {\mathbb R}.$
The coefficient matrix
\\[2.5ex]
\mbox{}\hfill
$
A\colon  t\to 
\left\|\!\!
\begin{array}{ccc}
\alpha_3^{}(t) & \alpha_3^{}(t)-\alpha_1^{}(t) & \alpha_1^{}(t)-\alpha_3^{}(t) 
\\[1.5ex]
\alpha_3^{}(t)-\alpha_2^{}(t) & \alpha_3^{}(t) & \alpha_2^{}(t)-\alpha_3^{}(t)
\\[1.5ex]
\alpha_3^{}(t)-\alpha_2^{}(t) & \alpha_3^{}(t)-\alpha_1^{}(t) & \alpha_1^{}(t)+\alpha_2^{}(t)-\alpha_3^{}(t)
\end{array}
\!\!\right\|
$
\
for all 
$
t\in  J
\hfill
$
\\[2.5ex]
such that 
the transposed matrix $B\colon t\to A^{T}(t)$ for all $t\in J$
\vspace{0.5ex}
has the real eigenvectors
$\nu^1=(0,{}-1,1),\ \nu^2=({}-1,0,1),$ and $\nu^3=(1,1,{}-1)$
\vspace{0.75ex}
corresponding to the eigenfunctions
$\lambda_{\xi}^{}\colon t\to \alpha_{\xi}^{}(t)$ for all 
$t\in J,\ \xi=1,2, 3,$ respectively.
\vspace{0.75ex}

Therefore the differential system (2.7) is algebraic reducible.
\vspace{0.35ex}

By Theorem 2.1, we can build the first integrals on the domain 
$J\times {\mathbb R}^3$ of system (2.7)
\\[1.75ex]
\mbox{}\hfill
$
\displaystyle
F_{1}^{}\colon (t, x)\to\, 
(x_3^{}-x_2^{})\exp\biggl({}-\int\limits_{t_0^{}}^{t} \alpha_1^{}(\zeta)\,d\zeta\biggr)\, - \,
\int\limits_{t_0^{}}^{t} \bigl(f_3^{}(\zeta)-f_2^{}(\zeta)\bigr)
\exp\biggl({}-\int\limits_{\zeta_0^{}}^{\zeta} \alpha_1^{}(\theta)\,d\theta\biggr)\,d\zeta,
\hfill
$
\\[2ex]
\mbox{}\hfill
$
\displaystyle
F_{2}^{}\colon (t, x)\to\, 
(x_3^{}-x_1^{})\exp\biggl({}-\int\limits_{t_0^{}}^{t} \alpha_2^{}(\zeta)\,d\zeta\biggr)\ -
\ \int\limits_{t_0^{}}^{t} \bigl(f_3^{}(\zeta)-f_1^{}(\zeta)\bigr)
\exp\biggl({}-\int\limits_{\zeta_0^{}}^{\zeta} \alpha_2^{}(\theta)\,d\theta\biggr)\,d\zeta,
\hfill
$
\\[2ex]
\mbox{}\hfill
$
\displaystyle
F_{3}^{}\colon (t, x)\to 
(x_1^{}+x_2^{}-x_3^{})\exp\biggl(\!\!{}-\int\limits_{t_0^{}}^{t}\!\! \alpha_3^{}(\zeta)\,d\zeta\biggr) -
\int\limits_{t_0^{}}^{t}\!\! \bigl(f_1^{}(\zeta)+f_2^{}(\zeta)-f_3^{}(\zeta)\bigr)
\exp\biggl(\!{}-\int\limits_{\zeta_0^{}}^{\zeta}\!\! \alpha_3^{}(\theta)\,d\theta\biggr)\,d\zeta,
\hfill
$
\\[2ex]
where $t_0^{}$ and $\zeta_0^{}$ are fixed points from the interval $J.$
\vspace{0.5ex}

The functionally independent first integrals $F_1^{},\ F_2^{},$ and $F_3^{}$ 
\vspace{0.25ex}
are an integral basis on the domain $J\times {\mathbb R}^3$ 
of the algebraic reducible differential system (2.7).
\vspace{0.75ex}

{\bf Theorem 2.2.}
{\it
Let 
$\nu={\stackrel{*}{\nu}}+\widetilde{\nu}\,i\ 
({\stackrel{*}{\nu}}={\rm Re}\,\nu,\ \widetilde{\nu}={\rm Im}\,\nu\ne 0)$ 
be a complex eigenvector of the matrix $B\colon t\to B(t)$ for all $t\in J$ 
\vspace{0.35ex}
corresponding to the eigenfunction
$\lambda\colon t\to {\stackrel{*}{\lambda}}(t)+\widetilde{\lambda}(t)\;\!i$ for all 
$t\in J\ 
({\stackrel{*}{\lambda}}\colon t\to {\rm Re}\,\lambda(t), \ 
\widetilde{\lambda}\colon t\to {\rm Im}\,\lambda(t)$ for all $t\in J).$
\vspace{0.5ex}
Then the algebraic reducible differential system {\rm (2.1)} has the first integrals  
\\[1.75ex]
\mbox{}\hfill                                           % (2.8)
$
\displaystyle
F_{\tau}^{}\colon (t,x)\to 
\gamma_{\tau}^{}(t, x)-
\int\limits_{t_0^{}}^{t} \gamma_{\tau}^{}(\zeta, f(\zeta))\,d\zeta
$
\ for all 
$
(t,x)\in J\times {\mathbb R}^n,
\quad
\tau=1, \ \tau=2,
$
\hfill {\rm(2.8)}
\\[2ex]
where $t_0^{}$ is a fixed point from the interval $J,$ the scalar functions
\\[1.75ex]
\mbox{}\hfill                                     
$
\displaystyle
\gamma_{1}^{}\colon (t, x)\to
\biggl({\stackrel{*}{\nu}}x\;\!\cos\!\int\limits_{t_0^{}}^{t} \widetilde{\lambda}(\zeta)\;\!d\zeta  +
\widetilde{\nu}x\;\!\sin\!\int\limits_{t_0^{}}^{t}  \widetilde{\lambda}(\zeta)\;\!d\zeta \biggr)
\!\exp\biggl(\!{}-\int\limits_{t_0^{}}^{t}  {\stackrel{*}{\lambda}}(\zeta)\;\!d\zeta\!\biggr)\!
$
for all 
$
(t,x)\in J\times {\mathbb R}^n,
\hfill
$
\\[2.25ex]
\mbox{}\hfill                                     
$
\displaystyle
\gamma_{2}^{}\colon (t, x)\to
\biggl(\widetilde{\nu}x\;\!\cos\!\int\limits_{t_0^{}}^{t} \widetilde{\lambda}(\zeta)\;\!d\zeta  -
{\stackrel{*}{\nu}}x\;\!\sin\!\int\limits_{t_0^{}}^{t}  \widetilde{\lambda}(\zeta)\;\!d\zeta \biggr)
\!\exp\biggl(\!{}-\int\limits_{t_0^{}}^{t}  {\stackrel{*}{\lambda}}(\zeta)\;\!d\zeta\!\biggr)\!
$
for all 
$
(t,x)\in J\times {\mathbb R}^n.
\hfill
$
}
\\[1.75ex]
\indent
{\sl Proof}.  
Formally using Theorem 2.1, we have 
the complex-valued function (2.6) is a first integral 
of the algebraic reducible system (2.1). Then the real and imaginary parts 
of this complex-valued first integral are the real-valued first integrals (2.8) of 
system (2.1). \k
\vspace{0.5ex}

\newpage

{\bf Example 2.2.}
The algebraic reducible differential system
\\[1.75ex]
\mbox{}\hfill                                                                   % (2.9)
$
\dfrac{dx_1^{}}{dt} ={}-\tanh t\,x_1^{}-\cosh t\,(x_2^{}-x_3^{})-3,
\hfill 
$
\\[2ex]
\mbox{}\hfill
$
\dfrac{dx_2^{}}{dt} = \Bigl(\dfrac{1}{t}-1-\tanh t\Bigr) x_1^{} +
\Bigl(1-\dfrac{1}{t}-\cosh t\Bigr) x_2^{} +\cosh t\,x_3^{} +
2t^3+\sin t- e^t,
$
\hfill {\rm (2.9)}
\\[2.25ex]
\mbox{}\hfill
$
\dfrac{dx_3^{}}{dt} = 
\Bigl(\dfrac{1}{t}-1-\tanh t-\cosh t\Bigr) x_1^{} +
\Bigl(1-\dfrac{1}{t}+\tanh t-\cosh t\Bigr) x_2^{} +(\cosh t-\tanh t)\,x_3^{} \ +
\hfill 
$
\\[2ex]
\mbox{}\hfill
$
+\ 2t^3+\sin t- e^t-\sinh^2 t
\hfill 
$
\\[2.25ex]
has the eigenvectors  
\vspace{0.75ex}
$\nu^1=({}-1,1,0), \ \nu^2=(1,i,{}-i),$ and $\nu^3=(1,{}-i,i)$  
corresponding to the eigenfunctions
\vspace{0.75ex}
$
\lambda_1^{}\colon t\to 1-\dfrac{1}{t}\;\!,\  
\lambda_2^{}\colon t\to {}-\tanh t+\cosh t\,i,$ and 
$\lambda_3^{}\colon t\to {}-\tanh t-\cosh t\,i$ for all $t\in J,$ 
respectively,
\vspace{0.5ex}
where $J$ is an interval from the set $\{t\colon t\ne 0\}.$  

Using the eigenvector $\nu^1=({}-1,1,0)$ and the 
\vspace{0.35ex}
corresponding eigenfunction $\lambda_1^{}\colon t\to 1-\dfrac{1}{t}$ for all $t\in J,$
we can build (by Theorem 2.1) the first integral
\\[1.75ex]
\mbox{}\hfill
$
F_1^{}\colon (t,x)\to\,
te^{{}-t}(x_2^{}-x_1^{})+\dfrac{t^2}{2}+
\Bigl(2t^4+8t^3+24t^2+45t+45+\dfrac{1}{2}\,\cos t+ \dfrac{1}{2}\,t(\cos t+\sin t)\Bigr)e^{{}-t}
\hfill 
$
\\[1.75ex]
on the domain $J\times {\mathbb R}^3$ of the algebraic reducible differential system (2.9). 
\vspace{0.5ex}

Using the eigenvector $\nu^2=(1,i,{}-i)$ and the 
\vspace{0.35ex}
eigenfunction $\lambda_2^{}\colon t\to {}-\tanh t+\cosh t\,i$ for all $t\in J,$
we can build (by Theorem 2.2) the first integrals of system (2.9) 
\\[2ex]
\mbox{}\quad \ 
$
F_2^{}\colon (t,x)\to\
\cosh t\bigl(\cos(\sinh t)\, x_1^{}+\sin(\sinh t)\, (x_2^{}-x_3^{})\bigr) \ +
\hfill 
$
\\[2ex]
\mbox{}\hfill
$
+\ (\sinh^2 t-2)\cos(\sinh t) + (3-2\sinh t)\sin(\sinh t)
$
\ for all 
$
J\times {\mathbb R}^3,
$
\\[2.5ex]
\mbox{}\quad \ 
$
F_3^{}\colon (t,x)\to\
\cosh t\bigl({}-\sin(\sinh t)\, x_1^{}+\cos(\sinh t)\, (x_2^{}-x_3^{})\bigr) \ +
\hfill 
$
\\[2ex]
\mbox{}\hfill
$
+\ (3-2\sinh t)\cos(\sinh t) + (2-\sinh^2 t)\sin(\sinh t) 
$
\ for all 
$
J\times {\mathbb R}^3.
$
\\[2ex]
\indent
The functionally independent first integrals $F_1^{},\ F_2^{},$ and $F_3^{}$ 
\vspace{0.25ex}
are an integral basis on the domain $J\times {\mathbb R}^3$ 
of the algebraic reducible differential system (2.9).
\vspace{0.75ex}

In the case of the homogeneous algebraic reducible differential system (2.2), 
we have the following statements (Corollary 2.1 and Theorem 2.3).
\vspace{0.75ex}

{\bf Corollary 2.1.}
{\it
Under the conditions of Theorem {\rm 2.1}, we have the scalar function
\\[1.5ex]
\mbox{}\hfill                                            
$
\displaystyle
F\colon (t,x)\to\, \nu x\;\!\exp\biggl({}-\int\limits_{t_0^{}}^{t}\! \lambda(\zeta)\;\!d\zeta\biggr)
$ 
\ for all 
$
(t,x)\in J\times {\mathbb R}^n,
\quad
t_0^{}\in J, 
\hfill 
$
\\[1.5ex]
is a first integral of the homogeneous algebraic reducible differential system {\rm (2.2)}. 
}
\vspace{1.25ex}

{\bf Theorem 2.3.}
{\it
Let the assumptions of Theorem  {\rm 2.2} hold, then the scalar functions  
\\[1.5ex]
\mbox{}\hfill                                            
$
\displaystyle
F_1^{}\colon (t,x)\to 
\Bigl(\bigl({\stackrel{*}{\nu}}x\bigr)^2+\bigl(\widetilde{\nu}x\bigr)^2\Bigr)
\exp\biggl({}-2\int\limits_{t_0^{}}^{t}  {\stackrel{*}{\lambda}}(\zeta)\,d\zeta\biggr)
$
\ for all 
$
(t,x)\in J\times{\mathbb R}^n
\hfill 
$
\\[0.25ex]
and
\\[1ex]
\mbox{}\hfill                                            
$
\displaystyle
F_2^{}\colon (t,x)\to \,
\arctan\dfrac{\widetilde{\nu}x}{{\stackrel{*}{\nu}}x} \ - \
\int\limits_{t_0^{}}^{t} \widetilde{\lambda}(\zeta)\,d\zeta
$ 
\ for all 
$
(t,x)\in  J\times {\mathscr X},
\quad
t_0^{}\in J,
\quad
{\mathscr X}\subset\bigl\{x\colon {\stackrel{*}{\nu}}x\ne 0\bigr\},
\hfill 
$
\\[1ex]
are first integrals of the homogeneous algebraic reducible differential system {\rm (2.2)}. 
}
\vspace{0.75ex}

{\sl Proof}.
From Lemma 2.2, we get  
\\[1.25ex]
\mbox{}\hfill                                            
$
\displaystyle
{\frak A}F_1^{}(t,x)= 
\exp\biggl(\!\!{}-2\int\limits_{t_0^{}}^{t}\!  {\stackrel{*}{\lambda}}(\zeta)\,d\zeta\!\biggr)\,
{\frak A}\Bigl(\!\bigl({\stackrel{*}{\nu}}x\bigr)^2+\bigl(\widetilde{\nu}x\bigr)^2\Bigr)  +
\Bigl(\!\bigl({\stackrel{*}{\nu}}x\bigr)^2+\bigl(\widetilde{\nu}x\bigr)^2\Bigr)\, 
{\frak A}\exp\biggl(\!\!{}-2\int\limits_{t_0^{}}^{t}\!  {\stackrel{*}{\lambda}}(\zeta)\,d\zeta\!\biggr) =
\hfill 
$
\\[1.5ex]
\mbox{}\hfill                                            
$
\displaystyle
=2{\stackrel{*}{\lambda}}(t)
\Bigl(\!\bigl({\stackrel{*}{\nu}}x\bigr)^2+\bigl(\widetilde{\nu}x\bigr)^2\Bigr)
\exp\biggl(\!\!{}-2\int\limits_{t_0^{}}^{t}\!  {\stackrel{*}{\lambda}}(\zeta)\,d\zeta\biggr)-
2{\stackrel{*}{\lambda}}(t)
\Bigl(\!\bigl({\stackrel{*}{\nu}}x\bigr)^2+\bigl(\widetilde{\nu}x\bigr)^2\Bigr)
\exp\biggl(\!\!{}-2\int\limits_{t_0^{}}^{t}\!  {\stackrel{*}{\lambda}}(\zeta)\,d\zeta\biggr)=0
\hfill 
$
\\[1.5ex]
\mbox{}\hfill                                            
for all 
$
(t,x)\in J\times{\mathbb R}^n,
\hfill 
$
\\[2ex]
\mbox{}\hfill                                            
$
\displaystyle
{\frak A}F_2^{}(t,x)= 
{\frak A}\arctan\dfrac{\widetilde{\nu}x}{{\stackrel{*}{\nu}}x} \ - \
{\frak A}\int\limits_{t_0^{}}^{t} \widetilde{\lambda}(\zeta)\,d\zeta=
\widetilde{\lambda}(t)- \widetilde{\lambda}(t)=0
$
\ for all 
$
(t,x)\in  J\times {\mathscr X}.
\ \k
\hfill 
$
\\[2.75ex]
\indent
{\bf Example 2.3.}
Consider the algebraic reducible differential system
\\[2.25ex]
\mbox{}\hfill
$
\dfrac{dx_1^{}}{dt} =\alpha_2^{}(t)\,x_1^{}+\alpha_3^{}(t)\,(x_3^{}-x_2^{}),
\ \ \ 
\dfrac{dx_2^{}}{dt} = (\alpha_2^{}(t)-\alpha_1^{}(t))\,x_1^{} +(\alpha_1^{}(t)-\alpha_3^{}(t))\,x_2^{} +
\alpha_3^{}(t)\,x_3^{},
\hfill 
$
\\[1.5ex]
\mbox{}\hfill {\rm (2.10)}
\\[0.25ex]
\mbox{}\hfill
$
\dfrac{dx_3^{}}{dt} = (\alpha_2^{}(t)-\alpha_1^{}(t)-\alpha_3^{}(t))\,x_1^{} + 
(\alpha_1^{}(t)-\alpha_2^{}(t)-\alpha_3^{}(t))\,x_2^{} + (\alpha_2^{}(t)+\alpha_3^{}(t))\,x_3^{},
\hfill 
$
\\[2.5ex]
where the scalar functions 
\vspace{0.75ex}
$\alpha_{\xi}^{}\colon J\to {\mathbb R},\ \xi=1, 2, 3,$ 
are continuous on an interval $J\subset {\mathbb R}.$

The system (2.10) has the eigenvectors
\vspace{0.5ex}
$\nu^1=({}-1,1,0), \ \nu^2=(1,i,{}-i),$ and $\nu^3=(1,{}-i,i)$
corresponding to the eigenfunctions 
\vspace{0.75ex}
$\lambda_1^{}\colon t\to \alpha_1^{}(t)$ for all $t\in J,\
\lambda_2^{}\colon t\to \alpha_2^{}(t)+\alpha_3^{}(t)\,i$ 
for all $t\in J,$ and 
$\lambda_3^{}\colon t\to \alpha_2^{}(t)-\alpha_3^{}(t)\,i$ for all $t\in J,$
respectively.
\vspace{0.75ex}

Using the eigenvector $\nu^1=({}-1,1,0)$ and the 
\vspace{0.5ex}
corresponding eigenfunction $\lambda_1^{}\colon t\to \alpha_1^{}(t)$ for all $t\in J,$
we can build (by Corollary 2.1) the first integral
\\[1.75ex]
\mbox{}\hfill
$
\displaystyle
F_1^{}\colon (t,x)\to\,
(x_2^{}-x_1^{})\exp\biggl({}-\int\limits_{t_0^{}}^{t} \alpha_1^{}(\zeta)\,d\zeta\biggr)
$
\ for all 
$
(t,x)\in J\times {\mathbb R}^3
\hfill 
$
\\[1.75ex]
of the algebraic reducible differential system (2.10). 
\vspace{0.5ex}

Using the eigenvector $\nu^2=(1,i,{}-i)$ and the 
\vspace{0.5ex}
corresponding complex-valued eigenfunction 
$\lambda_2^{}\colon t\to \alpha_2^{}(t)+\alpha_3^{}(t)\,i$ for all $t\in J,$
\vspace{0.5ex}
we can construct (by Theorem 2.3) the first integrals 
of the algebraic reducible differential system (2.10)
\\[2ex]
\mbox{}\hfill
$
\displaystyle
F_{2}^{}\colon (t, x)\to\, 
\bigl(x_1^2+(x_2^{}-x_3^{})^2\bigr)
\exp\biggl({}-2\int\limits_{t_0^{}}^{t}  \alpha_2^{}(\zeta)\,d\zeta\biggr)
$
\ for all 
$
(t,x)\in J\times {\mathbb R}^3
\hfill
$
\\[1.5ex]
and
\\[1.5ex]
\mbox{}\hfill
$
\displaystyle
F_{3}^{}\colon (t, x)\to\
\arctan\dfrac{x_2^{}-x_3^{}}{x_1}\,- \,
\int\limits_{t_0^{}}^{t} \alpha_3^{}(\zeta)\,d\zeta
$
\ for all 
$
(t,x)\in J\times {\mathscr X},
\hfill
$
\\[2ex]
where ${\mathscr X}$ is a domain from the set $\{x\colon x_1^{}\ne 0\}\subset {\mathbb R}^3.$
\vspace{0.75ex}

The functionally independent first integrals $F_1^{},\ F_2^{},$ and $F_3^{}$ 
\vspace{0.5ex}
are an integral basis on the domain $J\times {\mathscr X}$ 
of the algebraic reducible differential system (2.10).

\newpage

\mbox{}
\\[-1.75ex]
\centerline{                               %{П. 2, Параграф 2}
\bf 2.2.
Triangular systems
}
\\[1.5ex]
\indent
Consider an triangular linear nonhomogeneous differential system of the $n\!$-th order
\\[2ex]
\mbox{}\hfill                                          % (2.11)
$
\displaystyle
\dfrac{dx_i^{}}{dt} =\sum\limits_{j=i}^n a_{ij}^{}(t)\,x_j^{}+f_i^{}(t),
\quad
i=1,\ldots, n,
$
\hfill (2.11)
\\[2ex]
where $a_{ij}^{}\colon J\to {\mathbb R}$ and $f_i^{}\colon J\to {\mathbb R}$
\vspace{0.5ex}
are continuous functions on an interval $J\subset {\mathbb R}.$ 

The corresponding homogeneous system of the 
nonhomogeneous system (2.11) is
\\[1.75ex]
\mbox{}\hfill                                          % (2.12)
$
\displaystyle
\dfrac{dx_i^{}}{dt} =\sum\limits_{j=i}^n a_{ij}^{}(t)\,x_j^{},
\quad
i=1,\ldots, n.
$
\hfill (2.12)
\\[1.75ex]
\indent
By Perron's theorem [56], any linear differential system can be transformed to an 
upper triangular system by an orthogonal transformation\footnote[1]{
The transformation $x=U(t)\;\!y$ is called {\it orthogonal} if
\vspace{0.35ex}
the matrix $U\colon t\to U(t)$ for all $t\in J$ is orthogonal, i.e., 
$U(t)\;\!U^{T}\!(t)=U^{T}\!(t)\;\!U(t)=E$ for all $t\in J,$ where $E$ is 
the identity matrix.
}\!.\!
A triangular system is in\-te\-g\-ra\-ted in closed form and 
its fundamental matrix can be chosen triangular [8, pp. 32 -- 33].

Using Theorem 2.4 and Corollary 2.2, we can build integral bases for 
the triangular differential systems (2.11) and (2.12).
\vspace{0.75ex}
We now state the main results of this subsection.

{\bf Theorem 2.4.}
{\it 
First integrals of system {\rm (2.11)} are the functions
\\[1.5ex]
\mbox{}\hfill
$
\displaystyle
F_{{}_{\scriptstyle \tau}}^{}\colon (t,x)\to \
x_{{}_{\scriptstyle n+1-\tau}}\varphi_{{}_{\scriptstyle n+1-\tau}}(t)\, -\, 
\sum\limits_{\xi=1}^{\tau-1}A_{{}_{\scriptstyle \tau\xi}}(t)\;\!F_{{}_{\scriptstyle \xi}}(t,x)-
B_{{}_{\scriptstyle \tau}}(t) 
\hfill
$
\\[0.35ex]
\mbox{}\hfill {\rm(2.13)}
\\[0.35ex]
\mbox{}\hfill
for all 
$
(t,x)\in J\times {\mathbb R}^n,
\ \ 
\tau=1,\ldots, n,
\hfill
$
\\[1.75ex]
where the scalar functions
\\[2ex]
\mbox{}\hfill
$
\displaystyle
A_{{}_{\scriptstyle \tau\xi}}\colon t\to 
\int\limits_{t_0^{}}^{t}
\biggl(\ \sum\limits_{k=1}^{\tau-\xi}a_{{}_{\scriptstyle n+1-\tau,n+1+k-\tau}}(\zeta)\,
\psi_{{}_{\scriptstyle n+1+k-\tau}}(\zeta)\,
A_{{}_{\scriptstyle \tau-k,\xi}}(\zeta)\biggr)\varphi_{{}_{\scriptstyle n+1-\tau}}(\zeta)\,d\zeta,
\hfill
$
\\[-0.5ex]
\mbox{}\hfill {\rm (2.14)}
\\[1.5ex]
\mbox{}\hfill
$
\xi=1,\ldots, \tau-1, 
\  
\tau=1,\ldots, n,
\quad 
A_{\zeta\zeta}\colon t\to 1,
\ \zeta=1,\ldots, n-1,
\ \
A_{10}\colon t\to 0
$
\ for all $t\in J,
\hfill
$
\\[2.75ex]
\mbox{}\hfill
$
\displaystyle
B_{{}_{\scriptstyle \tau}}\colon t\to 
\int\limits_{t_0^{}}^{t}
\biggl(f_{{}_{\scriptstyle n+1-\tau}}(\zeta)+ \sum\limits_{\xi=1}^{\tau-1}
a_{{}_{\scriptstyle n+1-\tau,n+1-\tau+\xi}}(\zeta)\,
\psi_{{}_{\scriptstyle n+1-\tau+\xi}}(\zeta)
\ B_{{}_{\scriptstyle \tau-\xi}}(\zeta)\biggr)\varphi_{{}_{\scriptstyle n+1-\tau}}(\zeta)\,d\zeta,
\hfill
$
\\[2.25ex]
\mbox{}\hfill
$
\displaystyle
\varphi_{{}_{\scriptstyle \tau}}\colon t\to 
\exp\biggl({}-\int\limits_{t_0^{}}^{t} a_{{}_{\scriptstyle \tau\tau}}(\zeta)\,d\zeta\biggr),
\ 
\psi_{{}_{\scriptstyle \tau}}\colon t\to 
\exp \int\limits_{t_0^{}}^{t} a_{{}_{\scriptstyle \tau\tau}}(\zeta)\,d\zeta
$
for all 
$
t\in J,
\ 
\tau=1,\ldots, n.
$
\hfill {\rm (2.15)}
\\[2ex]
}
\indent
{\sl Proof}.
The system (2.11) is induced the linear differential ope\-ra\-tor of first order
\\[1.5ex]
\mbox{}\hfill
$
\displaystyle
{\frak L}(t,x)=\partial_{t}^{}+
\sum\limits_{i=1}^n\;\!\Bigl(\, \sum\limits_{j=i}^na_{ij}^{}(t)\,x_j^{}+f_i^{}(t)\!\Bigr)\,\partial_{x_i}^{}
$ 
\ for all 
$
(t,x)\in J\times {\mathbb R}^n.
\hfill
$
\\[1.5ex]
\indent
If $\tau=1,$ then the scalar function 
\vspace{0.5ex}
$
F_1^{}\colon (t,x)\to x_{n}^{}\varphi_{n}^{}(t)-B_1^{}(t)
$ 
for all 
$
(t,x)\in J\times {\mathbb R}^n
$
is a first integral of the differential system (2.11):
\\[1.5ex]
\mbox{}\hfill
$
{\frak L}\,F_{{}_{\scriptstyle 1}}(t,x)=
x_{{}_{\scriptstyle n}}\partial_{t}^{}\,\varphi_{{}_{\scriptstyle n}}(t) +
\bigl(a_{{}_{\scriptstyle nn}}(t)\;\!x_{{}_{\scriptstyle n}}\!+\!f_{n}(t)\!\bigr)
\partial_{x_n}\bigl(x_{n}\varphi_{{}_{\scriptstyle n}}(t)\!\bigr)- \partial_{t}B_1(t)=
{}-a_{{}_{\scriptstyle nn}}(t)\;\!\varphi_{{}_{\scriptstyle n}}(t)\;\!x_{n}\ +
\hfill
$
\\[1.5ex]
\mbox{}\hfill
$
\displaystyle
+\
\bigl((a_{{}_{\scriptstyle nn}}(t)\;\!x_{{}_{\scriptstyle n}}+f_{n}^{}(t)\bigr)\varphi_{{}_{\scriptstyle n}}(t)-
\partial_{t}^{}\int\limits_{t_0^{}}^{t}\! f_{n}(\zeta)\;\!\varphi_{{}_{\scriptstyle n}}(\zeta)\;\! d\zeta=0
$
\ for all 
$
(t,x)\in J\times {\mathbb R}^n.
\hfill
$
\\[1.5ex]
\indent 
Let $\tau=2.$ Then the Lie derivative of function $F_2^{}$ by virtue of system (2.11) is equal to
\\[2ex]
\mbox{}\hfill
$
\displaystyle
{\frak L}\,F_{{}_{\scriptstyle 2}}(t,x)=
{\frak L}\bigl(
x_{{}_{\scriptstyle n-1}}\varphi_{{}_{\scriptstyle n-1}}(t)-
A_{{}_{\scriptstyle 21}}(t)F_{{}_{\scriptstyle 1}}(t,x)-B_{{}_{\scriptstyle 2}}(t)
\bigr)=
\hfill
$
\\[2.25ex]
\mbox{}\hfill
$
\displaystyle
={}- a_{{}_{\scriptstyle n-1,n-1}}(t)\,\varphi_{{}_{\scriptstyle n-1}}(t)\,x_{{}_{\scriptstyle n-1}} +
\bigl(a_{{}_{\scriptstyle n-1,n-1}}(t)\,x_{{}_{\scriptstyle n-1}} +
a_{{}_{\scriptstyle n-1,n}}(t)\,x_{{}_{\scriptstyle n}}+f_{n-1}(t)\bigr)\,\varphi_{{}_{\scriptstyle n-1}}(t) \ -
\hfill
$
\\[1.5ex]
\mbox{}\hfill
$
\displaystyle
-\, \partial_{t}\biggl(
\int\limits_{t_0^{}}^{t}\!\!  
a_{{}_{\scriptstyle n-1,n}}\!(\zeta)\;\!\psi_{{}_{\scriptstyle n}}\!(\zeta)\;\!
\varphi_{{}_{\scriptstyle n-1}}\!(\zeta)d\zeta\!\biggr)\! 
F_{{}_{\scriptstyle 1}}\!(t,x) -\;\!
\partial_{t}\biggl(\int\limits_{t_0^{}}^{t}\!\! 
\bigl(f_{n-1}(\zeta)+
 a_{{}_{\scriptstyle n-1,n}}\!(\zeta)\;\!\psi_{{}_{\scriptstyle n}}\!(\zeta) B_{{}_{\scriptstyle 1}}\!(\zeta)\!\bigr)
\varphi_{{}_{\scriptstyle n-1}}\!(\zeta)d\zeta\!\biggr)\!\! =
\hfill
$
\\[1.5ex]
\mbox{}\hfill
$
=
\bigl(a_{{}_{\scriptstyle n-1,n}}\!(t)\;\!x_{{}_{\scriptstyle n}}+f_{n-1}\!(t)\bigr)
\varphi_{{}_{\scriptstyle n-1}}\!(t) 
-
a_{{}_{\scriptstyle n-1,n}}\!(t)\;\!\psi_{{}_{\scriptstyle n}}\!(t)\;\!\varphi_{{}_{\scriptstyle n-1}}\!(t)\;\! 
F_{{}_{\scriptstyle 1}}(t,x)\ -
\hfill
$
\\[2.5ex]
\mbox{}\hfill
$
-\ \bigl(f_{{}_{\scriptstyle n-1}}(t)+
a_{{}_{\scriptstyle n-1,n}}(t)\,\psi_{{}_{\scriptstyle n}}(t)\,
B_{{}_{\scriptstyle 1}}(t)\bigr)\,\varphi_{{}_{\scriptstyle n-1}}(t) =
\hfill
$
\\[2.25ex]
\mbox{}\hfill
$
=a_{{}_{\scriptstyle n-1,n}}(t)\psi_{{}_{\scriptstyle n}}(t)\varphi_{{}_{\scriptstyle n-1}}(t) 
\bigl(x_{{}_{\scriptstyle n}}\varphi_{{}_{\scriptstyle n}}(t) - B_{{}_{\scriptstyle 1}}(t)-
F_{{}_{\scriptstyle 1}}(t,x)\bigr)=0
$
\ for all 
$
(t,x)\in J\times {\mathbb R}^n.
\hfill
$
\\[2ex]
\indent
Therefore the function $F_2^{}\colon J\times {\mathbb R}^n\to {\mathbb R}$ is
a first integral of system (2.11).
\vspace{0.5ex}

Suppose the functions
\vspace{0.35ex}
$F_{\tau}^{}\colon J\times {\mathbb R}^n\to {\mathbb R},\ \tau=1,\ldots, \varepsilon-1,$
are first integrals of the triangular differential system (2.11). 
Then, for $\tau=\varepsilon,$ we have
\\[1.75ex]
\mbox{}\qquad
$
\displaystyle
{\frak L}\,F_{{}_{\scriptstyle \varepsilon}}(t,x)=
{\frak L}\,x_{{}_{\scriptstyle n+1-\varepsilon}}^{}\, \varphi_{{}_{\scriptstyle n+1-\varepsilon}}^{}(t) \,+\,
x_{{}_{\scriptstyle n+1-\varepsilon}}^{}\, {\frak L}\;\! \varphi_{{}_{\scriptstyle n+1-\varepsilon}}(t) \ -
\hfill
$
\\[2ex]
\mbox{}\hfill
$
\displaystyle
-\ \sum\limits_{\xi=1}^{\varepsilon-1}{\frak L}\;\!A_{{}_{\scriptstyle \varepsilon\xi}}(t)\, 
F_{{}_{\scriptstyle \xi}}(t,x) \ - \ 
\sum\limits_{\xi=1}^{\varepsilon-1}A_{{}_{\scriptstyle \varepsilon\xi}}(t)\, 
{\frak L}\;\! F_{{}_{\scriptstyle \xi}}(t,x)-
{\frak L}\;\!B_{{}_{\scriptstyle \varepsilon}}(t) =
\quad\mbox{}
$
\\[1.75ex]
\mbox{}
$
\displaystyle
= \!\!
\biggl(\sum\limits_{\xi=0}^{\varepsilon-1}\!
a_{{}_{\scriptstyle n+1-\varepsilon, n+1-\varepsilon+\xi}}(t)x_{{}_{\scriptstyle n+1-\varepsilon+\xi}}\!+
f_{{}_{\scriptstyle n+1-\varepsilon}}(t)\!\!\biggr)\!\varphi_{{}_{\scriptstyle n+1-\varepsilon}}\!(t) -
a_{{}_{\scriptstyle n+1-\varepsilon, n+1-\varepsilon}}\!(t)x_{{}_{\scriptstyle n+1-\varepsilon}}
\varphi_{{}_{\scriptstyle n+1-\varepsilon}}\!(t) -
\hfill
$
\\[1.75ex]
\mbox{}\hfill
$
-\ \displaystyle
\sum\limits_{\xi=1}^{\varepsilon-1}\biggl(\, 
\sum\limits_{k=1}^{\varepsilon-\xi}\!\!
a_{{}_{\scriptstyle n+1-\varepsilon, n+1-\varepsilon+k}}(t)\,\psi_{{}_{\scriptstyle n+1-\varepsilon+k}}(t)
A_{{}_{\scriptstyle \varepsilon-k,\xi}}(t)\!\!\biggr) \varphi_{{}_{\scriptstyle n+1-\varepsilon}}(t)
F_{{}_{\scriptstyle \xi}}(t,x)\ -  
\hfill
$
\\[2ex]
\mbox{}\hfill
$
\displaystyle
- \ \biggl(f_{{}_{\scriptstyle n+1-\varepsilon}}(t) + 
\sum\limits_{\xi=1}^{\varepsilon-1}
a_{{}_{\scriptstyle n+1-\varepsilon,n+1-\varepsilon+\xi}}(t)\,
\psi_{{}_{\scriptstyle n+1-\varepsilon+\xi}}(t)
B_{{}_{\scriptstyle \varepsilon-\xi}}(t)\biggr)\varphi_{{}_{\scriptstyle n+1-\varepsilon}}(t) =
\hfill
$
\\[2ex]
\mbox{}\hfill
$
\displaystyle
=\!\!\biggl(\sum\limits_{\xi=1}^{\varepsilon-1}\!
a_{{}_{\scriptstyle n+1-\varepsilon, n+1-\varepsilon+\xi}}\!(t)x_{{}_{\scriptstyle n+1-\varepsilon+\xi}}\!-
\sum\limits_{\xi=1}^{\varepsilon-1}
\biggl(\! a_{{}_{\scriptstyle n+1-\varepsilon, n+1-\varepsilon+\xi}}(t)
\psi_{{}_{\scriptstyle n+1-\varepsilon+\xi}}\!(t)\!
\sum\limits_{k=1}^{\varepsilon-\xi}\!
A_{{}_{\scriptstyle \varepsilon-\xi,k}}\!(t)F_{{}_{\scriptstyle k}}\!(t,x)\!\!\biggr)\! -
\hfill
$
\\[2ex]
\mbox{}\hfill
$
\displaystyle
-\ \sum\limits_{\xi=1}^{\varepsilon-1}
a_{{}_{\scriptstyle n+1-\varepsilon,n+1-\varepsilon+\xi}}(t)\,
\psi_{{}_{\scriptstyle n+1-\varepsilon+\xi}}(t)
B_{{}_{\scriptstyle \varepsilon-\xi}}(t)
\biggr)
\varphi_{{}_{\scriptstyle n+1-\varepsilon}}(t) =
\hfill
$
\\[2ex]
\mbox{}\hfill
$
\displaystyle
=\sum\limits_{\xi=1}^{\varepsilon-1}
a_{{}_{\scriptstyle n+1-\varepsilon, n+1-\varepsilon+\xi}}(t)\,\psi_{{}_{\scriptstyle n+1-\varepsilon+\xi}}(t)
\biggl(x_{{}_{\scriptstyle n+1-\varepsilon+\xi}}\varphi_{{}_{\scriptstyle n+1-\varepsilon-\xi}}(t) -
\sum\limits_{k=1}^{\varepsilon-\xi-1}
A_{{}_{\scriptstyle \varepsilon-\xi,k}}(t)F_{{}_{\scriptstyle k}}(t,x) \ -
\hfill
$
\\[2ex]
\mbox{}\hfill
$
\displaystyle
-\  B_{{}_{\scriptstyle \varepsilon-\xi}}(t) -
A_{{}_{\scriptstyle \varepsilon-\xi,\varepsilon-\xi}}(t)F_{{}_{\scriptstyle \varepsilon-\xi}}(t,x)\biggr)
\varphi_{{}_{\scriptstyle n+1-\varepsilon}}(t) =0
$
\ for all 
$
(t,x)\in J\times {\mathbb R}^n.
\hfill
$
\\[1.75ex]
\indent
This yields that if $\tau=\varepsilon,$ then 
\vspace{0.35ex}
$F_{\varepsilon}^{}\colon J\times {\mathbb R}^n\to {\mathbb R}$ is a 
first integral of system (2.11).

Thus the functions (2.13) are functionally independent 
first integrals of system (2.11). \k

{\bf Example 2.4.}
The linear nonhomogeneous differential system
\\[2ex]
\mbox{}\hfill
$
\dfrac{dx_1^{}}{d\;\!t}=
{}-\dfrac{1}{t}\,x_1^{}+
\dfrac{\sqrt{2}}{2}\,(6-e^{{}-t})x_2^{}-
\dfrac{\sqrt{2}}{2}\,(6+e^{{}-t})x_3^{}+
8t^5+4t+2(t^3+3t^2+6t+6)e^{{}-t},
\hfill
$
\\[3ex]
\mbox{}\hfill
$
\dfrac{dx_2^{}}{d\;\!t}=
\dfrac{1}{2}\,\Bigl(1+\dfrac{2}{t}-5e^t\Bigr)x_2^{}+
\dfrac{1}{2}\,\Bigl(1-\dfrac{2}{t}+5e^t\Bigr)x_3^{}+
2\sqrt{2}\,(e^t-t^3),
$
\hfill (2.16)
\\[3ex]
\mbox{}\hfill
$
\dfrac{dx_3^{}}{d\;\!t}=
\dfrac{1}{2}\,\Bigl(1-\dfrac{2}{t}-5e^t\Bigr)x_2^{}+
\dfrac{1}{2}\,\Bigl(1+\dfrac{2}{t}+5e^t\Bigr)x_3^{}+
2\sqrt{2}\,e^t
\hfill
$
\\[2.25ex]
can be reduced by the orthogonal transformation
\\[1.5ex]
\mbox{}\hfill
$
x_1^{}=y_1^{},
\quad 
x_2^{}=\dfrac{\sqrt{2}}{2}\,(y_2^{}-y_3^{}),
\quad
x_3^{}=\dfrac{\sqrt{2}}{2}\,(y_2^{}+y_3^{})
\hfill
$ 
\\[1.75ex]
to the triangular differential system
\\[2.25ex]
\mbox{}\hfill
$
\dfrac{dy_1^{}}{d\;\!t}=
{}-\dfrac{1}{t}\,y_1^{}-e^{{}-t}y_2^{}-6y_3^{}+
8t^5+4t+2(t^3+3t^2+6t+6)e^{{}-t},
\hfill
$
\\[0.5ex]
\mbox{}\hfill (2.17)
\\[0.5ex]
\mbox{}\hfill
$
\dfrac{dy_2^{}}{d\;\!t}=
y_2^{}+5e^ty_3^{}-2t^3+ 4e^t,
\qquad
\dfrac{dy_3^{}}{d\;\!t}=
\dfrac{2}{t}\,y_3^{}+2t^3.
\hfill
$
\\[2.5ex]
\indent
By Theorem 2.4, using the functions
\\[2ex]
\mbox{}\hfill
$
\varphi_1^{}\colon t\to t,
\ \
\varphi_2^{}\colon t\to e^{{}-t},
\ \
\varphi_3^{}\colon t\to \dfrac{1}{t^2}\,,
\ \
\psi_1^{}\colon t\to \dfrac{1}{t}\,,
\ \
\psi_2^{}\colon t\to e^{t},
\ \
\psi_3^{}\colon t\to t^2
$
for all 
$
t\in J,
\hfill
$
\\[2.5ex]
\mbox{}\hfill
$
A_{11}^{}\colon t\to 1,
\ \
A_{21}^{}\colon t\to \dfrac{5}{3}\,t^3,
\ \
A_{31}^{}\colon t\to \dfrac{1}{6}\,(9-2t)\;\!t^4,
\ \
A_{32}^{}\colon t\to {}-\dfrac{1}{2}\,t^2
$
\ for all 
$
t\in J,
\hfill
$
\\[2.75ex]
\mbox{}\hfill
$
B_{1}^{}\colon t\to t^2,
\ \
B_{2}^{}\colon t\to t^5+4t+2(t^3+3t^2+6t+6)e^{{}-t},
\ \
B_{3}^{}\colon t\to (t-1)\;\!t^6
$
\ for all 
$
t\in J,
\hfill
$
\\[2ex]
where $J$ 
\vspace{0.35ex}
is an interval from the set $\{t\colon t\ne 0\},$ we can build the 
functionally independent first integrals on the domain $J\times {\mathbb R}^3$
of  the triangular differential system (2.17):
\\[2ex]
\mbox{}\hfill
$
F_{1}^{}\colon (t,y)\to \dfrac{1}{t^2}\,y_3^{}-t^2,
\quad
F_{2}^{}\colon (t,y)\to e^{{}-t}y_2^{}-\dfrac{5}{3}\,ty_3^{}+
\dfrac{2}{3}\,t^5-4t -2(t^3+3t^2+6t+6)e^{{}-t},
\hfill
$
\\[2.75ex]
\mbox{}\hfill
$
F_{3}^{}\colon (t,y)\to ty_1^{}+\dfrac{1}{2}\,t^2e^{{}-t}y_2^{}+\dfrac{1}{2}\,(3-t)t^2y_3^{}-
\dfrac{1}{2}\,t^3(2t^4+t^3+4)-t^2(t^3+3t^2+6t+6)e^{{}-t}.
\hfill
$
\\[2ex]
%\mbox{}\hfill
%$
%\forall (t,y)\in J\times {\mathbb R}^3,
%\quad
%J\subset \{t\colon t\ne 0\}.
%\hfill
%$
%\\[2ex]
\indent
Now using the inverse transformation 
\vspace{1ex}
$y_1^{}=x_1^{},\, y_2^{}=\dfrac{\sqrt{2}}{2}\,(x_2^{}+x_3^{}),\,
y_3^{}=\dfrac{\sqrt{2}}{2}\,(x_3^{}-x_2^{}),$
we obtain the first integrals of the differential system (2.16):
\\[2ex]
\mbox{}\hfill
$
\widetilde{F}_{1}^{}\colon (t,x)\to\ 
\dfrac{\sqrt{2}}{2t^2}\,(x_3^{}-x_2^{})-t^2
$
\ for all 
$
(t,x)\in J\times {\mathbb R}^3,
\hfill
$
\\[2.25ex]
\mbox{}\hfill
$
\widetilde{F}_{2}^{}\colon (t,x)\to \
\dfrac{\sqrt{2}}{6}\,(5t+3e^{{}-t})x_2^{}-
\dfrac{\sqrt{2}}{6}\,(5t-3e^{{}-t})x_3^{}+
\dfrac{2}{3}\,t^5-4t -2(t^3+3t^2+6t+6)e^{{}-t},
\hfill
$
\\[2.25ex]
\mbox{}\hfill
$
\widetilde{F}_{3}^{}\colon \ 
(t,x)\to tx_1^{}+
\dfrac{\sqrt{2}}{4}\,t^2(t-3+e^{{}-t})x_2^{}+
\dfrac{\sqrt{2}}{4}\,t^2(3-t+e^{{}-t})x_3^{} \ -
\hfill
$
\\[1.75ex]
\mbox{}\hfill
$
- \ \dfrac{1}{2}\,t^3(2t^4+t^3+4)-t^2(t^3+3t^2+6t+6)e^{{}-t}
$
\ for all 
$
(t,x)\in J\times {\mathbb R}^3,
\quad
J\subset \{t\colon t\ne 0\}.
\hfill
$
\\[1.75ex]
\indent
The functionally independent first integrals 
\vspace{0.25ex}
$\widetilde{F}_{1},\ \widetilde{F}_{2},$ and $\widetilde{F}_{3}$ 
are an integral basis on any domain $J\times {\mathbb R}^3$ of 
the linear differential system (2.16). 
\vspace{0.75ex}

From Theorem 2.4 with $f_{\tau}^{}(t)=0$ for all $t\in J,\ \tau=1,\ldots, n,$
we have the following  
\vspace{1ex}

{\bf Corollary 2.2.}
{\it 
First integrals of system {\rm (2.12)} are the functions
\\[1.75ex]
\mbox{}\hfill
$
\displaystyle
F_{{}_{\scriptstyle \tau}}\colon (t,x)\to \,
x_{{}_{\scriptstyle n+1-\tau}}\varphi_{{}_{\scriptstyle n+1-\tau}}(t)-
\sum\limits_{\xi=1}^{\tau-1}A_{{}_{\scriptstyle \tau\xi}}(t)F_{{}_{\scriptstyle \xi}}(t,x)
$
\ for all 
$
(t,x)\in J\times {\mathbb R}^n,
\ \
\tau=1,\ldots, n,
\hfill
$
\\[2ex]
where the functions 
\vspace{0.5ex}
$A_{{}_{\scriptstyle \tau\xi}}\colon J\to {\mathbb R},\ \xi=1,\ldots,\tau,\ \tau=1,\ldots, n,$ and 
$\varphi_{{}_{\scriptstyle \tau}}\colon J\to {\mathbb R},\ \tau=1,\ldots, n,$ 
are given by the formulas {\rm (2.14)} and {\rm (2.15)}, respectively.
}
\vspace{1.25ex}

{\bf Example 2.5.}
The second-order linear homogeneous differential system 
\\[1.75ex]
%\mbox{}\hfill
$
\dfrac{dx_1^{}}{dt} =\bigl(\cos^2t\, a(t)-\cos t\sin t\, b(t)+\sin^2 t\, c(t)\bigr)\,x_1^{} +
\bigl(\cos t\sin t (a(t)- c(t))+\cos^2 t\, b(t)-1\bigr)\,x_2^{},
\hfill
$
\\[0.75ex]
\mbox{}\hfill {\rm (2.18)}
\\[0.5ex]
%\mbox{}\hfill
$
\dfrac{dx_2^{}}{dt} =\bigl(\cos t\sin t (a(t)- c(t))-\sin^2 t\, b(t)+1\bigr)\,x_1^{} +
\bigl(\sin^2t\, a(t)+\cos t\sin t\, b(t)+\cos^2 t\, c(t)\bigr)\,x_2^{},
\hfill
$
\\[2.25ex]
where $a\colon J\to {\mathbb R},\ b\colon J\to {\mathbb R},$ and 
\vspace{0.35ex}
$c\colon J\to {\mathbb R}$ are continuous functions on an interval $J\subset {\mathbb R},$  
can be reduced by the orthogonal transformation
\vspace{0.35ex}
$x_1^{}=\cos t\, y_1^{}-\sin t\, y_2^{},\, x_2^{}=\sin t\, y_1^{}+\cos t\, y_2^{}$ 
to the triangular differential system
\\[2.25ex]
\mbox{}\hfill
$
\dfrac{dy_1^{}}{d\;\!t}=
a(t)y_1^{}+b(t)y_2^{},
\qquad
\dfrac{dy_2^{}}{d\;\!t}=c(t)y_2^{}.
$
\hfill (2.19)
\\[2.25ex]
\indent
By Corollary 2.2, the scalar functions
\\[2ex]
\mbox{}\hfill
$
\displaystyle
F_1^{}\colon (t,y_1^{},y_2^{})\to\  y_2^{}
\exp\biggl({}-\int\limits_{t_0^{}}^{t} c(\zeta)\,d\zeta\biggr)
$
\ for all 
$
(t,y_1^{},y_2^{})\in J\times {\mathbb R}^2
\hfill
$ 
\\[0.75ex]
and
\\[0.75ex]
\mbox{}\hfill
$
\displaystyle
F_2^{}\colon (t,y_1^{},y_2^{})\to\ 
 y_1^{}\exp\biggl({}-\int\limits_{t_0^{}}^{t} a(\zeta)\,d\zeta\biggr) \, -\ 
\int\limits_{t_0^{}}^{t} b(\zeta)
\exp\biggl(\, \int\limits_{\zeta_0^{}}^{\zeta} (c(\theta)-a(\theta))d\theta\biggr)\, d\zeta\, F_1^{}(t,y_1^{},y_2^{})
\hfill
$
\\[2ex]
\mbox{}\hfill
for all 
$
(t,y_1^{},y_2^{})\in J\times {\mathbb R}^2,
\quad 
t_0^{},\, \zeta_0^{}\in J,
\hfill
$
\\[1.75ex]
are first integrals of the triangular differential system (2.19).
\vspace{0.5ex}

Using the inverse transformation
\vspace{0.35ex}
$
y_1^{}=\cos t\, x_1^{}+\sin t\, x_2^{},
\
y_2^{}={}-\sin t\, x_1^{}+\cos t\, x_2^{},
$
we get the integral basis of the linear differential system (2.18):
\\[2ex]
\mbox{}\hfill
$
\displaystyle
\widetilde{F}_1^{}\colon (t,x_1^{},x_2^{})\to \
(\cos t\, x_2^{}-\sin t\, x_1^{})\exp\biggl({}-
\int\limits_{t_0^{}}^{t}\! c(\zeta)\,d\zeta\biggr)
$ 
\ for all 
$
(t,x_1^{},x_2^{})\in J\times {\mathbb R}^2,
\hfill
$
\\[1.5ex]
\mbox{}\hfill
$
\displaystyle
\widetilde{F}_2^{}\colon (t,x_1^{},x_2^{})\to \
(\cos t\, x_1^{}+\sin t\, x_2^{})\exp\biggl({}-\int\limits_{t_0^{}}^{t}\! a(\zeta)\,d\zeta\biggr) \ -
\hfill
$
\\[0.5ex]
\mbox{}\hfill
$
\displaystyle
-\, \int\limits_{t_0^{}}^{t}\! b(\zeta)
\exp\biggl(\,\int\limits_{\zeta_0^{}}^{\zeta}\! (c(\theta)-a(\theta))\,d\theta\!\biggr)\,d\zeta\, 
\widetilde{F}_1^{}(t,x_1^{},x_2^{})
$
\ for all 
$
(t,x_1^{},x_2^{})\in J\times {\mathbb R}^2.
\hfill
$
\\

\newpage

\mbox{}
\\[-2.25ex]
\centerline{                               %{П. 3, Параграф 2}
\bf 2.3.
The Lappo-Danilevskii systems
}
\\[1.25ex]
\indent
{\bf 2.3.1. Linear homogeneous differential system}
\\[0.75ex]
\indent
Consider an nonautonomous linear homogeneous differential system
\\[1.5ex]
\mbox{}\hfill                                          % (2.20)
$
\displaystyle
\dfrac{dx}{dt}=A(t)\,x,
\ \quad 
x\in {\mathbb R}^n,
\qquad 
A\colon t\to 
\sum\limits_{j=1}^{m} \alpha_j^{}(t)A_j^{}
$ 
\ for all 
$
t\in J,
$
\hfill (2.20)
\\[1.5ex]
where continuous functions $\alpha_j^{}\colon J\to {\mathbb R}, \ j=1,\ldots, m,$ 
\vspace{0.35ex}
are lineary independent on an interval $J\subset {\mathbb R},$ 
and $A_j^{},\ j=1,\ldots, m,$ are a commuting family of real constant $n\times n$ matrices: 
\\[1.5ex]
\mbox{}\hfill
$
A_j^{}A_k^{}=A_k^{}A_j^{},
\quad 
j=1,\ldots, m, 
\ \
k=1,\ldots, m.
\hfill
$
\\[1.75ex]
\indent
The differential system (2.20) is a Lappo-Danilevskii system [8, pp. 34 -- 36], 
\vspace{0.25ex}
i.e. the co\-ef\-fi\-ci\-ent matrix $A$ of system (2.20) is commutative with its integral [15]: 
\\[1.5ex]
\mbox{}\hfill                                        
$
\displaystyle
A(t)\int\limits_{t_0}^{t}A(\tau)\,d\tau=\int\limits_{t_0}^{t}A(\tau)\,d\tau \, A(t)
$
\ for all 
$
t\in J,
\quad 
t_0^{}\in J.
\hfill
$
\\[1.5ex]
\indent
In this subsection we study the problem of building first integrals for the Lappo-Danilevskii 
differential system (2.20). Using the approaches of constructing 
first integrals for linear partial differential systems [25], 
we also solve the N.P. Erugin problem of the existence of autonomous first integrals [54, p. 469]
for the Lappo-Danilevskii system (2.20).
\vspace{0.5ex}

{\sl Partial integrals}.
A complex-valued linear homogeneous function
\\[1.5ex]
\mbox{}\hfill                                        % (2.21)
$
p \colon x\to \nu x
$ 
\ for all 
$
x\in {\mathbb R}^{n}, 
\quad
\nu\in {\mathbb C}^n,
$
\hfill (2.21)
\\[1.5ex]
is a partial integral of the system (2.20) if and only if 
\\[1.75ex]
\mbox{}\hfill                                          % (2.22)
$
{\frak A} p(x) =  \lambda(t)\;\! p(x)
$
\ for all 
$
(t, x)\in J\times {\mathbb R}^{n}, 
$
\hfill (2.22)
\\[1.75ex]
where the linear differential operator
\vspace{0.5ex}
${\frak A}(t,x)=\partial_{t}^{}+A(t)\;\!x\,\partial_{x}^{}$ for all $(t,x)\in J\times\R^{n},$
the function $\lambda\colon J\to {\mathbb C}.$
The identity (2.22) is equivalent to the linear homogeneous system
\\[1.25ex]
\mbox{}\hfill                                        
$
\displaystyle
\biggl(\, \sum\limits_{j=1}^{m} \alpha_j^{}(t)B_j^{} - \lambda(t) E \biggr)\, \nu = 0
$
\ for all 
$
t\in J,
$ 
\hfill (2.23)
\\[1.25ex]
where $E$ is the identity matrix of order $n,$ 
\vspace{0.35ex}
the matrices $B_j^{}$ are the transpose of the matrices
$A_j^{},\ j=1,\ldots, m,$ respectively.
\vspace{0.35ex}
Since $A_j^{},\ j=1,\ldots, m,$ are a commuting family of matrices,
we see that there exists a relation [60, pp. 203 -- 207] between 
\vspace{0.25ex}
common eigenvectors and eigenvalues of the matrices $B_j^{},\, j=1,\ldots, m.$
\vspace{0.5ex}

The following basic statements (Lemmas 2.3 and 2.4) are base for the method 
of building first integrals of the Lappo-Danilevskii differential system (2.20). 
The proof of Lemmas 2.3 and 2.4 is similar to that one in Lemmas 2.1 and 2.2. 
\vspace{0.5ex}

{\bf Lemma 2.3.}
{\it 
Suppose $\nu$ is a common real eigenvector of the matrices $B_{j}^{}$
corresponding to the eigenvalues $\lambda^j,\ j=1,\ldots, m.$ 
\vspace{0.35ex}
Then the function {\rm (2.21)}  is a linear partial integral of system {\rm (2.20)}, 
where the scalar function $\lambda\colon J\to {\mathbb R}$ in the identity {\rm (2.22)} is given by
\\[1.5ex]
\mbox{}\hfill
$
\displaystyle
\lambda\colon t\to\  
\sum\limits_{j=1}^{m} \lambda^j\;\!\alpha_j^{}(t)
$
\ for all 
$
t\in J.
\hfill
$
}
\\[1.5ex]
\indent
{\sl Proof}. 
If $\nu$ is a common real eigenvector of the matrices $B_{j}^{}$
corresponding to the eigen\-va\-lues $\lambda^j,\ j=1,\ldots, m,$ then 
$\nu$ is a common real eigenvector of the matrix
$
B\colon t\to \sum\limits_{j=1}^{m} \alpha_j^{}(t)\;\!B_j^{}
$
for all $t\in J$
corresponding to the eigenfunction
$
\lambda\colon t\to \sum\limits_{j=1}^{m} \lambda^j\alpha_j^{}(t)
$
for all $t\in J.$

Therefore $\nu$ is a solution to the functional system (2.23).
This yields that the identity (2.22) is satisfied. 
Consequently the linear function (2.21) is a partial integral 
of the Lappo-Danilevskii homogeneous differential system {\rm (2.20)}. \k
\vspace{0.5ex}

{\bf Lemma 2.4.}
{\it
Suppose 
\vspace{0.35ex}
$\nu\!={\stackrel{*}{\nu}}+\widetilde{\nu}\,i\ 
({\stackrel{*}{\nu}}\!=\!{\rm Re}\,\nu,\ \widetilde{\nu}\!=\!{\rm Im}\,\nu\ne 0)$ 
is a common complex eigenvec\-tor of the matrices $\!B_j^{}\!$ 
corresponding to the eigenvalues
\vspace{0.5ex}
$\lambda^j\!=\!{\stackrel{*}{\lambda}}{}^j+\widetilde{\lambda}{}^j\,i\ 
({\stackrel{*}{\lambda}}{}^j\!={\rm Re}\,\lambda^j,\ \widetilde{\lambda}{}^j\!={\rm Im}\,\lambda^j),$ 
$j=1,\ldots, m.$ 
Then the Lie derivatives of the scalar functions 
\\[1.35ex]
\mbox{}\hfill                                            
$
P\colon x\to 
({\stackrel{*}{\nu}}x)^2+(\widetilde{\nu}x)^2
$ 
\ for all 
$
x\in {\mathbb R}^n
$
\ \ and \ \
$
\displaystyle
\psi\colon x\to \
\arctan\dfrac{\widetilde{\nu}x}{{\stackrel{*}{\nu}}x}
$
\ for all 
$
x\in {\mathscr X}
\hfill
$
\\[1.35ex]
by virtue of system {\rm (2.20)} are equal to
\\[1.5ex]
\mbox{}\hfill                                            
$
\displaystyle
{\frak A}\;\!P(x)=\;\! 
2\sum\limits_{j=1}^{m}{\stackrel{*}{\lambda}}{}^j\,\alpha_j^{}(t)\;\!P(x)
$  
\ for all 
$
(t,x)\in J\times {\mathbb R}^n
\hfill
$
\\[0.5ex]
and
\\[0.5ex]
\mbox{}\hfill
$
\displaystyle
{\frak A}\;\!\psi(x)=\,
\sum\limits_{j=1}^{m}\;\!
\widetilde{\lambda}{}^j\,\alpha_j^{}(t)
$ 
\ for all 
$
(t,x)\in J\times {\mathscr X},
\hfill 
$
\\[1ex]
where ${\mathscr X}$ is a domain from the set 
\vspace{0.5ex}
$\bigl\{x\colon {\stackrel{*}{\nu}}x\ne 0\bigr\}\subset {\mathbb R}^n.$ 
}

{\sl Proof}. 
Formally using Lemma 2.3, we obtain 
the complex-valued function (2.21) is a par\-ti\-al integral 
of the differential system (2.20) and the following identity holds
\\[1.25ex]
\mbox{}\hfill
$
{\frak A}\;\!
\bigl({\stackrel{*}{\nu}}x+i\,\widetilde{\nu} x\bigr)=
\bigl({\stackrel{*}{\lambda}}(t)+i\,\widetilde{\lambda}(t)\bigr)
\bigl({\stackrel{*}{\nu}}x+i\,\widetilde{\nu} x\bigr)
$
\ for all 
$
(t,x)\in J\times {\mathbb R}^{n}.
\hfill
$
\\[1.35ex]
\indent
This complex identity is equivalent to the real system of identities
\\[1.5ex]
\mbox{}\hfill
$
{\frak A}\;\!{\stackrel{*}{\nu}} x =
{\stackrel{*}{\lambda}}(t)\,{\stackrel{*}{\nu}} x -
\widetilde{\lambda}(t)\,\widetilde{\nu} x,
\quad \
{\frak A}\;\!\widetilde{\nu} x=
{\stackrel{*}{\lambda}}(t)\,\widetilde{\nu} x+
\widetilde{\lambda}(t)\,{\stackrel{*}{\nu}} x
$
\ for all 
$
(t,x)\in J\times {\mathbb R}^{n},
\hfill
$
\\[1.5ex]
where the scalar functions  
\vspace{0.35ex}
$
{\stackrel{*}{\lambda}}\colon t\to
\sum\limits_{j=1}^{m}{\stackrel{*}{\lambda}}{}^j\;\!\alpha_j^{}(t)
$
for all $t\in J, \ 
\widetilde{\lambda}\colon t\to 
\sum\limits_{j=1}^{m}
\widetilde{\lambda}{}^j\;\!\alpha_j^{}(t)
$
for all $t\in J.$

Using this system of identities, we have
\\[1.35ex]
\mbox{}\hfill                        
$
{\frak A}\;\!P(x)=
{\frak A}\bigl(
({\stackrel{*}{\nu}} x)^2+(\widetilde{\nu} x)^2\bigr)=
2\;\!{\stackrel{*}{\nu}} x\,{\frak A}\;\!{\stackrel{*}{\nu}} x +
2\;\!\widetilde{\nu} x\,{\frak A}\;\!\widetilde{\nu} x=
2\;\!{\stackrel{*}{\nu}} x\,
\bigl(\;\!{\stackrel{*}{\lambda}}(t)\,{\stackrel{*}{\nu}} x -
\widetilde{\lambda}(t)\,\widetilde{\nu} x\bigr)  \ + 
\hfill                        
$
\\[1.75ex]
\mbox{}\hfill                        
$
+\ 2\;\!\widetilde{\nu} x\,
\bigl(\;\!{\stackrel{*}{\lambda}}(t)\,\widetilde{\nu} x+
\widetilde{\lambda}(t)\,{\stackrel{*}{\nu}} x\bigr)=
2\;\!{\stackrel{*}{\lambda}}(t)\;\!
\bigl(({\stackrel{*}{\nu}} x)^2+(\widetilde{\nu} x)^2\bigr)=
2\;\!{\stackrel{*}{\lambda}}(t)\,P(x)
$
\ for all 
$
(t,x)\in J\times {\mathbb R}^{n},
\hfill                        
$
\\[2ex]
\mbox{}\hfill                        
$
{\frak A}\,\psi(x)\, =\, 
\dfrac{
{\stackrel{*}{\nu}} x\;\! 
\bigl(\;\!{\stackrel{*}{\lambda}}(t)\,\widetilde{\nu} x+
\widetilde{\lambda}(t)\,{\stackrel{*}{\nu}} x\bigr)  -\,
\widetilde{\nu} x\;\! 
\bigl(\;\!{\stackrel{*}{\lambda}}(t)\,{\stackrel{*}{\nu}} x -
\widetilde{\lambda}(t)\,\widetilde{\nu} x\bigr)}
{({\stackrel{*}{\nu}} x)^2+(\widetilde{\nu} x)^2}=
\widetilde{\lambda}(t)
$
for all 
$
(t,x)\in J\times {\mathscr X}.
\ \k
\hfill                        
$
\\[2ex]
\indent
{\sl Nonautonomous first integrals}.
Using Theorems 2.5, 2.6, and 2.7, we can obtain first in\-teg\-rals of the 
the Lappo-Danilevskii homogeneous differential system (2.20). 
\vspace{0.35ex}

{\bf Theorem 2.5.}
{\it
Let the assumptions of Lemma {\rm 2.3} hold, then the scalar function  
\\[1ex]
\mbox{}\hfill                                            % (2.24)
$
\displaystyle
F\colon (t,x)\to\, 
\nu x\;\!
\exp\biggl({}-\int\limits_{t_0^{}}^{t}\sum\limits_{j=1}^m \lambda^j\,\alpha_j^{}(\tau)\,d\tau\biggr)
$
\ for all 
$
(t,x)\in J\times\R^n,
\quad
t_0^{}\in J,
$
\hfill {\rm (2.24)}
\\[1ex]
is a first integral on the domain $J\times\R^n$ of 
the Lappo-Danilevskii system {\rm(2.20)}.
}
\vspace{0.35ex}

{\sl Proof}.
From Lemma 2.3, we get
\\[1ex]
\mbox{}\hfill                                           
$
\displaystyle
{\frak A}F(t,x)=
\exp\biggl({}-\int\limits_{t_0^{}}^{t}
\sum\limits_{j=1}^m \lambda^j\,\alpha_j^{}(\tau)\;\! d\tau\!\biggr)\, {\frak A}\;\!\nu x \, +\, 
\nu x\, {\frak A}\exp\biggl({}-\int\limits_{t_0^{}}^{t}
\sum\limits_{j=1}^m \lambda^j\,\alpha_j^{}(\tau)\;\!d\tau\biggr) =
\hfill 
$
\\[1.25ex]
\mbox{}\hfill                                           
$
\displaystyle
=\;\!\sum\limits_{j=1}^{m} \lambda^j\alpha_j^{}(t)\, F(t,x) \,-\, 
F(t,x)\, \partial_{{}_{\scriptstyle t}}^{}\int\limits_{t_0^{}}^{t} 
\sum\limits_{j=1}^m \lambda^j\,\alpha_j^{}(\tau)\;\!d\tau =0
$ 
\ for all 
$
(t,x)\in J\times\R^n.\ \k
\hfill 
$
\\[2ex]
\indent
{\bf Example 2.6.}
%\marginpar{\No\ 09\,018}
The linear homogeneous differential system
\\[2.25ex]
\mbox{}\hfill                               % (2.25)              
$
\dfrac{dx_1^{}}{dt}=(t+2\sin t)\, x_1^{}+\sin t\, x_2^{},
\qquad
\dfrac{dx_2^{}}{dt}=\sin t\, x_1^{}+ t\, x_2^{}
$
\hfill (2.25)
\\[3.75ex]
has the coefficient matrix
\vspace{1.5ex}
$
A\colon t\to 
\left\|\!\!
\begin{array}{cc}
t+2\sin t & \sin t
\\[0.5ex]
\sin t & t
\end{array}
\!\!\right\|
$
such that 
$
A(t)=t\;\! A_1^{} +\sin t\;\! A_2^{}
$ 
for all $t\in {\mathbb R},$
where the constant matrices
$
A_1^{}=
\left\|\!\!
\begin{array}{cc}
1 & 0
\\[0.5ex]
0 & 1
\end{array}
\!\!\right\|
$
and
\vspace{0.75ex}
$
A_2^{}=
\left\|\!\!
\begin{array}{cc}
2 & 1
\\[0.5ex]
1 & 0
\end{array}
\!\!\right\|.
$
Since the matrices $A_1^{}$ and $A_2^{}$ are commutative, 
we see that the system (2.25) is a  Lappo-Danilevskii system. 
\vspace{0.75ex}

The matrices 
\vspace{0.5ex}
$
B_1^{}=A_1^{T}=A_1^{}
$
and
$
B_2^{}=A_2^{T}=A_2^{}
$
have the common real eigenvectors $\nu^{1}=\bigl(1-\sqrt{2}, 1\bigr)$ and 
$\nu^{2}=\bigl(1+\sqrt{2}, 1\bigr)$ corresponding to the eigenvalues
\vspace{0.35ex}
$\lambda_1^1=1,\ \lambda_1^2=1-\sqrt{2}$ and 
$\lambda_2^1=1,\ \lambda_2^2=1+\sqrt{2}\,,$ respectively.
\vspace{0.75ex}

Using the functions 
\vspace{0.5ex}
$\alpha_1^{}\colon t\to t$ for all $t\in {\mathbb R}$ and
$\alpha_2^{}\colon t\to \sin t$ for all $t\in {\mathbb R},$ 
we can build (by Theorem 2.5) the first integrals of 
the Lappo-Danilevskii system (2.25)
\\[1.5ex]
\mbox{}\hfill
$
F_1^{}\colon (t,x_1^{},x_2^{})\to\ 
\bigl((1-\sqrt{2}\,)\;\!x_1^{}+x_2^{}\bigr)\, 
\exp\Bigl({}-\dfrac{t^2}{2}+ (1-\sqrt{2}\,)\cos t\Bigr)
$
\ for all 
$
(t,x_1^{},x_2^{})\in {\mathbb R}^3
\hfill
$
\\[1.75ex]
and
\\[1.25ex]
\mbox{}\hfill
$
F_2^{}\colon (t,x_1^{},x_2^{})\to\ 
\bigl((1+\sqrt{2}\,)\;\!x_1^{}+x_2^{}\bigr)\, 
\exp\Bigl({}-\dfrac{t^2}{2}+ (1+\sqrt{2}\,)\cos t\Bigr)
$
\ for all 
$
(t,x_1^{},x_2^{})\in {\mathbb R}^3.
\hfill
$
\\[2ex]
\indent
The functionally independent first integrals $F_1^{}$ and $F_2^{}$ are 
\vspace{0.35ex}
an integral basis of the Lappo-Danilevskii differential system (2.25) on space ${\mathbb R}^3.$
\vspace{0.75ex}

{\bf Theorem 2.6.}
\vspace{0.25ex}
{\it
Let the assumptions of Lemma  {\rm 2.4} hold, then first integrals of 
the Lappo-Danilevskii differential  system {\rm (2.20)} are the scalar functions  
\\[1.5ex]
\mbox{}\hfill                                            
$
\displaystyle
F_1^{}\colon (t,x)\to\, 
\Bigl(\bigl({\stackrel{*}{\nu}}x\bigr)^2+\bigl(\widetilde{\nu}x\bigr)^2\Bigr)
\exp\biggl({}-2\int\limits_{t_0^{}}^{t}\, \sum\limits_{j=1}^{m}\;\!
{\stackrel{*}{\lambda}}{}^j\;\!\alpha_j^{}(\tau)\;\!d\tau\biggr)
$
\ for all 
$
(t,x)\in J\times {\mathbb R}^n
\hfill 
$
\\[0.75ex]
and
\\[1.5ex]
\mbox{}\hfill                                            
$
\displaystyle
F_2^{}\colon (t,x)\to \
\arctan\dfrac{\widetilde{\nu}x}{{\stackrel{*}{\nu}}x} \, -\, 
\int\limits_{t_0^{}}^{t} \,
\sum\limits_{j=1}^{m}\;\!
\widetilde{\lambda}{}^j\;\!\alpha_j^{}(\tau)\;\!d\tau
$
\ for all 
$
(t,x)\in  J\times {\mathscr X},
\hfill 
$
\\[1.5ex]
where 
\vspace{1ex}
$t_0^{}$ is a fixed point from the interval $J,$ and 
a domain ${\mathscr X}$ from the set $\bigl\{x\colon {\stackrel{*}{\nu}}x\ne 0\bigr\}.$
}

{\sl Proof}. 
Taking into account Lemma 2.4, we obtain
\\[2ex]
\mbox{}\hfill
$
\displaystyle
{\frak A}\;\!F_1^{}(t,x)=
\exp\biggl({}-2\int\limits_{t_0^{}}^{t}\, \sum\limits_{j=1}^{m}\;\!
{\stackrel{*}{\lambda}}{}^j\;\!\alpha_j^{}(\tau)\;\!d\tau\biggr)\, 
{\frak A}\;\!\bigl( ({\stackrel{*}{\nu}} x)^2+(\widetilde{\nu} x)^2\bigr) \ + 
\hfill
$
\\[2.25ex]
\mbox{}\hfill
$
\displaystyle
+\ \bigl( ({\stackrel{*}{\nu}} x)^2+(\widetilde{\nu} x)^2\bigr)
\ {\frak A}\exp\biggl({}-2\int\limits_{t_0^{}}^{t}\, \sum\limits_{j=1}^{m}\;\!
{\stackrel{*}{\lambda}}{}^j\;\!\alpha_j^{}(\tau)\;\!d\tau\biggr) \,=\ 
2\,\sum\limits_{j=1}^{m}{\stackrel{*}{\lambda}}{}^j\,\alpha_j^{}(t)\;\!F_1^{}(t,x)\ + 
\hfill
$
\\[2.25ex]
\mbox{}\hfill
$
\displaystyle
+\ \bigl( ({\stackrel{*}{\nu}} x)^2+(\widetilde{\nu} x)^2\bigr)
\,\partial_{{}_{\scriptstyle t}}^{}\exp\biggl({}-2\int\limits_{t_0^{}}^{t}\, \sum\limits_{j=1}^{m}\;\!
{\stackrel{*}{\lambda}}{}^j\;\!\alpha_j^{}(\tau)\;\!d\tau\biggr) =0
$ 
\ for all 
$
(t,x)\in J\times {\mathbb R}^n,
\hfill
$
\\[2.25ex]
\mbox{}\hfill                                            
$
\displaystyle
{\frak A}\;\!F_2^{}(t,x)= \,
{\frak A}\;\!\arctan\dfrac{\widetilde{\nu} x}{{\stackrel{*}{\nu}} x}\, - \,
{\frak A}\;\!\int\limits_{t_0^{}}^{t} \,
\sum\limits_{j=1}^{m}
\widetilde{\lambda}{}^j\;\!\alpha_j^{}(\tau)\;\!d\tau \, = 
\hfill
$
\\[1.75ex]
\mbox{}\hfill                                            
$
\displaystyle
=\ \sum\limits_{j=1}^{m}
\widetilde{\lambda}{}^j\;\!\alpha_j^{}(t)
\, -\,
\partial_{t}^{}\;\!\int\limits_{t_0^{}}^{t} \,
\sum\limits_{j=1}^{m}
\widetilde{\lambda}{}^j\;\!\alpha_j^{}(\tau)\;\!d\tau=0
$
\ for all 
$
(t,x)\in J\times {\mathscr X}.
\hfill
$
\\[1.75ex]
\indent
Thus the functions $F_1^{}$ and $F_2^{}$ are first integrals of system (2.20). \k 
\vspace{1ex}

{\bf Example 2.7.}
%\marginpar{\No\ 09\,017}
Consider the second order Lappo-Danilevskii differential  system
\\[2ex]
\mbox{}\hfill                                             % (2.26)
$
\dfrac{dx_1^{}}{dt}=a\cos\omega t\, x_1^{}+b\sin\omega t\, x_2^{},
\qquad
\dfrac{dx_2^{}}{dt}={}-d\sin\omega t\, x_1^{}+a\cos\omega t\, x_2^{},
$
\hfill (2.26)
\\[2.25ex]
where $a,\ b,$ and $d$ are positive constants, and
$\omega$ is an non-zero real number.
\vspace{1ex}

The matrices 
$
B_1^{}=
\left\|\!\!
\begin{array}{cc}
a & 0
\\[0.5ex]
0 & a
\end{array}
\!\!\right\|
$
and
\vspace{1ex}
$
B_2^{}=
\left\|\!\!
\begin{array}{cc}
0 & {}-d
\\[0.5ex]
b & 0
\end{array}
\!\!\right\|
$
have the common complex eigenvectors
\vspace{0.5ex}
$
\nu^{1}=\bigl(\sqrt{d}, \sqrt{b}\,i\bigr)
$
and
$
\nu^{2}=\bigl(\sqrt{d}, {}-\sqrt{b}\,i\bigr)
$
associated with the eigenvalues
\vspace{0.35ex}
$\lambda^1_1=a\ \lambda_1^2={}-\sqrt{b\,d}\,i$ and 
$\lambda_2^1=a,\ \lambda_2^2=\sqrt{b\,d}\,i,$ respectively.
\vspace{0.35ex}

Using the numbers
\vspace{0.35ex}
${\stackrel{*}{\lambda}}{}^1_1=a,\ {\stackrel{*}{\lambda}}{}^2_1=0,\ 
\widetilde{\lambda}{}^1_1=0,\ \widetilde{\lambda}{}^2_1={}-\sqrt{b\,d},$ and 
the scalar functions 
$\alpha_1^{}\colon t\to \cos\omega t,\ \, 
\alpha_2^{}\colon t\to \sin\omega t$ for all $t\in {\mathbb R},$
\vspace{0.25ex}
we can find (by Theorem 2.6) the basis of first integrals 
for the Lappo-Danilevskii differential system (2.26)
\\[2ex]
\mbox{}\hfill
$
F_1^{}\colon (t,x_1^{},x_2^{})\to\  
\bigl(dx_1^2+bx_2^2\bigr) 
\exp\Bigl({}-\dfrac{2a}{\omega}\ \sin\omega t\Bigr)
$
\ for all 
$
(t,x_1^{},x_2^{})\in {\mathbb R}^3,
\hfill
$
\\[2.75ex]
\mbox{}\hfill
$
F_2^{}\colon (t,x_1^{},x_2^{})\to\ 
\arctan\dfrac{\sqrt{b}\,x_2^{}}{\sqrt{d}\,x_1^{}}\, -\, 
\dfrac{\sqrt{b\,d}}{\omega}\ \cos\omega t
$
\ for all 
$
(t,x_1^{},x_2^{})\in \Omega
\hfill
$
\\[2ex]
on any domain $\Omega$ from the set 
$\{(t,x_1^{},x_2^{})\colon x_1^{}\ne 0\}\subset {\mathbb R}^3.$ 
\vspace{1.25ex}

{\bf Example 2.8.}
The linear homogeneous differential system 
\\[2ex]
\mbox{}\hfill                            
$
\dfrac{dx_1^{}}{dt}=3(\tanh t+t^2)\, x_1^{}+\tanh t\, (x_2^{}-x_3^{}),
\qquad
\dfrac{dx_2^{}}{dt}={}-2\tanh t\, x_1^{}+3t^2\, x_2^{},
\hfill
$
\\[0.75ex]
\mbox{}\hfill (2.27)
\\[0.25ex]
\mbox{}\hfill
$
\dfrac{dx_3^{}}{dt}=\tanh t\, (x_1^{}+x_2^{})+(\tanh t+3t^2)\, x_3^{}
\hfill
$
\\[2.25ex]
is the Lappo-Danilevskii system of form (2.20) with the constant matrices
\\[2ex]
\mbox{}\hfill
$
B_1^{}=
\left\|\!\!
\begin{array}{ccc}
1 & 0 & 0
\\[0.5ex]
0 & 1 & 0
\\[0.5ex]
0 & 0 & 1
\end{array}
\!\!\right\|
$
\ \ and \ \ 
$
B_2^{}=
\left\|\!\!
\begin{array}{rrr}
{}-2 & 2 & {}-1
\\[0.5ex]
{}-1 & 1 & {}-1
\\[0.5ex]
1 & 0 & 0
\end{array}
\!\!\right\|,
\hfill
$
\\[2ex]
and the scalar functions
\vspace{0.75ex}
$\alpha_1^{}\colon t \to \tanh t+3t^2$ for all $t\in {\mathbb R}, \
\alpha_2^{}\colon t\to {}-\tanh t$ for all $t\in {\mathbb R}.$ 

Using the linearly independent common eigenvectors
\\[1.5ex]
\mbox{}\hfill
$
\nu^1=(1, 0,{}-1),
\quad 
\nu^2=(i, i, 1),
\quad 
\nu^3=({}-i,{}-i, 1)
\hfill
$ 
\\[1.5ex]
of the matrices $B_1^{},\ B_2^{}$ and the corresponding eigenvalues 
\\[1.5ex]
\mbox{}\hfill
$
\lambda^1_1=1,\ \ \lambda_1^2={}-1,
\qquad
\lambda_2^1=1,\ \ \lambda_2^2=i,
\qquad
\lambda_3^1=1, \ \ \lambda^2_3={}-i,
\hfill
$
\\[1.5ex]
we can construct the integral basis of the Lappo-Danilevskii differential system (2.27)  
\\[2ex]
\mbox{}\hfill
$
F_1^{}\colon (t,x)\to\  
\dfrac{1}{e^{\;\!t^{{}^{3}}}\!\cosh^2 t}\, (x_1^{}-x_3^{})
$
\ for all 
$
(t,x)\in {\mathbb R}^4
$ \
(by Theorem 2.5),
\hfill\mbox{}
\\[2.5ex]
\mbox{}\hfill
$
F_2^{}\colon (t,x)\to\  
\dfrac{1}{e^{\;\!2t^{{}^{3}}}\!\cosh^2 t}\, 
\bigl(x_3^2+(x_1^{}+x_2^{})^2\bigr)
$
\ for all 
$
(t,x)\in {\mathbb R}^4
$ \
(by Theorem 2.6),
\hfill\mbox{}
\\[2.5ex]
\mbox{}\hfill
$
F_3^{}\colon (t,x)\to\ 
 \arctan\dfrac{x_1^{}+x_2^{}}{x_3^{}}\, +\, 
\ln\cosh t
$ 
\ for all 
$
(t,x)\in \Omega
$
\ (by Theorem 2.6)
\hfill\mbox{}
\\[2ex]
on any domain $\Omega$ from the set $\{(t,x)\colon x_3^{}\ne 0\}\subset {\mathbb R}^4.$ 
\vspace{1.25ex}

{\bf Lemma 2.5.}                                           
{\it
Let the following conditions hold}:
\vspace{0.35ex}

(i) 
{\it
$\nu^0$ is a common eigenvector of $B_j^{}$ corresponding 
\vspace{0.35ex}
to the eigenvalues $\lambda^j,\ j=1,\ldots, m;$
}

(ii) 
{\it
$\nu^{\theta},\, \theta=1,\ldots, s-1,$ are
\vspace{0.25ex}
generalized eigenvectors of the matrix $B_{\zeta}^{}$
corresponding to the eigenvalue $\lambda^{\zeta}$ with
elementary divisor of multiplicity $s\geq 2;$ 
}
\vspace{0.35ex}

(iii) 
{\it the linear differential system $\dfrac{dx}{dt}=A_{\zeta}^{}\;\!x$ 
has no the first integrals of the form
\\[2ex]
\mbox{}\hfill                                                             % (2.28)
$
\displaystyle
F_{{}_{\scriptstyle j\theta}}^{\,\zeta}\colon x\to\ 
{\frak a}_j^{}\;\!\Psi_{\theta}^{\zeta}(x) 
$
\ for all 
$
x\in {\mathscr X}, 
\quad
j=1,\ldots, m, \ j\ne\zeta,
\quad 
\theta=1,\ldots, s-1.
$
\hfill {\rm (2.28)}
\\[1.75ex]
Then, we claim that  
\\[1.75ex]
\mbox{}\hfill                                                           % (2.29)
$
{\frak a}_{\zeta}^{}\Psi_{1}^{\zeta}(x)=1 
$
\ for all 
$
x\in {\mathscr X},
\quad
{\frak a}_{\zeta}^{}\Psi_{\theta}^{\zeta}(x)=0
$
\ for all 
$
x\in {\mathscr X},
\ \, 
\theta=2,\ldots, s-1,
$
\hfill {\rm (2.29)}
\\[1.25ex]
and
\\[1.25ex]
\mbox{}\hfill                                                     % (2.30)
$
{\frak a}_{j}^{}\Psi_{\theta}^{\zeta}(x)=\mu_{\theta}^{j\zeta}={\rm const}
$
\ for all 
$
x\in {\mathscr X}, 
\ \ j=1,\ldots, m,\ j\ne \zeta,
\ \ \theta=1,\ldots, s-1,
$
\hfill {\rm (2.30)}
\\[2ex]
where the linear differential operators
\vspace{0.25ex}
${\frak a}_{j}^{}(x)=A_j^{}x\;\!\partial_{x}$ for all $x\in {\mathbb R}^n,$ the scalar
functions $\Psi_{\theta}^{\zeta}\colon {\mathscr X}\to {\mathbb R}$ 
are the solution to the functional system}
\\[1.35ex]
\mbox{}\hfill                                  % (2.31)                                
$
\nu^{\,\theta}x=
{\displaystyle \sum\limits_{\rho=1}^{\theta} }
\binom{\theta-1}{\rho-1}\Psi_{\rho}^{\zeta}(x)\, \nu^{\theta-\rho}x
$
\ for all 
$
x\in {\mathscr X},
\ \ 
\theta=1,\ldots, s-1,
\quad
{\mathscr X}\subset \{x\colon \nu^0x\ne 0\}.
$
\hfill {\rm(2.31)}
\\[1.5ex]
\indent
{\sl Proof.}  
\vspace{0.35ex}
From the proof of Theorem 1.8 it follows that the identities (2.29) 
for the linear differential system $\dfrac{dx}{dt}=A_{\zeta}^{}x$ are true.
Since the linear differential operators of first order ${\frak a}_j^{},\ j=1,\ldots, m,$ are 
commutative and the system $\dfrac{dx}{dt}=A_{\zeta}^{}x$ hasn't 
\vspace{0.35ex}
the first integrals (2.28), we see that 
the system of identities (2.30) is satisfied.

Thus there exist the functionally independent functions 
\vspace{0.5ex}
$
\Psi_{\theta}^{\zeta}\colon {\mathscr X}\to {\mathbb R},
\ \theta=1,\ldots, s-1
$ 
(the solution to the system (2.31)) and these functions satisfies (2.29) and (2.30).  \k
\vspace{0.75ex}

{\bf Theorem 2.7.}                                           
{\it
Under the conditions of Lemma {\rm 2.5}, we get 
first integrals of the Lappo-Da\-ni\-lev\-s\-kii differential system {\rm (2.20)} are the scalar functions
\\[1.5ex]
\mbox{}\hfill                        % (2.32)                  
$
\displaystyle
F_{\theta}^{}\colon (t, x)\to\ 
\Psi_{\theta}^{\zeta}(x) \, -\,
\int\limits_{t_0^{}}^{t}\, \sum\limits_{j=1}^{m}\mu_{\theta }^{j\zeta}\,\alpha_j^{}(\tau)\,d\tau 
$
\ for all 
$
(t,x)\in J\times {\mathscr X}, 
\quad
\theta=1,\ldots, s-1,
$
\hfill {\rm (2.32)}
\\[1.5ex]
where  
\vspace{0.75ex}
$t_0^{}$ is a fixed point from the interval $J,$ a domain 
${\mathscr X}$ from the set $\bigl\{x\colon \nu^0 x\ne 0\bigr\}\subset {\mathbb R}^n.$
}

{\sl Proof}.
From the identities (2.29) and (2.30), on the domain $J\times {\mathscr X}$ we get the following 
\\[1.5ex]
\mbox{}\hfill                                          
$
\displaystyle
{\frak A}\;\! F_{\theta}^{}(t, x)=
{}- \partial_{{}_{\scriptstyle t}}^{}
\int\limits_{t_0^{}}^{t}\, \sum\limits_{j=1}^{m}\mu_{\theta }^{j\zeta}\,\alpha_j^{}(\tau)\,d\tau \, +\, 
\sum\limits_{j=1}^{m}\alpha_j^{}(t)\,{\frak a}_j^{}\Psi_{\theta}^{\zeta}(x) = 0,
\quad
\theta=1,\ldots, s-1.
\hfill 
$
\\[1.5ex]
\indent
Therefore the functions (2.32) are functionally independent first integrals
on the domain $J\times {\mathscr X}$ of the 
Lappo-Da\-ni\-lev\-s\-kii differential system (2.20). \k
\vspace{0.75ex}

{\bf Example 2.9.}
%\marginpar{\No\ 09\,016}
Lappo-Da\-ni\-lev\-s\-kii differential systems 
with non-diagonal coefficient mat\-ri\-ces of the second order
are the linear differential systems of the form [58]
\\[2ex]
\mbox{}\hfill                                               % (2.33)
$
\dfrac{dx_1^{}}{dt} =\bigl(\alpha_1^{}(t)+b_1^{}\,\alpha_2^{}(t)\bigr)\,x_1^{}+\alpha_2^{}(t)\,x_2^{},
\qquad 
\dfrac{dx_2^{}}{dt} = b_2^{}\,\alpha_2^{}(t)\,x_1^{} +\alpha_1^{}(t)\,x_2^{},
$
\hfill {\rm (2.33)}
\\[2.25ex]
where the functions $\alpha_1^{}\colon J\to {\mathbb R}$ and 
$\alpha_2^{}\colon J\to {\mathbb R}$ are continuous,
\vspace{0.35ex}
$b_1^{}$ and $b_2^{}$ are real numbers.

Consider the number $D=b_1^2+4b_2^{}.$ 
We have three possible cases.
\vspace{0.35ex}

Let $D>0.$ Then, using the common real eigenvectors
\vspace{0.35ex}
$\nu^{1}=(1, {}-\lambda_2^2), \ \nu^{2}=(1, {}-\lambda_1^2)$ 
and the corresponding eigenvalues 
\vspace{0.35ex}
$\lambda_1^1=1,\ \lambda_{1}^{2}=\bigl(b_1-\sqrt{D}\,\bigr)/2,\
\lambda_2^1=1,\ \lambda_{2}^{2}=\bigl(b_1+\sqrt{D}\,\bigr)/2,$
we can build (by Theorem 2.5) the integral basis of system (2.33)
\\[1.25ex]
\mbox{}\hfill
$
\displaystyle
F_k^{}\colon (t,x_1^{},x_2^{})\to\, \bigl(x_1^{}-\lambda_{3-k}^2\,x_2^{}\bigr) 
\exp\biggl({}-\int\limits_{t_0^{}}^{t}\bigl(\alpha_1^{}(\tau)+\lambda_k^2\,\alpha_2^{}(\tau)\bigr)\,d\tau\biggr)
\hfill
$
\\[1.5ex]
\mbox{}\hfill
for all 
$
(t,x_1,x_2)\in J\times {\mathbb R}^2,
\quad 
t_0^{}\in J,
\quad 
k=1,2.
\hfill
$
\\[1.5ex]
\indent
Let $D<0.$ Then, using the common complex eigenvector
$\nu^{1}=(\lambda^2_1, 1)$ 
and the corresponding eigenvalues 
\vspace{0.5ex}
$\lambda_1^1=1,\ \lambda^{2}_{1}={\stackrel{*}{\lambda}}-\widetilde{\lambda}\,i,$
where ${\stackrel{*}{\lambda}}=b_1^{}/2, \ \widetilde{\lambda}=\sqrt{{}-D}/2,$
we can construct (by Theorem 2.6) the basis of first integrals for system (2.33)
\\[1.5ex]
\mbox{}\hfill
$
\displaystyle
F_1^{}\colon (t,x_1^{},x_2^{})\to 
\bigl(({\stackrel{*}{\lambda}}\;\!x_1^{}+x_2^{})^2+
(\widetilde{\lambda}\;\!x_1^{})^2\bigr)
\exp\biggl(\!\!{}-2\int\limits_{t_0^{}}^t\!
\bigl(\alpha_1^{}(\tau)+{\stackrel{*}{\lambda}}\,\alpha_2^{}(\tau)\bigr)\,d\tau\biggr)
\hfill
$
\\[1.5ex]
\mbox{}\hfill
for all 
$
(t,x_1^{},x_2^{})\in J\times {\mathbb R}^2,
\quad 
t_0^{}\in J,
\hfill
$
\\[2.25ex]
\mbox{}\hfill
$
\displaystyle
F_2\colon (t,x_1^{},x_2^{})\to \
\arctan\dfrac{\widetilde{\lambda}\,x_1^{}}{{\stackrel{*}{\lambda}}\,x_1^{}+x_2^{}}\ +\
\int\limits_{t_0^{}}^t \widetilde{\lambda}\,\alpha_2^{}(\tau)\,d\tau
$
\ for all 
$
(t,x_1^{},x_2^{})\in J\times {\mathscr X},
\hfill
$
\\[1.5ex]
where ${\mathscr X}$ is a domain from the set 
$\{(x_1,x_2)\colon {\stackrel{*}{\lambda}}\,x_1^{}+x_2^{}\ne 0\}\subset {\mathbb R}^2.$
\vspace{0.5ex}

Let $D=0.$ Then, using the common real eigenvector
\vspace{0.35ex}
$\nu^{0}=(\lambda^2_1, 1),$ the 1-st order real generalized eigenvector $\nu^{1}=(1, 0),$ 
\vspace{0.35ex}
and the corresponding eigenvalue $\lambda^{2}_{1}=b_1^{}/2$ with elementary divisor of 
multiplicity $s=2,$ we can find (by Theorems 2.5 and 2.7) the functionally independent 
first integrals of system (2.33)
\\[1.5ex]
\mbox{}\hfill
$
\displaystyle
F_1^{}\colon (t,x_1^{},x_2^{})\to \
\bigl(\lambda^{2}_1\,x_1^{}+x_2^{}\bigr) 
\exp\biggl({}-\int\limits_{t_0^{}}^{t}
\bigl(\alpha_1^{}(\tau)+\lambda^2_1\,\alpha_2^{}(\tau)\bigr)\,d\tau\biggr)
$ 
\ for all 
$
(t,x_1^{},x_2^{})\in J\times {\mathbb R}^2,
\hfill
$
\\[1.5ex]
\mbox{}\hfill
$
\displaystyle
F_2^{}\colon (t,x_1^{},x_2^{})\to \
\dfrac{x_1^{}}{\lambda^{2}_1\,x_1^{}+x_2^{}} -
\int\limits_{t_0^{}}^{t}\alpha_2^{}(\tau)\,d\tau
$
\ for all 
$
(t,x_1^{},x_2^{})\in J\times {\mathscr X},
\quad 
t_0^{}\in J,
\hfill
$
\\[1.5ex]
where ${\mathscr X}$ is a domain from the set 
$\{(x_1^{},x_2^{})\colon \lambda^2_1\,x_1^{}+x_2^{}\ne 0\}\subset {\mathbb R}^2.$
\vspace{1ex}

{\it In the complex case}, from the complex-valued first integrals (2.32) of 
the Lappo-Da\-ni\-lev\-s\-kii differential system (2.20), 
we obtain the real first integrals 
\\[1ex]
\mbox{}\hfill                                         
$
\displaystyle
F_{\theta}^{\,1}\colon (t, x)\to\ 
{\rm Re}\,\Psi_{\theta}^{\zeta}(x) \, -\, 
\int\limits_{t_0^{}}^{t}\,\sum\limits_{j=1}^{m}{\rm Re}\,\mu_{\theta}^{j\zeta}\,\alpha_j^{}(\tau)\;\!d\tau 
$
\ for all 
$(t,x)\in J\times {\mathscr X}, 
\quad
\theta=1,\ldots, s-1,
\hfill                                         
$
\\[1ex] 
\mbox{}\hfill                                         
$
\displaystyle
F_{\theta}^{\,2}\colon (t, x)\to \
{\rm Im}\,\Psi_{\theta}^{\zeta}(x)\, -\,
\int\limits_{t_0^{}}^{t}\, \sum\limits_{j=1}^{m}{\rm Im}\,\mu_{\theta}^{j\zeta}\,\alpha_j^{}(\tau)\;\!d\tau 
$
\ for all 
$
(t,x)\in J\times {\mathscr X}, 
\quad
\theta=1,\ldots, s-1.
\hfill 
$
\\[2ex]
\indent
{\sl Autonomous first integrals}.
The base of building autonomous first integrals for 
the Lappo-Da\-ni\-lev\-s\-kii differential system (2.20) is
the following basic proposition. 
\vspace{0.5ex}

{\bf Lemma 2.6.}
{\it
A function $F\colon {\mathscr X}\to {\mathbb R}$ is an autonomous first integral of
the Lappo-Da\-ni\-lev\-s\-kii differential system {\rm (2.20)} if and only if
this function is a first integral of the linear homogeneous system of partial 
differential equations}
\\[1.5ex]
\mbox{}\hfill                                         % (2.34)
$
{\frak a}_j^{}(x)\;\!y=0,
\quad 
j=1,\ldots, m,
$
\hfill (2.34)
\\[1.5ex]
{\it 
where the linear differential operators of first order}
\vspace{0.75ex}
${\frak a}_{j}^{}(x)=A_j^{}\;\!x\,\partial_{x}^{}$ for all $x\in {\mathbb R}^n.$ 

{\sl Proof} [54, p. 470]. 
A function $F\colon {\mathscr X}\to {\mathbb R}$ is an autonomous first integral of
the Lappo-Da\-ni\-lev\-s\-kii differential system {\rm (2.20)} if and only if 
the following identity hold
\\[1.5ex]
\mbox{}\hfill                                        
$
\displaystyle
\sum\limits_{j=1}^{m} \alpha_j^{}(t)\;\!A_j^{}\;\!x\,\partial_x^{} F(x)=0
$
\ for all 
$
(t,x)\in J\times {\mathbb R}^n.
\hfill 
$
\\[1.5ex]
\indent
Since the functions $\alpha_j^{}\colon J\to {\mathbb R}, \ j=1,\ldots, m,$ 
are linearly independent on the interval $J,$
we see that this identity is equivalent to the system of identities
\\[1.5ex]
\mbox{}\hfill                                        
$
\displaystyle
A_j^{}\;\!x\,\partial_x^{} F(x)=0
$
\ for all 
$
x\in {\mathbb R}^n,
\quad 
j=1,\ldots, m.
\hfill 
$
\\[1.5ex]
\indent
This yields that the function $F\colon {\mathscr X}\to {\mathbb R}$ 
is a first integral of system (2.34). \k
\vspace{0.5ex}

By Theorems 2.8, 2.9, 2.10, and 2.11, we can construct autonomous first integrals of the 
Lappo-Da\-ni\-lev\-s\-kii homogeneous differential system (2.20). 
\vspace{0.5ex}

{\bf Theorem 2.8.}
{\it
Suppose $\nu^{k}$ are common real eigenvectors of the matrices $B_{j}^{}$
cor\-res\-pon\-ding to the eigenvalues
\vspace{0.35ex}
$\lambda_k^j, \ j=1,\ldots, m,\ k=1,\ldots, m+1.$
Then the Lappo-Da\-ni\-lev\-s\-kii dif\-fe\-ren\-ti\-al system {\rm (2.20)} has 
the autonomous first integral}
\\[1.5ex]
\mbox{}\hfill                                            %(2.35)
$
\displaystyle
F\colon  x\to\ 
\prod\limits_{k=1}^{m+1}\,\bigl|\nu^{k} x\bigr|^{h_{k}^{}}
$
\ for all 
$
x\in {\mathscr X},
\quad
{\mathscr X}\subset {\rm D}F,
$
\hfill (2.35)
\\[1.75ex]
{\it
where the real numbers $h_{k}^{},\ k=1,\ldots, m+1,$ are 
an nontrivial solution to the system}
\\[1.5ex]
\mbox{}\hfill
$
\displaystyle
\sum\limits_{k=1}^{m+1}\,\lambda_{k}^{j}\;\!h_{k}^{} = 0, 
\quad
j=1,\ldots, m.
\hfill
$
\\[1.5ex]
\indent
{\sl Proof}.
First note that the Lappo-Da\-ni\-lev\-s\-kii system {\rm (2.20)} is 
induced the linear differential systems
$\dfrac{dx}{dt}=A_j^{}\;\!x$ with the operators
${\frak a}_{j}^{}(x)=A_j^{}\;\!x\,\partial_{x}^{}$ for all $x\in {\mathbb R}^n,\ j=1,\ldots, m.$ 
\vspace{0.75ex}

If $\nu^{k}$ are common real eigenvectors of the matrices $B_{j}^{}$
cor\-res\-pon\-ding to the eigenvalues
\vspace{0.75ex}
$\lambda_{k}^{j},\ j=1,\ldots, m,\ k=1,\ldots, m+1,$ 
then 
the linear homogeneous functions (by Lemma 1.1) 
$
p_k^{}\colon x\to\nu^{k}x
$  
for all 
$
x\in {\mathbb R}^{n},
\ k=1,\ldots, m+1,
$
are partial integrals of the linear autonomous differential systems
$\dfrac{dx}{dt}=A_j^{}\;\!x,\ j=1,\ldots, m,$ and the following system of identities hold
\\[2ex]
\mbox{}\hfill                                            % (2.36)
$
{\frak a}_{j}^{}\;\!
p_{k}^{}(x) =\lambda_{k}^{j}\;\!p_{k}^{}(x) 
$
\ for all 
$
x\in {\mathbb R}^{n},
\quad 
j=1,\ldots, m,
\ \ 
k=1,\ldots, m+1.
$
\hfill (2.36)
\\[2ex]
\indent
We obviously have 
\\[2ex]
\mbox{}\hfill
$
\displaystyle
{\frak a}_{j}^{}\;\! F(x)=
\prod\limits_{k=1}^{m+1}
\bigl|p_{k}^{}(x)\bigr|^{h_{k}^{}-1}
\sum\limits_{k=1}^{m+1}\,{\rm sgn}\,p_{k}^{}(x)\;\!
h_{k}^{}
\prod\limits_{l=1,l\ne k}^{m+1}\bigl|p_{l}^{}(x)\bigr|\,
{\frak a}_{j}^{}\;\!p_{k}^{}(x)
$
for all 
$
x \in  {\mathscr X},
\
j\!=\!1,\ldots, m.
\hfill
$
\\[2ex]
\indent
Now taking into account the identities (2.36), we obtain
\\[2ex]
\mbox{}\hfill
$
\displaystyle
{\frak a}_{j}^{}\;\!F(x) =
\sum\limits_{k=1}^{m+1}\,\lambda_{k}^{j}\;\! h_{k}^{}\,F(x)
$
\ for all 
$
x\in {\mathscr X},
\quad 
j=1,\ldots, m.
\hfill
$
\\[2ex]
\indent
If the real numbers $h_{k}^{},\ k=1,\ldots, m+1,$ 
are an nontrivial solution to the linear homogeneous system
$
\sum\limits_{k=1}^{m+1}\,\lambda_{k}^{j}\;\!h_{k}^{} = 0, \
j=1,\ldots, m,
$
then the scalar function (2.35) is a first integral for 
the linear homogeneous system 
of partial differential equations (2.34).
\vspace{0.25ex}

By Lemma 2.6, the function (2.35) is an autonomous first integral
\vspace{0.75ex}
of system {\rm (2.20)}. \k

{\bf Example 2.10.}
Consider the third-order Lappo-Da\-ni\-lev\-s\-kii dif\-fe\-ren\-ti\-al system
\\[1.75ex]
\mbox{}\hfill                            
$
\dfrac{dx_1^{}}{dt}=e^t\cos e^t\, x_1^{}+2\sinh t\, (x_2^{}+x_3^{}),
\qquad
\dfrac{dx_2^{}}{dt}=\sinh t\, x_1^{}+e^t\cos e^t\, x_3^{},
\hfill                            
$
\\[0.35ex]
\mbox{}\hfill (2.37)                           
\\[0.1ex]
\mbox{}\hfill                            
$
\dfrac{dx_3^{}}{dt}=\sinh t\, x_1^{}+e^t\cos e^t\, x_2^{}.
\hfill 
$
\\[2.5ex]
\indent
The matrices
$
B_1^{}=
\left\|\!\!
\begin{array}{ccc}
0 & 1 & 1
\\[0.35ex]
2 & 0 & 0
\\[0.35ex]
2 & 0 & 0
\end{array}
\!\!\right\|
$
\vspace{1ex}
and 
$
B_2^{}=
\left\|\!\!
\begin{array}{ccc}
1 & 0 & 0
\\[0.35ex]
0 & 0 & 1
\\[0.35ex]
0 & 1 & 0
\end{array}
\!\!\right\|
$
have the common real eigenvectors
$\nu^1=(1, {}-1,{}-1),\ \nu^2=(0, 1, {}-1),$ and $\nu^3=(1, 1, 1)$ 
\vspace{0.5ex}
corresponding to the real eigenvalues
$
\lambda^1_1={}-2,\ \lambda_1^2=1,\
\lambda_2^1=0,\ \lambda_2^2={}-1,
$
and
$  
\lambda_3^1=2, \ \lambda^2_3=1,
$ 
respectively.
\vspace{0.5ex}

From the linear homogeneous system 
\vspace{0.35ex}
${}-2h_1^{}+2h_3^{}=0,\ h_1^{}-h_2^{}+h_3^{}=0,$ we have, for example, 
$h_1^{}=1,\ h_2^{}=2,\ h_3^{}=1.$
By Theorem 2.8, the scalar function
\\[2ex]
\mbox{}\hfill
$
F\colon (x_1^{},x_2^{},x_3^{})\to\  
(x^2_1-(x_2^{}-x_3^{})^2)(x_2^{}-x_3^{})^2
$
\ for all 
$
(x_1^{},x_2^{},x_3^{})\in {\mathbb R}^3
\hfill
$ 
\\[2ex]
is an autonomous first integral on space ${\mathbb R}^3$ of 
the Lappo-Da\-ni\-lev\-s\-kii system (2.37).
\vspace{0.35ex}

By Theorem 2.5, using the functions  
\vspace{0.75ex}
$\alpha_1^{}\colon t \to \sinh t$ and
$\alpha_2^{}\colon t\to e^t\cos e^t$ for all $t\in {\mathbb R},$ 
we can build the nonautonomous first integrals of system (2.37)  
\\[2ex]
\mbox{}\hfill
$
F_1^{}\colon (t,x_1^{},x_2^{},x_3^{})\to\  
e^{{}^{\;\!\scriptstyle 2\cosh t-\sin e^t}}
(x_1^{}-x_2^{}-x_3^{}),
\quad
F_2^{}\colon (t,x_1^{},x_2^{},x_3^{})\to\  
e^{{}^{\;\!\scriptstyle \sin e^t}}(x_2^{}-x_3^{}),
\hfill\mbox{}
$ 
\\[1.5ex]
and
\\[0.5ex]
\mbox{}\hfill
$
F_3^{}\colon (t,x_1^{},x_2^{},x_3^{})\to\  
e^{{}^{\;\!\scriptstyle {}-2\cosh t-\sin e^t}}
(x_1^{}+x_2^{}+x_3^{})
$
\ for all 
$
(t,x_1^{},x_2^{},x_3^{})\in {\mathbb R}^4.
\hfill\mbox{}
$ 
\\[2ex]
\indent
Thus every set of the functionally independent first integrals 
\vspace{0.5ex}
$\{F,\, F_1^{},\, F_2^{}\}, \ \{F,\, F_1^{}, \, F_3^{}\},$ 
$\{F,\, F_2^{},\, F_3^{}\},$ and $\{F_1^{},\, F_2^{},\, F_3^{}\}$ is an 
integral basis on the space ${\mathbb R}^4$ of system (2.37).
\vspace{1.25ex}

{\bf Corollary\! 2.3.}\!                                       %  2.3
{\it
Let $\nu^{k}\!$ be real common eigenvectors of  the matrices $B_{j}^{}$
cor\-res\-pon\-ding to the eigenvalues 
\vspace{0.5ex}
$\lambda_k^j, \ j=1,\ldots, m,\ k= 1,\ldots, m+1.$
Then an autonomous first integral of the Lappo-Da\-ni\-lev\-s\-kii 
differential system {\rm (2.20)} is the scalar function}
\\[1.75ex]
\mbox{}\hfill
$
\displaystyle
F_{_{\scriptstyle 12\ldots m(m+1)}}\colon x\to \
\prod\limits_{k=1}^{m}\,
\bigl|\nu^{k} x\bigr|^{{}^{\scriptstyle \triangle_{k}}}\;\!
\bigl|\nu^{m+1} x\bigr|^{{}^{\scriptstyle {}-\triangle}} 
$
\ for all 
$
x\in {\mathscr X},
\hfill
$
\\[2ex]
{\it 
where 
\vspace{0.75ex}
the determinants $\triangle_{k},\, k=1,\ldots, m$ 
are obtained by replacing the $k\!$-th column of the determinant $\triangle = \big|\lambda_{k}^{j}\big|$  by  
${\rm colon}\!\left(\lambda_{m+1}^{1},\ldots,\lambda_{m+1}^{m}\right)\!,$ respectively.
}
\vspace{1.25ex}

{\bf Example 2.11.}
Consider the fourth order Lappo-Da\-ni\-lev\-s\-kii differential system
\\[1.5ex]
\mbox{}\hfill                            
$
\dfrac{dx}{dt}=\alpha_1^{}(t)A_1^{}\;\!x+\alpha_2^{}(t)A_2^{}\;\!x,
\quad
x\in {\mathbb R}^4,
$
\mbox{}\hfill (2.38)
\\[1.5ex]
with linearly independent continuous functions 
\vspace{1ex}
$\alpha_1^{}\colon J\to {\mathbb R},\ \alpha_2^{}\colon J\to {\mathbb R},$
and the commutative matrices
$
A_1^{} = \left\|\!\!
\begin{array}{rrrr}
  0&  6&  {}-1&  3
\\
  1& {}-7& 1&  {}-5
\\
  2& {}-6&  3&  {}-3
\\
 {}-2& 10& {}-2&   8
\end{array}\!\! \right\|,
\
A_2^{} = \left\|\!\!
\begin{array}{rrrr}
1&  4&  0& 2
\\
1& {}-4& 1& {}-3
\\
0&  {}-4&  1& {}-2
\\
{}-2& 6& {}-2& 5
\end{array}\!\! \right\|
$

The matrices $B_1^{}=A_1^{T}$ and $B_2^{}=A_2^{T},$ 
\vspace{0.25ex}
where $T$ denotes the matrix transpose, have 
the linearly independent common real eigenvectors
\\[1.25ex]
\mbox{}\hfill
$
\nu^{\;\!1} = ( 0, 2, 0, 1 ), 
\quad 
\nu^{\,2} = ( 2, 2, 1, 1 ),
\quad 
\nu^{\,3} = ( 1, 0, 1, 0 ), 
\quad 
\nu^{\,4} = ({}-1, 1, {}-1, 1)
\hfill
$
\\[1.25ex]
corresponding to the real eigenvalues
\\[1.25ex]
\mbox{}\hfill
$
\lambda_1^1 = {}-2, \ \ 
\lambda_2^1 = 1, \ \ 
\lambda_3^1 = 2, \ \ 
\lambda_4^1 =3,
\quad \ 
\lambda_1^2 ={}-1,\ \ 
\lambda_2^2=1, \ \
\lambda_3^2 = 1,\ \
\lambda_4^2 = 2.
\hfill
$
\\[1.25ex]
\indent
The determinants 
\\[2ex]
\mbox{}\hfill
$
\triangle =
\left|\!\!
\begin{array}{rr}
{}-2 &  1
\\
{}-1  & 1
\end{array}\!\!\right| = {}-1,
\qquad
\triangle_{11}^{} =
\left|\!\!
\begin{array}{rr}
2 &  1
\\
1  & 1
\end{array}\!\!\right| = 1,
\qquad
\triangle_{21}^{} =
\left|\!\!
\begin{array}{rr}
{}-2 &  2
\\
{}-1  & 1
\end{array}\!\!\right| = 0,
\hfill
$
\\[2.25ex]
\mbox{}\hfill
$
\triangle_{12}^{} =
\left|\!\!
\begin{array}{rr}
3 &  1
\\
2  & 1
\end{array}\!\!\right| =1,
\qquad
\triangle_{22}^{} =
\left|\!\!
\begin{array}{rr}
{}-2 & 3
\\
{}-1  & 2
\end{array}\!\!\right| ={}-1.
\hfill
$
\\[1.5ex]
\indent
By Corollary 2.3, the scalar functions
\\[1.5ex]
\mbox{}\hfill                                   
$
F_{123}^{}\colon x\to\,
(2x_{2}^{} + x_{4}^{})(x_{1}^{} + x_{3}^{})
$
\ for all 
$
x\in {\mathbb R}^4
\hfill
$
\\[0.75ex]
and
\\[1ex]
\mbox{}\hfill                                   
$
F_{124}^{}\colon x\to \
\dfrac{(2x_{2}^{}+x_{4}^{})({}-x_{1}^{}+x_{2}^{}-x_{3}^{}+x_{4}^{})}
{2x_{1}^{}+2x_{2}^{}+x_{3}^{}+x_{4}^{}}
$
\ for all 
$
x\in {\mathscr X},
\hfill
$
\\[2ex]
where ${\mathscr X}$ is a domain from the set 
\vspace{0.35ex}
$\{x\colon 2x_{1}^{}+2x_{2}^{}+x_{3}^{}+x_{4}^{}\ne 0\}\subset {\mathbb R}^4,$
are autonomous first integrals of the Lappo-Da\-ni\-lev\-s\-kii differential system (2.38).
\vspace{0.75ex}

{\bf Theorem 2.9.}                                    % 2.9
{\it
Suppose  
\vspace{0.25ex}
$\nu^{k}\! =\! \overset{*}{\nu}\;\!{}^{k}+\widetilde{\nu}\;\!{}^{k}\,i\ 
(\overset{*}{\nu}\;\!{}^{k}\!=\!{\rm Re}\,\nu^{k},\, \widetilde{\nu}\;\!{}^{k}\!=\!{\rm Im}\,\nu^{k}),\, 
k\!=\!1,\ldots, s,\, 
s\!\leq\! (m\!+\!1)/2$ 
{\rm(}this set hasn't complex conjugate vectors{\rm)}, 
and $\nu^{\theta},\ \theta=s+1,\ldots, m+1-s,$ are 
common complex and real linearly independent eigenvectors 
\vspace{0.35ex}
of the matrices $B_{j}^{},\,  j\!=\!1,\ldots, m.$  
Then the Lappo-Da\-ni\-lev\-s\-kii dif\-fe\-ren\-ti\-al system {\rm (2.20)} has 
the autonomous first integral}
\\[2ex]
\mbox{}\hfill                                       % (2.39)
$
\displaystyle
F\colon x \to\ 
\prod\limits_{k=1}^{s}
\bigl(P_k^{}(x)\bigr)^{^{\scriptstyle \stackrel{*}h_{_k}}}
\exp\bigl({}-2\,\widetilde{h}_{_{\scriptstyle k}}
\varphi_{_{\scriptstyle k}}(x)\bigr)
\prod\limits_{\theta=s+1}^{m+1-s}\,
\bigl|\nu^{\theta} x\bigr|^{^{\scriptstyle h_{_{\scriptsize\theta}}}}
$
\ for all 
$
x\in {\mathscr X},
\ \ 
{\mathscr X}\subset {\mathbb R}^n,
$
\hfill {\rm (2.39)}
\\[2ex]
{\it where the scalar functions
$
P_k^{}\colon x\to
\bigl(\overset{*}{\nu}\;\!{}^{k}\;\! x\bigr)^{2} +
\bigl(\widetilde{\nu}\,{}^{k}\;\! x\bigr)^{2}$ for all $x\in {\mathbb R}^n,
\ \
\varphi_{_{\scriptstyle k}}\colon x\to 
\arctan\dfrac{\widetilde{\nu}\,{}^{k}\;\! x}
{\overset{*}{\nu}\;\!{}^{k}\;\! x}
$
for all $x\in {\mathscr X},
$
the real numbers
\vspace{0.35ex}
$
\overset{*}{h}_{k}^{},\ \widetilde{h}_{k}^{},\ h_{\theta}^{}$
are an nontrivial solution to the system
\\[1.5ex]
\mbox{}\hfill                  % (2.40)
$
\displaystyle
2\,\sum\limits_{k=1}^{s}\,
\bigl(\overset{*}{\lambda}{}_{_{\scriptstyle k}}^{j}\,
\overset{*}{h}_{_{\scriptstyle k}} -
\widetilde{\lambda}{}_{_{\scriptstyle k}}^{j}\,
\widetilde{h}_{_{\scriptstyle k}}\bigr) \, +\, 
\sum\limits_{\theta=s+1}^{m+1-s}
\lambda_{_{\scriptstyle \theta}}^{j}\,h_{_{\scriptstyle \theta}} = 0,
\quad 
j=1,\ldots, m,
$
\hfill {\rm (2.40)}
\\[1.5ex]
the numbers
\vspace{0.75ex}
$
\lambda_{k}^{j} =\overset{*}{\lambda}{}_{k}^{j} +\widetilde{\lambda}{}_{k}^{j}\,i\ 
(\overset{*}{\lambda}{}_{k}^{j}={\rm Re}\,\lambda_{k}^{j},\ 
\widetilde{\lambda}{}_{k}^{j}={\rm Im}\,\lambda_{k}^{j})$ 
and $\lambda_{\theta}^{j}$ 
are eigenvalues of the matrices $B_{j}^{}$
corresponding to the eigenvectors $\nu^{k}$ and $\nu^{\theta},$ respectively.
}
\vspace{0.75ex}

{\sl Proof.} 
By Lemma 1.1, Properties 1.4 and 1.5, the linear homogeneous functions 
\\[1.25ex]
\mbox{}\hfill
$
p_\xi^{}\colon x\to\nu^{\xi}x
$  
\ for all 
$
x\in {\mathbb R}^{n},
\quad 
\xi=1,\ldots, m+1,
\hfill
$
\\[1.25ex]
are partial integrals of the linear differential systems $\dfrac{dx}{dt}=A_j^{}\;\!x,\ j=1,\ldots, m,$  
\vspace{0.5ex}
and the following system of identities hold
\\[1.5ex]
\mbox{}\hfill                                         %(2.41)
$
{\frak a}_{j}^{}\;\!P_k^{}(x) =
2\;\!\overset{*}{\lambda}{}_{k}^{j}P_k^{}(x)
$
\ for all 
$
x\in {\mathbb R}^{n},
\quad
j=1,\ldots, m, 
\ \
k=1,\ldots, s,
\hfill
$
\\[2ex]
\mbox{}\hfill                                   
$
{\frak a}_{j}^{}\;\!\varphi_k^{}(x) =
\widetilde{\lambda}{}_{k}^{j}
$
\ for all 
$
x\in {\mathscr X},
\quad
j=1,\ldots, m, 
\ \
k=1,\ldots, s,
$
\hfill (2.41)
\\[2ex]
\mbox{}\hfill                                   
$
{\frak a}_{ j}^{}\;\!\nu^{\theta} x =
\lambda_{\theta}^{j}\;\!\nu^{\theta} x
$
\ for all 
$
x\in {\mathbb R}^{n},
\quad
j=1,\ldots, m, 
\quad
\theta=s+1,\ldots, m+1-s.
\hfill
$
\\[1.75ex]
\indent
Since the function $F\colon {\mathscr X}\to {\mathbb R}$ is given by (2.39), it follows that
\\[2ex]
\mbox{}\hfill
$
\displaystyle
{\frak a}_{j}^{}\;\! F(x) =
\biggl(\
\prod\limits_{k=1}^{s}
\bigl(P_{_{\scriptstyle k}}(x)\bigr)
^{^{\scriptstyle \stackrel{*}h_{_k}-1}}
\exp\bigl({}-2\,\widetilde{h}_{_{\scriptstyle k}}
\varphi_{_{\scriptstyle k}}(x)\bigr)
\ \sum\limits_{k=1}^{s}\,\overset{*}{h}_{_{\scriptstyle k}}\
\prod\limits_{l=1, l\ne k}^{s}\!P_{_{\scriptstyle l}}(x)\,
{\frak a}_{j}^{}\;\! P_{k}^{}(x)\ +
\hfill
$
\\[2.25ex]
\mbox{}\hfill
$
\displaystyle
+ \
\prod\limits_{k=1}^{s}
\bigl(P_{_{\scriptstyle k}}(x)\bigr)
^{^{\scriptstyle \stackrel{*}h_{_k}}}
\exp\bigl({}-2\,\widetilde{h}_{_{\scriptstyle k}}
\varphi_{_{\scriptstyle k}}(x)\bigr)\,
\sum\limits_{k=1}^{s}\,
{\frak a}_{j}^{}
\bigl({}-2\,\widetilde{h}_{_{\scriptstyle k}}
\varphi_{_{\scriptstyle k}}(x)\bigr)\biggr)\,
\prod\limits_{\theta=s+1}^{m+1-s}
\bigl|\nu^{\theta} x\bigr|^{^{\scriptstyle h_{_{\scriptsize\theta}}}} \ +
\hfill
$
\\[3ex]
\mbox{}\hfill
$
\displaystyle
+ \
\prod\limits_{k=1}^{s}
\bigl(P_{_{\scriptstyle k}}(x)\bigr)
^{^{\scriptstyle \stackrel{*}h_{_k}}}
\exp\bigl({}-2\,\widetilde{h}_{_{\scriptstyle k}}
\varphi_{_{\scriptstyle k}}(x)\bigr)
\prod\limits_{\theta=s+1}^{m+1-s}
\bigl|\nu^{\theta} x\bigr|
^{^{\scriptstyle h_{_{\scriptsize\theta}}-1}}
\sum\limits_{\theta=s+1}^{m+1-s}\,
{\rm sgn}\;\!\bigl(\nu^{\theta} x\bigr)\;\!
h_{_{\scriptstyle \theta}}
\prod\limits_{l=s+1,l\ne \theta}^{m+1-s}
\bigl|\nu^{l} x\bigr|\,
{\frak a}_{ j}^{}\bigl(\nu^{\theta} x\bigr)
\hfill
$
\\[2.75ex]
\mbox{}\hfill
for all 
$
x\in {\mathscr X},
\quad 
j=1,\ldots, m.
\hfill
$
\\[1.75ex]
\indent
Using the system of identities (2.41), we get
\\[2ex]
\mbox{}\hfill
$
\displaystyle
{\frak a}_{j}^{}\;\!F(x) =
\biggl(\
2\sum\limits_{k=1}^{s}\,
\bigl(\overset{*}{\lambda}{}_{_{\scriptstyle k}}^{j}\,
\overset{*}{h}_{_{\scriptstyle k}} -
\widetilde{\lambda}{}_{_{\scriptstyle k}}^{j}\,
\widetilde{h}_{_{\scriptstyle k}}\bigr) +
\sum\limits_{\theta=s+1}^{m+1-s}
\lambda_{_{\scriptstyle \theta}}^{j}\,h_{_{\scriptstyle \theta}}\biggr) F(x)
$
\ for all 
$
x\in {\mathscr X},
\quad 
j=1,\ldots, m.
\hfill
$
\\[1.75ex]
\indent
If the numbers
\vspace{0.35ex}
$
\overset{*}{h}_{k}^{},\ \widetilde{h}_{k}^{},\
k=1,\ldots, s,
$ 
and 
$
h_{\theta}^{},\ \theta=s+1,\ldots, m+1-s,
$
are an nontrivial real solution to the system (2.40), 
then (by Lemma 2.6) the scalar function (2.39) is 
an autonomous first integral of the Lappo-Da\-ni\-lev\-s\-kii 
dif\-fe\-ren\-ti\-al system {\rm (2.20)}. \k
\vspace{1ex}

{\bf Example 2.12.}
Consider the fourth order Lappo-Da\-ni\-lev\-s\-kii differential system
\\[2ex]
\mbox{}\hfill       % (2.42)                     
$
\dfrac{dx}{dt}=\alpha_1^{}(t)A_1^{}\;\!x+\alpha_2^{}(t)A_2^{}\;\!x,
\quad
x\in {\mathbb R}^4,
$
\mbox{}\hfill (2.42)
\\[2ex]
where linearly independent functions 
\vspace{1ex}
$\alpha_1^{}\colon J\to {\mathbb R}$ and $\alpha_2^{}\colon J\to {\mathbb R}$ are 
continuous, the matrices 
$
A_1^{} = \left\|\!\!
\begin{array}{rrrr}
  {}-1&  2&  {}-2&  0
\\
 {}-6& 5& {}-4&  2
\\
 {}-3& 2&  {}-2&  1
\\
 2& {}-1& 2&  0
\end{array}\!\! \right\|
$
and
$
A_2^{} = \left\|\!\!
\begin{array}{rrrr}
{}-4&  2&  2& 4
\\
{}-5& 1& 6& 6
\\
0&  {}-1&  3& 1
\\
{}-3& 2& 0& 3
\end{array}\!\! \right\|
$
commute.
\vspace{1ex}

The matrices $B_1^{}=A_1^{T}$ and $B_2^{}=A_2^{T},$ 
\vspace{0.35ex}
where $T$ denotes the matrix transpose, have 
the common eigenvectors
\vspace{0.5ex}
$
\nu^{\;\!1} = (1+i, {}-i, {}-1+i, {}-1 ), \ \nu^{\,2} = (1-i,\, i, {}-1-i, {}-1),
$
$ 
\nu^{\;\!3} = ( 1, {}-1, 2, 0 ),\ \nu^{\;\!4} = ({}-1, 0, 2, 2)
$
\vspace{0.75ex}
corresponding to the eigenvalues 
$
\lambda_1^1 = 1+i, \
\lambda_2^1 = 1-i,$ 
$\lambda_3^1 = {}-1, \ \lambda_4^1 =1,
$
and
$
\lambda_1^2 =i,\ \lambda_2^2={}-i, \
\lambda_3^2 = 1,\ \lambda_4^2 = 2.
$
\vspace{0.75ex}

The linear homogeneous systems
\\[1.25ex]
\mbox{}\hfill
$
2{\stackrel{*}{h}}_{11}^{}-2\widetilde{h}_{11}^{}-h_{31}^{} = 0,
\ \
{}-2\widetilde{h}_{11}^{}+h_{31}^{} = 0
$
\ \ and \ \ 
$
2{\stackrel{*}{h}}_{12}^{}-2\widetilde{h}_{12}^{}+h_{42}^{} = 0,
\ \
{}-2\widetilde{h}_{12}^{}+2h_{42}^{} = 0
\hfill
$
\\[1.5ex]
have, for example, the solutions
${\stackrel{*}{h}}_{11}^{}=2, \ \widetilde{h}_{11}^{}=1, \
h_{31}^{}=2,$ and ${\stackrel{*}{h}}_{12}^{}=1, \ \widetilde{h}_{12}^{}=2, \ h_{42}^{}=2.$ 
\vspace{0.75ex}

By Theorem 2.9, the scalar functions 
\\[1.75ex]
\mbox{}\hfill                                   
$
F_{1}^{}\colon x\to\,
(x_{1}^{}-x_{2}^{}+2x_{3}^{})
\bigl( (x_{1}^{}-x_{3}^{}-x_{4}^{})^{2} + (x_{1}^{}-x_{2}^{}+x_{3}^{})^{2}\,\bigr)  
\exp \biggl({}-\arctan\dfrac{x_{1}^{}-x_{2}^{}+x_{3}^{}}{x_{1}^{}-x_{3}^{}-x_{4}^{}}\biggr)
\hfill
$
\\[1.25ex]
and
\\[1.25ex]
\mbox{}\hfill                                   
$
F_{2}^{}\colon x\to \,
({}-x_{1}^{}+2x_{3}^{}+2x_{4}^{})^{2}
\bigl( ( x_{1}^{}-x_{3}^{}-x_{4}^{} )^{2} + ( x_{1}^{}-x_{2}^{}+x_{3}^{} )^{2}\,\bigr)
\exp\biggl({}-4\arctan\dfrac{x_{1}^{}-x_{2}^{}+x_{3}^{}}{x_{1}^{}-x_{3}^{}-x_{4}^{}}\biggr)
\hfill
$
\\[2ex]
are autonomous first integrals 
\vspace{0.35ex}
of the Lappo-Da\-ni\-lev\-s\-kii dif\-fe\-ren\-ti\-al system {\rm (2.42)}
on any domain  ${\mathscr X}$ from the set 
$\{x\colon x_{1}^{}-x_{3}^{}-x_{4}^{}\ne 0\}\subset {\mathbb R}^4.$
%\vspace{0.75ex}

\newpage

{\bf Example 2.13} (continuation of Example 2.8).
From the linear homogeneous system
\vspace{0.25ex}
$
h_1^{}+2{\stackrel{*}{h}}_{2}^{} = 0,
\
{}-h_1^{}-2\widetilde{h}_{2}^{} = 0
$
it follows that, for example, $h_1^{}={}-2,\ {\stackrel{*}{h}}_{2}=\widetilde{h}_{2}^{} =1.$
\vspace{0.35ex}

Then, using the common eigenvectors
\vspace{0.35ex}
$\nu^1=(1, 0,{}-1),\ \nu^2=(i, i, 1)$ 
and the corresponding eigenvalues
\vspace{0.35ex}
$\lambda^1_1=\lambda_2^1=1,\ \lambda_1^2={}-1,\,\lambda_2^2=i,$
we can build (by Theorem 2.9) the autonomous first integral 
\vspace{0.35ex}
of the Lappo-Da\-ni\-lev\-s\-kii dif\-fe\-ren\-ti\-al system {\rm (2.27)}
\\[1.75ex]
\mbox{}\hfill                                   
$
F\colon (x_1^{},x_2^{},x_3^{})\to\
\dfrac{(x_{1}^{}+x_{2}^{})^2+x_{3}^{2}}{(x_{1}^{}-x_{3}^{})^2}\,
\exp \biggl({}-2\arctan\dfrac{x_{1}^{}+x_{2}^{}}{x_{3}^{}}\biggr)
$
\ for all 
$
(x_1^{},x_2^{},x_3^{})\in {\mathscr X},
\hfill
$
\\[2.25ex]
where ${\mathscr X}$ is a domain from the set
$\{(x_1^{},x_2^{},x_3^{})\colon x_{1}^{}\ne x_{3}^{},\, x_{3}^{}\ne 0\}.$
\vspace{1ex}

{\bf Theorem 2.10.}                                        % 2.10
{\it
Suppose  
\vspace{0.25ex}
$\nu^{\tau}\! = {\stackrel{*}{\nu}\!{}^{\tau}}+\;\!\widetilde{\nu}\,{}^{\tau}i,\
\nu^{s+\tau}\! = {\stackrel{*}{\nu}\!{}^{\tau}} -\;\!\widetilde{\nu}\,{}^{\tau}i,\
\tau\!=1,\ldots, s,\, s\leq m/2,\ \nu^{2s+1} =$
$={\stackrel{*}{\nu}\!{}^{2s+1}} +\widetilde{\nu}\,{}^{2s+1}\,i,$
\vspace{0.35ex}
and $\nu^{\theta},\ \theta=2s+2,\ldots, m+1,$ are 
common complex and real linearly independent eigenvectors 
\vspace{0.5ex}
of the matrices $B_{j}^{},\,  j\!=\!1,\ldots, m.$  
Then the Lappo-Da\-ni\-lev\-s\-kii dif\-fe\-ren\-ti\-al system {\rm (2.20)} has 
the autonomous first integrals}
\\[2.25ex]
\mbox{}\hfill                                        % (2.43)
$
\displaystyle
F_{1}^{}\colon  x\to \
\prod\limits_{\tau=1}^{s}
\bigl(P_{_{\scriptstyle \tau}}(x)\bigr)^{
^{\scriptstyle \stackrel{*}{h}_
{_{\scriptsize \tau}}+\stackrel{*}{h}_{_{\scriptsize s+\tau}}}}
\exp\Bigl({}-2\;\!\bigl(\;\!
\widetilde{h}_{_{\scriptstyle \tau}}-
\widetilde{h}_{_{\scriptstyle s+\tau}}\bigr)\;\!
\varphi_{_{\scriptstyle \tau}}(x)\Bigr)\,  \cdot
\hfill
$
\\
\mbox{}\hfill {\rm(2.43)}
\\[1ex]
\mbox{}\hfill
$
\displaystyle
\cdot\,
\bigl(P_{_{\scriptstyle 2s+1}}(x)\bigr)^{
^{\scriptstyle\stackrel{*}{h}_{_{\scriptsize 2s+1}}}}
\exp\Bigl({}-2\;\!
\widetilde{h}_{_{\scriptstyle 2s+1}}\;\!
\varphi_{_{\scriptstyle 2s+1}}(x)\Bigr)
\prod\limits_{\theta=2s+2}^{m+1}\!\!
\bigl(\nu^{\theta} x\bigr)^{^{\scriptstyle 2
\overset{*}{h}_{_{\scriptsize \theta}}}} 
$
\ for all 
$x\in {\mathscr X}
\hfill
$
\\[2.75ex]
and
\\[2.25ex]
\mbox{}\hfill                                        % (2.44)
$
\displaystyle
F_{2}^{}\colon  x\to \
\prod\limits_{\tau=1}^{s}
\bigl(P_{_{\scriptstyle \tau}}(x)\bigr)^
{^{\scriptstyle
\widetilde{h}_{_{\scriptsize \tau}}+\widetilde{h}_{_{\scriptsize s+\tau}}}}
\exp\Bigl(2\;\!\bigl(\;\!
\overset{*}{h}_{_{\scriptstyle \tau}}-
\overset{*}{h}_{_{\scriptstyle s+\tau}}\bigr)\;\!
\varphi_{_{\scriptstyle \tau}}(x) \Bigr)\, \cdot
\hfill
$
\\
\mbox{}\hfill {\rm(2.44)}
\\[1ex]
\mbox{}\hfill
$
\displaystyle
\cdot\,
\bigl(P_{_{\scriptstyle2s+1}}(x)\bigr)^{^{\scriptstyle \widetilde{h}_{_{\scriptsize 2s+1}}}}
\exp\Bigl(2\;\!\overset{*}{h}_{_
{\scriptstyle 2s+1}}\;\!\varphi_{_{\scriptstyle 2s+1}}(x)\Bigr)
\prod\limits_{\theta=2s+2}^{m+1}
\bigl(\nu^{\theta} x\bigr)^{^{\scriptstyle 2\widetilde{h}_
{_{\scriptsize \theta }} }}
$
\ for all 
$x\in {\mathscr X},
\quad
{\mathscr X}\subset {\rm D}F_1^{}\cap{\rm D}F_2^{}\,,
\hfill
$
\\[2.5ex]
{\it where the scalar functions
\\[1.5ex]
\mbox{}\hfill
$
P_\tau^{}\colon x\to
\bigl(\overset{*}{\nu}\;\!{}^{\tau} x\bigr)^{2}\! +
\bigl(\widetilde{\nu}\,{}^{\tau} x\bigr)^{2},
\quad
\varphi_{_{\scriptstyle \tau}}\colon x\to 
\arctan\dfrac{\widetilde{\nu}\,{}^{\tau} x}{\overset{*}{\nu}\;\!{}^{\tau} x}
$
\ for all 
$
x\in {\mathscr X},
\quad 
\tau=1,2,\ldots, s, 2s+1,
\hfill
$
\\[2ex]
the numbers
$
h_{k}^{} = {\overset{*}{h}}_{k}^{} + \widetilde{h}_{k}^{}\,i,
\ k=1,\ldots, m+1,
$
are an nontrivial solution to the system
\\[1.5ex]
\mbox{}\hfill
$
\displaystyle
\sum\limits_{k=1}^{m+1}\, \lambda_{k}^{j}\;\!h_{k}^{} = 0, 
\quad
j=1,\ldots, m,
\hfill
$
\\[1.5ex]
and $\lambda_{k}^{j}$ are eigenvalues of the matrices $B_{j}$ 
corresponding to the eigenvectors $\nu^{k}.$
}
\vspace{1ex}

{\sl Proof.}
We form two complex-valued functions of real variables
\\[2ex]
\mbox{}\hfill
$
\displaystyle
\overset{*}{F}\colon x\to \
\prod\limits_{k=1}^{2s}\bigl(\nu^{k} x\bigl)^{\!^{\scriptstyle
h_k^{}}}
\bigl(\nu^{2s+1} x\bigr)^{\!^{\scriptstyle h_{2s+1}^{}}}
\prod\limits_{\theta=2s+2}^{m+1}
\bigl(\nu^{\theta} x\bigr)^{\!^{\scriptstyle h_{\theta}^{}}} 
$
\ for all 
$
x\in {\mathscr X}
\hfill
$
\\[1.5ex]
and
\\[1.5ex]
\mbox{}\hfill
$
\displaystyle
\overset{**}{F}\colon x\to \
\prod\limits_{k=1}^{2s}\bigl( \nu^{k} x\bigr)^{\!^{\scriptstyle
l_k^{}}}
\bigl(\;\!\overline{\nu^{2s+1}}\, x\bigr)^{\!^{\scriptstyle
l_{2s+1}^{}}}
\prod\limits_{\theta=2s+2}^{m+1}
\bigl(\nu^{\theta} x\bigr)^{\!^{\scriptstyle l_{\theta}^{}}} 
$
\ for all 
$
x\in {\mathscr X},
\hfill
$
\\[1.75ex]
where $h_{k}^{},\ l_{k}^{},\, k=1,\ldots, m+1,$ are some complex numbers, 
a domain ${\mathscr X}\subset {\mathbb R}^n.$  

We have
\\[1.5ex]
\mbox{}\hfill
$
\displaystyle
{\frak a}_{j}^{}\overset{*}{F}(x) =
\sum\limits_{k=1}^{m+1}
\, \lambda_{k}^{j}\,h_{k}^{}\,
\overset{*}{F}(x)
$
\ for all 
$
x\in {\mathscr X},
\quad
j=1,\ldots, m,
\hfill
$
\\[2ex]
\mbox{}\hfill
$
\displaystyle
{\frak a}_{j}\overset{**}{F}(x)  =
\biggl(\,
\sum\limits_{k=1}^{2s}\;\!
\lambda_{_{\scriptstyle k}}^{j}\;\!l_{_{\scriptstyle k}} +
\overline{\lambda_
{_{\scriptstyle 2s+1}}^{j}}\;l_{_{\scriptstyle 2s+1}} +
\!\sum\limits_{\theta=2s+2}^{m+1}\!\lambda_{_{\scriptstyle \theta}}^{j}\;\!
l_{_{\scriptstyle \theta}}
\biggr)\;\!
\overset{**}{F}(x)
$
\ for all 
$
x\in {\mathscr X},
\quad
j=1,\ldots, m.
\hfill
$
\\[1.5ex]
\indent
Let $h_{k}^{}\! =\overset{*}{h}_{k}^{}+\widetilde{h}_{k}^{}\,i,\ k=1,\ldots, m+1,$ 
be an nontrivial solution to the linear system
\\[1.5ex]
\mbox{}\hfill
$
\displaystyle
\sum\limits_{k=1}^{m+1}\lambda_{k}^{j}\;\!h_{k}^{} =  0, 
\quad 
j=1,\ldots, m.
\hfill
$
\\[1.5ex]
Then 
\vspace{0.5ex}
$
l_{_{\scriptstyle k}}\!=
\overset{*}{h}_{_{\scriptstyle s+k}} -\;\!
\widetilde{h}_{_{\scriptstyle s+k}}\,i, \
\vspace{0.5ex}
l_{_{\scriptstyle s+k}} \!=
\overset{*}{h}_{_{\scriptstyle k}} -\;\!
\widetilde{h}_{_{\scriptstyle k}}\,i,\
k=1,\ldots, s,
\ \,
l_{_{\scriptstyle 2s+1}}\! =
\overset{*}{h}_{_{\scriptstyle 2s+1}} -\;\!
\widetilde{h}_{_{\scriptstyle 2s+1}}\,i,\ \,
l_{_{\scriptstyle \theta}}\! =
\overset{*}{h}_{_{\scriptstyle \theta}} -\;\!
\widetilde{h}_{_{\scriptstyle \theta}}\,i,
$
$\theta=2s+2,\ldots, m+1$
is an nontrivial solution to the linear system 
\\[1.5ex]
\mbox{}\hfill
$
\displaystyle
\sum\limits_{k=1}^{2s}\lambda_{k}^{j}\,l_{k}^{} +
\overline{\lambda_{2s+1}^{j}}\,l_{2s+1}^{} +
\sum\limits_{\theta=2s+2}^{m+1}
\lambda_{\theta}^{j}\,l_{\theta}^{} = 0,
\quad 
j=1,\ldots, m,
\hfill
$
\\[1.5ex]
and the functions ${\stackrel{*}{F}}\colon {\mathscr X}\to {\mathbb C},\
{\stackrel{**}{F}}\colon {\mathscr X}\to {\mathbb C}$ 
\vspace{0.5ex}
are autonomous first integrals of the system (2.20).

Since 
\vspace{0.35ex}
$
F_{1}^{}(x) =\overset{*}{F}(x)\;\!\overset{**}{F}(x)
$
for all 
$x\in {\mathscr X}$
and
$
F_{2}^{}(x) = \overset{**}{F}{}^{\;\!^{\scriptstyle i}}(x)\;\!
\overset{*}{F}{}^{\;\!^{\scriptstyle{}- i}}(x)
$
for all 
$x\in {\mathscr X},$
we see that the scalar functions (2.43) and (2.44) are 
\vspace{0.35ex}
autonomous first integrals (by Lemma 2.6) 
of the Lappo-Da\-ni\-lev\-s\-kii differential system (2.20). \k
\vspace{1ex}

{\bf Example 2.14.}
The small oscillations of gyrocompass in the precession gyroscope theory
are described by the following fourth-order linear differential system [61]
\\[2.5ex]
\mbox{}\hfill                                             % (2.45)
$
\begin{array}{l}
\dfrac{dx_1^{}}{dt} =\omega_0^{}\;\! x_2^{} +\Omega(t)\;\! x_4^{},
\qquad \ \ \ \ \,
\dfrac{dx_2^{}}{dt} = {}-\omega_0^{}\;\! x_1^{} +\Omega(t)\;\! x_3^{},
\\[3.75ex]
\dfrac{dx_3^{}}{dt} = {}-\Omega(t)\;\! x_2^{} -\omega_0^{}\;\! x_4^{},
\qquad \  
\dfrac{dx_4^{}}{dt} ={}-\Omega(t)\;\! x_1^{} +\omega_0^{}\;\! x_3^{},
\end{array}
$
\hfill {\rm (2.45)}
\\[2.5ex]
where a continuous on an interval $J\subset {\mathbb R}$ function $\Omega\colon J\to {\mathbb R}$ 
is the projection of the absolute angular velocity of the gyroscope's sensitive element 
\vspace{0.35ex}
onto the direction of geocentric vertical line,
the constant $\omega_0^{}=\sqrt{g/R}\ (g$ is the acceleration due to gravity,
\vspace{0.35ex}
$R$ is the Earth's radius).

The system (2.45) is a 
Lappo-Da\-ni\-lev\-s\-kii differential system [8, pp. 37 -- 38].
\vspace{0.35ex}

By Theorem 2.10, using the 
\vspace{0.35ex}
linearly independent common eigenvectors
$\nu^{1}=(1, i, {}-1, i),$ 
$\nu^{2}\!=\!(i, {}-1, i, 1),\ \nu^{3}\!=\!(1, {}-i, {}-1, {}-i),\ \nu^{4}\!=\!({}-i, {}-1, {}-i, 1)$ 
\vspace{0.5ex}
and the corresponding eigen- va\-lues
$\lambda_1^1=\lambda_2^1={}-\omega_0\,i, \ 
\lambda_3^1=\lambda_4^1=\omega_0\,i,\ 
\lambda_1^2={}-i,\ \lambda_2^2=\lambda_3^2=i,\ \lambda_4^2={}-i$ of the matrices
\\[2ex]
\mbox{}\hfill
$
B_1^{}=
\left\|\!\!
\begin{array}{cccc}
0 & {}-\omega_0^{} & 0 & 0
\\[0.35ex]
\omega_0^{} & 0 & 0 & 0
\\[0.35ex]
0 & 0 & 0 & \omega_0^{}
\\[0.35ex]
0 & 0 & {}-\omega_0^{} & 0
\end{array}
\!\!\right\|
\qquad
\text{and}
\qquad
B_2^{}=
\left\|\!\!
\begin{array}{ccrr}
0 & 0 & 0 & {}-1
\\[0.35ex]
0 & 0 & {}-1 & 0
\\[0.35ex]
0 & 1 & 0 & 0
\\[0.35ex]
1 & 0 & 0 & 0
\end{array}
\!\!\right\|,
\hfill
$
\\[2.25ex]
we can build the autonomous first integrals  
of the Lappo-Da\-ni\-lev\-s\-kii differential system (2.45)
\\[2ex]
\mbox{}\hfill
$
F_1^{}\colon x\to \ (x_1^{}-x_3^{})^2+(x_2^{}+x_4^{})^2,
\quad
F_2^{}\colon x\to\  ({}-x_2^{}+x_4^{})^2+(x_1^{}+x_3^{})^2
$
\ for all 
$
x\in {\mathbb R}^4.
\hfill
$
\\[2.25ex]
\indent
By Theorem 2.6, using the functions 
\vspace{0.5ex}
$\alpha_1^{}\colon t\to 1$ for all $t\in J,\ 
\alpha_2^{}\colon t\to \Omega(t)$ for all $t\in J,$ 
the numbers
\vspace{0.75ex}
${\stackrel{*}{\lambda}}{}^1_1= {\stackrel{*}{\lambda}}{}^2_1=0,\ 
\widetilde{\lambda}{}^1_1={}-\omega_0,\ \widetilde{\lambda}{}^2_1={}-1,$  
${\stackrel{*}{\lambda}}{}^1_2= {\stackrel{*}{\lambda}}{}^2_2=0,\ 
\widetilde{\lambda}{}^1_2={}-\omega_0,\ \widetilde{\lambda}{}^2_2=1,$ 
and the common eigenvectors $\nu^1$ and $\nu^2,$ 
we can construct the first integrals of system (2.45)
\\[2ex]
\mbox{}\hfill
$
\displaystyle
F_3^{}\colon (t,x)\to \arctan\dfrac{x_2^{}+x_4^{}}{x_1^{}-x_3^{}} + 
\int\limits_{t_0^{}}^{t} \bigl(\omega_0^{}+\Omega(\tau)\bigr)\;\! d\tau
$
for all 
$
(t,x)\!\in\! J\!\times\! {\mathscr X}_1^{},
\,
{\mathscr X}_1^{}\!\subset\! \{x\colon x_1^{}-x_3^{}\ne 0\},
\hfill
$
\\[1.25ex]
and
\\[1.25ex]
\mbox{}\hfill
$
\displaystyle
F_4^{}\colon (t,x)\to  \arctan\!\dfrac{x_1^{}+x_3^{}}{-\;\!x_2^{}\!+x_4^{}} + 
\int\limits_{t_0^{}}^{t} \bigl(\omega_0^{}-\Omega(\tau)\bigr)\;\! d\tau
$
for all 
$
(t,x)\!\in\! J\!\times\! {\mathscr X}_2^{},
\,
{\mathscr X}_2^{}\!\subset\! \{x\colon x_2^{}\!-x_4^{}\!\ne 0\}.
\hfill
$
\\[2ex]
\indent
The functionally independent first integrals $F_1^{},\ldots, F_4^{}$ 
\vspace{0.35ex}
are an integral basis of 
the Lappo-Da\-ni\-lev\-s\-kii system
(2.45) on any domain $J\times {\mathscr X},$
\vspace{0.75ex}
where ${\mathscr X}\subset \{x\colon x_1^{}-x_3^{}\ne 0,\ x_2^{}-x_4^{}\ne 0\}.$

{\bf Example 2.15.}
The Lappo-Da\-ni\-lev\-s\-kii differential system
\\[2ex]
\mbox{}\hfill                                             % (2.46)
$
\begin{array}{c}
\dfrac{dx_1^{}}{dt} =\alpha_1^{}(t)\;\! x_1^{} +\alpha_2^{}(t)\;\! x_3^{},
\qquad \ \ \
\dfrac{dx_2^{}}{dt} ={}-\alpha_1^{}(t)\;\! x_2^{} +\alpha_2^{}(t)\;\! x_4^{},
\\[3.75ex]
\dfrac{dx_3^{}}{dt} ={}-\alpha_2^{}(t)\;\! x_1^{} +\alpha_1^{}(t)\;\! x_3^{},
\qquad
\dfrac{dx_4^{}}{dt} ={}-\alpha_2^{}(t)\;\! x_2^{} -\alpha_1^{}(t)\;\! x_4^{}
\end{array}
$
\hfill {\rm (2.46)}
\\[2.5ex]
with linealy independent on an interval $J\subset {\mathbb R}$ functions
\vspace{0.75ex}
$\alpha_1^{}\colon J\to {\mathbb R}$ and $\alpha_2^{}\colon J\to {\mathbb R}$ 
such that the matrices
\vspace{0.75ex}
$
B_1^{} =  \left\|\!\!
\begin{array}{rrrr}
 1 &    0 & 0 &    0
\\
 0 & -\;\!1 & 0 &    0
\\
 0 &    0 & 1 &    0
\\
 0 &    0 & 0 & -\;\!1
\end{array}
\!\!
\right\|
$
and
$
B_2^{} =  \left\|\!\!
\begin{array}{rrrr}
 0 & 0 & -\;\!1 &    0
\\
 0 & 0 &    0 & -\;\!1
\\
 1 & 0 &    0 &    0
\\
 0 & 1 &    0 &    0
\end{array}
\!\!
\right\|
$
have the linealy independent eigenvectors
\vspace{0.5ex}
$
\nu^{1} = (0, {}-i, 0, 1),\
\nu^{2} = (0, i, 0, 1),\
\nu^{3}\! =\! ( -\;\! i, 0, 1, 0), \,
\nu^{4}\! =\! (i, 0, 1, 0)
$
corresponding to the eigenvalues 
\vspace{0.5ex}
$
\lambda_1^1\!=\lambda_2^1\!=-\;\!1,\
\lambda_3^1=\lambda_4^1=1,\ 
\lambda_1^2=\lambda_3^2={}-i,\
\lambda_2^2=\lambda_4^2=i.
$

From the linear system
\vspace{0.35ex}
$
{}-h_{1}^{} - h_{2}^{} + h_{3}^{} = 0, \
{}-ih_{1}^{} + ih_{2}^{} - ih_{3}^{} = 0
$
it follows that, for example, $h_{1}^{}=0,\ h_{2}^{}=h_{3}^{}=1.$
\vspace{0.25ex}
Then, by Theorem 2.10, the  Lappo-Da\-ni\-lev\-s\-kii differential system (2.46) 
has the autonomous first integrals
\\[2ex]
\mbox{}\hfill
$
F_1^{}\colon x\to\,
(x_1^2+x_3^2)(x_2^2+x_4^2)
$
\ for all 
$
x\in {\mathbb R}^4,
\qquad
F_2^{}\colon x\to\
\dfrac{x_1^{} x_2^{} + x_3^{} x_4^{}}{x_1^{} x_4^{} - x_2^{} x_3^{}}
$
\ for all 
$
x\in{\mathscr X},
\hfill
$
\\[2ex]
where ${\mathscr X}$ is any domain from the set 
\vspace{0.75ex}
$\{x\colon x_1^{} x_4^{} - x_2^{} x_3^{}\ne 0\}\subset {\mathbb R}^4.$

{\bf Example 2.16.}
Consider the fifth-order Lappo-Da\-ni\-lev\-s\-kii differential system 
\\[1.75ex]
\mbox{}\hfill       % (2.47)                     
$
\dfrac{dx}{dt}=\alpha_1^{}(t)A_1^{}\;\!x+\alpha_2^{}(t)A_2^{}\;\!x+\alpha_3^{}(t)A_3^{}\;\!x,
\quad
x\in {\mathbb R}^5,
$
\mbox{}\hfill (2.47)
\\[1.75ex]
where linearly independent functions $\alpha_j^{}\colon J\to {\mathbb R},\ j=1,2,3$ are
continuous, the matrices
\\[2ex]
\mbox{}\hfill
$
B_1^{} =A_1^{T}=
\left\|\!\!
\begin{array}{rrrrr}
 4&  {}-2&   2& {}-7&  11
\\
0&   2&   2& {}-7&   7
\\
 0&   4&   0& {}-2&   2
\\
{}-4&   8&  {}-8&  6&  {}-2
\\
{}-4&   4&  {}-4&  4&   0
\end{array}\!\! \right\|,
\qquad
B_2^{} =A_2^{T}=
\left\|\!\!
\begin{array}{rrrrr}
 4& 0&  0& {}-5& 7
\\
 0& 4&  0& {}-5& 5
\\
 0& 4&  0& 0& 0
\\
{}-2& 6& {}-6& 6& {}-2
\\
{}-2& 2& {}-2& 2& 2
\end{array}\!\!\right\|,
\hfill
$
\\[2.5ex]
\mbox{}\hfill
$
B_3^{} =  A_3^{T}=
\left\|\!\!
\begin{array}{rrrrr}
 0&  6&  {}-6& {}-9& 17
\\
 0&  6& {}-18&  3& {}-3
\\
 0& 12& {}-24& 18& {}-18
\\
{}-8& 20& {}-20& 14& {}-14
\\
{}-8&  8&  {}-8&  8& {}-8
\end{array}\!\!\right\|
$
\quad
$(T$ denotes the matrix transpose).
\hfill\mbox{}
\\[2ex]
\indent
Using the linearly independent common eigenvectors
\vspace{0.5ex}
$
\nu^{1}\! =\! (1, 0, 0, i, i),\,
\nu^{2}\! =\! (1, 0, 0, -i, -i),
$
$
\nu^{3}\! =(1+2i, 1+2i, 2, 2, 0),\
\nu^{4}\! = (1-2i, 1-2i, 2, 2, 0), \
\nu^{5}\! = (0, 1, 1, 0, 0)
$
\vspace{0.5ex}
of the matrices $B_1^{},$ $B_2^{}, B_3^{}$
\vspace{0.5ex}
and the corresponding eigenvalues
$
\lambda_1^1=4+4i, \
\lambda_2^1=4-4i, \
\lambda_3^1=4i,   \
\lambda_4^1={}-4i,\
\lambda_{5}^{1}=4,
$ 
$
\lambda_1^2=4+2i,\
\lambda_2^2=4-2i,\
\lambda_3^2=2+4i,\
\lambda_4^2=2-4i, \
\lambda_{5}^{2}=4, 
$
\vspace{0.5ex}
$
\lambda_1^3=8i,   \ 
\lambda_2^3={}-8i, \
\lambda_3^3=12i,
$
$
\lambda_4^3={}-12i,\
\lambda_{5}^{2}={}-12,
$
\vspace{0.5ex}
we can construct (by Theorem 2.10) the autonomous first integrals of 
the Lappo-Da\-ni\-lev\-s\-kii differential system (2.47)
\\[1ex]
\mbox{}\hfill  
$
F_{1}^{}\colon x\to 
\dfrac{\left( x_{1}^{2} + (x_{4}^{}+x_{5}^{})^{2}\,\right)^{2}}
{(x_{2}^{}\!+\!x_{3}^{})^2
\left((x_{1}^{}\!+\!x_{2}^{}\!+\!2x_{3}^{}\!+\!2x_{4}^{})^{2} +4(x_{1}^{}\!+\!x_{2}^{})^{2}\right)}\,
\exp\!\left(-\;\!2\arctan\dfrac{2(x_{1}+x_{2})}{x_{1}\!+\!x_{2}\!+\!2x_{3}\!+\!2x_{4}}\right)
\hfill
$
\\[1.5ex]
and
\\[1.5ex]
\mbox{}
$
F_{2}^{}\colon x\to \
\dfrac{\left( ( x_{1}^{}+x_{2}^{}+2x_{3}^{}+2x_{4}^{} )^{2} +
4(x_{1}^{}+x_{2}^{})^{2}\,\right)^{5}}{\left(x_{2}^{}+x_{3}^{}\right)^{10}}\ \cdot
\hfill  
$
\\[2ex]
\mbox{}\hfill  
$
\cdot \
\exp\left( 12\arctan\dfrac{x_{4}^{}+x_{5}^{}}{x_{1}^{}}
\,-\,10\arctan\dfrac{2(x_{1}^{}+x_{2})}{x_{1}^{}+x_{2}^{}+2x_{3}^{}+2x_{4}^{}} \right)
$
\\[2.25ex]
on any domain ${\mathscr X}$ from the set 
$\{x\colon x_{2}^{}+x_{3}^{}\ne 0,\ 
x_{1}^{}+x_{2}^{}+2x_{3}^{}+2x_{4}^{}\ne 0,\ x_1^{}\ne 0\}\subset {\mathbb R}^5.$
\vspace{1.5ex}

{\bf Theorem 2.11.}
{\it
Let the following conditions hold}:
\vspace{0.35ex}

(i) 
{\it
$\nu^{0l}$ are linearly independent common real eigenvectors of the matrices
$B_{j}^{}$ corresponding 
to the eigenvalues $\lambda_{\,l}^{j},\ l=1,\ldots, r,\ j=1,\ldots, m;$
}
\vspace{0.5ex}

(ii) 
{\it
$\!\nu^{\theta l}\!,\ \theta\!=\!1,\ldots, s_{l}^{}\!-\!1,\!$
\vspace{0.35ex}
are linearly independent real generalized eigenvectors 
of the matrix $B_{\zeta}^{}$
corresponding to the eigenvalue $\lambda_{\, l}^{\zeta}$
with elementary divisor of multiplicity $s_{l}^{},\ l=1,\ldots, r,$ 
such that $\sum\limits_{l=1}^{r}s_{l}^{}\geq m+1;$
}
\\[-0.75ex]
\indent
(iii) 
{\it the linear differential system $\dfrac{dx}{dt}=A_{\zeta}^{}\;\!x$ 
has no the first integrals of the form
\\[2ex]
\mbox{}\hfill                                                            
$
\displaystyle
F_{{}_{\scriptstyle j\;\!\theta l}}^{\,\zeta}\colon x\to\ 
{\frak a}_j^{}\;\!\Psi_{\theta l}^{\zeta}(x) 
$ 
\ for all 
$
x\in {\mathscr X}, 
\quad
j=1,\ldots, m, \ \, j\ne\zeta, 
\quad 
\theta=1,\ldots, s_{l}^{}-1,\ \, l=1,\ldots, r.
\hfill 
$
\\[2ex]
Then autonomous first integrals of 
the Lappo-Da\-ni\-lev\-s\-kii system {\rm (2.20)} are the functions
\\[2ex]
\mbox{}\hfill                                       % (2.48)
$
\displaystyle
F\colon x\to
\prod\limits_{\xi=1}^{k}\bigl(\nu^{0\xi} x\bigr)^{h_{0 \xi}^{}}
\exp\sum\limits_{q=1}^{\varepsilon_{\xi}}h_{q\xi}^{}\;\! \Psi_{q\xi}^{\;\!\zeta}(x)
$ 
\ for all 
$
x\in {\mathscr X},
\quad 
{\mathscr  X}\subset {\rm D}F,
$
\hfill {\rm (2.48)}
\\[2ex]
where the scalar functions
$\Psi_{q\xi}^{\;\!\zeta}\colon {\mathscr X}\to {\mathbb R}$ 
are the solution to the functional system
\\[1.5ex]
\mbox{}\hfill                                            %(2.49)
$
\displaystyle
\nu^{q\;\!\xi} x = \sum_{\rho=1}^{q}
{\textstyle\binom{q-1}{\rho-1}}\,\Psi_{\rho\;\!\xi}^{\;\!\zeta}(x)\,
\nu^{\;\!q-\rho{,}\,\xi} x
$
\ for all 
$
x\in {\mathscr X},
\quad 
q=1,\ldots, \varepsilon_{\xi}^{},
\quad 
\xi =1,\ldots, k,
$
\hfill {\rm (2.49)}
\\[2ex]
and $\sum\limits_{\xi=1}^{k}\varepsilon_{\xi}^{}=m-k+1, \
\varepsilon_{\xi}^{}\leq s_{\xi}^{}-1, \ \xi=1,\ldots, k,\ k\leq r.$ 
Also, here 
\\[1.75ex]
\mbox{}\hfill                             % (2.50)
$
{\frak  a}_{j}^{}\;\!\Psi_{q\xi}^{\zeta}(x) = \mu_{q\xi}^{j\zeta},
\ \ \mu_{q\xi}^{j\zeta} ={\rm const},
\ \ j=1,\ldots, m,
\ \ q=1,\ldots, \varepsilon_{\xi}^{},
\ \ \xi=1,\ldots, k,
$
\hfill {\rm (2.50)}
\\[2.25ex]
the numbers $h_{q\xi}^{},\ q=0,\ldots, \varepsilon_{\xi}^{},
\ \xi=1,\ldots, k,$ are an nontrivial solution to the system
\\[1.5ex]
\mbox{}\hfill                                   % (2.51)
$
\displaystyle
\sum\limits_{\xi=1}^{k}
\bigl(\lambda_{\xi}^{j}\,h_{0\xi}^{}
+ \sum\limits_{q=1}^{\varepsilon_{\xi}^{}}
\mu_{q\xi}^{j\zeta}\,h_{q\xi}^{}\big) = 0.
\quad 
j=1,\ldots, m,
$
\hfill {\rm(2.51)}
\\[1.5ex]
}
\indent
{\sl Proof.} 
\vspace{0.35ex}
The functional system (2.49) with any fixed index $\xi\in \{1,\ldots, k\}$
has the determinant 
$
\bigl(\nu^{0\xi} x\bigr)^{\varepsilon_{\xi}^{}}
$
\vspace{0.75ex}
for all $x\in {\mathbb R}^{n}$
such that 
$\bigl(\nu^{0\xi} x\bigr)^{\varepsilon_{\xi}^{}}\ne 0$
for all $x\in {\mathscr X},$ where 
a domain ${\mathscr X}\subset \{x\colon \nu^{0\xi}x\ne 0,\ \xi=1,\ldots, k\}.$ 
\vspace{0.75ex}
Therefore for any fixed index $\xi\in \{1,\ldots, k\}$
there exists the solution 
$\Psi_{\rho\;\!\xi}^{\;\!\zeta}\colon {\mathscr X}\to {\mathbb R},\
q=1,\ldots, \varepsilon_{\xi}^{},$
\vspace{0.75ex}
on the domain ${\mathscr X}$ of the functional system (2.49).

For any fixed index $\xi\in \{1,\ldots, k\}$
\vspace{0.25ex}
from the results of Lemma 2.5 it follows that the 
identities (2.50) are satisfied.

Consequently there exist $\sum\limits_{\xi=1}^{k}\varepsilon_{\xi}^{}$ 
functionally independent on the domain ${\mathscr X}$ scalar functions
$
\Psi_{_{\scriptstyle q \xi}}^{\,\zeta}\colon {\mathscr X}\to {\mathbb R},
\ q=1,\ldots, \varepsilon_{\xi}^{},\ \xi=1,\ldots, k,
$
such that the identities (2.50) hold.
 \vspace{0.75ex}

Since the functions
\vspace{0.25ex}
$
p_\xi^{}\colon x\to\nu^{0\xi}x
$  
for all 
$
x\in {\mathbb R}^{n},
\ \xi=1,\ldots, k,
$
are partial integrals of 
the linear differential systems $\dfrac{dx}{dt}=A_j^{}\;\!x,\ j=1,\ldots, m,$ we have 
\\[2ex]
\mbox{}\hfill
$
\displaystyle
{\frak a}_{j}^{}\;\! F(x) =
\biggl(\, \sum\limits_{\xi=1}^{k}\;\!
\bigl(\lambda_{\xi}^{j}\,h_{0\xi}^{} +
\sum\limits_{q=1}^{\varepsilon_{\xi}^{}}\mu_{q\xi}^{j\zeta}\,
h_{q\xi}^{}\bigr)\!\biggr)\;\! F(x)
$
\ for all 
$
x\in {\mathscr X},
\quad 
j=1,\ldots, m.
\hfill
$
\\[2ex]
\indent
If the numbers $h_{q\xi}^{},\ q=0,\ldots, \varepsilon_{\xi}^{},\ \xi=1,\ldots, k,$ 
\vspace{0.5ex}
are an nontrivial solution to the linear homogeneous system (2.51),
then the function (2.48) is a first integral of system (2.34).
\vspace{0.35ex}

By Lemma 2.6, the scalar function (2.48) is autonomous first integral on the domain
${\mathscr X}$ of the Lappo-Da\-ni\-lev\-s\-kii differential system (2.20). \k
\vspace{1ex}

{\bf Example 2.17.}
Consider the fourth-order Lappo-Da\-ni\-lev\-s\-kii differential system
\\[1.75ex]
\mbox{}\hfill       % (2.52)                     
$
\dfrac{dx}{dt}=\alpha_1^{}(t)\,A_1^{}\;\!x+\alpha_2^{}(t)\,A_2^{}\;\!x,
\quad
x\in {\mathbb R}^4,
$
\mbox{}\hfill (2.52)
\\[1.75ex]
where linearly independent functions 
\vspace{0.25ex}
$\alpha_1^{}\colon J\to {\mathbb R}$ and $\alpha_2^{}\colon J\to {\mathbb R}$
are continuous on an interval $J\subset {\mathbb R},$ 
the commuting matrices
\\[2ex]
\mbox{}\hfill
$
A_1^{}\! =\!
\left\|\!\!
\begin{array}{rcrr}
 0&  1&  0& 0
\\
 0&   2&   {}-1& {}-1
\\
 1&   0&   0& {}-1
\\
{}-1&  0&  2&  2
\end{array}\!\! \right\|
$
\ \ and \ \ 
$
A_2^{}\! =\!  \left\|\!\!
\begin{array}{rcrc}
 2&  0&  {}-1& 0
\\
 {}-1&  2& 0&  1
\\
 {}-1& 0& 3& 1
\\
0& 1& {}-3& 1
\end{array}\!\!\right\|\!.
\hfill
$
\\[2ex]
\indent
The matrix $B_1^{}=A_1^{T},$ where $T$ denotes the matrix transpose, 
\vspace{0.5ex}
has the eigenvalue $\lambda_1^1=1$ with 
elementary divisor $(\lambda^1-1)^4$
\vspace{0.5ex}
corresponding to the eigenvector $\nu^{01}=({}-1,1,{}-1,0)$ and to the 
\vspace{0.5ex}
generalized eigenvectors
$\nu^{11}=(1,0,{}-1,{}-1),$ 
$\nu^{21}=(1,{}-1,3,0),\, \nu^{31}=({}-3,0,9,9).$

From the functional system (2.49), we get the scalar functions
\\[2.25ex]
\mbox{}\hfill
$
\Psi_{11}^1\colon x\to\  \dfrac{x_1^{}-x_3^{}-x_4^{}}{{}-x_1^{}+x_2^{}-x_3^{}}
$
\ for all 
$
x\in {\mathscr X},
\hfill
$
\\[2.5ex]
\mbox{}\hfill
$
\Psi_{21}^1\colon x\to\
\dfrac{({}-x_1^{}+x_2^{}-x_3^{})(x_1^{}-x_2^{}+3x_3^{})-(x_1^{}-x_3^{}-x_4^{})^2 }
{({}-x_1^{}+x_2^{}-x_3^{})^{2}}
$
\ for all 
$
x\in {\mathscr X},
\hfill
$
\\[2.5ex]
\mbox{}\hfill
$
\Psi_{31}^1\colon x\to\
\dfrac{1}{({}-x_1^{}+x_2^{}-x_3^{})^{3}}\,
\bigl(({}-3x_1^{}+9x_3^{}+9x_4^{})({}-x_1^{}+x_2^{}-x_3^{})^2 \ -
\hfill
$
\\[2.15ex]
\mbox{}\hfill
$
- \ 3({}-x_1^{}+x_2^{}-x_3^{})(x_1^{}-x_3^{}-x_4^{})
(x_1^{}-x_2^{}+3x_3^{}) +2(x_1^{}-x_3^{}-x_4^{})^3\bigr)
$
\ for all 
$
x\in {\mathscr X},
\hfill
$
\\[2.5ex]
where ${\mathscr X}$ is any domain from the space 
$\{ x\colon x_1^{}-x_2^{}+x_3^{}\ne 0\}$ of the space ${\mathbb R}^4.$
\vspace{0.75ex}

The Lappo-Da\-ni\-lev\-s\-kii differential system (2.52) 
is induced on the space ${\mathbb R}^4$ the linear differential operators of first order
\\[1.5ex]
\mbox{}\hfill
$
{\frak a}_{1}^{}(x)=
x_2^{}\;\!\partial_{x_1^{}}^{} + (2x_2^{}-x_3^{}-x_4^{})\;\!\partial_{x_2^{}}^{}  +
(x_1^{}-x_4^{})\;\partial_{x_3^{}}^{}  +  ({}-x_1^{}+2x_3^{}+2x_4^{})\;\!\partial_{x_4^{}}^{}
\hfill
$
\\[1ex]
and 
\\[1.25ex]
\mbox{}\hfill
$
{\frak  a}_{2}^{}(x)  =
(2x_1^{}-x_3^{})\;\!\partial_{x_1^{}}^{} +
({}-x_1^{}+2x_2^{}+x_4^{})\,\partial_{x_2^{}}^{} +
({}-x_1^{}+3x_3^{}+x_4^{})\,\partial_{x_3^{}}^{} + (x_2^{}-3x_3^{}+x_4^{})\,\partial_{x_4^{}}^{}.
\hfill
$
\\[-1.75ex]

\newpage 

Using the identities 
\\[1.5ex]
\mbox{}\hfill
$
{\frak a}_{1}^{}\nu^{01} x =\nu^{01}x
$
\ for all 
$
x\in {\mathbb R}^4,
\quad \ 
{\frak a}_{2}^{}\;\!\nu^{01}x = 2\;\!\nu^{01}x
$
\ for all 
$
x\in {\mathbb R}^4,
\hfill
$
\\[2ex]
\mbox{}\hfill
$
{\frak a}_{1}^{} \Psi_{11}^1(x)=1,
\quad \ 
{\frak a}_{1}^{} \Psi_{21}^1(x)= 0,
\quad \ 
{\frak a}_{1}^{} \Psi_{31}^1(x)=0
$
\ for all 
$
x\in {\mathscr X},
\hfill
$
\\[2.25ex]
\mbox{}\hfill
$
{\frak a}_{2}^{} \Psi_{11}^1(x) =  {}-1,
\quad  \
{\frak a}_{2}^{} \Psi_{21}^1(x)=0,
\quad  \
{\frak a}_{2}^{} \Psi_{31}^1(x)=6
$
\ for all 
$
x\in {\mathscr X},
\hfill
$
\\[1.75ex]
we obtain the linear homogeneous systems 
\\[1.75ex]
\mbox{}\hfill
$
h_{11}^{}+ h_{21}^{} = 0,
\ \ 
2h_{11}^{}- h_{21}^{} = 0
$
\ \ and \ \ 
$
h_{12}^{}+ h_{22}^{} = 0,
\ \ 
2h_{12}^{} - h_{22}^{} +6h_{32}^{} = 0.
\hfill
$
\\[1.75ex]
From these systems it follows that
\vspace{0.5ex}
$h_{11}^{}\!=h_{21}^{}\!=0,\, h_{31}^{}\!=1,$ and 
$h_{12}^{}\!=2,\,  h_{22}^{}\!={}-2,\, h_{32}^{}\!={}-1.$

Then, by Theorem 2.11, the scalar functions
\\[2ex]
\mbox{}\hfill                                       
$
F_{1}^{}\colon x\to \Psi_{21}^1(x)
$ 
\ and \
$
F_{2}^{}\colon x\to\, ({}-x_{1}^{}+x_{2}^{}-x_{3}^{})^{2}
\exp\bigl({}-2\;\!\Psi_{11}^1(x)-\Psi_{31}^1(x)\bigr)
$
\ for all 
$
x\in {\mathscr X}
\hfill
$
\\[2ex]
are autonomous first integrals of the Lappo-Da\-ni\-lev\-s\-kii differential system (2.52).
\vspace{1ex}

If the matrices $B_{j}^{}$ have some  
complex common eigenvectors $\nu^{0\,l}$ corresponding to the 
eigenvalues $\lambda_{\,l}^{\zeta}$ with elementary divisors $s_{l}^{},$
then the proof of Theorem 2.11 is also true. 
\vspace{0.5ex}

Let the set $V$ of $m+1$ functions be given by
\\[1.25ex]
\mbox{}\hfill
$
\displaystyle
V\!=\!\Bigl\{\nu^{0\xi}x\ \text{for all}\, x\in{\mathbb R}^n,\, 
\Psi_{q\xi}^{\zeta}(x)\ \text{for all}\, x\in{\mathscr X}\colon  
q\!=1,\ldots,\varepsilon_{\xi}^{},\, \xi\!=1,\ldots, k,\,
\sum\limits_{\xi=1}^{k}\varepsilon_{\xi}^{}=m-k+1\!\Bigr\}. 
\hfill
$
\\[1.25ex]
\indent
In the complex case, we shall have two logical possibilities:

1. Any function from the set $V$ has 
the complex conjugate function in the set $V.$
\vspace{0.25ex}

2.\! At least one function from the set $\!V\!$ has not the complex conjugate function 
\vspace{0.5ex}
in the $\!V.\!\!$

In these cases the system (2.20) has 
\vspace{0.75ex}
the following autonomous first integrals.

{\sl Case}\ 1. The Lappo-Da\-ni\-lev\-s\-kii differential system (2.20) 
has the autonomous first integral
\\[1.75ex]
\mbox{}\hfill
$
\displaystyle
F\colon x\to \
\prod\limits_{\xi=1}^{k_1}
\Bigl(\bigl({\stackrel{*}{\nu}}{}^{0\xi} x\bigr)^{2} +
\bigl(\,{\stackrel{\sim}{\nu}}\,{}^{0\xi} x\bigr)^{2}\,
\Bigr)^{{\stackrel{*}{h}}_{0\xi}}
\exp\Bigl({}-2\,{\stackrel{\sim}{h}}_{0\xi}^{}
\arctan\dfrac{{\stackrel{\sim}{\nu}}\,{}^{0\xi} x}
{{\stackrel{*}{\nu}}{}^{0\xi} x} \ +
\hfill
$
\\[2ex]
\mbox{}\hfill
$
\displaystyle
+ \ 2\, \sum\limits_{q=1}^{\varepsilon_{\xi}}
\bigl(\,{\stackrel{*}{h}}_{q\xi}^{}\,{\stackrel{*}{\Psi}}{}_{q\xi}^{\,\zeta}(x) -
{\stackrel{\sim}{h}}_{q\xi}^{}\,{\stackrel{\sim}{\Psi}}{}_{q\xi}^{\,\zeta}(x)
\bigr)\Bigr)\, 
\prod\limits_{\rho=1}^{k_2}\;\!
\bigl|\nu^{0\rho} x\bigr|^{h_{0\rho}^{}}
\exp\sum\limits_{q=1}^{\varepsilon_{\rho}}
h_{q\rho}^{}\,\Psi_{q\rho}^{\,\zeta}(x)
$ 
\ for all 
$
x\in {\mathscr X}\subset {\rm D}F,
\hfill
$
\\[1.75ex]
where 
\vspace{0.75ex}
$\stackrel{*}{h}_{q\xi}^{},\ \stackrel{\sim}{h}_{q\xi}^{},$ and
$h_{q\rho}^{},\ q=0,\ldots, \varepsilon_k^{},\ k=\xi$
or $k=\rho, \ \xi=1,\ldots, k_1^{},$ $\rho=1,\ldots, k_2^{},$
are a real nontrivial solution to
the linear homogeneous system 
\\[1.75ex]
\mbox{}\hfill
$
\displaystyle
2\sum\limits_{\xi=1}^{k_1}
\Bigl( \bigl(
{\stackrel{*}{\lambda}}{}_{\xi}^{j}\,{\stackrel{*}{h}}_{0\xi}^{} -
{\stackrel{\sim}{\lambda}}{}_{\xi}^{j}\,{\stackrel{\sim}{h}}_{0\xi}^{}
\bigr) \,+
\sum\limits_{q=1}^{\varepsilon_{\xi}}
\bigl(
{\stackrel{*}{\mu}}{}_{q\xi}^{j\zeta}\,{\stackrel{*}{h}}_{q\xi}^{} -
{\stackrel{\sim}{\mu}}{}_{q\xi}^{j\zeta}\,{\stackrel{\sim}{h}}_{q\xi}^{}
\bigr)
\!\Bigr) \,+\,
\sum\limits_{\rho=1}^{k_2}
\Bigl(
\lambda_{\rho}^{j}\,h_{0\rho} +
\sum\limits_{q=1}^{\varepsilon_{\rho}}
\mu_{q\rho}^{j\zeta}\,h_{q\rho}^{} \Bigr)\! =\! 0,
\ j=1,\ldots, m.
\hfill
$
\\[1.5ex]
\indent
Here $\nu^{0\xi}={\stackrel{*}{\nu}}{}^{0\xi} +\widetilde{\nu}\,{}^{0\xi}\,i$ are complex 
common eigenvectors of the matrices $B_j^{}$ corresponding to the eigenvalues
\vspace{0.35ex}
$\lambda_{\xi}^{j} ={\stackrel{*}{\lambda}}{}_{\xi}^{j} +\widetilde{\lambda}{}_{\xi}^{j}\,i,\ 
j=1,\ldots, m,\ \xi=1,\ldots, k_1^{},$ respectively; $\nu^{0\rho}$ are real 
common eigenvectors of $B_j^{}$ corresponding to the eigenvalues
\vspace{0.35ex}
$\lambda_{\rho}^{j}, \ j=1,\ldots, m,\ \rho=1,\ldots, k_2^{},$ respectively;
the functions 
$\Psi_{q\xi}^{\,\zeta}={\stackrel{*}{\Psi}}{}_{q\xi}^{\,\zeta} +
{\stackrel{\sim}{\Psi}}{}_{q\xi}^{\,\zeta}\,i$ and $\Psi_{q\rho}^{\,\zeta}$
\vspace{0.75ex}
are the solution to the functional system (2.49); 
the real numbers 
\\[1.75ex]
\mbox{}\hfill
$
{\stackrel{*}{\mu}}{}_{q\xi}^{j\zeta} =
{\rm Re}\,\bigl({\frak a}_{j}^{}\, \Psi_{q\xi}^{\,\zeta}(x)\bigr),
\quad
{\stackrel{\sim}{\mu}}{}_{q\xi}^{j\zeta} =
{\rm Im}\,\bigl({\frak a}_{j}^{}\, \Psi_{q\xi}^{\,\zeta}(x)\bigr),
\quad
{\mu}_{q\rho}^{j\zeta}={\frak a}_{j}^{}\, \Psi_{q\rho}^{\,\zeta}(x)
$
\ for all 
$
x\in {\mathscr X},
\hfill
$
\\[2.5ex]
\mbox{}\hfill
$
q=1,\ldots, \varepsilon_k^{},\ \ k=\xi
$
\, or \,
$
k=\rho,
\ \ \xi=1,\ldots, k_1^{},
\ \ \rho=1,\ldots, k_2^{},
\ \ \ j=1,\ldots, m;
\hfill
$
\\[1.75ex]
the natural numbers $\varepsilon_{\xi}^{}$ and $\varepsilon_{\rho}^{}$ such that 
\vspace{0.35ex}
$
2\sum\limits_{\xi=1}^{k_1^{}}\varepsilon_{\xi}
+ \sum\limits_{\rho=1}^{k_2^{}}
\varepsilon_{\rho}^{} = m-2k_{1}^{}-k_{2}^{}+1$
with
$
2k_1^{}+k_2^{}\leq r,$ 
$
\varepsilon_{\xi}^{}\leq s_{\xi}^{}-1, \ \xi=1,\ldots, k_1^{}, \
\varepsilon_{\rho}^{}\leq s_{\rho}^{}-1,\ \rho=1,\ldots, k_2^{},$
\vspace{0.5ex}
where $k_1^{}\!$ is an number of pairs of 
complex-conjugate common eigen\-vec\-tors  
of  the matrices $B_j^{},$ and  
\vspace{0.35ex}
$k_2^{}$ is an number of real common eigenvectors of the matrices $B_j^{}.$
\vspace{1.25ex}

{\bf Example 2.18.}
Consider the sixth-order Lappo-Da\-ni\-lev\-s\-kii differential system 
\\[2ex]
\mbox{}\hfill       % (2.53)                     
$
\dfrac{dx}{dt}=\alpha_1^{}(t)A_1^{}\;\!x+\alpha_2^{}(t)A_2^{}\;\!x+\alpha_3^{}(t)A_3^{}\;\!x,
\quad
x\in {\mathbb R}^6,
$
\mbox{}\hfill (2.53)
\\[2.25ex]
where $\alpha_j^{}\colon J\to {\mathbb R},\ j=1,2,3,$ 
\vspace{0.5ex}
are linearly independent continuous functions on an interval $J\subset {\mathbb R},$
the commuting matrices
\\[2.25ex]
\mbox{}\hfill
$
B_1^{}\! =A_1^{T}\! =
\left\| \!\!
\begin{array}{rrrrrr}
   3 &-1 & -3 & 3 & 5 &-2 \\
  -4 & 3 &  5 &-6 &-5 & 5 \\
   4 &-3 & -5 & 4 & 8 &-4 \\
   1 & 0 & -1 & 4 & 3 &-3 \\
   0 &-2 & -2 &-1 & 3 & 1 \\
   2 &-3 & -4 & 5 & 6 &-2
\end{array}
\right\|, 
\quad
B_2^{}\! = A_2^{T}\! =
\left\| \!\!
\begin{array}{rrrrrr}
 0 &   1 & 0 & 2 & 1 &-2 \\
-4 &   3 & 6 &-6 &-6 & 4 \\
 2 &   0 &-2 & 2 & 3 &-3 \\
 1 &   0 &-1 & 3 & 2 &-3 \\
-1 &   1 & 2 &-4 &-2 & 2 \\
 1 &  -1 &-1 & 2 & 2 &-2
\end{array}
\right\|, 
\hfill
$
\\[2.75ex]
\mbox{}\hfill
and \
$
B_3^{} =  A_3^{T} =
\left\|\!\!
\begin{array}{rrrrrr}
-\,3 &   2 &   3 &-\,3 &-\,3 &   2 \\
   2 &-\,3 &-\,3 &   2 &   3 &-\,1 \\
-\,4 &   3 &   5 &-\,6 &-\,6 &   4 \\
-\,3 &   3 &   4 &-\,4 &-\,4 &   2 \\
   0 &-\,1 &   0 &-\,1 &-\,1 &   1 \\
-\,2 &   2 &   2 &-\,1 &-\,2 &   0
\end{array}\!\! \right\|
$
\ $(T$ denotes the matrix transpose).
\hfill\mbox{}
\\[2ex]
\indent
Using the eigenvalue $\lambda_{1}^{1}=1+2i$ with elementary divisor $(\lambda^1-1-2i)^3$
\vspace{0.35ex}
of the matrix $B_1^{}$ and the corresponding complex eigenvector
\vspace{0.35ex}
$\nu^{01}=(1, 0, 1+i, 1, i, 1),$ the generalized eigen\-vec\-tor of the 1-st order
\vspace{0.5ex}
$\nu^{11}=(1, 1+i, 0, 0, i, i),$ the generalized eigen\-vec\-tor of the 2-nd order
$\nu^{21}=(2+2i, 0, 2+2i, 0, 2i, 2i),$ 
we can construct the scalar functions
\\[2ex]
\mbox{}\hfill
$
{\stackrel{*}{\Psi}}{}_{11}^{\,1}\colon x\to\
\dfrac{(x_1^{}+x_2^{})(x_1^{}+x_3^{}+x_4^{}+x_6^{})\,+\,(x_3^{}+x_5^{})(x_2^{}+x_5^{}+x_6^{})}
{(x_1^{}+x_3^{}+x_4^{}+x_6^{})^2\,+\,(x_3^{}+x_5^{})^2}
$
\ for all  
$
x\in {\mathscr X},
\hfill
$
\\[2.75ex]
\mbox{}\hfill
$
\widetilde{\Psi}{}_{11}^{\,1}\colon x\to\
\dfrac{(x_1^{}+x_3^{}+x_4^{}+x_6^{})(x_2^{}+x_5^{}+x_6^{})\,-\,(x_1^{}+x_2^{})(x_3^{}+x_5^{})}
{(x_1^{}+x_3^{}+x_4^{}+x_6^{})^2\,+\,(x_3^{}+x_5^{})^2}
$
\ for all  
$
x\in {\mathscr X},
\hfill
$
\\[2.75ex]
\mbox{}\hfill
$
{\stackrel{*}{\Psi}}{}_{21}^{\,1}\colon x\to\
\dfrac{1}{(x_1^{}+x_3^{}+x_4^{}+x_6^{})^2\,+\,(x_3^{}+x_5^{})^2}\,
\Bigl(
2\bigl( (x_1^{}+x_3^{}+x_4^{}+x_6^{})^2 + (x_3^{}+x_5^{})^2 \bigr)\, \cdot
\hfill
$
\\[2.25ex]
\mbox{}\hfill
$
\cdot \, \bigl(
(x_1^{}+x_3^{})(x_1^{}+x_3^{}+x_4^{}+x_6^{})\,+\,(x_3^{}+x_5^{})(x_1^{}+x_3^{}+x_5^{}+x_6^{})
\bigr) \ - \
\bigl((x_1^{}+x_2^{})(x_1^{}+x_3^{}+x_4^{}+x_6^{})\ +
\hfill
$
\\[2.25ex]
\mbox{}\hfill
$
+\,(x_3^{}+x_5^{})(x_2^{}+x_5^{}+x_6^{})\bigr)^2\ + \
\bigl((x_1^{}+x_3^{}+x_4^{}+x_6^{})(x_2^{}+x_5^{}+x_6^{})\,-\,(x_1^{}+x_2^{})(x_3^{}+x_5^{})\bigr)^2\,
\Big),
\hfill
$
\\[2.75ex]
\mbox{}\hfill
$
\widetilde{\Psi}{}_{21}^{\,1}\colon x\to\
\dfrac{2}{(x_1^{}+x_3^{}+x_4^{}+x_6^{})^2\,+\,(x_3^{}+x_5^{})^2}\,
\Big(
\bigl((x_1^{}+x_3^{}+x_4^{}+x_6^{})^2  +  (x_3^{}+x_5^{})^2 \bigr)\,\cdot
\hfill
$
\\[2.25ex]
\mbox{}\hfill
$
\cdot\, 
\bigl((x_1^{}+x_3^{}+x_4^{}+x_6^{})(x_1^{}+x_3^{}+x_5^{}+x_6^{})\,-\,(x_1^{}+x_3^{})(x_3^{}+x_5^{})\bigr) 
\, +\, \bigl( (x_1^{}+x_2^{})(x_1^{}+x_3^{}+x_4^{}+x_6^{})\ +
\hfill
$
\\[2.25ex]
\mbox{}\hfill
$
+\,  (x_3^{}\!+x_5^{})(x_2^{}\!+x_5^{}\!+x_6^{})\bigr)
\bigl((x_3^{}\!+x_5^{})(x_1^{}\!+x_2^{}) - (x_1^{}\!+x_3^{}\!+x_4^{}+x_6^{})(x_2^{}\!+x_5^{}\!+x_6^{})
\bigr)\!\Bigr)\!
$
for all 
$\!x\in\! {\mathscr X},
\hfill
$
\\[2ex]
where ${\mathscr X}$ is any domain from the set 
\vspace{0.5ex}
$\{x\colon x_1^{}+x_3^{}+x_4^{}+x_6^{}\ne 0\}$ of the space ${\mathbb R}^6.$

\newpage

Then, the Lappo-Da\-ni\-lev\-s\-kii differential system (2.53) 
has the autonomous first integrals
\\[1.75ex]
\mbox{}\hfill
$
F_{1}^{}\colon x\to\
P(x)\exp\Bigl({}- 4\arctan\dfrac{x_3^{}+x_5^{}}{x_1^{}+x_3^{}+x_4^{}+x_6^{}} +
6{\stackrel{*}{\Psi}}\!{}_{11}^{\,1}(x) + 2\widetilde{\Psi}\!{}_{11}^{\,1}(x)\Bigr),
\hfill
$
\\[2.5ex]
\mbox{}\hfill
$
F_{2}^{}\colon x\to\
P^2(x)\exp\Bigl({} -2\arctan\dfrac{x_3^{}+x_5^{}}{x_1^{}+x_3^{}+x_4^{}+x_6^{}} +
{\stackrel{*}{\Psi}}{}_{21}^{\,1}(x) - \widetilde{\Psi}{}_{21}^{\,1}(x)\Bigr),
\hfill
$
\\[1.5ex]
and
\\[1ex]
\mbox{}\hfill
$
F_{3}^{}\colon x\to\
2\;\!\widetilde{\Psi}{}_{11}^{\,1}(x)\, -\, 2\;\!{\stackrel{*}{\Psi}}{}_{21}^{\,1}(x)
\, -\, \widetilde{\Psi}{}_{21}^{\,1}(x)
$
\ for all $x\in {\mathscr X},
\hfill
$
\\[2ex]
where the polynomial
\vspace{1.25ex}
$P\colon x\to (x_1^{}+x_3^{}+x_4^{}+x_6^{})^2 + (x_3^{}+x_5^{})^2$ for all $x\in {\mathbb R}^6.$

{\sl Case}\ 2A. 
\vspace{0.35ex}
A common complex eigenvector of the matrices $B_{j}, \ j=1,\ldots, m,$ 
hasn't the com\-p\-lex conjugate vector.
Autonomous first integrals of the system (2.20) are the functions
\\[1.5ex]
\mbox{}\hfill
$
\displaystyle
F_{1}^{}\colon  x\to \
\prod\limits_{\xi=1}^{k_1^{}}
\bigl(P_{\xi}^{}(x)\bigr)^{{\stackrel{*}{h}}_{0\xi}^{}+
{\stackrel{*}{h}}_{0,(k_1^{}+\xi)}^{}}
\exp\Bigl({}-2\bigl(\,\widetilde{h}_{0\xi}^{}-
\widetilde{h}_{0,(k_1^{}+\xi)}^{}\bigr)\varphi_{\xi}^{}(x)\ +
\hfill
$
\\[1.5ex]
\mbox{}\hfill
$
\displaystyle
+\ 2\sum\limits_{q=1}^{\varepsilon_{\xi}^{}}
\Bigl(
\bigl(\,
{\stackrel{*}{h}}_{q\xi}^{}+{\stackrel{*}{h}}_{q,(k_1^{}+\xi)}^{}
\bigr)\,{\stackrel{*}{\Psi}}{}_{q\xi}^{\,\zeta}(x) +
\bigl(\,
\widetilde{h}_{q,(k_1^{}+\xi)}^{}-\widetilde{h}_{q\xi}^{}
\bigr)\,\widetilde{\Psi}{}_{q\xi}^{\,\zeta}(x)
\Bigr)\,\Bigr)\, \cdot
\hfill
$
\\[1.5ex]
\mbox{}\hfill
$
\displaystyle
\cdot\, 
\bigl(P_{2k_1^{}+1}^{}(x)\bigr)^{{\stackrel{*}{h}}_{0,(2k_1^{}+1)}^{}}
\exp\Bigl({}-2\;\!
\widetilde{h}_{0,(2k_1^{}+1)}^{}\;\!\varphi_{2k_1^{}+1}^{}(x)\Bigr)
\prod\limits_{\rho=1}^{k_2^{}}
\bigl(\nu^{0\rho} x\bigr)^{2\;\!{\stackrel{*}{h}}_{0\rho}^{}}\;\!
\exp\Bigl( 2\sum\limits_{q=1}^{\varepsilon_{\rho}^{}}
{\stackrel{*}{h}}_{q\rho}^{}\;\!\Psi_{q\rho}^{\,\zeta}(x)\Bigr)
\hfill
$
\\[1ex]
and
\\[1ex]
\mbox{}\hfill
$
\displaystyle
F_{2}^{}\colon x\to \
\prod\limits_{\xi=1}^{k_1^{}}
\bigl(P_{\xi}^{}(x)\bigr)^{\widetilde{h}_{0\xi}^{}+
\widetilde{h}_{0,(k_1^{}+\xi)}^{}}
\exp\Big( 2\bigl(\, {\stackrel{*}{h}}_{0\xi}^{}-
{\stackrel{*}{h}}_{0,(k_1^{}+\xi)}^{}\bigr)\,\varphi_{\xi}^{}(x)\ +
\hfill
$
\\[1.5ex]
\mbox{}\hfill
$
\displaystyle
+\ 2\sum\limits_{q=1}^{\varepsilon_{\xi}^{}}
\Bigl(
\bigl(\,
\widetilde{h}_{q\xi}^{}+\widetilde{h}_{q,(k_1^{}+\xi)}^{}
\bigr)\,{\stackrel{*}{\Psi}}{}_{q\xi}^{\,\zeta}(x) +
\bigl(\,
{\stackrel{*}{h}}_{q\xi}^{}-{\stackrel{*}{h}}_{q,(k_1^{}+\xi)}^{}
\bigr)\,\widetilde{\Psi}{}_{q\xi}^{\,\zeta}(x)
\Bigr)
\Bigr)\, \cdot
\hfill
$
\\[1.5ex]
\mbox{}\hfill
$
\displaystyle
\cdot\,
\bigl(P_{2k_1^{}+1}^{}(x)\bigr)^{\widetilde{h}_{0,(2k_1^{}+1)}^{}}
\exp\Bigl( 2\;\!{\stackrel{*}{h}}_{0,(2k_1^{}+1)}^{}\;\!\varphi_{2k_1^{}+1}^{}(x)
\bigr)
\prod\limits_{\rho=1}^{k_2^{}}
\bigl( \nu^{0\rho} x \bigr)^{2\;\!\widetilde{h}_{0\rho}^{}}
\exp\Bigl( 2\sum\limits_{q=1}^{\varepsilon_{\rho}^{}}
\widetilde{h}_{q\rho}^{}\;\!\Psi_{q\rho}^{\,\zeta}(x)
\Bigr)
$
\\[2.5ex]
\mbox{}\hfill
for all 
$
x\in {\mathscr X},
\quad
{\mathscr X}\subset {\rm D}F_1^{}\cap {\rm D}F_2^{},
\hfill
$
\\[2.25ex]
where the polynomials
\vspace{0.5ex}
$
P_{\xi}^{}\colon x \to
\bigl({\stackrel{*}{\nu}}{}^{0\xi} x \bigr)^{2} +
\bigl(\,\widetilde{\nu}\,{}^{0\xi}x \bigr)^{2}
$
for all $x\in {\mathbb R}^n,$
the scalar functions
$
\varphi_{\xi}^{}\colon x\to
\arctan\dfrac{\widetilde{\nu}\,{}^{0\xi} x}
{{\stackrel{*}{\nu}}{}^{0\xi} x}
$
for all  
$
x\in {\mathscr X},
\ \xi=1,\ldots, k_1^{},
\ \xi=2k_1^{}+1,
$
\vspace{0.5ex}
the numbers
$h_{q\xi}^{}={\stackrel{*}{h}}_{q\xi}^{}+\widetilde{h}_{q\xi}^{}\,i,$
$h_{q\rho}^{}={\stackrel{*}{h}}_{q\rho}^{}+\widetilde{h}_{q\rho}^{}\,i,\
q=0,\ldots, \varepsilon_k^{},\ k=\xi$ or
\vspace{1ex}
$k=\rho, \ \xi=1,\ldots, 2k_1^{}+1,\ \rho=1,\ldots, k_2^{},$
are an nontrivial solution to the linear homogeneous system 
\\[1.25ex]
\centerline{
$
\displaystyle
\sum\limits_{\xi=1}^{2k_1^{}}
\Bigl(
\lambda_{\xi}^{j}\;\!h_{0\xi}^{} +
\sum\limits_{q=1}^{\varepsilon_{\xi}^{}}
\mu_{q\xi}^{j\zeta}\;\!h_{q\xi}^{}
\Bigr)\, +\,
\lambda_{2k_1^{}+1}^{j}\;\!h_{0,(2k_1^{}+1)}^{} +
\sum\limits_{\rho=1}^{k_2^{}}
\Bigl(
\lambda_{\rho}^{j}\;\!h_{0\rho}^{} +
\sum\limits_{q=1}^{\varepsilon_{\rho}^{}}
\mu_{q\rho}^{j\zeta}\;\!h_{q\rho}^{} \Bigr)\! =  0,
\ \ 
j=1,\ldots, m.
$
}
\\[1.5ex]
\indent
Here $\nu^{0\xi}={\stackrel{*}{\nu}}{}^{\,0\xi}\, +\widetilde{\nu}{}^{\,0\xi}\,i, \
\vspace{0.35ex}
\nu^{0,(k_1+\xi)}=\overline{\nu^{0\xi}},$ and 
$\nu^{0,(2k_1+1)}={\stackrel{*}{\nu}}{}^{\,0,(2k_1+1)}\,+\widetilde{\nu}{}^{\,\, 0,(2k_1+1)}\,i$
are com\-p\-lex common eigenvectors of the matrices $B_j^{}$ corresponding 
to the eigenvalues
\\[1.5ex]
\mbox{}\hfill
$
\lambda_{\xi}^{j}= {\stackrel{*}{\lambda}}{}_{\xi}^{j} +\widetilde{\lambda}{}_{\xi}^{j}\,i, 
\ \ 
\lambda_{k_1+\xi}^{j}=\overline{\lambda_{\xi}^{j}},\
\xi=1,\ldots, k_1^{},
$
and 
$
\lambda_{2k_1+1}^{j}={\stackrel{*}{\lambda}}{}_{2k_1+1}^{j}\,+
\widetilde{\lambda}{}_{2k_1+1}^{j}\,i,\ j=1,\ldots, m;
\hfill
$  
\\[1.5ex]
$\nu^{0\rho}$ are real common eigenvectors of 
the matrices $B_j^{}$ corresponding to the real eigenvalues
$\lambda_{\rho}^{j},\ j=1,\ldots,m,\ \rho=1,\ldots, k_2,$ respectively;
\vspace{0.75ex}
the functions  $\Psi_{q\xi}^{}= {\stackrel{*}{\Psi}}{}_{q\xi}^{\,\zeta} +
\widetilde{\Psi}{}_{q\xi}^{\,\zeta}\;\!i$ and $\Psi_{q\rho}^{\,\zeta}$
are the solution to the functional system (2.49);
the numbers 
\\[1.5ex]
\mbox{}\hfill
$
{\mu}_{q\xi}^{j\zeta}={\frak a}_{j}^{}\;\!\Psi_{q\xi}^{\,\zeta}(x),
\quad
{\stackrel{*}{\mu}}{}_{q\xi}^{j\zeta} =\, {\rm Re}\,{\mu}_{q\xi}^{j\zeta},
\quad
\widetilde{\mu}{}_{q\xi}^{j\zeta} =\, {\rm Im}\,{\mu}_{q\xi}^{j\zeta},
\quad
\mu_{q\rho}^{j\zeta} ={\frak a}_{j}^{}\;\!\Psi_{q\rho}^{\,\zeta}(x)
$
\ for all 
$
x\in {\mathscr X},
\hfill
$
\\[2.5ex]
\mbox{}\hfill
$
q=1,\ldots, \varepsilon_k^{},
\ \ k=\xi
$ 
\, or \,
$
k=\rho,\ \ \
\xi=1,\ldots, 2k_1^{},
 \ \
\rho=1,\ldots, k_2^{},
\ \
j=1,\ldots, m.
\hfill
$
\\[1.75ex]
the natural numbers $\varepsilon_{\xi}^{}$ and $\varepsilon_{\rho}^{}$ such that
\vspace{0.35ex}
$
2\sum\limits_{\xi=1}^{k_1^{}} \varepsilon_{\xi}^{} +
\sum\limits_{\rho=1}^{k_2^{}} \varepsilon_{\rho}^{} =m-2k_1^{}-k_2^{}$
with
$2k_1^{}+1+k_2^{} \leq r,$
$\varepsilon_{\xi}^{}\leq s_{\xi}^{}-1,\
\xi=1,\ldots, k_1^{},\ 
\varepsilon_{\rho}^{}\leq s_{\rho}^{}-1,\ \rho=1,\ldots, k_2^{},$
\vspace{0.5ex}
where $k_1^{}$ is an number of  complex common eigen\-vec\-tors 
{\rm(}this set hasn't complex conjugate vectors{\rm)} 
\vspace{0.35ex}
of  the matrices $B_j^{},$ and  
$k_2^{}$ is an number of real common eigenvectors of the matrices $B_j^{}.$
\vspace{1.25ex}

{\bf Example 2.19.}
Consider the sixth-order Lappo-Da\-ni\-lev\-s\-kii differential system 
\\[0.75ex]
\mbox{}\hfill       % (2.54)                     
$
\displaystyle
\dfrac{dx}{dt}\, =\,
\sum\limits_{j=1}^{4}\alpha_j^{}(t)A_j^{}\;\!x,
\quad
x\in {\mathbb R}^6,
\quad
t\in J\subset {\mathbb R},
$
\mbox{}\hfill (2.54)
\\[1ex]
with linearly independent continuous functions 
$\alpha_j^{}\colon J\to {\mathbb R},\ j=1,\ldots, 4,$ and the matrices
\\[1.5ex]
\centerline{
$
A_1^{} =
\left\| \!\!
\begin{array}{rrrrrr}
   1 &{}-2 & 2 & 0 & 1 & 1 \\
  0 & 2 & {}-2 & 0 &{}-2 & {}-2 \\
   0 & 3 & {}-2 & 0 & {}-2 &{}-2 \\
  0 & {}-4 & 0 & 2 & {}-2 & 2 \\
  2 &{}-3 & 4 & 2 & 2 & 4 \\
   {}-1 & 3 & {}-2 & {}-2 & 1 &{}-1
\end{array}
\right\|, 
\quad
A_2^{} = \left\| \!\!
\begin{array}{rrrrrr}
 0 &  2 & 0 & 0 & 1 & 1 \\
{}-1 & {}-3 & 0 & 0 &{}-1 & {}-1 \\
{}-1  & {}-3 &{}-1 & 0 & {}-2 &{}-2 \\
 2 &  2 & 0 & 1 & 0 & 4 \\
3 &  3 & 2 & 2 & 1 & 4 \\
{}-2 &  {}-1 &{}-2 & {}-2 & {}-1 &{}-4
\end{array}\!\! \right\|,
$
}
\\[2.5ex]
\centerline{
$
A_3^{} =  \left\|\!\!
\begin{array}{rrrrrr}
3 &  0 &  0 & 0 &{}-1 &  {}-1 \\
{}-1 & 2 & 0 &  0 &  1 & 1 \\
{}-2 & {} -1 & 2 & 0 & 1 & 1 \\
1 &  2 &  {}-2 & 1 & {}-1 & {}-1 \\
2 & 1 & {}-1 & 0 & 0 & {}-1 \\
{}-1 & {}-1 & 1 & 0 & 1 & 2
\end{array}\!\! \right\|,
\ \
\text{and}
\ \
A_4^{} =  \left\|\!\!
\begin{array}{rrrrrr}
1 &  {}-2 &  4 & 0 & 2 &  2 \\
{}-2 &{}-1 &{}-4 &  0 & {}-4 &{}-4 \\
{}-3 &  2 &  {}-7 & 0 &{}-5 & {}-5 \\
3 &  {}-4 & 10 & 2 & 7 & 7 \\
3 &{}-2 &  9 & 0 & 7 & 5 \\
1 & 4 & {}-5 & 0 &{}-4 & {}-2
\end{array}\!\! \right\|.
$
}
\\[2ex]
\indent
The system (2.54) has the eigenvalue 
\vspace{0.5ex}
$\lambda_{1}^{1}=1+i$ with elementary divisor $(\lambda^1-1-i)^2$  
and the eigenvalue $\lambda_{2}^1=2i$ with elementary divisor $\lambda^1-2i.$
\vspace{0.5ex}

Since the eigenvalue $\lambda_{1}^{1}=1+i$ 
\vspace{0.65ex}
corresponding to the eigen\-vec\-tor $\nu^{01}=(1, 1+i, 0, 0, i, i)$ 
of the matrices $B_j^{}=A_j^{T},\ j=1,\ldots, 4\ (T$ denotes the matrix transpose) and 
\vspace{0.5ex}
to the 1-st order generalized eigenvector  
$\nu^{11}=(1+i, 0, 1+i, 0, i, i)$ of the matrix $B_1^{},$ we have
\\[2ex]
\mbox{}\hfill
$
{\stackrel{*}{\Psi}}{}_{11}^{\,1} \colon x\to\ 
\dfrac{(x_1^{}+x_2^{})(x_1^{}+x_3^{})+(x_2^{}+x_5^{}+x_6^{})(x_1^{}+x_3^{}+x_5^{}+x_6^{})}{P_1^{}(x)}
$ 
\ for all 
$
x\in {\mathscr X},
\hfill
$
\\[2.75ex]
\mbox{}\hfill
$
\widetilde{\Psi}{}_{11}^{\,1} \colon x\to\ 
\dfrac{(x_1^{}+x_2^{})(x_1^{}+x_3^{}+x_5^{}+x_6^{})-(x_1^{}+x_3^{})(x_2^{}+x_5^{}+x_6^{})}{P_1^{}(x)}
$ 
\ for all 
$
x\in {\mathscr X},
\hfill
$
\\[2ex]
where ${\mathscr X}$ is a domain from the set 
\vspace{1ex}
$\bigl\{ x\colon x_1^{}+x_2^{}\ne 0,\ x_1^{}+x_3^{}+x_4^{}+x_6^{}\ne 0\bigr\}$ 
of the space ${\mathbb R}^6,$
the polynomial
$P_1^{}\colon x\to (x_1^{}+x_2^{})^2+(x_2^{}+x_5^{}+x_6^{})^2$ for all $x\in {\mathbb R}^6.$
\vspace{1ex}

Using the common eigenvector 
\vspace{0.65ex}
$\nu^{02}=(1, 0, 1+i, 1,i, 1)$ of the matrices $B_1^{},\ldots, B_4^{}$
cor\-res\-pon\-ding to the simple eigenvalue $\lambda_{2}^1=2i,$ we can build 
\vspace{0.35ex}
autonomous first integrals of the Lappo-Da\-ni\-lev\-s\-kii differential system (2.54) 
\\[1.75ex]
\mbox{}\hfill                                      
$
F_1^{}\colon x\to\, P_1^{}(x)\bigl(P_2^{}(x)\bigr)^2
\exp\bigl({}-10\varphi_1^{}(x) + 8\;\!{\stackrel{*}{\Psi}}{}_{11}^{\,1}(x) +
6\;\!\widetilde{\Psi}{}_{11}^{\,1}(x)\bigr)
$ 
\ for all 
$
x\in {\mathscr X},
\hfill 
$
\\[2ex]
\mbox{}\hfill                                    
$
F_2^{}\colon x\to\,
\bigl(P_1^{}(x)\bigr)^3 
\exp\bigl({}-10\varphi_1^{}(x) - 4\varphi_2^{}(x) + 12\;\!{\stackrel{*}{\Psi}}{}_{11}^{\,1}(x) +
14\;\!\widetilde{\Psi}{}_{11}^{\,1}(x)\bigr)
$
\ for all 
$
x\in {\mathscr X},
\hfill
$
\\[2.25ex]
where the polynomial
\vspace{1.25ex}
$
P_2^{}\colon x\to (x_1^{}+x_3^{}+x_4^{}+x_6^{})^2+(x_3^{}+x_5^{})^2
$
for all 
$
x\in {\mathbb R}^6,
$
the scalar functions 
$
\varphi_1^{}\colon x\to
\arctan\dfrac{x_2^{}+x_5^{}+x_6^{}}{x_1^{}+x_2^{}}\,, \ \
\varphi_2^{}\colon x\to
\arctan\dfrac{x_3^{}+x_5^{}}{x_1^{}+x_3^{}+x_4^{}+x_6^{}}
$
\vspace{1.5ex}
for all 
$
x\in {\mathscr  X}.
$

{\sl Case} 2B. 
\vspace{0.75ex}
Suppose a function $\Psi_{\beta\gamma}^{\,\zeta},\
\gamma\in \{1,\ldots,k_1^{}\},\
\beta\in \{1,\ldots,\varepsilon_{\gamma}^{}\}$
hasn't the complex conjugate fun\-c\-tion.
Then the system (2.20) has the autonomous first integrals 
\\[2ex]
\mbox{}\hfill
$
\displaystyle
F_{1}^{}\colon x\to \
\prod\limits_{\xi=1}^{k_1^{}}
\bigl(P_{\xi}^{}(x)\bigr)^{{\stackrel{*}{h}}_{0\xi}^{}+
{\stackrel{*}{h}}_{0,(k_1^{}+\xi)}^{}}
\exp\Bigl({}-2\bigl(\,\widetilde{h}_{0\xi}^{}-
\widetilde{h}_{0,(k_1^{}+\xi)}^{}\bigr)\,\varphi_{\xi}^{}(x) \ +
\hfill
$
\\[1.75ex]
\centerline{
$
\displaystyle
+\ 2\sum\limits_{q=1}^{\varepsilon_{\xi}^{}}
(1-\delta_{q\beta}\delta_{\xi\gamma})\Bigl(
\bigl(\,
{\stackrel{*}{h}}_{q\xi}^{}+{\stackrel{*}{h}}_{q,(k_1^{}+\xi)}^{}
\bigr)\,{\stackrel{*}{\Psi}}{}_{q\xi}^{\,\zeta}(x) +
\bigl(\,
\widetilde{h}_{q,(k_1^{}+\xi)}^{}-\widetilde{h}_{q\xi}^{}
\bigr)\,\widetilde{\Psi}{}_{q\xi}^{\,\zeta}(x)
\Bigr) \ +
$
}
\\[1.75ex]
\mbox{}\hfill
$
\displaystyle
+ \ 2\,\Bigl(\,{\stackrel{*}{h}}_{\beta\gamma}\,
{\stackrel{*}{\Psi}}{}_{\beta\gamma}^{\,\zeta}(x)  -
\widetilde{h}_{\beta\gamma}^{}\,
\widetilde{\Psi}{}_{\beta\gamma}^{\,\zeta}(x)\Bigr)
\Bigr)\
\prod\limits_{\rho=1}^{k_2^{}}
\bigl(\nu^{0\rho} x \bigr)^{2\,{\stackrel{*}{h}}_{0\rho}^{}}
\exp\biggl( 2\sum\limits_{q=1}^{\varepsilon_{\rho}^{}}
{\stackrel{*}{h}}_{q\rho}^{}\, \Psi_{q\rho}^{\,\zeta}(x)\biggr)
$ 
\ for all 
$
x\in {\mathscr X}
\hfill
$
\\[1.75ex]
and
\\[2ex]
\mbox{}\hfill
$
\displaystyle
F_{2}^{}\colon x\to \
\prod\limits_{\xi=1}^{k_1^{}}
\bigl(P_{\xi}^{}(x)\bigr)^{\widetilde{h}_{0\xi}^{}+
\widetilde{h}_{0,(k_1^{}+\xi)}^{}}
\exp\Bigl( 2\bigl(\, {\stackrel{*}{h}}_{0\xi}^{}-
{\stackrel{*}{h}}_{0,(k_1^{}+\xi)}^{}\bigr)\,\varphi_{\xi}^{}(x)\ +
\hfill
$
\\[1.75ex]
\centerline{
$
\displaystyle
+\ 2\sum\limits_{q=1}^{\varepsilon_{\xi}^{}}
(1-\delta_{q\beta}^{}\delta_{\xi\gamma}^{})\Bigl(
\bigl(\,
\widetilde{h}_{q\xi}^{} +
\widetilde{h}_{q,(k_1^{}+\xi)}^{}\bigr)\,
{\stackrel{*}{\Psi}}{}_{q\xi}^{\,\zeta}(x) +
\bigl(\,
{\stackrel{*}{h}}_{q\xi}^{}-{\stackrel{*}{h}}_{q,(k_1^{}+\xi)}^{}
\bigr)\,
\widetilde{\Psi}_{q\xi}^{\,\zeta}(x)
\Bigr) \ +
$
}
\\[1.75ex]
\mbox{}\hfill
$
\displaystyle
+\ 2\bigl(\, {\stackrel{*}{h}}_{\beta\gamma}^{}\,
\widetilde{\Psi}{}_{\beta\gamma}^{\,\zeta}(x) +
\widetilde{h}_{\beta\gamma}^{}\,
{\stackrel{*}{\Psi}}{}_{\beta\gamma}^{\,\zeta}(x)\bigr)
\Bigr)\ 
\prod\limits_{\rho=1}^{k_2^{}}
\bigl(\nu^{0\rho}x\bigr)^{2\;\!\widetilde{h}_{0\rho}^{}}
\exp\biggl( 2\sum\limits_{q=1}^{\varepsilon_{\rho}^{}}
\widetilde{h}_{q\rho}^{}\;\!\Psi_{q\rho}^{\,\zeta}(x)\biggr)
$ 
\ for all 
$
x\in {\mathscr X},
\hfill
$
\\[2ex]
where ${\mathscr X}$ is a domain from the set ${\rm D}F_1^{}\cap {\rm D}F_2^{},$
\vspace{0.5ex}
the polynomials
$
P_{\xi}^{}\colon x \to
\bigl({\stackrel{*}{\nu}}{}^{0\xi} x \bigr)^{2}  +
\bigl(\,\widetilde{\nu}\;\!{}^{0\xi} x\bigr)^{2}\;\;
$
for all $x\in {\mathbb R}^n,$
\vspace{0.5ex}
the scalar functions
$
\varphi_{\xi}^{}\colon x\to
\arctan\dfrac{\widetilde{\nu}\,{}^{0\xi} x}{{\stackrel{*}{\nu}}{}^{0\xi} x}
$
for all $x\in {\mathscr X},
\ \xi=1,\ldots, k_1^{},
$
the numbers
$
h_{q\xi}^{}\!={\stackrel{*}{h}}_{q\xi}^{}+\widetilde{h}_{q\xi}^{}\,i,\ 
h_{q\rho}^{}\!={\stackrel{*}{h}}_{q\rho}^{}+\widetilde{h}_{q\rho}^{}\,i,
\ q=0,\ldots, \varepsilon_k^{},\ k=\xi$ 
or
\vspace{1ex}
$k=\rho, \
\xi=1,\ldots, 2k_1^{},$
$
\rho=1,\ldots, k_2^{},$
are an nontrivial solution to the linear homogeneous system 
\\[1.5ex]
\centerline{
$
\displaystyle
\sum\limits_{\xi=1}^{2k_1^{}}
\Bigl(
\lambda_{\xi}^{j}\;\!h_{0\xi}^{}+
\sum\limits_{q=1}^{\varepsilon_{\xi}^{}}
\mu_{q\xi}^{j\zeta}\;\!h_{q\xi}^{}\Bigr) -\,
\mu_{\beta,(k_1^{}+\gamma)}^{j\zeta}\;\!h_{\beta,(k_1^{}+\gamma)} +
\sum\limits_{\rho=1}^{k_2^{}}
\Bigl(
\lambda_{\rho}^{j}\;\!h_{0\rho}^{} +
\sum\limits_{q=1}^{\varepsilon_{\rho}^{}}
\mu_{q\rho}^{j\zeta}\;\! h_{q\rho}^{}\Bigr)\!  = 0,
\ j=1,\ldots, m,
$
}
\\[1.5ex]
and $\delta$ is the Kronecker symbol.
\vspace{0.75ex}

Here 
\vspace{0.35ex}
$\nu^{0\xi}={\stackrel{*}{\nu}}{}^{\,0\xi}\, +\widetilde{\nu}{}^{\,0\xi}\,i, \
\nu^{0,(k_1+\xi)}=\overline{\nu^{0\xi}}$
are complex common eigenvectors of the matrices $B_j^{}$ corresponding to the eigenvalues
\vspace{0.5ex}
$\lambda_{\xi}^{j}={\stackrel{*}{\lambda}}{}_{\xi}^{j} +\widetilde{\lambda}{}_{\xi}^{j}\,i,\ 
\lambda_{k_1+\xi}^{j}=\overline{\lambda_{\xi}^{j}},\ j=1,\ldots, m,\ \xi=1,\ldots, k_1^{},$ respectively; 
$\nu^{0\rho}$ are real common eigenvectors of the matrices $B_j^{}$ corresponding 
to the eig\-en\-va\-lu\-es
\vspace{0.5ex}
$\lambda_{\rho}^{j},\ j=1,\ldots,m,\ \rho=1,\ldots, k_2^{},$ respectively;
the scalar functions 
$
\Psi_{q\xi}^{\,\zeta}={\stackrel{*}{\Psi}}{}_{q\xi}^{\,\zeta} +\widetilde{\Psi}{}_{q\xi}^{\,\zeta}\;\!i$
and $\Psi_{q\rho}^{\,\zeta}$
are the solution to the system (2.49); the numbers 
\\[2ex]
\mbox{}\hfill
$
{\mu}_{q\xi}^{j\zeta} = {\frak a}_{j}^{}\;\! \Psi_{q\xi}^{\,\zeta}(x),
\quad
{\stackrel{*}{\mu}}{}_{q\xi}^{j\zeta} =\,
{\rm Re}\,{\mu}_{q\xi}^{j\zeta},
\quad
\widetilde{\mu}{}_{q\xi}^{j\zeta} =\,
{\rm Im}\,{\mu}_{q\xi}^{j\zeta},
\quad
\mu_{q\rho}^{j\zeta} = {\frak a}_{j}^{}\;\!\Psi_{q\rho}^{\,\zeta}(x)
$
\ for all 
$
x\in {\mathscr X},
\hfill
$
\\[2.75ex]
\mbox{}\hfill
$
q=1,\ldots, \varepsilon_k^{},
\ \ k=\xi
$ 
\, or \, 
$
k=\rho,\ \ 
\xi=1,\ldots, 2k_1^{},
\ \ \rho=1,\ldots, k_2^{},
\ \ \ j=1,\ldots, m;
\hfill
$
\\[2ex]
the natural numbers $\varepsilon_{\xi}^{}$ and $\varepsilon_{\rho}^{}$ such that
\vspace{0.35ex}
$
2\sum\limits_{\xi=1}^{k_1^{}} \varepsilon_{\xi}^{} +
\sum\limits_{\rho=1}^{k_2^{}} \varepsilon_{\rho} =m-2k_{1}^{}-k_{2}^{}+2
$
with
$2k_1^{}+k_2^{} \leq r,$ 
$\varepsilon_{\xi}^{}\leq s_{\xi}^{}-1,\
\xi=1,\ldots, k_1^{},\ \varepsilon_{\rho}\leq s_{\rho}^{}-1, \ \rho=1,\ldots, k_2^{},$
\vspace{0.5ex}
where $k_1^{}$ is an number of complex common eigen\-vec\-tors 
{\rm(}this set hasn't complex conjugate vectors{\rm)} 
\vspace{0.25ex}
of  the matrices $B_j^{},$ and  
$k_2^{}$ is an number of real common eigenvectors of the matrices $B_j^{}.$
\vspace{1.25ex}

{\bf Example 2.20.}
Consider the sixth-order Lappo-Da\-ni\-lev\-s\-kii differential system 
\\[2ex]
\mbox{}\hfill
$
\dfrac{dx}{dt}=
\alpha_1^{}(t)
\left\| \!\!
\begin{array}{rrrrrr}
   3\! &-4\! & 4\! & 1\! & 0\! & 2 \\
  -1\! & 3\! & -3\! & 0\! &-2\! &-3 \\
  -3\! & 5\! & -5\! & -1\! & -2\! &-4 \\
   3\! & -6\! & 4\! & 4\! & -1\! & 5 \\
   5\! &-5\! & 8\! & 3\! & 3\! & 6 \\
  -2\! & 5\! & -4\! & -3\! & 1\! &-2
\end{array}
\right\| x\;\! +\;\!
\alpha_2^{}(t)
\left\| \!\!
\begin{array}{rrrrrr}
 0\! &  -4\! & 2\! & 1\! & -1\! & 1 \\
 1\! &  3\! & 0\! & 0\! & 1\! & -1 \\
 0\! &  6\! &-2\! & -1\! & 2\! &-1 \\
 2\! & -6\! & 2\! & 3\! & -4\! & 2 \\
 1\! &  -6\! & 3\! & 2\! &-2\! & 2 \\
 -2\! &  4\! &-3\! & -3\! & 2\! &-2
\end{array}
\right\| x,
\ x\!\in\! {\mathbb R}^{6},
\hfill
$
\\[2.25ex]
where linearly independent functions
\vspace{0.5ex}
$\alpha_1^{}\colon J\to {\mathbb R}$ and $\alpha_2^{}\colon J\to {\mathbb R}$ 
are continuous.

Using the complex eigenvalue 
\vspace{0.5ex}
$\lambda_{1}^{1}=1+2i$ with elementary divisor $(\lambda^1\!-1-2i)^3$ of the 
matrix $B_1^{}=A_1^{T}\ (T$ denotes the matrix transpose) 
\vspace{0.5ex}
corresponding to the common eigen\-vec\-tor $\nu^{01}=(1, 0, 1+i, 1, i, 1),$
\vspace{0.5ex}
to the 1-st order generalized eigenvector 
$\nu^{11}=(1, 1+i, 0, 0, i, i),$ and 
to the 2-nd order generalized eigenvector 
\vspace{0.35ex}
$\nu^{21}=(2+2i, 0, 2+2i, 0, 2i, 2i),$ we can construct the 
autonomous first integrals of this  Lappo-Da\-ni\-lev\-s\-kii differential system 
\\[1.75ex]
\mbox{}\hfill                                      
$
F_1^{}\colon x\to\
P(x)\exp\bigl({} -\varphi(x)- \widetilde{\Psi}{}_{11}^{\,}(x)\bigr),
\quad
F_2^{}\colon x\to\
P(x)\exp\bigl({}-2\varphi(x) + 2\,{\stackrel{*}{\Psi}}{}_{11}^{\,1}(x)\bigr),
\hfill 
$
\\[2.25ex]
\mbox{}\hfill                                          
$
F_3^{}\colon x\to\,
P^2(x)\exp\big({}-2\varphi(x)- \widetilde{\Psi}_{21}^{\,1}(x)\bigr),
$
\ and \
$
F_4^{}\colon x\to\ {\stackrel{*}{\Psi}}{}_{21}^{\,}(x)
$ 
\ for all 
$
x\in {\mathscr X},
\hfill 
$
\\[2.5ex]
where the polynomial  
\vspace{1.25ex}
$P\colon x\to (x_1^{}+x_3^{}+x_4^{}+x_6^{})^2 + (x_3^{}+x_5^{})^2
$
for all $x\in {\mathbb R}^6,$
the scalar functions 
$
\varphi\colon x\to \arctan\dfrac{x_3^{}+x_5^{}}{x_1^{}+x_3^{}+x_4^{}+x_6^{}} 
$
for all 
$
x\in {\mathscr X},
$
\\[2.5ex]
\mbox{}\hfill
$
{\stackrel{*}{\Psi}}{}_{11}^{\,1}\colon x\to\
\dfrac{(x_1^{}+x_2^{})(x_1^{}+x_3^{}+x_4^{}+x_6^{})\,+\,(x_3^{}+x_5^{})(x_2^{}+x_5^{}+x_6^{})}
{(x_1^{}+x_3^{}+x_4^{}+x_6^{})^2\,+\,(x_3^{}+x_5^{})^2}
$
\ for all  
$
x\in {\mathscr X},
\hfill
$
\\[2.75ex]
\mbox{}\hfill
$
\widetilde{\Psi}{}_{11}^{\,1}\colon x\to\
\dfrac{(x_1^{}+x_3^{}+x_4^{}+x_6^{})(x_2^{}+x_5^{}+x_6^{})\,-\,(x_1^{}+x_2^{})(x_3^{}+x_5^{})}
{(x_1^{}+x_3^{}+x_4^{}+x_6^{})^2\,+\,(x_3^{}+x_5^{})^2}
$
\ for all  
$
x\in {\mathscr X},
\hfill
$
\\[2.75ex]
\mbox{}\hfill
$
{\stackrel{*}{\Psi}}{}_{21}^{\,1}\colon x\to\
\dfrac{1}{(x_1^{}+x_3^{}+x_4^{}+x_6^{})^2\,+\,(x_3^{}+x_5^{})^2}\,
\Bigl(
2\bigl( (x_1^{}+x_3^{}+x_4^{}+x_6^{})^2 + (x_3^{}+x_5^{})^2 \bigr)\, \cdot
\hfill
$
\\[2.25ex]
\mbox{}\hfill
$
\cdot \, \bigl(
(x_1^{}+x_3^{})(x_1^{}+x_3^{}+x_4^{}+x_6^{})\,+\,(x_3^{}+x_5^{})(x_1^{}+x_3^{}+x_5^{}+x_6^{})
\bigr) \ - \
\bigl((x_1^{}+x_2^{})(x_1^{}+x_3^{}+x_4^{}+x_6^{})\ +
\hfill
$
\\[2.25ex]
\mbox{}\hfill
$
+\,(x_3^{}+x_5^{})(x_2^{}+x_5^{}+x_6^{})\bigr)^2\ + \
\bigl((x_1^{}+x_3^{}+x_4^{}+x_6^{})(x_2^{}+x_5^{}+x_6^{})\,-\,(x_1^{}+x_2^{})(x_3^{}+x_5^{})\bigr)^2\,
\Big),
\hfill
$
\\[2.75ex]
\mbox{}\hfill
$
\widetilde{\Psi}{}_{21}^{\,1}\colon x\to\
\dfrac{2}{(x_1^{}+x_3^{}+x_4^{}+x_6^{})^2\,+\,(x_3^{}+x_5^{})^2}\,
\Big(
\bigl((x_1^{}+x_3^{}+x_4^{}+x_6^{})^2  +  (x_3^{}+x_5^{})^2 \bigr)\,\cdot
\hfill
$
\\[2.25ex]
\mbox{}\hfill
$
\cdot\, 
\bigl((x_1^{}+x_3^{}+x_4^{}+x_6^{})(x_1^{}+x_3^{}+x_5^{}+x_6^{})\,-\,(x_1^{}+x_3^{})(x_3^{}+x_5^{})\bigr) 
\, +\, \bigl( (x_1^{}+x_2^{})(x_1^{}+x_3^{}+x_4^{}+x_6^{})\ +
\hfill
$
\\[2.25ex]
\mbox{}\hfill
$
+\,  (x_3^{}\!+x_5^{})(x_2^{}\!+x_5^{}\!+x_6^{})\bigr)
\bigl((x_3^{}\!+x_5^{})(x_1^{}\!+x_2^{}) - (x_1^{}\!+x_3^{}\!+x_4^{}+x_6^{})(x_2^{}\!+x_5^{}\!+x_6^{})
\bigr)\!\Bigr)\!
$
for all 
$\!x\in\! {\mathscr X},
\hfill
$
\\[2ex]
where ${\mathscr X}$ is any domain from the set 
\vspace{0.5ex}
$\{x\colon x_1^{}+x_3^{}+x_4^{}+x_6^{}\ne 0\}$ of the space ${\mathbb R}^6.$

\newpage

\mbox{}
\\[-1.75ex]
\indent
{\bf 2.3.2. Linear nonhomogeneous differential system}
\\[1ex]
\indent
Consider an nonhomogeneous Lappo-Da\-ni\-lev\-s\-kii differential system
\\[1.5ex]
\mbox{}\hfill                                          % (2.55)
$
\displaystyle
\dfrac{dx}{dt}=\sum\limits_{j=1}^{m} \alpha_j^{}(t)A_j^{}\,x+f(t),
\ \quad 
x\in {\mathbb R}^n,
\quad 
t\in J,
\quad 
f\in C(J),
$
\hfill (2.55)
\\[1.25ex]
where linearly independent on an interval $J\subset {\mathbb R}$ functions 
$\alpha_j^{}\colon J\to {\mathbb R}$ are continuous,
real constant $n\times n$ matrices $A_j^{}$ such that 
$A_j^{}A_k^{}=A_k^{}A_j^{},\ j =1,\ldots, m,\ k=1,\ldots, m.$
\vspace{0.35ex}

The cor\-res\-pon\-ding homogeneous differential system of system (2.55) is the system (2.20).

Using Theorems 2.12, 2.13, and 2.14, we can build first integrals of the 
nonhomogeneous Lappo-Da\-ni\-lev\-s\-kii differential system (2.55). 
\vspace{0.5ex}

{\bf Theorem 2.12.}
{\it
Suppose that the conditions of Lemma {\rm 2.3} hold. 
Then the Lappo-Da\-ni\-lev\-s\-kii differential system {\rm(2.55)} has the first integral
\\[1.5ex]
\mbox{}\hfill                                            % (2.56)
$
\displaystyle
F\colon (t,x)\to \
\nu x\,\varphi (t) -
\int\limits_{t_0^{}}^{t}\nu f(\tau)\,\varphi (\tau)\, d\tau
$
\ for all 
$
(t,x)\in J\times {\mathbb R}^n,
$
\hfill {\rm (2.56)}
\\[1.5ex]
where $t_0^{}$ is a fixed point from the interval $J,$ the exponential function}
\\[1.5ex]
\mbox{}\hfill                                          % (2.57)
$
\displaystyle
\varphi\colon t\to \
\exp\biggl({}-\int\limits_{t_0^{}}^{t}\sum\limits_{j=1}^m \lambda^j\;\!\alpha_j(\tau)\, d\tau\biggr)
$
\ for all 
$
t\in J.
$
\hfill {\rm (2.57)}
\\[1.5ex]
\indent
{\sl Proof}. 
By Lemma 2.3, it follows that 
the Lie derivative of the function (2.56) by virtue of 
the Lappo-Da\-ni\-lev\-s\-kii differential system (2.55) is equal to
\\[1.5ex]
\mbox{}\hfill
$
\displaystyle
{\frak B}F(t,x)\!=
\nu x\;\!\partial_t^{}\varphi(t)- \partial_t^{}\int\limits_{t_0^{}}^{t}\!\nu f(\tau)\,\varphi (\tau)\, d\tau
+\bigl({\frak A}\;\!\nu x + f(t)\;\!\partial_{x}^{}\nu x\bigr)\varphi(t)\!=0
$
for all 
$\!
(t,x)\!\in\! J\!\times\! {\mathbb R}^n\!,
\hfill
$
\\[2ex]
where the linear differential operator
\vspace{0.5ex}
$
{\frak B}(t,x)={\frak A}(t,x)+f(t)\;\!\partial_{x}^{}
$
for all $(t,x)\in J\times {\mathbb R}^{n}$
is induced by the Lappo-Da\-ni\-lev\-s\-kii differential system (2.55). \k
\vspace{1.25ex}

{\bf Example 2.21.}
The second-order Lappo-Da\-ni\-lev\-s\-kii differential system
\\[2ex]
\mbox{}\hfill                               % (2.58)
$
\dfrac{dx_1^{}}{dt}={}-3\Bigl(\dfrac{1}{t}+t^2\Bigr)\, x_1^{}-\dfrac{2}{t}\, x_2^{} +
t\bigl(\sin t+2t\;\!e^{t^2}\bigr)\;\!e^{{}-t^3},
\hfill
$
\\[0.5ex]
\mbox{}\hfill (2.58)
\\[0.5ex]
\mbox{}\hfill 
$
\dfrac{dx_2^{}}{dt}=\dfrac{4}{t}\, x_1^{}+3\Bigl(\dfrac{1}{t}-t^2\Bigr)\, x_2^{} +
t\Bigl(\cosh t+6\arctan\dfrac{t}{3}\Bigr)\;\!e^{{}-t^3}
\hfill
$
\\[2ex]
such that the coefficient matrix
\vspace{1ex}
$A(t)=3t^2A_1^{}+\Bigl(\dfrac{1}{t}+6t^2\Bigr)A_2^{}$ for all $t\in J\subset \{t\colon t\ne 0\},$
where the constant matrices
\vspace{1ex}
$
A_1^{} = \left\|\!\!
\begin{array}{rr}
5& 4
\\
{}-8& {}-7
\end{array}\!\!
\right\|
$
and 
$
A_2^{} = \left\|\!\!
\begin{array}{rr}
{}-3 & {}-2
\\
4 & 3
\end{array}\!\!
\right\|\,.
$
\vspace{0.5ex}

The matrices $B_1^{}=A_1^{T}$ and $B_2^{}=A_2^{T},$ 
\vspace{0.75ex}
where $T$ denotes the matrix transpose, have the eigenvalues
\vspace{0.5ex}
$\lambda^1_1=1,\ \lambda_2^1={}-3,$ and $\lambda_1^2={}-1,\,\lambda_2^2=1$
corresponding to the linearly independent common real eigenvectors
$\nu^1=(2, 1),\ \nu^2=(1, 1).$
\vspace{0.5ex}

Using the scalar functions 
\\[1.5ex]
\mbox{}\hfill
$
\alpha_1^{}\colon t \to\, 3t^2
$
\ for all 
$
t\in J,
\quad \
\alpha_2^{}\colon t\to\ \dfrac{1}{t}+6t^2
$
\ for all 
$
t\in J,
\hfill
$ 
\\[2.5ex]
\mbox{}\hfill 
$
f_1^{}\colon t \to \,
t\bigl(\sin t+2t\;\!e^{t^2}\bigr)\;\!e^{{}-t^3},
\quad
f_2^{}\colon t \to \
t\Bigl(\cosh t+6\arctan\dfrac{t}{3}\Bigr)\;\!e^{{}-t^3}
$
\ for all 
$
t\in J,
\hfill
$
\\[2ex]
\mbox{}\hfill
$
\displaystyle
\varphi_1^{}(t)=
\exp\biggl({}-\int\bigl(\lambda_1^1\;\!\alpha_1(t)+\lambda_1^2\;\!\alpha_2^{}(t)\bigr)\;\!dt\biggr)=  
\exp\int\Bigl(\dfrac{1}{t}+3t^2\Bigr)\;\!dt =  t\exp t^3
$
for all 
$
t\in J,
\hfill
$ 
\\[2ex]
\mbox{}\hfill
$
\displaystyle
\varphi_2^{}(t)=
\exp\biggl({}-\int\bigl(\lambda_2^1\;\!\alpha_1(t)+\lambda_2^2\;\!\alpha_2^{}(t)\bigr)\;\!dt\biggr)=  
\exp\int\Bigl(3t^2-\dfrac{1}{t}\Bigr)\;\!dt =  \dfrac{1}{t}\,\exp t^3
$
for all 
$
t\in J,
\hfill
$ 
\\[2ex]
\mbox{}\hfill
$
\displaystyle
I_1^{}(t)=
\int\bigl(2f_1^{}(t)+f_2^{}(t)\bigr)\varphi_1^{}(t)\;\!dt=  
\int\Bigl(
2t^2\sin t+4t^3e^{t^2}+t^2\cosh t+6t^2\arctan\dfrac{t}{3}\Bigr)\;\!dt=  
{}-3t^2\ +
\hfill
$ 
\\[2ex]
\mbox{}\hfill
$
+\ 4t\sin t-2(t^2-2)\cos t +2(t^2-1)e^{t^2}-
2t\cosh t+(t^2+2)\sinh t  +2t^3\arctan\dfrac{t}{3}+27\ln(t^2+9),
\hfill
$ 
\\[2ex]
\mbox{}\hfill
$
\displaystyle
I_2^{}(t)=
\int\bigl(f_1^{}(t)+f_2^{}(t)\bigr)\varphi_2^{}(t)\;\!dt=  
\int\Bigl(
\sin t+2te^{t^2}+\cosh t+6\arctan\dfrac{t}{3}\Bigr)\;\!dt=  
\hfill
$ 
\\[2ex]
\mbox{}\hfill
$
={}-\cos t +e^{t^2}+\sinh t+6t\arctan\dfrac{t}{3}-9\ln(t^2+9)
$
\ for all 
$
t\in J,
\hfill
$ 
\\[2ex]
we can build (by Theorem 2.12) first integrals of 
the Lappo-Da\-ni\-lev\-s\-kii system (2.58)
\\[2ex]
\mbox{}\hfill
$
F_1^{}\colon (t,x_1^{},x_2^{})\to\  
t\;\!e^{t^3}(2x_1^{}+x_2^{})+3t^2- 4t\sin t+2(t^2-2)\cos t -2(t^2-1)e^{t^2}\ +
\hfill
$ 
\\[2ex]
\mbox{}\hfill
$
+\ 2t\cosh t-(t^2+2)\sinh t  - 2t^3\arctan\dfrac{t}{3}-27\ln(t^2+9)
$ 
\ for all 
$
(t,x_1^{},x_2^{})\in J\times {\mathbb R}^2
\hfill
$ 
\\[1.5ex]
and 
\\[1.5ex]
\mbox{}\hfill
$
F_2^{}\colon (t,x_1^{},x_2^{})\to  \
\dfrac{1}{t}\;\!e^{t^3}(x_1^{}+x_2^{})+
\cos t -e^{t^2}-\sinh t  -6t\arctan\dfrac{t}{3}+9\ln(t^2+9)
\hfill
$
\\[2ex]
\mbox{}\hfill
for all 
$
(t,x_1^{},x_2^{})\in J\times {\mathbb R}^2.
\hfill
$
\\[1.5ex]
\indent
The functionally independent first integrals $F_1^{}$ and $F_2^{}$ are an integral basis of
\vspace{0.35ex}
the Lappo-Da\-ni\-lev\-s\-kii differential system (2.58) on any domain $J\times {\mathbb R}^2.$
\vspace{0.75ex}

{\bf Theorem 2.13.}
{\it
Suppose that the conditions of Lemma {\rm 2.4} are satisfied. 
Then the Lappo-Da\-ni\-lev\-s\-kii differential system {\rm (2.55)} has the first integrals
\\[1.25ex]
\mbox{}\hfill                                           % (2.59)
$
\displaystyle
F_{\rho}^{}\colon (t,x)\to\  
\gamma_{\rho}^{}(t, x)- \int\limits_{t_0^{}}^{t} \gamma_{\rho}^{}(t, f(t))\;\!dt
$
\ for all 
$
(t,x)\in J\times {\mathbb R}^n,
\quad
\rho=1, 2,
$
\hfill {\rm (2.59)}
\\[1.25ex]
where $t_0^{}$ is a fixed point from the interval $J,$ the scalar functions 
\\[1.25ex]
\mbox{}\hfill                                     
$
\displaystyle
\gamma_1^{}\colon (t,x)\to
\biggl({\stackrel{*}{\nu}}x
\cos\int\limits_{t_0^{}}^{t}\!
\sum\limits_{j=1}^{m}\widetilde{\lambda}{}^j\;\!\alpha_j^{}(\tau)\;\!d\tau  +
\widetilde{\nu}x\sin\int\limits_{t_0^{}}^{t}\!
\sum\limits_{j=1}^{m}\widetilde{\lambda}{}^j\;\!\alpha_j^{}(\tau)\;\!d\tau\biggr)
\exp\biggl(\!\!{}-\int\limits_{t_0^{}}^{t}\! 
\sum\limits_{j=1}^{m}{\stackrel{*}{\lambda}}{}^j\;\!\alpha_j^{}(\tau)\;\!d\tau\!\biggr),
\hfill
$
\\[2ex]
\mbox{}\hfill                                     
$
\displaystyle
\gamma_2^{}\colon (t,x)\to
\biggl(\widetilde{\nu}x
\cos\int\limits_{t_0^{}}^{t}\! 
\sum\limits_{j=1}^{m}\widetilde{\lambda}{}^j\;\!\alpha_j^{}(\tau)\;\!d\tau  -
{\stackrel{*}{\nu}}x\sin\int\limits_{t_0^{}}^{t}\! 
\sum\limits_{j=1}^{m}\widetilde{\lambda}{}^j\;\!\alpha_j^{}(\tau)\;\!d\tau\biggr)
\exp\biggl(\!\!{}-\int\limits_{t_0^{}}^{t}\! 
 \sum\limits_{j=1}^{m}{\stackrel{*}{\lambda}}{}^j\;\!\alpha_j^{}(\tau)\;\!d\tau\!\biggr).
\hfill
$
}
\\[1.5ex]
\indent
{\sl Proof}.
The idea of the proof of Theorem 2.13 is analogous to the proof of 
Corollary 1.6 or Theorem 2.2.
Formally using Theorem 2.12, we obtain 
the complex-valued function (2.56) is a first integral 
of system (2.55). Then the real and imaginary parts 
of this complex-valued first integral are the real first integrals (2.59) of 
the Lappo-Da\-ni\-lev\-s\-kii system (2.55). \k
\vspace{0.75ex}

{\bf Example 2.22.}
Consider the second-order Erugin system of differential equations
(right-hand side of system satisfy the Cauchy --- Riemann equations) [62; 48, pp. 152 -- 153] 
\\[1.75ex]
\mbox{}\hfill                               % (2.60)
$
\dfrac{dx_1^{}}{dt}=\alpha_1^{}(t)\, x_1^{}+\alpha_2^{}(t)\, x_2^{}+f_1^{}(t),
\qquad
\dfrac{dx_2^{}}{dt}={}-\alpha_2^{}(t)\, x_1^{} + \alpha_1^{}(t)\, x_2^{}+f_2^{}(t),
$
\hfill (2.60)
\\[2ex]
where functions  
\vspace{0.35ex}
$\alpha_j^{}\colon J\to {\mathbb R}$ and  $f_j^{}\colon J\to {\mathbb R},\ j=1,2$ 
are continuous on an interval $J\subset {\mathbb R}.$

Since $\nu^1=(1, i)$ and $\nu^2=(1, {}-i)$ are common eigenvectors of the matrices
\vspace{1ex}
$B_1^{}=
\left\|\!\!
\begin{array}{cc}
1 & 0
\\
0 & 1
\end{array}\!\!
\right\|
$ 
and 
$B_2^{}=
\left\|\!\!
\begin{array}{cc}
0 & {}-1
\\
1 & 0
\end{array}\!\!
\right\|
$
corresponding to the eigenvalues 
\vspace{1ex}
$
\lambda^1_1=1,\
\lambda_1^2={}-i,
$
and 
$
\lambda_2^1=1,\ \,\lambda_2^2=i,
$ 
respectively,
we see that
the Erugin system (2.60) has the integral basis (by Theorem 2.13)
\\[2ex]
\mbox{}\hfill
$
\displaystyle
F_\rho^{}\colon (t,x_1^{},x_2^{})\to
\gamma_{\rho}^{}(t,x_1^{},x_2^{})  -
\int\limits_{t_0^{}}^{t}\gamma_{\rho}^{}(\tau,f_1^{}(\tau),f_2^{}(\tau))\;\!d\tau
$
for all 
$
(t,x_1^{},x_2^{})\in J\times {\mathbb R}^2,
\  
\rho=1,2,
\hfill
$
\\[1.5ex]
where $t_0^{}$ is a fixed point from the interval $J,$ the scalar functions 
on the domain $J\times {\mathbb R}^2$
\\[1.5ex]
\mbox{}\hfill
$
\displaystyle
\gamma_1^{}\colon (t,x_1^{},x_2^{})\to \
\biggl(\!x_1^{}\cos\int\limits_{t_0^{}}^{t}\!\!\alpha_2^{}(t)\;\!dt \;\!- \;\!
x_2^{}\sin\int\limits_{t_0^{}}^{t}\!\!\alpha_2^{}(t)\;\!dt\!\biggr)
\exp\biggl(\!\!{}-\int\limits_{t_0^{}}^{t}\!\!\alpha_1^{}(t)\;\!dt\!\biggr),
\hfill
$ 
\\[2ex]
\mbox{}\hfill
$
\displaystyle
\gamma_2^{}\colon (t,x_1^{},x_2^{})\to \
\biggl(\!x_1^{}\sin\int\limits_{t_0^{}}^{t}\!\!\alpha_2^{}(t)\;\!dt \;\!+\;\!
x_2^{}\cos\int\limits_{t_0^{}}^{t}\!\!\alpha_2^{}(t)\;\!dt\!\biggr)
\exp\biggl(\!\!{}-\int\limits_{t_0^{}}^{t}\!\!\alpha_1^{}(t)\;\!dt\!\biggr).
\hfill
$ 
\\[1.5ex]
\indent
For example, if 
\vspace{0.75ex}
$\alpha_1^{}\colon t\to \dfrac{1}{t}\,,\ \alpha_2^{}\colon t\to 2t,\ 
f_1^{}\colon t\to t,\ f_2^{}\colon t\to 2t^2$ for all $t\in J\!\subset\! \{t\colon t\!\ne\! 0\},$
then the Erugin system (2.60) has the functionally independent first integrals
\\[1.5ex]
\mbox{}\hfill
$
\displaystyle
F_1^{}\colon (t,x_1^{},x_2^{})\to  \
\dfrac{1}{t}\,(x_1^{}\cos t^2- x_2^{}\sin t^2)  -\cos t^2-
\int\limits_{t_0^{}}^{t}\cos\tau^2\;\!d\tau
$
\ for all 
$
(t,x_1^{},x_2^{})\in J\times {\mathbb R}^2
\hfill
$ 
\\[1ex]
and
\\[1ex]
\mbox{}\hfill
$
\displaystyle
F_2^{}\colon (t,x_1^{},x_2^{})\to  \
\dfrac{1}{t}\,(x_1^{}\sin t^2+ x_2^{}\cos t^2)  -\sin t^2-
\int\limits_{t_0^{}}^{t}\sin\tau^2\;\!d\tau
$ 
\ for all 
$
(t,x_1^{},x_2^{})\in J\times {\mathbb R}^2.
\hfill
$ 
\\[1.5ex]
Note that these first integrals of  the Erugin system (2.60) are nonelementary functions.
\vspace{1.25ex}

{\bf Lemma 2.7.}
{\it 
Under the conditions of Lemma {\rm 2.5}, we have
\\[2ex]
\mbox{}\hfill                                  % (2.61)                                
$
{\frak a}_j^{}\;\!\nu^{\,\theta}x=
{\displaystyle \sum\limits_{\rho=0}^{\theta} }
\binom{\theta}{\rho}\,\mu_{\rho}^{j\zeta}\,\nu^{\,\theta-\rho}x
$
\ for all 
$
x\in {\mathscr X},
\quad
j=1,\ldots, m, 
\quad 
\theta=1,\ldots, s-1,
$
\hfill {\rm(2.61)}
\\[2ex]
where the numbers 
$\mu_{0}^{j\zeta}=\lambda^{j}, \ \mu_{\theta}^{j\zeta}={\frak a}_j^{}\Psi_{\theta}^{\zeta}(x), \  
j=1,\ldots, m, \ \theta=1,\ldots, s-1.$ 
}
\vspace{1.25ex}

{\sl Proof}. The proof of Lemma 2.7 is by induction on $s.$
\vspace{0.35ex}

Let $s=2.$ Using the functional system (2.31), we get 
\\[1.5ex]
\mbox{}\hfill
$
\nu^1x=\Psi_1^{\zeta}(x)\;\!\nu^{0}x$ for all $x\in {\mathscr X}.
\hfill
$
\\[1.5ex]
\indent
Then, taking into account Lemmas 2.3 and 2.5, we obtain 
\\[1.5ex]
\mbox{}\hfill
$
{\frak a}_j^{}\;\!\nu^1x=\mu_{0}^{j\zeta}\nu^1x+\mu_{1}^{j\zeta}\nu^0x
$
\ for all 
$
x\in {\mathscr X},
\quad 
j=1,\ldots, m.
\hfill
$
\\[1.5ex]
\indent
Therefore the system of identities (2.61)  for $s=2$ is true.
\vspace{0.35ex}

Suppose that the assertion of Lemma 2.7 is valid for $s=\varepsilon.$
\vspace{0.25ex}
Then, from the functional system (2.31) with $s=\varepsilon+1$ and $\theta=\varepsilon,$
we have
\\[1.5ex]
\mbox{}\hfill                               
$
{\frak a}_j^{}\;\!\nu^{\varepsilon}x\!=\!
{\displaystyle \sum\limits_{\rho=1}^{\varepsilon} }
\binom{\varepsilon-1}{\rho-1} \mu_{\rho}^{j\zeta}\;\! \nu^{\,\varepsilon-\rho}x  +
{\displaystyle \sum\limits_{\rho=1}^{\varepsilon} }
\binom{\varepsilon-1}{\rho-1} \Psi_{\rho}^{\zeta}(x)\;\! 
{\displaystyle \sum\limits_{\beta=0}^{\varepsilon-\rho} }
\binom{\varepsilon-\rho}{\beta} \mu_{\beta}^{j\zeta}\;\! \nu^{\,\varepsilon-\rho-\beta}x
$  
for all 
$\!
x\!\in\! {\mathscr X},
\, j\!=\!1,\ldots, m.
\hfill
$
\\[1.5ex]
\indent
By the induction hypothesis, so that
\\[1.5ex]
\mbox{}\hfill                               
$
{\frak a}_j^{}\;\!\nu^{\varepsilon}x =
{\displaystyle \sum\limits_{\rho=1}^{\varepsilon} }
\binom{\varepsilon-1}{\rho-1}\,\mu_{\rho}^{j\zeta}\, \nu^{\,\varepsilon-\rho}x +
{\displaystyle \sum\limits_{\delta=0}^{\varepsilon-1} }
\binom{\varepsilon-1}{\delta}\,\mu_{\delta }^{j\zeta}\
{\displaystyle \sum\limits_{\eta=1}^{\varepsilon-\delta} }
\binom{\varepsilon-\delta-1}{\eta-1}\,\Psi_{\eta}^{\zeta}(x)\, 
\nu^{(\varepsilon-\delta)-\eta}x=
\hfill
$
\\[2ex]
\mbox{}\hfill                               
$
={\displaystyle \sum\limits_{\rho=1}^{\varepsilon} }
\binom{\varepsilon-1}{\rho-1}\,\mu_{\rho}^{j\zeta}\, \nu^{\,\varepsilon-\rho}x +
{\displaystyle \sum\limits_{\delta=0}^{\varepsilon-1} }
\binom{\varepsilon-1}{\delta}\,\mu_{\delta}^{j\zeta}\, \nu^{\varepsilon-\delta}x=
{\displaystyle \sum\limits_{\rho=0}^{\varepsilon} }
\binom{\varepsilon}{\rho}\,\mu_{\rho}^{j\zeta}\, \nu^{\varepsilon-\rho}x
$ 
\ for all 
$
x\in {\mathscr X},
\ \ j=1,\ldots, m.
\hfill
$
\\[1.5ex]
\indent
Consequently the statement (2.61) for $s=\varepsilon+1$ is true. 
Thus by the principle of induction, the system of identities (2.61) is true 
for every natural number $s\geq 2.$ \k
\vspace{0.75ex}

{\bf Theorem 2.14.}
{\it 
Suppose the system {\rm (2.55)} satisfies the conditions of Lemma {\rm 2.5.}
\vspace{0.25ex}
Then the Lappo-Da\-ni\-lev\-s\-kii system {\rm (2.55)} 
on the domain $J\times {\mathscr X}$ has the first integrals
\\[1.5ex]
\mbox{}\hfill                                  % (2.62)                                
$
F_{\theta}^{}\colon (t,x)\to\ 
\nu^{\theta}x\;\! \varphi(t) \, -\, 
{\displaystyle \sum\limits_{\rho=1}^{\theta} }\,
K_{\rho-1}^{\theta}(t)\, F_{\rho-1}^{}(t,x) \, -\, C_{\theta}^{}(t),
\quad \
\theta = 0,\ldots, s-1,
\hfill 
$
\mbox{}\hfill {\rm (2.62)}
\\[1.75ex]
where ${\mathscr X}$ is a domain from the set $\{x\colon \nu^{0}x\ne 0\}\subset {\mathbb R}^n,$ 
\vspace{0.5ex}
the exponential function $\varphi\colon J\to {\mathbb R}$ is given by 
the formula {\rm (2.57)}, the scalar functions
\\[1.5ex]
\mbox{}\hfill
$
K_{\rho-1}^{\theta}\colon t\to \
{\displaystyle \int\limits_{t_0^{}}^{t}} 
\biggl(\!\binom{\theta}{\rho-1}\;\!\mu_{\theta-\rho+1}(\tau)+
{\displaystyle \sum\limits_{\eta=1}^{\theta-\rho}}\binom{\theta}{\eta}\,
\mu_{\eta}^{\zeta}(\tau)\;\! K_{\rho-1}^{\theta-\eta}(\tau)\biggr)d\tau,
\ \ 
\rho=1,\ldots, \theta,
\ \theta=1,\ldots, s-1,
\hfill
$
\\[1.75ex]
\mbox{}\hfill
$
C_{\theta}^{}\colon t\to\, 
{\displaystyle \int\limits_{t_0^{}}^{t}} 
\biggl(\nu^{\theta}f(\tau)\;\! \varphi(\tau)+
{\displaystyle \sum\limits_{\rho=1}^{\theta}}
\binom{\theta}{\rho}\;\!\mu_{\rho}^{\zeta}(\tau)\;\! C_{\theta-\rho}^{}(\tau)\biggr)\;\!d\tau,
\quad 
\theta=0,\ldots, s-1,
$ 
\ for all 
$
t\in J,
\hfill
$
\\[1.75ex]
\mbox{}\hfill
$
\displaystyle 
\mu_{\theta}^{\zeta}\colon t\to 
\sum\limits_{j=1}^{m}\mu_{\theta}^{j\zeta}\;\!\alpha_j^{}(t)
$  
for all 
$
t\in J
\
\bigl(\mu_{\theta}^{j\zeta}={\frak a}_j^{}\Psi_{\theta}^{\zeta}(x), \,  j=1,\ldots, m\bigr),
\ \theta=1,\ldots, s-1,
\ t_0^{}\in J.
\hfill
$
\\[1.5ex]
}
\indent
{\sl Proof}. The proof is by induction on $s.$
The case $s=1$ was considered in Theorem 2.12.
\vspace{0.25ex}

Suppose $s=2.$ Using Theorem 2.12 and Lemma 2.7, we obtain  
\\[2ex]
\mbox{}\hfill                                                 
$
{\frak B}\;\! F_1^{}(t,x)= 
\mu_1^{\zeta}(t)\bigl(\nu^0x\;\!\varphi(t)-C_{0}^{}(t)-F_{0}^{}(t,x)\bigr)=0
$
\ for all 
$
(t,x)\in J\times {\mathscr X}.
\hfill 
$
\\[2ex]
\indent
Therefore the scalar function
$F_1^{}\colon J\times {\mathscr X}\to {\mathbb R}$ 
is a first integral of system (2.55).
\vspace{0.5ex}

Suppose that the assertion of Theorem 2.14 is valid for 
\vspace{0.75ex}
$s=\varepsilon,$ i.e., the scalar functions
$F_{\theta}\colon \widetilde{\mathcal T}\times {\mathscr X}\to {\mathbb R},
\ \theta=1,\ldots,\varepsilon-1,$
\vspace{0.35ex}
are first integrals of the  Lappo-Da\-ni\-lev\-s\-kii system (2.55). 
Then, from Lemma 2.7, we get on the domain $J\times {\mathscr X}$
\\[1.5ex]
\mbox{}\hfill                                 
$
{\frak B} F_{\varepsilon}^{}(t,x)=
{\displaystyle \sum\limits_{\rho=1}^{\varepsilon} }
\binom{\varepsilon}{\rho}\,\mu_{\rho}^{\zeta}(t)
\biggl(
\Bigl( \nu^{\,\varepsilon-\rho}x\, \varphi(t) -
{\displaystyle \sum\limits_{\eta=1}^{\varepsilon-\rho} }
K_{\eta-1}^{\varepsilon-\rho}(t)\, F_{\eta-1}^{}(t,x) -
C_{\varepsilon-\rho}^{}(t)\Bigr)
-F_{\varepsilon-\rho}^{}(t,x)\biggr)=0.
\hfill                                 
$
\\[1.5ex]
\indent
Consequently if $s=\varepsilon+1,$ then the scalar function
\vspace{0.35ex}
$F_{\varepsilon}^{}\colon J\times {\mathscr X}\to {\mathbb R}$ 
is a first integral on the domain $J\times {\mathscr X}$
of the Lappo-Da\-ni\-lev\-s\-kii differential system (2.55).
\vspace{0.25ex}

Now from the method of building scalar functions (2.62) it follows that 
the Lappo-Da\-ni\-lev\-s\-kii differential system (2.55) 
has the functionally independent first integrals (2.62). \k
\vspace{0.75ex}

{\bf Example 2.23.}
The third-order Lappo-Da\-ni\-lev\-s\-kii differential system
\\[1.5ex]
\mbox{}\hfill                               % (2.63)
$
\dfrac{dx_1^{}}{dt}=t\;\! x_1^{}-(2t+1)\;\! x_2^{}+(t+1)\;\!x_3^{}+\dfrac{1}{2}\,,
\quad
\dfrac{dx_2^{}}{dt}=(2t+1)\;\! x_1^{}- x_2^{}+(2t+1)\;\!x_3^{}-2t^2,
\hfill
$
\\[1ex]
\mbox{}\hfill (2.63)
\\[0.25ex]
\mbox{}\hfill
$
\dfrac{dx_3^{}}{dt}={}-(t+1)\;\! x_1^{}+(2t+1)\;\! x_2^{}-(t+2)\;\!x_3^{}+t
\hfill
$
\\[2.5ex]
such that the coefficient matrix
$A(t)=A_1^{}+t\;\!A_2^{}$ for all $t\in {\mathbb R},$ 
where the constant matrix
\\[1.5ex]
\mbox{}\hfill
$
A_1^{}=
\left\|\!\!
\begin{array}{rrr}
0 & {}-1 & 1
\\
1 & {}-1 & 1
\\
{}-1 & 1 & {}-2
\end{array}\!\!
\right\|
$ 
\ \, and \ \,
$A_2^{}=
\left\|\!\!
\begin{array}{rrr}
1 & {}-2 & 1
\\
2 & 0 & 2
\\
{}-1 & 2 & {}-1
\end{array}\!\!
\right\|.
\hfill
$ 
\\[1.5ex]
\indent
The matrix $B_1^{}=A_1^{T},$ where $T$ denotes the matrix transpose,
\vspace{0.35ex}
has triple eigenvalue $\lambda_1^{1}={}-1.$ 
The rank of the matrix $B_1^{}-\lambda_1^{1}\;\!E$ is equal $2.$ 
\vspace{0.5ex}
Therefore the eigenvalue $\lambda_1^{1}={}-1$ has $\varkappa_1^{}=3-2=1$ 
elementary divisor $(\lambda^1+1)^3.$
\vspace{0.5ex}

The matrices $B_1^{}$ and $B_2^{}=A_2^{T}$ have
\vspace{0.5ex}
the common real eigenvector $\nu^0=(1,0,1)$ 
cor\-re\-s\-pon\-ding to the eigenvalues $\lambda_1^{1}={}-1$ and $\lambda_1^2=0.$ 
\vspace{0.5ex}
Also, the matrix $B_1^{}$ has the  1-st order generalized eigenvector
\vspace{0.35ex}
$\nu^1=(0,1,0)$ and the 2-nd order generalized eigenvector $\nu^2=(0,2,2)$ 
corresponding to the eigenvalue $\lambda_1^{1}={}-1.$
\vspace{0.35ex}

Now using the common eigenvector $\nu^0$ and the 
\vspace{0.35ex}
generalized eigenvectors $\nu^1,\ \nu^2,$ we can build the functions
$
\nu^0x\colon x\to\, x_1^{}+x_3^{},
\ \,
\nu^1x\colon x\to\, x_2^{},
\ \,
\nu^2x\colon x\to\, 2(x_2^{}+x_3^{})
$
for all 
$
x\in {\mathbb R}^3,
$
\\[2.75ex]
\mbox{}\hfill
$
\Psi_1^{1}\colon x\to\ \dfrac{x_2^{}}{x_1^{}+x_3^{}}\,,
\ \
\Psi_2^{1}\colon x\to\, \dfrac{2(x_1^{}+x_3^{})(x_2^{}+x_3^{})-x_2^{2}}{(x_1^{}+x_3^{})^2} 
$
\ for all 
$
x\in {\mathscr X}\subset \{x\colon x_1^{}+x_3^{}\ne 0\}.
\hfill
$
\\[2.25ex]
\indent
The real numbers 
\vspace{0.75ex}
$
\mu_1^1={\frak a}_1^{}\;\!\Psi_1^1(x)=1,
\ \, 
\mu_2^1={\frak a}_1^{}\;\!\Psi_2^1(x)=0,
$
and
$ 
\mu_1^2={\frak a}_2^{}\;\!\Psi_1^1(x)=2,
$
$ 
\mu_2^2={\frak a}_2^{}\;\!\Psi_2^1(x)=2, 
$
where the linear differential operators of first order
\\[1.5ex]
\mbox{}\hfill
$
{\frak a}_1^{}(x)=({}-x_2^{}+x_3^{})\;\!\partial_{x_1^{}}^{}+
(x_1^{}-x_2^{}+x_3^{})\;\!\partial_{x_2^{}}^{}-
(x_1^{}-x_2^{}+2x_3^{})\;\!\partial_{x_3^{}}^{}
$ 
\ for all 
$
x\in {\mathbb R}^3,
\hfill
$
\\[2ex]
\mbox{}\hfill
$
{\frak a}_2^{}(x)=(x_1^{}-2x_2^{}+x_3^{})\;\!\partial_{x_1^{}}^{}+
2(x_1^{}+x_3^{})\;\!\partial_{x_2^{}}^{}-
(x_1^{}-2x_2^{}+x_3^{})\;\!\partial_{x_3^{}}^{}
$ 
\ for all 
$
x\in {\mathbb R}^3.
\hfill
$
\\[1.75ex]
%are induced by the Lappo-Da\-ni\-lev\-s\-kii differential system (2.63).
%\vspace{0.35ex}
\indent
Using the scalar functions 
\\[1.5ex]
\mbox{}\hfill
$
\alpha_1^{}\colon t\to 1,
\ \ 
\alpha_2^{}\colon t\to t,
\ \
f_1^{}\colon t\to \dfrac{1}{2}\,,
\quad 
f_2^{}\colon t\to {}-2\;\!t^2,
\ \ 
f_3^{}\colon t\to\, t,
\ \ 
\varphi\colon t\to e^t 
$
\ for all 
$
t\in {\mathbb R}, 
\hfill
$
\\[2.5ex]
\mbox{}\hfill
$
C_0^{}\colon t\to \Bigl(t-\dfrac{1}{2}\Bigr)e^t,
\ \
\mu_1^{1}\colon t\to 2t+1, 
\ \
K_0^1\colon t\to t(t+1), 
\ \
C_1^{}\colon t\to {}-\dfrac{1}{2}\,e^t, 
\ \
\mu_2^{1}\colon t\to 2t, 
\hfill
$
\\[2.5ex]
\mbox{}\hfill
$
K_0^2\colon t\to t^2(t^2+2t+2),
\ \ 
K_1^2\colon t\to 2t(t+1), 
\ \
C_2^{}\colon t\to {}-(2t^2-3t+4)\;\!e^t
$
\ for all 
$
t\in {\mathbb R}, 
\hfill
$
\\[2ex]
we can build (by Theorem 2.14) the basis of first integrals on the space ${\mathbb R}^4$
for system (2.63)
\\[2ex]
\mbox{}\hfill
$
F_0^{}\colon (t,x)\to \
\Bigl(x_1^{}+x_3^{}-t+\dfrac{1}{2}\Bigr)e^{\;\!t},
\quad \
F_1^{}\colon (t,x)\to \
\Bigl(x_2^{}+\dfrac{1}{2}\Bigr)e^{\;\!t}-t(t+1)F_0^{}(t,x),
\hfill
$
\\[3ex]
\mbox{}\hfill
$
F_2^{}\colon (t,x)\to \
(2x_2^{}+2x_3^{}+2t^2-3t+4)e^{t} -t^2(t^2+2t+2)F_0^{}(t,x)-
2t(t+1)F_1^{}(t,x).
\hfill
$
\\[2ex]
\indent
The proof of Theorem 2.14 is true both for 
real eigenvectors (common and generalized) of  the matrix $B_{\zeta}^{},$ 
and for 
\vspace{0.25ex}
complex eigenvectors (common and generalized) of the matrix $B_{\zeta}^{}.$ 

{\it In the complex case}, from the complex-valued first integrals (2.62) of 
the Lappo-Da\-ni\-lev\-s\-kii differential system (2.55), 
we obtain the real first integrals 
\\[1.5ex]
\mbox{}\hfill
$
F_{\theta}^{\;\!1}\colon (t,x)\to\, {\rm Re}\,F_{\theta}^{}(t,x),
\ \
F_{\theta}^{\;\!2}\colon (t,x)\to\, {\rm Im}\,F_{\theta}^{}(t,x)
$
for all 
$
(t,x)\in J\times {\mathscr X},
\
\theta=0,\ldots, s-1.
\hfill
$
\\[2ex]
\indent
{\bf Example 2.24.}
The fourth-order Lappo-Da\-ni\-lev\-s\-kii differential system
\\[2ex]
\mbox{}\hfill                               % (2.64)
$
\dfrac{dx_1^{}}{dt}=\dfrac{1}{t}\;\! (x_1^{}- x_2^{})-x_4^{}+t\cos t,
\qquad
\dfrac{dx_2^{}}{dt}={}- x_1^{}+\dfrac{1}{t}\, x_2^{}-x_3^{}-t\sin t,
\hfill
$
\\[1.25ex]
\mbox{}\hfill (2.64)
\\[0.5ex]
\mbox{}\hfill
$
\dfrac{dx_3^{}}{dt}=\Bigl(1+\dfrac{1}{t}\Bigr)\;\!x_2^{}+\dfrac{1}{t}\;\! x_3^{}+x_4^{}+t^2,
\quad
\dfrac{dx_4^{}}{dt}=\Bigl(3-\dfrac{1}{t}\Bigr)\;\!x_1^{}+
\Bigl(2-\dfrac{1}{t}\Bigr)\;\! x_3^{}+\dfrac{1}{4}\, x_4^{}+2t
\hfill
$
\\[2.5ex]
has the coefficient matrix of the form
\vspace{1.25ex}
$A(t)=A_1^{}+\dfrac{1}{t}\,A_2^{}$ for all $t\in J\subset \{t\colon t\ne 0\},$ 
where the constant matrix 
$A_1^{}=
\left\|\!\!
\begin{array}{rrrr}
0 & 0 & 0 & {}-1
\\
{}-1 & 0 & {}-1 & 0
\\
0 & 1 & 0 & 1
\\
3 & 0 & 2 & 0
\end{array}\!\!
\right\|
$ 
and
\vspace{1.25ex}
$A_2^{}=
\left\|\!\!
\begin{array}{rrrr}
1 & {}-1 & 0 & 0
\\
0 & 1 & 0 & 0
\\
0 & 1 & 1 & 0
\\
{}-1 & 0 & {}-1 & 1
\end{array}\!\!
\right\|.
$ 

The matrix $B_1^{}=A_1^{T},$ 
\vspace{0.5ex}
where $T$ denotes the matrix transpose, has two double complex eigenvalues
\vspace{0.5ex}
$\lambda_1^{1}=i$ and  $\lambda_2^{1}={}-i.$
The rank of the matrix $B_1^{}-\lambda_1^{1}\;\!E$ is equal $3.$ 
Therefore the eigenvalue $\lambda_1^{1}=i$ has $\varkappa_1^{}=4-3=1$ 
elementary divisor $(\lambda^1-i)^2.$ 
\vspace{0.5ex}

Thus the matrix $B_1^{}$ has two elementary divisors $(\lambda^1-i)^2$ 
and $(\lambda^1+i)^2.$ 
\vspace{0.5ex}

The matrices $B_1^{}$ and $B_2^{}=A_2^{T}$ 
\vspace{0.5ex}
have the common complex eigenvector $\nu^0=(1,{}-i,1,0)$ 
corresponding to the eigenvalues 
\vspace{0.5ex}
$\lambda_1^{1}=i$ and $\lambda_1^2=1.$ 
Also, the matrix $B_1^{}$ has the generalized eigenvector of the 1-st order $\nu^1=({}-i, 1,0, 1).$ 
\vspace{0.5ex}

From the eigenvector $\nu^0$ and  the generalized eigenvector $\nu^1,$ we obtain the function
\\[2ex]
\mbox{}\hfill
$
\Psi_1^{1}\colon x\to\ \dfrac{{}-ix_1^{}+x_2^{}+x_4^{}}{x_1^{}-ix_2^{}+x_3^{}}
$
\ for all 
$
x\in {\mathscr X}\subset \{x\colon (x_1^{}+x_3^{})^2+x_2^2\ne 0\}.
\hfill
$
\\[2.25ex]
\indent
The numbers 
$
\mu_0^1=i,\
\mu_0^2=1,\
\mu_1^1={\frak a}_1^{}\;\!\Psi_1^1(x)=1,\ 
\mu_1^2={\frak a}_2^{}\;\!\Psi_1^1(x)={}-1,
$
where 
\\[1.75ex]
\mbox{}\hfill
$
{\frak a}_1^{}(x)={}-x_4\;\!\partial_{x_1^{}}^{}-
(x_1^{}+x_3^{})\;\!\partial_{x_2^{}}^{}+
(x_2^{}+x_4^{})\;\!\partial_{x_3^{}}^{}+
(3x_1^{}+2x_3^{})\;\!\partial_{x_4^{}}^{}
$ 
\ for all 
$
x\in {\mathbb R}^4,
\hfill
$
\\[2ex]
\mbox{}\hfill
$
{\frak a}_2^{}(x)=(x_1^{}-x_2^{})\;\!\partial_{x_1^{}}^{}+
x_2^{}\;\!\partial_{x_2^{}}^{}+
(x_2^{}+x_3^{})\;\!\partial_{x_3^{}}^{}-
(x_1^{}+x_3^{}-x_4^{})\;\!\partial_{x_4^{}}^{}
$ 
\ for all 
$
x\in {\mathbb R}^4.
\hfill
$
\\[1.75ex]
\indent
Using the scalar functions 
\\[1.5ex]
\mbox{}\hfill
$
\alpha_1^{}\colon t\to 1,
\ \
\alpha_2^{}\colon t\to\, \dfrac{1}{t}\;\!, 
\ \
f_1^{}\colon t\to\, t\cos t,
\ \
f_2^{}\colon t\to {}-t\sin t,
\ \
f_3^{}\colon t\to\, t^2, 
\ \
f_4^{}\colon t\to\, 2t,
\hfill
$
\\[1.5ex]
\mbox{}\hfill
$
\varphi\colon t\to \dfrac{1}{t}\,(\cos t-i\sin t),
\quad 
C_0^{}\colon t\to t+\cos t+t\sin t-i(\sin t-t\cos t),
\quad 
\mu_1^{1}\colon t\to 1-\dfrac{1}{t}\,, \
\hfill
$
\\[2ex]
\mbox{}\hfill
$
\displaystyle
K_0^1\colon t\to t-\ln |t|, 
\quad
C_1^{}\colon t\to {}-t+\dfrac{t^2}{2}+\cos t+4\sin t+\dfrac{1}{2}\,\cos 2t -t\cos t-
\int\dfrac{\cos t}{t}\ dt \ +
\hfill
$
\\[2.5ex]
\mbox{}\hfill
$
\displaystyle
+\ i\;\!\Bigl(
4\cos t-\sin t-\dfrac{1}{2}\sin 2t+t\sin t +\int\dfrac{\sin t}{t}\ dt\Bigr)
$
\ for all 
$
t\in J, 
\hfill
$
\\[2ex]
we can build the integral basis of the Lappo-Da\-ni\-lev\-s\-kii differential system (2.64)
\\[2ex]
\mbox{}\hfill
$
F_0^{\;\!1}\colon (t,x)\to\  
\dfrac{1}{t}\,\bigl(\cos t\;\!(x_1^{}+x_3^{})-\sin t\;\!x_2^{}\bigr)-t-\cos t-t\sin t
$ 
\ for all 
$
(t,x)\in J\times {\mathbb R}^4,
\hfill
$
\\[2.5ex]
\mbox{}\hfill
$
F_0^{\;\!2}\colon (t,x)\to \ 
{}-\dfrac{1}{t}\,\bigl(\cos t\;\!x_2^{}+\sin t\;\!(x_1^{}+x_3^{})\bigr)+\sin t-t\cos t
$ 
\ for all 
$
(t,x)\in J\times {\mathbb R}^4,
\hfill
$
\\[2.75ex]
\mbox{}\hfill
$
F_1^{\;\!1}\colon (t,x)\to\  
\dfrac{1}{t}\,\bigl(\cos t\;\!(x_2^{}+x_4^{})-\sin t\;\!x_1^{}\bigr)-(t-\ln |t|)\;\!F_0^{\;\!1}(t,x)\ +
\hfill
$
\\[2.25ex]
\mbox{}\hfill
$
\displaystyle
+\ t-\dfrac{t^2}{2}-\cos t-4\sin t-\dfrac{1}{2}\,\cos 2t +t\cos t+
\int\dfrac{\cos t}{t}\ dt 
$ 
\ for all 
$
(t,x)\in J\times {\mathbb R}^4,
\hfill
$
\\[2.75ex]
\mbox{}\hfill
$
F_1^{\;\!2}\colon (t,x)\to\  
\dfrac{1}{t}\,\bigl(\cos t\;\!x_1+\sin t\;\!(x_2^{}+x_4^{})\bigr)+(t-\ln |t|)\;\!F_0^{\;\!2}(t,x)\ +
\hfill
$
\\[2.25ex]
\mbox{}\hfill
$
\displaystyle
+\ 4\cos t-\sin t-\dfrac{1}{2}\,\sin 2t +t\sin t+
\int\dfrac{\sin t}{t}\ dt 
$
\ for all 
$
(t,x)\in J\times {\mathbb R}^4.
\hfill
$
\\[-3ex]

\newpage

\mbox{}
\\[-1.75ex]
\centerline{                               %{Параграф 3}
\large\bf
3. Integrals of reducible ordinary differential systems
}
\\[2ex]
\centerline{                               %{П. 1, Параграф 1}
\bf 3.1.
Linear homogeneous differential system
}
\\[1.5ex]
\indent
Consider a system of $n$ ordinary real linear differential equations
\\[1.75ex]
\mbox{}\hfill                                          % (3.1)
$
\dfrac{dx}{dt}=A(t)\,x,
$
\hfill (3.1)
\\[1.75ex]
where 
\vspace{0.25ex}
$x={\rm colon}(x_1,\ldots,x_n)$ from the arithmetical phase space ${\mathbb R}^n,$  
the real $n\times n$ coefficient matrix 
$A\colon t\to A(t)$ for all $t\in J$ is continuous on an interval $J\subset {\mathbb R}.$ 
\vspace{0.25ex}

One of the most efficient methods for investigation of linear 
nonautonomous differential systems is the method of reducibility. 
In this method linear nonautonomous systems are re\-du\-ced to linear 
systems with constant coefficients by various transformation groups [63; 8; 10].
The idea of this method is due to  A.M. Lyapunov [57].
He studied linear periodic systems and showed that there 
exists a transformation which does not change the character of the growth of solutions 
and reduces a system with periodic coefficients to a system with constant coefficients. 
Lyapunov called systems having this condition {\it reducible systems}.
The general the\-o\-ry of reducible systems was developed by N.P.Erugin in his article [58]. 

Later on, the notion of reducible system was given for 
systems of difference equations [64] and for
multidimensional differential systems 
(see [65 -- 69; 9, pp. 72 -- 82, 242 -- 246]).

Let $G$ be a multiplicative group of real continuously differentiable on the interval $J$ 
nonsingular matrices of order $n.$
The differential system (3.1) is called 
{\it reducible with respect to the nonautonomous transformation group} $G$ if 
there exist a constant matrix $B$ and a matrix ${\rm g}\in G$ such that 
the linear transformation $y={\rm g}(t)\;\!x$ 
reduces the nonautonomous differential system (3.1) to the system with constant coefficients
\\[1.75ex]
\mbox{}\hfill                                          % (3.2)
$
\dfrac{dy}{dt}=B\,y,
\quad
y\in {\mathbb R}^n.
$
\hfill (3.2)
\\[1.75ex]
\indent
In addition, the transformation matrix ${\rm g}$ such that
\\[1.75ex]
\mbox{}\hfill                                          % (3.3)
$
\dfrac{d\,{\rm g}(t)}{dt}=B\,{\rm g}(t)-\;\!{\rm g}(t)\;\!A(t)
$
\ for all 
$
t\in J.
$
\hfill (3.3)
\\[1.75ex]
\indent
By Theorems 3.1 -- 3.4, using eigenvectors and corresponding eigenvalues 
\vspace{0.25ex}
of the matrix $C=B^{T},$ where $T$ denotes the matrix transpose, and 
a transformation matrix ${\rm g}\in G,$ we can construct first integrals 
of the reducible system (3.1).\!
The following basic statements (Lem\-mas 3.1 and 3.2)
\vspace{0.5ex}
are base for the method of building integral basis of system (3.1).

{\bf Lemma 3.1.}
\vspace{0.25ex}
{\it
Suppose the system {\rm (3.1)} 
is reducible to the sys\-tem 
%with con\-s\-tant coefficients 
{\rm (3.2)} 
by a transformation matrix ${\rm g}\!\in\! G,\!$ 
\vspace{0.35ex}
and $\nu\!$ is a real eigenvector of the matrix $C\!$ corresponding to the 
eigen\-va\-lue $\!\lambda.$
Then the linear function
$
p\colon (t,x)\to \nu\, {\rm g}(t)\;\!x
$
for all 
$
(t,x)\in J\times {\mathbb R}^{n}
$
\vspace{0.5ex}
is a partial integral of the reducible differential system {\rm (3.1)} such that}
\\[1.25ex]
\mbox{}\hfill
$
{\frak A}\;\!p(t,x)=\lambda\;\!p(t,x)
$
\ for all 
$
(t,x)\in J\times {\mathbb R}^{n},
\hfill
$
\\[1.5ex]
{\it where the linear differential operator 
$
{\frak A}(t,x)=\partial_{t}^{}+A(t)\;\!x\,\partial_{x}^{}
$
for all} 
$(t,x)\in J\times {\mathbb R}^{n}.$
\vspace{0.75ex}

{\sl Indeed}, using the matrix identity (3.3), we get
\\[1.5ex]
\mbox{}\hfill
$
{\frak A}\;\!p(t,x)=\partial_{t}^{}\;\!p(t,x)+A(t)\;\!x\,\partial_{x}^{}\;\!p(t,x)=
\nu\, {\rm g}^{\prime}(t)x \, +\, A(t)\;\!x\,\nu\, {\rm g}(t)\,=
\hfill
$
\\[2ex]
\mbox{}\hfill
$
=\nu\bigl(B\,{\rm g}(t)-\;\!{\rm g}(t)\;\!A(t)\bigr)\;\!x \,+\,
A(t)\;\!x\,\nu\, {\rm g}(t) =
\nu B\,{\rm g}(t)\;\!x +
\bigl(A(t)x\,\nu{\rm g}(t)- \nu{\rm g}(t)\;\!A(t)x\bigr)=
\hfill
$
\\[2ex]
\mbox{}\hfill
$
=C\nu \,{\rm g}(t)\;\!x =
\lambda\, \nu\, {\rm g}(t)x=
\lambda\;\!p(t,x)
$
\ for all 
$
(t,x)\in J\times {\mathbb R}^{n}.
\ \k
\hfill
$
\\[2ex]
\indent
{\bf Lemma 3.2.}
\vspace{0.25ex}
{\it
Suppose the system {\rm (3.1)} is reducible to the sys\-tem 
with con\-s\-tant coefficients {\rm (3.2)}
by a transformation matrix ${\rm g}\in G,$ and
$\nu={\stackrel{*}{\nu}}+\widetilde{\nu}\,i\
({\rm Re}\,\nu={\stackrel{*}{\nu}},\ {\rm Im}\,\nu=\widetilde{\nu}\,)$
is a complex eigenvector of the matrix $C$ 
corresponding to the eigenvalue
\vspace{0.35ex}
$\lambda={\stackrel{*}{\lambda}}+\widetilde{\lambda}\,i\ 
({\rm Re}\,\lambda={\stackrel{*}{\lambda}},\ {\rm Im}\, \lambda=\widetilde{\lambda}\,).$
Then the Lie derivatives of the scalar functions 
\\[1.5ex]
\mbox{}\hfill                                            
$
P\colon (t,x)\to \
({\stackrel{*}{\nu}}\,{\rm g}(t)\;\!x)^2+(\widetilde{\nu}\,{\rm g}(t)\;\!x)^2
$
\ for all 
$
(t,x)\in J\times {\mathbb R}^n
\hfill
$
\\[1ex]
and
\\[1.5ex]
\mbox{}\hfill                                            
$
\displaystyle
\varphi\colon (t,x)\to \
\arctan\dfrac{\widetilde{\nu}\,{\rm g}(t)\;\!x}{{\stackrel{*}{\nu}}\,{\rm g}(t)\;\!x}
$ 
\ for all 
$
(t,x)\in \Lambda,
\ \
\Lambda\subset
\bigl\{(t,x)\colon t\in J,\ {\stackrel{*}{\nu}}\,{\rm g}(t)\;\!x\ne 0\bigr\}\subset J\times {\mathbb R}^n,
\hfill
$
\\[1.5ex]
by virtue of the reducible system {\rm (3.1)} are equal to
\\[1.5ex]
\mbox{}\hfill                                            
$
{\frak A}\;\! P(t,x)=2\;\!{\stackrel{*}{\lambda}}\;\!P(t,x)
$  
\ for all 
$
(t,x)\in J\times {\mathbb R}^n
$
\ \
and
\ \ 
$
\displaystyle
{\frak A}\;\!\varphi(t,x)=
\widetilde{\lambda}
$ 
\ for all 
$
(t,x)\in \Lambda.
\hfill 
$
}
\\[2ex]
\indent
{\sl Proof}. 
\vspace{0.35ex}
Formally using Lemma 3.1, we obtain 
the linear function 
$
p\colon (t,x)\to\, \nu\, {\rm g}(t)\;\!x
$
for all 
$
(t,x)\in J\times {\mathbb R}^{n}
$
\vspace{0.25ex}
is a complex-valued par\-ti\-al integral 
of the reducible linear differential system (3.1) and the following identity holds
\\[1.25ex]
\mbox{}\hfill
$
{\frak A}\;\!
\bigl({\stackrel{*}{\nu}}\, {\rm g}(t)x+i\,\widetilde{\nu}\, {\rm g}(t)x\bigr)=
\bigl({\stackrel{*}{\lambda}}+\widetilde{\lambda}\,i\bigr)
\bigl({\stackrel{*}{\nu}}\, {\rm g}(t)x+i\,\widetilde{\nu}\, {\rm g}(t)x\bigr)
$
\ for all 
$
(t,x)\in J\times {\mathbb R}^{n}.
\hfill
$
\\[1.75ex]
\indent
This complex identity is equivalent to the real system of identities
\\[1.5ex]
\mbox{}\hfill
$
{\frak A}\;\!{\stackrel{*}{\nu}}\, {\rm g}(t)\;\!x =\,
{\stackrel{*}{\lambda}}\,{\stackrel{*}{\nu}}\, {\rm g}(t)\;\!x -
\widetilde{\lambda}\,\widetilde{\nu}\, {\rm g}(t)\;\!x
$
\ for all 
$
(t,x)\in J\times {\mathbb R}^{n},
\hfill
$
\\[-0.15ex]
\mbox{}\hfill 
\\[-0.35ex]
\mbox{}\hfill
$
{\frak A}\;\!\widetilde{\nu}\, {\rm g}(t)\;\!x=\,
{\stackrel{*}{\lambda}}\,\widetilde{\nu}\, {\rm g}(t)\;\!x+
\widetilde{\lambda}\,{\stackrel{*}{\nu}}\, {\rm g}(t)\;\!x
$
\ for all 
$
(t,x)\in J\times {\mathbb R}^{n}.
\hfill 
$
\\[1.75ex]
\indent
Using this system of identities, we have
\\[1.75ex]
\mbox{}\hfill                        
$
{\frak A}\;\!P(t,x)=
{\frak A}\bigl(
({\stackrel{*}{\nu}}\,{\rm g}(t)\;\!x)^2+(\widetilde{\nu}\,{\rm g}(t)\;\!x)^2\bigr)=
2\;\!{\stackrel{*}{\nu}}\;\!{\rm g}(t)\;\!x\,{\frak A}\;\!{\stackrel{*}{\nu}}\;\!{\rm g}(t)\;\!x +
2\;\!\widetilde{\nu}\;\!{\rm g}(t)\;\!x\,{\frak A}\;\!\widetilde{\nu}\;\!{\rm g}(t)\;\!x=
\hfill                        
$
\\[1.75ex]
\mbox{}\hfill                        
$
=
2\;\!{\stackrel{*}{\nu}}\;\!{\rm g}(t)\;\!x\,
\bigl(\;\!{\stackrel{*}{\lambda}}\,{\stackrel{*}{\nu}}\, {\rm g}(t)\;\!x -
\widetilde{\lambda}\,\widetilde{\nu}\, {\rm g}(t)\;\!x\bigr)  +\, 
2\;\!\widetilde{\nu}\;\!{\rm g}(t)\;\!x\,
\bigl(\;\!{\stackrel{*}{\lambda}}\,\widetilde{\nu}\, {\rm g}(t)\;\!x+
\widetilde{\lambda}\,{\stackrel{*}{\nu}}\, {\rm g}(t)\;\!x\bigr)=
\hfill                        
$
\\[1.75ex]
\mbox{}\hfill                        
$
= 2\;\!{\stackrel{*}{\lambda}}\;\!
\bigl(({\stackrel{*}{\nu}}\,{\rm g}(t)\;\!x)^2+(\widetilde{\nu}\,{\rm g}(t)\;\!x)^2\bigr)=
2\;\!{\stackrel{*}{\lambda}}\,P(t,x)
$
\ for all 
$
(t,x)\in J\times {\mathbb R}^{n};
\hfill                        
$
\\[2.5ex]
\mbox{}\hfill                        
$
{\frak A}\,\varphi(t,x)\, =\, 
{\frak A}\,\arctan\dfrac{\widetilde{\nu}\,{\rm g}(t)\;\!x}{{\stackrel{*}{\nu}}\,{\rm g}(t)\;\!x} \, =\,
\dfrac{1}{1+\dfrac{(\widetilde{\nu}\,{\rm g}(t)\;\!x)^2}{({\stackrel{*}{\nu}}\,{\rm g}(t)\;\!x)^2}}
\cdot
\dfrac{{\stackrel{*}{\nu}}\;\!{\rm g}(t)\;\!x\, {\frak A}\;\!\widetilde{\nu}\;\!{\rm g}(t)\;\!x - 
\widetilde{\nu}\;\!{\rm g}(t)\;\!x\, {\frak A}\;\!{\stackrel{*}{\nu}}\;\!{\rm g}(t)\;\!x}
{({\stackrel{*}{\nu}}\,{\rm g}(t)\;\!x)^2}\ =
\hfill                        
$
\\[2.25ex]
\mbox{}\hfill                        
$
=
\dfrac{
{\stackrel{*}{\nu}}\;\!{\rm g}(t)\;\!x\;\! 
\bigl(\;\!{\stackrel{*}{\lambda}}\,\widetilde{\nu}\, {\rm g}(t)\;\!x+
\widetilde{\lambda}\,{\stackrel{*}{\nu}}\, {\rm g}(t)\;\!x\bigr)  -\,
\widetilde{\nu}\;\!{\rm g}(t)\;\!x\;\! 
\bigl(\;\!{\stackrel{*}{\lambda}}\,{\stackrel{*}{\nu}}\, {\rm g}(t)\;\!x -
\widetilde{\lambda}\,\widetilde{\nu}\, {\rm g}(t)\;\!x\bigr)}
{({\stackrel{*}{\nu}}\,{\rm g}(t)\;\!x)^2+(\widetilde{\nu}\,{\rm g}(t)\;\!x)^2}=
\widetilde{\lambda}
$
\ for all 
$
(t,x)\in \Lambda.
\ \k
\hfill                        
$
\\[2.25ex]
\indent
Now we can state the Theorems 3.1 -- 3.4 for building first integrals of system (3.1).
\vspace{1.25ex}

{\bf Theorem 3.1.}
\vspace{0.25ex}
{\it
Let the conditions of Lemma {\rm 3.1} be satisfied.
Then the scalar function}
\\[1.5ex]
\mbox{}\hfill
$
F\colon (t,x)\to\ \nu\, {\rm g}(t)\;\!x\,\exp({}-\lambda\,t)
$ 
\ for all 
$
(t,x)\in J\times {\mathbb R}^{n}
\hfill
$
\\[1.5ex]
{\it
is a first integral on the domain $J\times {\mathbb R}^{n}$ 
of the reducible differential system {\rm (3.1)}.
}
\vspace{1ex}

{\sl Proof}. From Lemma 3.1, we get
\\[1.75ex]
\mbox{}\hfill
$
{\frak A}\;\!F(t,x)=
{\frak A}\,\bigl(\nu\;\!{\rm g}(t)\;\!x\;\!\exp({}-\lambda\,t)\bigr)=
\exp({}-\lambda\,t)\, {\frak A}\;\!\nu\;\!{\rm g}(t)\;\!x+ 
\nu\;\!{\rm g}(t)\;\!x\,{\frak A}\exp({}-\lambda\,t)=
\hfill
$
\\[2ex]
\mbox{}\hfill
$
=
\lambda\,\nu\;\!{\rm g}(t)\;\!x\;\!\exp({}-\lambda\,t) \, +\, 
\nu\;\!{\rm g}(t)\;\!x\,\partial_{t}^{}\exp({}-\lambda\,t)=0
$ 
\ for all 
$
(t,x)\in J\times {\mathbb R}^{n}.
\hfill
$
\\[1.75ex]
\indent
Therefore the function $F\colon J\times {\mathbb R}^{n}\to {\mathbb R}$ 
\vspace{1ex}
is a first integral of the reducible system (3.1). \k

{\bf Example 3.1.}
The second-order nonautonomous differential system  [8, p. 155]
\\[2ex]
\mbox{}\hfill
$
\dfrac{dx_1^{}}{dt}=(t^2+t+2)\;\!x_1^{}-(t^3+t^2+t-1)\;\!x_2^{},
\quad \ \,
\dfrac{dx_2^{}}{dt}=(t+1)\;\!x_1^{}-(t^2+t-1)\;\!x_2^{}
$
\hfill (3.4)
\\[2ex]
is reducible to the system with constant coefficients
\vspace{1.5ex}
$
\dfrac{dy_1^{}}{dt}=y_1^{},
\ \
\dfrac{dy_2^{}}{dt}=2x_2^{}
$
by the polynomial group\footnote[1]{ 
$P(n)$ is the multiplicative group of $n\times n$ polynomial matrices with
nonzero constant determinants. 
}
\vspace{1ex}
$\!P(2)$ with the transformation matrix 
$
{\rm g}(t)=
\left\|\!\!
\begin{array}{rc}
-\, t & 1+t^2
\\
1 & -\, t
\end{array}
\!\!\right\|
$
for all 
$
t\in{\mathbb R}.
$

Using the real eigenvectors $\nu^1=(1, 0),\ \nu^2=(0,1)$ of the matrix  
\vspace{0.35ex}
$
C ={\rm diag}(1, 2)
$
and the corresponding eigenvalues $\lambda_1^{}=1, \ \lambda_2^{}=2,$
\vspace{0.5ex}
we can build (by Theorem 3.1) the basis of first integrals on space ${\mathbb R}^3$ 
for the reducible differential system (3.4)
\\[2ex]
\mbox{}\hfill
$
F_1^{}\colon (t,x_1^{},x_2^{})\to \
({}-t\;\!x_1^{}+(1+t^2)\;\!x_2^{})\,e^{{}-t},
\quad \
F_2^{}\colon (t,x_1^{},x_2^{})\to \
(x_1^{}-t\;\!x_2^{})\,e^{{}-2t}.
\hfill
$
\\[2.5ex]
\indent
{\bf Theorem 3.2.}
\vspace{0.25ex}
{\it
Let the conditions of Lemma {\rm 3.2} be satisfied.
Then the reducible differential system {\rm (3.1)} has the first integrals}
\\[1.5ex]
\mbox{}\hfill                                            
$
F_1^{}\colon (t,x)\to 
\bigl( ({\stackrel{*}{\nu}}\,{\rm g}(t)\;\!x)^2+(\widetilde{\nu}\,{\rm g}(t)\;\!x)^2\bigr)
\exp({}-2\;\! {\stackrel{*}{\lambda}}\,t)
$ 
\ for all 
$
(t,x)\in J\times {\mathbb R}^n
\hfill
$
\\[0.75ex]
and 
\\[1.5ex]
\mbox{}\hfill                                            
$
\displaystyle
F_2^{}\colon (t,x)\to\ 
\arctan\dfrac{\widetilde{\nu}\,{\rm g}(t)\;\!x}{{\stackrel{*}{\nu}}\,{\rm g}(t)\;\!x}\,-\,
\widetilde{\lambda}\,t
$
\ for all 
$
(t,x)\in \Lambda,
\quad 
\Lambda\subset\bigl\{(t,x)\colon t\in J,\ {\stackrel{*}{\nu}}\,{\rm g}(t)\;\!x\ne 0\bigr\}.
\hfill
$
\\[2ex]
\indent
{\sl Proof}. It follows from Lemma 3.2 that
\\[1.75ex]
\mbox{}\hfill
$
{\frak A}\;\!F_1^{}(t,x)=
\exp({}-2\;\! {\stackrel{*}{\lambda}}\,t)\, 
{\frak A}\;\!\bigl( ({\stackrel{*}{\nu}}\,{\rm g}(t)\;\!x)^2+(\widetilde{\nu}\,{\rm g}(t)\;\!x)^2\bigr) +\, 
\bigl( ({\stackrel{*}{\nu}}\,{\rm g}(t)\;\!x)^2+(\widetilde{\nu}\,{\rm g}(t)\;\!x)^2\bigr)
\,{\frak A}\exp({}-2\;\! {\stackrel{*}{\lambda}}\,t)=
\hfill
$
\\[2ex]
\mbox{}\hfill
$
=2\;\! {\stackrel{*}{\lambda}}\;\!F_1^{}(t,x)\, +\, 
\bigl( ({\stackrel{*}{\nu}}\,{\rm g}(t)\;\!x)^2+(\widetilde{\nu}\,{\rm g}(t)\;\!x)^2\bigr)
\,\partial_{t}^{}\exp({}-2\;\! {\stackrel{*}{\lambda}}\,t)=0
$ 
\ for all 
$
(t,x)\in J\times {\mathbb R}^n,
\hfill
$
\\[2.25ex]
\mbox{}\hfill                                            
$
\displaystyle
{\frak A}\;\!F_2^{}(t,x)= 
{\frak A}\;\!\arctan\dfrac{\widetilde{\nu}\,{\rm g}(t)\;\!x}{{\stackrel{*}{\nu}}\,{\rm g}(t)\;\!x}\, - \,
{\frak A}\;\!\bigl(\;\!\widetilde{\lambda}\,t\bigr) \, =\, 
\widetilde{\lambda}\, -\,
\partial_{t}^{}\bigl(\;\!\widetilde{\lambda}\,t\bigr)=0
$ 
\ for all 
$
(t,x)\in \Lambda.
\hfill
$
\\[2ex]
\indent
Therefore the scalar functions $F_1^{}\colon J\times {\mathbb R}^n\to {\mathbb R}$ and 
\vspace{0.25ex}
$F_2^{}\colon \Lambda\to {\mathbb R}$ are first integrals 
of the reducible linear homogeneous differential system (3.1). \k
\vspace{1.5ex}

{\bf Example 3.2.}
The linear differential system [70, pp. 125 -- 126]
\\[2.5ex]
\mbox{}\hfill
$
\dfrac{dx_1^{}}{dt}={}-x_2^{}+\sqrt{2}\,\cos 2t\;\!x_3^{}+\sqrt{2}\,\sin 2t\;\!x_4^{},
\quad 
\dfrac{dx_2^{}}{dt}=x_1^{}+\sqrt{2}\,\sin 2t\;\!x_3^{}-\sqrt{2}\,\cos 2t\;\!x_4^{},
\hfill
$
\\[0.35ex]
\mbox{}\hfill (3.5)
\\[0.35ex]
\mbox{}\hfill
$
\dfrac{dx_3^{}}{dt}={}-\sqrt{2}\,\sin 2t\;\!x_1^{}+\sqrt{2}\,\cos 2t\;\!x_2^{}-x_4^{},
\quad 
\dfrac{dx_4^{}}{dt}=\sqrt{2}\,\cos 2t\;\!x_1^{}+\sqrt{2}\,\sin 2t\;\!x_2^{}+x_3^{}
\hfill
$
\\[2.5ex]
is reducible by the group\footnote[2]{ 
$P_n^{}(\omega)$ is the multiplicative group of 
$\omega\!$-periodic invertible continuously differentiable square matrices of order $n.$
At the same time $P_n^{}(\omega)$ is a sub-group of the Lyapunov group $L(n).$

$L(n)$ is the multiplicative group of 
invertible continuously differentiable on $T=(0;{}+\infty)$ 
square matrices of order $n$ such that 
these matrices and their inverse matrices are bounded on the interval $T.$

}
\vspace{0.5ex}
$\!P_4^{}(2\pi).$ Indeed, the 
$2\pi\!$-periodic nondegenerate transformation
\\[1.5ex]
\mbox{}\hfill
$
y_1^{}=\cos t\;\!x_1^{}+\sin t\;\!x_2^{},\ 
y_2^{}={}-\sin t\;\!x_1^{}+\cos t\;\!x_2^{}, \
y_3^{}=\cos t\;\!x_3^{}+\sin t\;\!x_4^{},\ 
y_4^{}={}-\sin t\;\!x_3^{}+\cos t\;\!x_4^{}
\hfill
$
\\[1.75ex]
reduces the system (3.5) to the linear system with constant coefficients
\\[2.25ex]
\mbox{}\hfill
$
\dfrac{dy_1^{}}{dt}=\sqrt{2}\,y_3^{},
\qquad
\dfrac{dy_2^{}}{dt}={}-\sqrt{2}\,y_4^{},
\qquad 
\dfrac{dy_3^{}}{dt}=\sqrt{2}\,y_2^{},
\qquad  
\dfrac{dy_4^{}}{dt}=\sqrt{2}\,y_1^{}.
\hfill 
$
\\[2.5ex]
\indent
By Theorem 3.2, using the eigenvalues $\lambda_1^{}=1+i,\ \lambda_2^{}={}-1+i$
\vspace{1ex}
and the corresponding complex eigenvectors 
\vspace{0.5ex}
$
\nu^1=(1+i,1-i,\sqrt{2},\sqrt{2}\,i),\ 
\nu^2=(1+i,{}-1+i,{}-\sqrt{2}\,i,{}-\sqrt{2}\,)
$
of the matrix $C,$ we can build the basis of first integrals for the reducible system (3.5)

\newpage

\mbox{}
\\[-1.5ex]
\mbox{}\hfill
$
F_1^{}\colon (t,x)\to
\Bigl(
\bigl((\cos t-\sin t)\;\!x_1^{}+(\cos t+\sin t)\;\!x_2^{}+
\sqrt{2}\;\!\cos t\;\!x_3^{}+\sqrt{2}\;\!\sin t\;\!x_4^{}\bigr)^2 \ +
\hfill
$
\\[1.5ex]
\mbox{}\hfill
$
+\ \bigl((\cos t+\sin t)\;\!x_1^{}+(\sin t-\cos t)\;\!x_2^{}-
\sqrt{2}\;\!\sin t\;\!x_3^{}+\sqrt{2}\;\!\cos t\;\!x_4^{}\bigr)^2
\Bigr)e^{{}-2t}
$
for all 
$
(t,x)\in {\mathbb R}^5,
\hfill 
$
\\[2.25ex]
\mbox{}\hfill
$
F_2^{}\colon (t,x)\to\
\arctan\dfrac{(\cos t+\sin t)\;\!x_1^{}+(\sin t-\cos t)\;\!x_2^{}-
\sqrt{2}\;\!\sin t\;\! x_3^{}+\sqrt{2}\;\!\cos t\;\!x_4^{}}
{(\cos t-\sin t)\;\!x_1^{}+(\cos t+\sin t)\;\!x_2^{}+\sqrt{2}\;\!\cos t\;\!x_3^{}+
\sqrt{2}\;\!\sin t\;\!x_4^{}} \, -\,  t
\hfill 
$
\\[1.75ex]
\mbox{}\hfill
for all 
$
(t,x)\in \Lambda,
\hfill
$
\\[2ex]
\mbox{}\hfill
$
F_3^{}\colon (t,x)\to
\Bigl(
\bigl((\cos t+\sin t)\;\!x_1^{}+(\sin t-\cos t)\;\!x_2^{}+
\sqrt{2}\;\!\sin t\;\!x_3^{}-\sqrt{2}\;\!\cos t\;\!x_4^{}\bigr)^2 \ +
\hfill
$
\\[1.5ex]
\mbox{}\hfill
$
+\ \bigl((\cos t-\sin t)\;\!x_1^{}+(\cos t+\sin t)\;\!x_2^{}-
\sqrt{2}\;\!\cos t\;\!x_3^{}-\sqrt{2}\;\!\sin t\;\!x_4^{}\bigr)^2
\Bigr)e^{2t}
$
for all 
$
(t,x)\in {\mathbb R}^5,
\hfill 
$
\\[2.25ex]
\mbox{}\hfill
$
F_4^{}\colon (t,x)\to\
\arctan\dfrac{(\cos t-\sin t)\;\!x_1^{}+(\cos t+\sin t)\;\!x_2^{}-
\sqrt{2}\;\!\cos t\;\! x_3^{}-\sqrt{2}\;\!\sin t\;\!x_4^{}}
{(\cos t+\sin t)\;\!x_1^{}+(\sin t-\cos t)\;\!x_2^{}+\sqrt{2}\;\!\sin t\;\!x_3^{}-
\sqrt{2}\;\!\cos t\;\!x_4^{}} \,-\, t
\hfill 
$
\\[1.75ex]
\mbox{}\hfill
for all 
$
(t,x)\in \Lambda,
\hfill
$
\\[1.75ex]
where a domain 
\vspace{0.5ex}
$\Lambda\subset
\big\{(t,x)\colon 
(\cos t-\sin t)\;\!x_1^{}+(\cos t+\sin t)\;\!x_2^{}+
\sqrt{2}\;\!\cos t\;\!x_3^{}+\sqrt{2}\;\!\sin t\;\!x_4^{}\ne 0,
$
$
(\cos t+\sin t)\;\!x_1^{}+(\sin t-\cos t)\;\!x_2^{}+\sqrt{2}\;\!\sin t\;\!x_3^{}-
\sqrt{2}\;\!\cos t\;\!x_4^{}\ne 0\big\}
$
of the space ${\mathbb R}^5.$
\vspace{1.5ex}

{\bf Example 3.3.}
The third-order Euler differential system
\\[2ex]
\mbox{}\hfill
$
\dfrac{dx_1^{}}{dt}=x_2^{},
\qquad 
\dfrac{dx_2^{}}{dt}=x_3^{},
\qquad 
\dfrac{dx_3^{}}{dt}=\dfrac{1}{t^3}\, x_1^{}- 
\dfrac{1}{t^2}\, x_2^{}-\dfrac{2}{t}\, x_3^{}
$
\hfill (3.6)
\\[2.25ex]
is reducible to the linear system with constant coefficients
\\[2ex]
\mbox{}\hfill
$
\dfrac{dy_1^{}}{d\tau}=y_2^{},
\qquad
\dfrac{dy_2^{}}{d\tau}=y_2^{}+y_3^{},
\qquad
\dfrac{dy_3^{}}{d\tau}=y_1^{}-y_2^{}
\hfill
$
\\[2.25ex]
by the exponential group\footnote[1]{ 
${\rm Exp}(n)$ is the multiplicative group of 
invertible continuously differentiable on the interval $(0;{}+\infty)$
square matrices ${\rm g}$ of order $n$ such that
$\lim\limits_{t\to{}+\infty}\,\dfrac{1}{t}\, \|{\rm g}^{{}\pm 1}(t)\|=0.$

At the same time the Lyapunov group $L(n)$ is a sub-group of the exponential group ${\rm Exp}(n).$

}
%\vspace{0.35ex}
$\!\!{\rm Exp}(3)$ with the transformation 
\\[1.5ex]
\mbox{}\hfill
$
y_1^{}=x_1^{},
\quad \
y_2^{}=tx_2^{},
\quad \
y_3^{}=t^2x_3^{},
\quad \
t=e^{\tau}.
\hfill
$
\\[2.25ex]
\indent
The matrix 
\vspace{0.75ex}
$
C\! =\! \left\|\!\!
\begin{array}{rrr}
0& 0 &   1
\\
1 &  1 &\!  {}- 1
\\
0& 1&  0
\end{array}
\!\! \right\|
$
has the eigenvalues 
$\lambda_{1}^{}\!=\!1,\ \lambda_{2}^{}\!=i,$ and $\!\lambda_{3}^{}\!=\!{}-i$
cor\-res\-pon\-ding to the eigenvectors
\vspace{0.75ex}
$\nu^{1}\!=(1,1,1),\ \nu^{2}\!=(1,{}-1,i),$ and $\nu^{3}\!=(1,{}-1,{}-i).$

An integral basis of the reducible system (3.6) is the scalar functions 
\\[2.25ex]
\mbox{}\hfill
$
F_1^{}\colon (t,x)\to\ 
\dfrac{1}{t}\,x_1^{}+ x_2^{}+tx_3^{}
$
\ for all 
$
(t,x)\in (0;{}+\infty)\times {\mathbb R}^3
\quad
$
(by Theorem 3.1),
\hfill\mbox{}
\\[2.75ex]
\mbox{}\hfill
$
F_2^{}\colon (t,x)\to\ 
(x_1^{}-tx_2^{})^2+t^4x_3^{2}
$
\ for all 
$
(t,x)\in (0;{}+\infty)\times {\mathbb R}^3
\quad
$
(by Theorem 3.2),
\hfill\mbox{}
\\[2ex]
and 
\\[1.5ex]
\mbox{}\hfill
$
F_3^{}\colon (t,x)\to\ 
\arctan\dfrac{t^2x_3^{}}{x_1^{}-tx_2^{}}\,-\,\ln t
$
\ for all 
$
(t,x)\in \Lambda
\quad
$
(by Theorem 3.2),
\hfill\mbox{}
\\[2ex]
where $\Lambda$ is any domain from the set 
$\{(t,x)\colon t>0,\  x_1^{}-tx_2^{}\ne 0\}.$
%\vspace{1.25ex}

\newpage

{\bf Theorem 3.3.}
\vspace{0.25ex}
{\it
Suppose the system {\rm (3.1)} is reducible to the system {\rm (3.2)} 
by a transformation matrix ${\rm g}\in G,$ and 
$\lambda$ is an eigenvalue with elementary divisor of multiplicity $m\geq 2$ of the matrix $C$
\vspace{0.35ex}
corresponding to a real eigenvector $\nu^{0}$ and to a 
real generalized eigenvector $\nu^{1}$ of the {\rm 1}-st order. 
Then the reducible system {\rm (3.1)} has the first integral}
\\[1.75ex]
\mbox{}\hfill                                            % (3.7)
$
F\colon (t,x)\to \
\dfrac{\nu^{1}{\rm g}(t)\;\!x}{\nu^{0}{\rm g}(t)\;\!x}\,-\,t
$ 
\ for all 
$
(t,x)\in \Lambda,
\quad 
\Lambda\subset \{(t,x)\colon t\in J,\ \nu^{0}{\rm g}(t)\;\!x\ne 0\}.
$
\hfill (3.7)
\\[2ex]
\indent
{\sl Proof}.
From the matrix identity (3.3) it follows that
\\[1.5ex]
\mbox{}\hfill
$
{\frak A}\;\!\nu^{1}{\rm g}(t)x=
\partial_{t}^{}\;\!\nu^{1}{\rm g}(t)x+A(t)\;\!x\,\partial_{x}^{}\;\!\nu^{1}{\rm g}(t)x=
\nu^{1}{\rm g}^{\prime}(t)x\, +\, A(t)\;\!x\,\nu^{1}{\rm g}(t)\,=
\hfill
$
\\[1.75ex]
\mbox{}\hfill
$
=\nu^{1}\bigl(B\,{\rm g}(t)-\;\!{\rm g}(t)\;\!A(t)\bigr)\;\!x \,+\,
A(t)\;\!x\,\nu^{1}{\rm g}(t) =
\nu^{1} B\,{\rm g}(t)\;\!x =
C\nu^{1} \,{\rm g}(t)\;\!x =
(\lambda\nu^1+\nu^0){\rm g}(t)\;\!x=
\hfill
$
\\[1.75ex]
\mbox{}\hfill
$
=\lambda\, \nu^{1}{\rm g}(t)x+\nu^{0}{\rm g}(t)x
$
\ for all 
$
(t,x)\in J\times {\mathbb R}^{n}.
\hfill
$
\\[1.75ex]
\indent
Using this identity and Lemma 3.1, we get on the domain $\Lambda$ 
\\[2ex]
\mbox{}\hfill                                            
$
{\frak A}\;\!F(t,x)=
{\frak A}\;\!\dfrac{\nu^{1}{\rm g}(t)x}{\nu^{0}{\rm g}(t)x}\,-\,{\frak A}\;\!t=
\dfrac{(\lambda\, \nu^{1}{\rm g}(t)x+\nu^{0}{\rm g}(t)x)\;\!\nu^{0}{\rm g}(t)x-
\lambda\, \nu^{0}{\rm g}(t)x\,\nu^{1}{\rm g}(t)x}{(\nu^{0}{\rm g}(t)x)^2}-1=0. \k
\hfill
$
\\[2.5ex]
\indent
{\it In the complex case},
if $\lambda$ is a complex eigenvalue, then
from the complex-valued first integral (3.7) of 
the reducible differential system (3.1), 
we obtain the real first integrals 
\\[2ex]
\mbox{}\hfill
$
F_1^{}\colon (t,x)\to\
\dfrac{{\stackrel{*}{\nu}}\;\!{}^{0}{\rm g}(t)\;\!x\; {\stackrel{*}{\nu}}{}^{1}{\rm g}(t)\;\!x \,+\,
\widetilde{\nu}{}^{\,0}{\rm g}(t)\;\!x\;\widetilde{\nu}{}^{\,1}{\rm g}(t)\;\!x}
{\bigl({\stackrel{*}{\nu}}\;\!{}^{0}{\rm g}(t)\;\!x\bigr)^{2} +
\bigl(\widetilde{\nu}{}^{\,0}{\rm g}(t)\;\!x\bigr)^{2}} \ - \, t
$
\ for all 
$
(t,x)\in\Lambda
\hfill
$
\\[1.25ex]
and
\\[1.25ex]
\mbox{}\hfill
$
F_2^{}\colon (t,x)\to\
\dfrac{{\stackrel{*}{\nu}}\;\!{}^{0}{\rm g}(t)\;\!x\; \widetilde{\nu}{}^{\,1}{\rm g}(t)\;\!x \,-\,
\widetilde{\nu}{}^{\,0}{\rm g}(t)\;\!x\; {\stackrel{*}{\nu}}{}^{1}{\rm g}(t)\;\!x}
{\bigl({\stackrel{*}{\nu}}\;\!{}^{0}{\rm g}(t)\;\!x\bigr)^{2} +
\bigl(\widetilde{\nu}{}^{\,0}{\rm g}(t)\;\!x\bigr)^{2}}
$ 
\ for all 
$
(t,x)\in\Lambda
\hfill
$
\\[2ex]
where $\Lambda$ is any domain from the set
\vspace{0.5ex}
$\bigl\{(t,x)\colon t\in J, \ \bigl({\stackrel{*}{\nu}}\;\! {}^0{\rm g}(t)\;\!x\bigr)^2 +
\bigl(\widetilde{\nu}{}^{\,0}{\rm g}(t)\;\!x\bigr)^2\ne 0\bigr\}$ of the space ${\mathbb R}^{n+1},$ 
the real vectors ${\stackrel{*}{\nu}}{}^{\,k}={\rm Re}\,\nu^{k},\ 
\widetilde{\nu}{}^{\,k}={\rm Im}\, \nu^{k},\ k=0, 1.$
\vspace{1ex}

{\bf Example 3.4.}
\vspace{0.25ex}
Linear Hamiltonian systems of second order 
reducible by orthogonal transformation group\footnote[1]{ 
$O(n)$ is the multiplicative group of continuously differentiable on an interval $T\subset {\mathbb R}$
squarte orthogonal matrices of order $n.$
At the same time if $T=(0;{}+\infty),$ then the orthogonal group of transformations
$O(n)$ is a sub-group of the Lyapunov group $L(n).$ 

}\!
are the linear differential systems of the form [8, pp. 142 -- 143]
\\[2.5ex]                             
\mbox{}\hfill                    % (3.8)
$
\dfrac{dx}{dt}=IA(t)\;\!x,
\quad
x\in {\mathbb R}^2,
\quad
A\colon t\to
\left\|\!\!
\begin{array}{cc}
\psi(t)+\beta(t)& \alpha(t)
\\[0.5ex]
\alpha(t) &  \psi(t)-\beta(t)
\end{array}
\!\! \right\|
$
\ for all 
$
t\in J,
$
\hfill (3.8)
\\[2.25ex]
where the simplectic matrix
\vspace{1ex}
$
I=
\left\|\!\!
\begin{array}{rc}
0& 1
\\
{}-1 &  0
\end{array}
\!\!\right\|,
$
the functions 
$
\alpha\colon t\to\, a\cos 2\varphi(t)-b\sin 2\varphi(t),
$
$
\beta\colon t\to\, a\sin 2\varphi(t)+b\cos 2\varphi(t), \
\psi\colon t\to\, \varphi^{\prime}(t)+c
$
\vspace{0.75ex}
for all 
$
t\in J
$
are continuously differentiable on the interval $J\subset {\mathbb R},$
and $a,\ b,\ c$ are some real numbers.
\vspace{0.35ex}

The orthogonal transformation
\vspace{0.5ex}
$
y_1^{}\!=\cos\varphi(t)\;\!x_1^{}-\,\sin\varphi(t)\;\!x_2^{}, \
y_2^{}=\sin\varphi(t)\;\!x_1^{}+\cos\varphi(t)\;\!x_2^{}
$
reduces the linear differential system (3.8) 
\vspace{0.5ex}
to the linear autonomous  Hamiltonian system
$
\dfrac{dy}{dt}=IB\;\!y
$
with the constant matrix
$
B=
\left\|\!\!
\begin{array}{cc}
c+b& a
\\
a &  c-b
\end{array}
\!\!\right\|.
$
\vspace{1ex}

Consider the real number $D=a^2+b^2-c^2.$ 
\vspace{0.35ex}
We have three possible cases for building first integrals of the 
reducible linear Hamiltonian system (3.8). 

\newpage

Let $D=0.$ 
\vspace{0.35ex}
Then, using the real eigenvector $\nu^0=(c+b, a),$  
the real 1-st order generalized eigenvector $\nu^1=(c+b,a-1),$ and 
the corresponding double eigenvalue $\lambda_1^{}=0,$
\vspace{0.35ex}
we can build (by Theorems 3.1 and 3.3) the integral basis on a domain $\Lambda$
of system (3.8)
\\[2ex]
\mbox{}\hfill
$
F_1^{}\colon (t,x)\to\, 
\bigl((c+b)\cos\varphi(t)+a\sin\varphi(t)\bigr)\;\!x_1^{} +
\bigl(a\cos\varphi(t)-(c+b)\sin\varphi(t)\bigr)\;\!x_2^{}\;\!,
\hfill
$
\\[2.25ex]
\mbox{}\hfill
$
F_2^{}\colon (t,x)\to 
\dfrac{\bigl((c\!+\!b)\cos\varphi(t)\!+\!(a\!-\!1)\sin\varphi(t)\bigr) x_1^{}\! +\!
\bigl((a\!-\!1)\cos\varphi(t)\!-\!(c\!+\!b)\sin\varphi(t)\bigr) x_2^{}\;\!}
{\bigl((c+b)\cos\varphi(t)+a\sin\varphi(t)\bigr)\;\!x_1^{} +
\bigl(a\cos\varphi(t)-(c+b)\sin\varphi(t)\bigr)\;\!x_2^{}}\,- \, t,
\hfill
$
\\[2ex]
where $\Lambda$ is any domain from the set
\\[1.5ex]
\mbox{}\hfill
$
M=\bigl\{(t,x)\colon t\in J,  \
\bigl((c+b)\cos\varphi(t)+a\sin\varphi(t)\bigr)x_1^{} +
\bigl(a\cos\varphi(t)-(c+b)\sin\varphi(t)\bigr)x_2^{}\ne 0\bigr\}.
\hfill
$
\\[1.75ex]
\indent
Let $D>0.$ 
\vspace{0.5ex}
Then, using the real eigenvectors $\nu^k=(c+b, a-\lambda_k^{}),\, k=1,2,$ 
and the corresponding real eigenvalues
$\lambda_1^{}=\sqrt{D},\ \lambda_2^{}={}-\sqrt{D},$
\vspace{0.35ex}
we can construct (by Theorem 3.1) the basis of first integrals for 
the Hamiltonian system (3.8)
\\[1.5ex]
\mbox{}\hfill
$
F_k^{}\!\colon (t,x)\!\to 
\Bigl(\!\bigl((c+b)\cos\varphi(t)+(a-\lambda_k^{})\sin\varphi(t)\bigr)x_1^{} +
\bigl((a-\lambda_k^{})\cos\varphi(t)-(c+b)\sin\varphi(t)\bigr)x_2^{}\!\Bigr)e^{{}-\lambda_kt}
\hfill
$
\\[1.5ex]
\mbox{}\hfill
for all 
$
(t,x)\in J\times {\mathbb R}^2,
\quad
k=1,2.
\hfill
$
\\[1.5ex]
\indent
Let $D<0.$ 
\vspace{0.75ex}
Then, using the eigenvectors $\nu^k=(c+b, a-\lambda_k^{}),\ k=1,2,$ 
and the cor\-res\-pon\-ding complex eigenvalues
\vspace{0.35ex}
$\lambda_1^{}\!=\sqrt{{}-D}\,i,\ \lambda_2^{}=-\,\sqrt{{}-D}\,i,$ 
we can find (by Theorem 3.2) the functionally independent 
first integrals of the Hamiltonian system (3.8)
\\[1.5ex]
\mbox{}\hfill
$
F_1^{}\colon (t,x)\to \,
\bigl( \bigl((c+b)\cos\varphi(t)+a\sin\varphi(t)\bigr)\;\!x_1^{} +
\bigl(a\cos\varphi(t)-(c+b)\sin\varphi(t)\bigr)\;\!x_2^{}\bigr)^2\ -
\hfill
$
\\[1.5ex]
\mbox{}\hfill
$
-\ D\;\!\bigl(\sin\varphi(t)\;\!x_1^{}+\cos\varphi(t)\;\!x_2^{}\bigr)^2
$
\ for all 
$
(t,x)\in J\times {\mathbb R}^2,
\hfill
$
\\[2.25ex]
\mbox{}\hfill
$
F_2^{}\!\colon (t,x)\!\to 
\arctan\dfrac{\sqrt{-\,D}\,\bigl(\sin\varphi(t)\;\!x_1^{} + \cos\varphi(t)\;\!x_2^{}\bigr)}
{\bigl((c+b)\cos\varphi(t)\!+\!a\sin\varphi(t)\bigr)x_1^{}\! +\!
\bigl(a\cos\varphi(t)\!-\!(c+b)\sin\varphi(t)\bigr)x_2^{}}\, +\sqrt{-\,D}\,t
\hfill
$
\\[2ex]
\mbox{}\hfill
for all 
$
(t,x)\in \Lambda\subset M,
\quad
M\subset {\mathbb R}^3.
\hfill
$
\\[2.25ex]
\indent
{\bf Theorem 3.4.}
\vspace{0.25ex}
{\it
Suppose the system {\rm (3.1)} is reducible to the system {\rm (3.2)} 
by a transformation matrix ${\rm g}\in G,$ and 
\vspace{0.25ex}
$\lambda$ is the eigenvalue with elementary divisor of multiplicity $m\geq 2$ 
of the matrix $C$ corresponding to a real eigenvector $\nu^{0}$ and to 
\vspace{0.35ex}
real generalized eigenvectors $\nu^{k},\ k=1,\ldots, m-1.$ 
Then the system {\rm (3.1)} has 
the fun\-c\-ti\-o\-na\-l\-ly independent first integrals
\\[2ex]
\mbox{}\hfill                                       % (3.9)
$
F_{\xi}^{}\colon (t,x)\to\, \Psi_{\xi}^{}(t,x)
$ 
\ for all 
$
(t,x)\in \Lambda,
\quad 
\xi=2,\ldots, m-1,
$
\hfill {\rm(3.9)}
\\[2.25ex]
where the functions 
$\Psi_{\xi}^{}\colon \Lambda\to {\mathbb R},\ \xi=2,\ldots, m-1,$ 
are the solution to the functional system
\\[1.5ex]
\mbox{}\hfill   % (3.10)
$
\displaystyle
\nu^{k}{\rm g}(t)x =
\sum\limits_{\tau=1}^{k}{\textstyle\binom{k-1}{\tau-1}}
\Psi_{\tau}^{}(t,x)\,\nu^{k-\tau}{\rm g}(t)x
$ 
\ for all 
$
(t,x)\in \Lambda,
\ \ k=1,\ldots, m-1,
$
\hfill {\rm (3.10)}
\\[1.5ex]
and $\Lambda$ is any domain from the set 
$\{(t,x)\colon t\in J,\ \nu^0{\rm g}(t)\;\!x\ne 0\}\subset J\times {\mathbb R}^n.$
}
\vspace{1ex}

{\sl Proof}.
From the matrix identity (3.3) 
and the notion of generalized eigenvectors for matrix (De\-fi\-ni\-ti\-on 1.1) it follows that
\\[1.5ex]
\mbox{}\hfill
$
{\frak A}\;\!\nu^{k}{\rm g}(t)\;\!x=
\partial_{t}^{}\;\!\nu^{k}{\rm g}(t)\;\!x+A(t)\;\!x\,\partial_{x}^{}\;\!\nu^{k}{\rm g}(t)\;\!x=
\nu^{k}{\rm g}^{\prime}(t)\;\!x\, +\, A(t)\;\!x\,\nu^{k}{\rm g}(t)\,=
\hfill
$
\\[2ex]
\mbox{}\hfill
$
=\nu^{k}\bigl(B\,{\rm g}(t)-\;\!{\rm g}(t)\;\!A(t)\bigr)\;\!x \,+\,
A(t)\;\!x\,\nu^{k}{\rm g}(t) =
\nu^{k} B\,{\rm g}(t)\;\!x =
C\nu^{k} \,{\rm g}(t)\;\!x =
(\lambda\nu^k+k\nu^{k-1}){\rm g}(t)\;\!x=
\hfill
$
\\[2ex]
\mbox{}\hfill
$
=\lambda\, \nu^{k}{\rm g}(t)\;\!x+k\nu^{k-1}{\rm g}(t)\;\!x
$
\ for all 
$
(t,x)\in J\times {\mathbb R}^{n},
\quad
k=1,\ldots, m-1.
\hfill
$
\\[1.5ex]
\indent
Hence using  Lemma 3.1, we obtain the system of identities 
\\[1.75ex]
\mbox{}\hfill                                     % (3.11)
$
{\frak A}\,\nu^0{\rm g}(t)\;\!x = \lambda\,\nu^0{\rm g}(t)\;\!x
$ 
\ for all 
$
(t,x)\in J\times {\mathbb R}^{n},
\hfill
$
\\[-0.35ex]
\mbox{}\hfill (3.11)
\\[-0.15ex]
\mbox{}\hfill
$
{\frak A}\,\nu^{k}{\rm g}(t)\;\!x =
\lambda\,\nu^{k}{\rm g}(t)\;\!x+k\,\nu^{k-1}{\rm g}(t)\;\!x
$
\ for all 
$
(t,x)\in J\times {\mathbb R}^{n},
\ \ 
k=1,\ldots, m-1.
\hfill
$
\\[1.75ex]
\indent
The functional system (3.10) has the determinant 
\vspace{0.5ex}
$(\nu^0{\rm g}(t)\;\!x)^{m-1}\ne 0$ for all $(t,x)\in\Lambda,$
where a domain $\Lambda\subset \{(t,x)\colon t\in J, \ \nu^0{\rm g}(t)\;\!x\ne 0\}.$ 
\vspace{0.5ex}
Therefore there exists the solution $\Psi_{\tau}^{},\ \tau=1,\ldots, m-1,$ 
on the domain $\Lambda$ of the functional system (3.10).
Let us show that
\\[2.5ex]
\mbox{}\hfill                                        % (3.12)
$
{\frak A}\Psi_{k}^{}(t,x)=
\left[\!
\begin{array}{lll}
1\! & \text{for all}\ \, (t,x)\in \Lambda, & k=1,
\\[1.25ex]
0\! & \text{for all}\ \, (t,x)\in \Lambda, & k=2,\ldots, m-1,
\end{array}
\right.
$
\hfill (3.12)
\\[2.5ex]
\indent
The proof of identities (3.12) is by induction on $m.$

For $m=2$ and $m=3,$ the assertions (3.12) follows from the identities (3.11).
%\vspace{0.35ex}

Assume that the identities (3.12) for $m=\varepsilon$ is true. 
Then, using the system of identities (3.11) and the identities (3.10) for 
$m=\varepsilon+1,\ m=\varepsilon,$ we get
\\[2ex]
\mbox{}\hfill
$
{\frak A}\,\nu^{\varepsilon}{\rm g}(t)\;\!x =
\lambda\ {\displaystyle \sum\limits_{\tau=1}^{\varepsilon}}
{\textstyle\binom{\varepsilon-1}{\tau-1}}\Psi_{\tau}^{}(t,x)\, \nu^{\varepsilon-\tau}{\rm g}(t)\;\!x \,+\,
(\varepsilon-1)\, {\displaystyle \sum\limits_{\tau=1}^{\varepsilon-1}}
{\textstyle\binom{\varepsilon-2}{\tau-1}}\Psi_{\tau}^{}(t,x)\, \nu^{\varepsilon-\tau-1}{\rm g}(t)\;\!x \ +
\hfill
$
\\[2.25ex]
\mbox{}\hfill
$
+\  \nu^{\varepsilon-1}{\rm g}(t)\;\!x \,+\, 
\nu^{0}{\rm g}(t)\;\!x\ {\frak A}\;\!\Psi_{\varepsilon}^{}(t,x)
$
\ for all 
$
(t,x)\in \Lambda.
\hfill
$
\\[2ex]
\indent
Now taking into account the system (3.10) with $k=\varepsilon-1$ and  $k=\varepsilon,$
\vspace{0.35ex}
the identity (3.11) with $k=\varepsilon,$ and $\nu^{0}{\rm g}(t)\;\!x\ne 0$ for all 
$(t,x)\in \Lambda,$ we have
\\[2ex]
\mbox{}\hfill
$
{\frak A}\, \Psi_{\varepsilon}^{}(t,x)=0
$
\ for all 
$
(t,x)\in \Lambda.
\hfill
$
\\[1.75ex]
\indent
This implies that the identities (3.12) for $m=\varepsilon+1$ are true.
So by the principle of mathematical induction, the statement (3.12) is true for every 
natural number $m\geq 2.$ 

Now from the method of building scalar functions (3.9) it follows that 
the reducible differential system (3.1) 
has the functionally independent first integrals (3.9). \k
\vspace{0.75ex}

The proof of Theorem 3.4 is true both for the case of the real eigenvalue $\lambda$ and 
for the case of the com\-p\-lex eigenvalue $\lambda\ ({\rm Im}\,\lambda\ne 0).$
{\it In the complex case}, from the complex-valued first integrals (3.9) of 
system (3.1), we obtain the real first integrals of system (3.1)
\\[1.75ex]
\mbox{}\hfill
$
F_{\xi}^{1}\colon (t,x)\to\, \mbox{Re}\,\Psi_\xi^{}(t,x),
\quad 
F_{\xi}^{2}\colon (t,x)\to\, \mbox{Im}\,\Psi_\xi^{}(t,x)
$
\ for all 
$
(t,x)\in \Lambda, 
\quad
\xi=2,\ldots, m-1,
\hfill
$
\\[1.5ex]
where $\Lambda$ is any domain from the set
$\bigl\{(t,x)\colon t\in J,\ \bigl(\,{\stackrel{*}{\nu}}{}^0{\rm g}(t)\;\!x\bigr)^2 +
\bigl(\,\widetilde{\nu}{}^{\,0}{\rm g}(t)\;\!x\bigr)^2\ne 0\bigr\}.$ 
\vspace{1.5ex}

{\bf Example 3.5.}
The third-order Euler differential system
\\[2ex]
\mbox{}\hfill                                               % (3.13)
$
\dfrac{dx_1^{}}{dt}=x_2^{},
\qquad 
\dfrac{dx_2^{}}{dt}=x_3^{},
\qquad 
\dfrac{dx_3^{}}{dt}=\dfrac{1}{t^3}\, x_1^{}- 
\dfrac{1}{t^2}\, x_2^{}
$
\hfill (3.13)
\\[2ex]
is reducible to the linear differential system with constant coefficients
\\[1.5ex]
\mbox{}\hfill
$
\dfrac{dy_1^{}}{d\tau}=y_2^{},
\qquad
\dfrac{dy_2^{}}{d\tau}=y_2^{}+y_3^{},
\qquad
\dfrac{dy_3^{}}{d\tau}=y_1^{}-y_2^{}+2y_3^{}
\hfill
$
\\[1.75ex]
by the exponential group ${\rm Exp}(3)$ with the transformation 
\\[1.5ex]
\mbox{}\hfill
$
y_1^{}=x_1^{},
\quad
y_2^{}=tx_2^{},
\quad
y_3^{}=t^2x_3^{},
\quad
t=e^{\tau}.
\hfill
$
\\[2.5ex]
\indent
The matrix 
$
C = \left\|\!\!
\begin{array}{ccr}
0& 0\!\!&   1
\\
1 &  1\!\!& \!\! {}- 1
\\
0& 1\!\!&  2
\end{array}
\!\! \right\|
$
has the triple eigenvalue $\lambda_{1}^{}\!=1$
\vspace{0.5ex}
with the elementary divisor $(\lambda-1)^3$ 
\vspace{0.5ex}
corresponding to the eigenvector
$\nu^{0}\!=(1,{}-1,1),$ to the generalized eigenvector of the 1-st order $\nu^{1}=(0, 0, 1),$ 
and to the generalized eigenvector of the 2-nd order $\nu^{2}=(0, 2, 0).$
\vspace{0.35ex}

A basis of first integrals for the reducible system (3.13) is the scalar functions 
\\[2.25ex]
\mbox{}\hfill
$
F_1^{}\colon (t,x)\to\ 
\dfrac{1}{t}\,x_1^{}- x_2^{}+tx_3^{}
$
\ for all 
$
(t,x)\in (0;{}+\infty)\times {\mathbb R}^3
\quad
$
(by Theorem 3.1),
\hfill\mbox{}
\\[2.5ex]
\mbox{}\hfill
$
F_2^{}\colon (t,x)\to\ 
\dfrac{t^2x_3^{}}{x_1^{}-t\;\!x_2^{}+t^2x_3^{}} -\ln t
$
\ for all 
$
(t,x)\in (0;{}+\infty)\times {\mathbb R}^3
\quad
$
(by Theorem 3.3),
\hfill\mbox{}
\\[2ex]
and (by Theorem 3.4)
\\[1.5ex]
\mbox{}\hfill
$
F_3^{}\colon (t,x)\to\ 
\dfrac{2\;\!t\;\!x_2^{}(x_1^{}-t\;\!x_2^{}+t^2x_3^{})-\;\!t^4\;\!x_3^{2}}{(x_1^{}-t\;\!x_2^{}+t^2x_3^{})^2}
$
\ for all 
$
(t,x)\in \Lambda,
\hfill
$
\\[2ex]
where $\Lambda$ is a domain from the set  
$\bigl\{(t,x)\colon t>0,\  x_1^{}-t\;\!x_2^{}+t^2x_3^{}\ne 0\bigr\}.$
\\[3.75ex]
\centerline{                               %{П. 2, Параграф 1}
\bf 3.2.
Linear nonhomogeneous differential system
}
\\[1.25ex]
\indent
A real linear nonhomogeneous differential system of the $n\!$-th order
\\[1.75ex]
\mbox{}\hfill                                          % (3.14)
$
\dfrac{dx}{dt}=A(t)\,x+f(t),
\quad
f\in C(J),
$
\hfill (3.14)
\\[1.75ex]
such that the system (3.1) is the corresponding homogeneous system of system (3.14).
\vspace{0.35ex}

The system (3.14) is called 
{\it reducible with respect to the nonautonomous transformation group} $G$ if 
there exist a constant matrix $B$ and a matrix ${\rm g}\in G$ such that 
the linear transformation $y={\rm g}(t)\;\!x$ 
reduces the system (3.14) to the system
\\[1.75ex]
\mbox{}\hfill                                          % (3.15)
$
\dfrac{dy}{dt}=B\,y+{\rm g}(t)f(t),
\quad
y\in {\mathbb R}^n.
$
\hfill (3.15)
\\[1.75ex]
\indent
Let us remark that if the nonhomogeneous system (3.14) is reduced to the system (3.15), then
homogeneous system (3.1) is reduced to the system with constant coefficients (3.2).

The reducible system (3.14) is induced the linear differential ope\-ra\-tor of first order
\\[1.75ex]
\mbox{}\hfill                                      
$
{\frak B}(t,x)=\partial_{t}^{}+
(A(t)\;\!x+f(t))\;\!\partial_{x}^{}=
{\frak A}(t,x)+f(t)\;\!\partial_{x}^{}
$
\ for all 
$
(t,x)\in J\times {\mathbb R}^n.
\hfill
$
\\[2ex]
\indent
{\bf 3.2.1. Case of simple elementary divisors}
\\[0.75ex]
\indent
If the matrix $C$ is diagonalizable, then using Theorem 3.5 and Corollary 3.1,
we can build first integrals of the reducible nonhomogeneous differential system (3.14).
\vspace{0.75ex}

{\bf Theorem 3.5.}
{\it
Let the system {\rm (3.14)} be reducible to the system {\rm (3.15)} 
by a transformation matrix $\!{\rm g}\in G,$  
let $\nu\!$ be a real eigenvector of the matrix $\!C\!$ corresponding to the eigenvalue $\!\lambda.$ 
Then a first integral of the reducible system {\rm (3.14)} is the scalar function 
\\[1.5ex]
\mbox{}\hfill                                            % (3.16)
$
\displaystyle
F\colon (t,x)\to\  \nu {\rm g}(t)x\;\!\exp({}-\lambda\;\! t)
-\int\limits_{t_0^{}}^t \nu{\rm g}(\zeta) f(\zeta)\exp({}-\lambda\;\! \zeta)\,d\zeta
$
\ for all 
$
(t,x)\in J\times {\mathbb R}^n,
$
\hfill {\rm (3.16)}
\\[1.5ex]
where $t_0^{}$ is a fixed point from the in\-ter\-val $J.$
}
\vspace{0.75ex}

{\sl Proof}.
From Lemma 3.1 it follows that on the domain $J\times {\mathbb R}^n\colon$
\\[1.5ex]
\mbox{}\hfill                                           
$
\displaystyle
{\frak B}\;\!F(t,x)=
{\frak B}\;\!\nu {\rm g}(t)x\;\!\exp({}-\lambda\;\! t)-
{\frak B}\;\!\int\limits_{t_0^{}}^t \nu{\rm g}(\zeta) f(\zeta)\exp({}-\lambda\;\! \zeta)\,d\zeta=
{}-\lambda\;\!\nu {\rm g}(t)x\exp({}-\lambda\;\! t)\ +
\hfill 
$
\\[1.5ex]
\mbox{}\hfill                                           
$
\displaystyle
+ \ {\frak A}\;\!\nu {\rm g}(t)x\exp({}-\lambda\;\! t)
+f(t)\;\!\partial_x^{}\nu {\rm g}(t)x\exp({}-\lambda\;\! t)
-\partial_t^{}\int\limits_{t_0^{}}^t \nu{\rm g}(\zeta) f(\zeta)\exp({}-\lambda\;\! \zeta)\,d\zeta=0.
\hfill 
$
\\[1.5ex]
\indent
Therefore the function (3.16) is a first integral of the reducible system (3.14). \k
\vspace{0.75ex}

{\bf Corollary 3.1.}
{\it
Let the system {\rm (3.14)} be reducible to the system {\rm (3.15)} 
by a transformation matrix ${\rm g}\in G,$ and
let $\nu={\stackrel{*}{\nu}}+\widetilde{\nu}\,i\
({\rm Re}\,\nu={\stackrel{*}{\nu}},\ {\rm Im}\,\nu=\widetilde{\nu}\,)$
be an eigenvector of the matrix $C$
corresponding to the complex eigenvalue 
$\lambda={\stackrel{*}{\lambda}}+\widetilde{\lambda}\,i$
$({\rm Re}\,\lambda\!=\!{\stackrel{*}{\lambda}},\, {\rm Im}\, \lambda\!=\!\widetilde{\lambda}\,).\!$
\vspace{0.35ex}
Then first integrals of the reducible system {\rm (3.14)} are the scalar functions 
\\[1.5ex]
\mbox{}\hfill                                           % (3.17)
$
\displaystyle
F_{{}_{\scriptstyle \theta}}\colon (t,x)\to\ 
\alpha_{{}_{\scriptstyle \theta}}(t, x)-
\int\limits_{t_0^{}}^{t} \alpha_{{}_{\scriptstyle \theta}}(\zeta, f(\zeta))\,d\zeta
$
\ for all 
$
(t,x)\in J\times\R^n,
\ \ \theta=1,2,
$
\hfill {\rm(3.17)}
\\[1.5ex]
where $t_0^{}$ is a fixed point from the in\-ter\-val $J,$ the functions
\\[1ex]
\mbox{}\hfill                                     
$
\displaystyle
\alpha_1^{}\colon (t,x)\to\, 
\bigl(\,{\stackrel{*}{\nu}}{\rm g}(t)\;\!x\;\!\cos\widetilde{\lambda}\,t +
\widetilde{\nu}{\rm g}(t)\;\!x\;\!\sin\widetilde{\lambda}\,t \bigr)\,
\exp\bigl({}-{\stackrel{*}{\lambda}}\,t\bigr)
$
\ for all 
$
(t,x)\in J\times {\mathbb R}^{n},
\hfill
$
\\[1.5ex]
\mbox{}\hfill                                     
$
\displaystyle
\alpha_{2}^{}\colon (t,x)\to\,
\bigl(\,\widetilde{\nu}{\rm g}(t)\;\!x\;\!\cos\widetilde{\lambda}\,t -
{\stackrel{*}{\nu}}{\rm g}(t)\;\!x\;\!\sin\widetilde{\lambda}\,t\bigr)\,
\exp\bigl({}-{\stackrel{*}{\lambda}}\,t\bigr)
$
\ for all 
$
(t,x)\in J\times {\mathbb R}^{n}.
\hfill
$
\\[1.75ex]
}
\indent
{\sl Proof}.  
Formally using Theorem 3.5, we obtain 
the complex-valued function (3.16) is a first integral 
of system (3.14). Then the real and imaginary parts 
of this complex-valued first integral are the real first integrals (3.17) of 
the reducible system (3.14). \k
\vspace{0.75ex}

{\bf Example 3.6.}
Consider the second-order linear nonhomogeneous differential system\footnote[1]{ 
Note that the linear homogeneous differential system corresponding to 
the linear nonhomogeneous differential system (3.18)
is a simplified differential system (has not nutations), which is
describing free oscillations of gyroscopic pendulum in 
twin gyrocompass [17, pp. 528 -- 529].

}\!\!
\\[2ex]                             
\mbox{}\hfill                    % (3.18)
$
\dfrac{dx_1^{}}{dt}=(a\sin\omega t\cos\omega t)\;\!x_1^{}-(b+a\cos^2\omega t)\;\!x_2^{}+f_1^{}(t),
\hfill
$
\\
\mbox{}\hfill (3.18)
\\
\mbox{}\hfill 
$
\dfrac{dx_2^{}}{dt}=(b+a\sin^2\omega t)\;\!x_1^{}-(a\sin\omega t\cos\omega t)\;\!x_2^{}+f_2^{}(t),
\hfill 
$
\\[2ex]
where $a,\, b,$ and $\omega$ are real numbers 
\vspace{0.35ex}
such that $\omega\ne b$ and $\omega\ne a+b,$ 
the scalar functions $f_1^{}\colon J\to{\mathbb R}$ and $f_2^{}\colon J\to{\mathbb R}$ are 
continuous on an interval $J\subset {\mathbb R}.$
\vspace{0.5ex}

The system (3.18) is reducible to the nonhomogeneous system with constant coefficients
\\[2ex]                             
\mbox{}\hfill                   
$
\dfrac{dy_1^{}}{dt}=(\omega -a-b)y_2^{}+f_1^{}(t)\cos\omega t+f_2^{}(t)\sin\omega t,
\ \ \
\dfrac{dy_2^{}}{dt}=(b-\omega)y_1^{}-f_1^{}(t)\sin\omega t+f_2^{}(t)\cos\omega t
\hfill
$
\\[2ex]
by the orthogonal group of transformations $O(J,2)$ with the transformation
\\[1.5ex]                             
\mbox{}\hfill                   
$
y_1^{}\!=\cos\omega t\;\!x_1^{}+\,\sin\omega t\;\!x_2^{},
\quad
y_2^{}={}-\sin\omega t\;\!x_1^{}+\cos\omega t\;\!x_2^{}
\hfill
$
\\[1.5ex]
\indent
Consider the real number $D=(b-\omega)(\omega-a-b).$ 
\vspace{0.35ex}
We have two possible cases for building first integrals of the 
reducible linear differential system (3.18). 
\vspace{0.35ex}

Let $D>0.$ 
\vspace{0.75ex}
Then, using the real eigenvectors
$\nu^k=(b-\omega, \lambda_k^{}),\ k=1,2,$ 
of the matrix 
\\[1.5ex]
\mbox{}\hfill
$
C=
\left\|\!\!
\begin{array}{cc}
0& b-\omega
\\
\omega-a-b &  0
\end{array}
\!\!\right\|,
\hfill
$
\\[1.5ex]
and the corresponding real eigenvalues
\vspace{0.35ex}
$\!\lambda_1^{}\!=\sqrt{D},\ \lambda_2^{}\! =-\,\sqrt{D},$
we can build (by Theorem~3.5) the basis of first integrals for 
the reducible system (3.18)
\\[2ex]
\mbox{}\hfill
$
F_k^{}\colon (t,x_1^{},x_2^{})\to 
\bigl(((b-\omega)\cos\omega t-\lambda_k^{}\sin\omega t)x_1^{} +
(\lambda_k^{}\cos\omega t+(b-\omega)\sin\omega t)x_2^{}\bigr)e^{{}-\lambda_kt} \ -
\hfill
$
\\[1.75ex]
\mbox{}\hfill
$
\displaystyle
-\ \int\limits_{t_0^{}}^{t}
\Bigl(\bigl((b-\omega)\cos\omega\zeta-\lambda_k^{}\sin\omega \zeta\bigr)f_1^{}(\zeta) +
\bigl(\lambda_k^{}\cos\omega\zeta+(b-\omega)\sin\omega\zeta\bigr)f_2^{}(\zeta)\Bigr)
e^{{}-\lambda_k^{}\zeta}\;\!d\zeta
\hfill
$
\\[1ex]
\mbox{}\hfill
for all 
$
(t,x_1^{},x_2^{})\in J\times {\mathbb R}^2,
\quad
k=1,2.
\hfill
$
\\[2ex]
\indent
Let $D<0.$ 
\vspace{0.75ex}
Then, using the eigenvectors $\nu^k=(\omega-b, \lambda_{3-k}^{}),\ k=1,2,$ 
and the cor\-res\-pon\-ding complex eigenvalues
\vspace{0.35ex}
$\lambda_1^{}\!=\sqrt{{}-D}\,i,\ \lambda_2^{}=-\,\sqrt{{}-D}\,i,$ 
we can (by Corollary 3.1) construct the functionally independent 
first integrals of the reducible system (3.18)
\\[1.5ex]
\mbox{}\hfill
$
\displaystyle
F_{{}_{\scriptstyle \theta}}\colon (t,x_1^{},x_2^{})\to 
\alpha_{{}_{\scriptstyle \theta}}(t, x_1^{},x_2^{})-
\int\limits_{t_0^{}}^{t}\! \alpha_{{}_{\scriptstyle \theta}}(\zeta, f_1^{}(\zeta),f_2^{}(\zeta))\;\!d\zeta
$
for all 
$\!
(t,x_1^{},x_2^{})\in J\times\R^2,
\ \theta=1,2,
\hfill
$
\\[1.5ex]
where $t_0^{}$ is a fixed point from the in\-ter\-val $J,$ the scalar functions
\\[1.5ex]
\mbox{}\hfill                                     
$
\displaystyle
\alpha_1^{}\colon (t,x_1^{},x_2^{})\to\, 
\bigl(\,(\omega-b)\cos\sqrt{{}-D}\ t\,\cos\omega t\, +
\,\sqrt{{}-D}\,\sin\sqrt{{}-D}\ t\, \sin\omega t\bigr)\;\!x_1^{} \ +
\hfill                                     
$
\\[2ex]
\mbox{}\hfill                                     
$
+\ \bigl(\,(\omega-b)\cos\sqrt{{}-D}\ t\, \sin\omega t \, -
\, \sqrt{{}-D}\, \sin\sqrt{{}-D}\ t\, \cos\omega t\bigr)\;\!x_2^{}
$
\ for all 
$
(t,x_1^{},x_2^{})\in J\times {\mathbb R}^{2},
\hfill
$
\\[2.15ex]
\mbox{}\hfill                                     
$
\displaystyle
\alpha_2^{}\colon (t,x_1^{},x_2^{})\to\, 
\bigl(\,\sqrt{{}-D}\,\cos\sqrt{{}-D}\ t\, \sin\omega t\,-\,
(\omega-b)\sin\sqrt{{}-D}\ t\,\cos\omega t\bigr)\;\!x_1^{} \ -
\hfill                                     
$
\\[2ex]
\mbox{}\hfill                                     
$
-\ \bigl(\,\sqrt{{}-D}\, \cos\sqrt{{}-D}\ t\, \cos\omega t\,+\,
(\omega-b)\sin\sqrt{{}-D}\ t\, \sin\omega t \bigr)\;\!x_2^{}
$
\ for all 
$
(t,x_1^{},x_2^{})\in J\times {\mathbb R}^{2}.
\hfill
$
\\[3.25ex]
\indent
{\bf 3.2.2. Case of multiple elementary divisors}
\\[0.95ex]
\indent
If the matrix $C$ has multiple elementary divisors, 
then using Theorem 3.6 and Co\-rol\-la\-ry~3.2,
we can build first integrals 
\vspace{0.75ex}
of the reducible nonhomogeneous differential system (3.14).

{\bf Theorem 3.6.}
\vspace{0.25ex}
{\it
Let the system {\rm (3.14)} be reducible to the system {\rm (3.15)} 
by a transformation matrix ${\rm g}\in G,$ and let 
\vspace{0.25ex}
$\lambda$ be the eigenvalue with elementary divisor of multiplicity $m\geq 2$ 
of the matrix $C$ corresponding to a real eigenvector $\nu^{0}$ and to real
\vspace{0.5ex}
generalized eigenvectors $\nu^{k},\ k=1,\ldots, m-1.\!$ 
Then the system {\rm (3.14)} has the fun\-c\-ti\-o\-na\-l\-ly independent first integrals
\\[1.75ex]
\mbox{}\hfill                                   % (3.19)
$
\displaystyle
F_{k+1}^{}\colon (t,x)\to \ 
\nu^k{\rm g}(t)\;\!x\;\!\exp({}-\lambda\;\! t) \, - \,  
\sum\limits_{\tau=0}^{k-1}\!
{\textstyle \binom{k}{\tau}}\;\! t^{k-\tau} F_{\tau+1}^{}(t,x) \, -\,  C_{k}^{}(t)
\hfill
$
\\
\mbox{}\hfill {\rm (3.19)}
\\
\mbox{}\hfill
for all 
$
(t,x)\in J\times {\mathbb R}^n,
\quad
k=1,\ldots, m-1,
\mbox{}\hfill
$
\\[2.25ex]
where the first integral {\rm(}by Theorem {\rm 3.5)}
\\[2ex]
\mbox{}\hfill
$
F_1^{}\colon (t,x)\to\
 \nu^0{\rm g}(t)\;\!x\;\!\exp({}-\lambda\;\! t)-C_0^{}(t)
$
\ for all 
$
(t,x)\in J\times {\mathbb R}^n,
\hfill
$
\\[1.5ex]
the scalar functions 
\\[2ex]
\mbox{}\hfill
$
\displaystyle
C_k^{}\colon t \to \
\int\limits_{t_0^{}}^{t} \bigl(\nu^k{\rm g}(\zeta)f(\zeta)\exp({}-\lambda\;\!\zeta) + 
k\,C_{k-1}^{}(\zeta)\bigr)\, d\zeta
$
\ for all 
$
t\in J, 
\quad 
k=0,\ldots, m-1,
\hfill
$
\\[1.5ex]
and  
\vspace{1ex}
$t_0^{}$ is a fixed point from the in\-ter\-val $J\subset {\mathbb R}.$
}

{\sl Proof}. 
The proof of Theorem 3.6 is by induction on $m.$
\vspace{0.25ex}

By the system of identities (3.11), it follows that
\\[1.75ex]
\mbox{}\hfill                    
$
{\frak B}\bigl(\nu^{\,\varepsilon}{\rm g}(t)\;\!x\;\!\exp({}-\lambda t)\bigr)=
{}-\lambda\nu^{\,\varepsilon}{\rm g}(t)\;\!x\,\exp({}-\lambda t)+
\bigl({\frak A}\;\!\nu^{\,\varepsilon}{\rm g}(t)\;\!x+f(t)\;\!\partial_x^{}\;\!\nu^{\,\varepsilon}{\rm g}(t)\;\!x\bigr)
\exp({}-\lambda t)=
\hfill                    
$
\\[2.25ex]
\mbox{}\hfill                    
$
={}-\lambda\nu^{\,\varepsilon}{\rm g}(t)\;\!x\,\exp({}-\lambda t)+
\bigl(\lambda\nu^{\,\varepsilon}{\rm g}(t)\;\!x+
\varepsilon\,\nu^{\,\varepsilon-1}{\rm g}(t)\;\!x\bigr)\exp({}-\lambda t)
+\nu^{\,\varepsilon}{\rm g}(t)f(t)\,\exp({}-\lambda t)=
\hfill                    
$
\\[2.25ex]
\mbox{}\hfill                    
$
=\bigl(\varepsilon\,\nu^{\,\varepsilon-1}{\rm g}(t)\;\!x+
\nu^{\,\varepsilon}{\rm g}(t)f(t)\bigr)\exp({}-\lambda t)
$
\ for all 
$
(t,x)\in J\times {\mathbb R}^n,
\ \ \varepsilon=1\ldots, m-1.
\hfill
$
\\[2ex]
\indent
Therefore, we have
\\[2ex]
\mbox{}\hfill                    % (3.20)
$
{\frak B}\bigl(\nu^{\,\varepsilon}{\rm g}(t)\;\!x\,\exp({}-\lambda t)\bigr)=
\bigl(\varepsilon\,\nu^{\,\varepsilon-1}{\rm g}(t)\;\!x+
\nu^{\,\varepsilon}{\rm g}(t)f(t)\bigr)\exp({}-\lambda t)
\hfill
$
\\[-0.15ex]
\mbox{}\hfill (3.20)
\\[-0.15ex]
\mbox{}\hfill
for all 
$
(t,x)\in J\times {\mathbb R}^n,
\quad
\varepsilon=1\ldots, m-1.
\hfill
$
\\[2ex]
\indent
Let $m=2.$ Using the identities (3.20) with $\varepsilon=1,$ we get
\\[2ex]
\mbox{}\hfill                                        
$
\displaystyle
{\frak B}\;\!F_2^{}(t,x)=
{\frak B}\;\!\bigl(
\nu^{1}{\rm g}(t)x\,\exp({}-\lambda t)-t\, F_1^{}(t,x) - C_1^{}(t)\bigr)=
\hfill                                        
$
\\[2ex]
\mbox{}\hfill                                        
$
\displaystyle
=\bigl(\nu^{0}{\rm g}(t)x+\nu^{1}{\rm g}(t)f(t)\bigr)\exp({}-\lambda t)- F_1^{}(t,x) -
\bigl(\nu^{1}{\rm g}(t)f(t)\,\exp({}-\lambda t) +C_0^{}(t)\bigr)=
\hfill                                        
$
\\[2ex]
\mbox{}\hfill                                        
$
\displaystyle
=\bigl(\nu^{0}{\rm g}(t)x\,\exp({}-\lambda t)- C_0^{}(t)\bigr) - F_1^{}(t,x)=0
$
\ for all 
$
(t,x)\in J\times{\mathbb R}^n.
\hfill
$
\\[2ex]
\indent
Hence the function
\vspace{0.35ex}
$F_2^{}\colon J\times{\mathbb R}^n\to {\mathbb R}$ 
is a first integral of the reducible system (3.14).

Assume that the functions (3.19) for $m=\mu$ are first integrals of system (3.14).
\vspace{0.25ex}
Then, from the identities (3.20) for the function $F_{\mu+1}^{}$ it follows that
\\[1.15ex]
\mbox{}\hfill                                  
$
\displaystyle
{\frak B}F_{\mu+1}^{}(t,x)= 
{\frak B}\Bigl(\nu^{\mu}{\rm g}(t)x\;\!\exp(-\lambda t) -  
\sum\limits_{\tau=0}^{\mu-1}
{\textstyle \binom{\mu}{\tau}}\;\! t^{\mu-\tau} F_{\tau+1}^{}(t,x)  -  C_{\mu}^{}(t)\Bigr)=
\hfill                              
$
\\[1.75ex]
\mbox{}\hfill                                  
$
\displaystyle
=\bigl(\;\!\mu\;\!\nu^{\mu-1}{\rm g}(t) x + \nu^{\mu}{\rm g}(t)f(t)\bigr)\;\!\exp({}-\lambda t) \, - \,
\mu\  \sum\limits_{\tau=0}^{\mu-2}
{\textstyle \binom{\mu-1}{\tau}}\, t^{\,{}^{\scriptstyle \mu-\tau-1}}\, F_{\tau+1}^{}(t,x) \ -
\hfill 
$
\\[1.75ex]
\mbox{}\hfill                                  
$
\displaystyle
- \ \mu\, F_{\mu}^{}(t,x)  \, -\,
\bigl(\;\!\nu^{\mu}{\rm g}(t)f(t)\,\exp({}-\lambda t) \;\! +\;\!
\mu\;\! C_{\mu-1}^{}(t)\bigr) \, =
\hfill 
$
\\[1.75ex]
\mbox{}\hfill                                  
$
\displaystyle
=\mu\Bigl(\! \nu^{\mu-1}{\rm g}(t) x-  
\sum\limits_{\tau=0}^{\mu-2}
{\textstyle \binom{\mu-1}{\tau}} t^{\,{}^{\scriptstyle \mu-\tau-1}} F_{\tau+1}^{}(t,x)  -
C_{\mu-1}^{}(t)\!\Bigr)\!
- \mu F_{\mu}(t,x) =0
$ 
for all 
$\!
(t,x)\!\in\! J\!\times\! \R^n\!.
\hfill 
$
\\[2ex]
\indent
This implies that the scalar function $F_{\mu+1}^{}\colon J\times{\mathbb R}^n\to {\mathbb R}$ 
\vspace{0.35ex}
(for $m=\mu+1)$ is a first integral 
of the reducible linear nonhomogeneous differential system (3.14).

So by the principle of mathematical induction, the scalar functions (3.19) are 
first integrals of the reducible system (3.14) for every natural number $m\geq 2.\ \k$
\vspace{0.75ex}

{\bf Example 3.7.}
Consider the linear nonhomogeneous differential system\footnote[1]{ 
Note that the linear homogeneous differential system corresponding to 
the nonhomogeneous system (3.21) is the 
linearization [8, pp. 139 -- 142] of differential equations, which  are describing
the motion of a symmetric balanced nonautonomous gyrostat 
with one point of attachment [71, pp. 219 -- 226; 72].

}\!\!
\\[1.75ex]                             
\mbox{}\hfill                    % (3.21)
$
\dfrac{dx_1^{}}{dt}=(\sigma +\psi(t))\;\!x_2^{}+\delta (\alpha\sin\varphi(t)-\beta\cos\varphi(t))\;\!x_3^{}+f_1^{}(t),
\hfill
$
\\[2ex]
\mbox{}\hfill 
$
\qquad\quad
\dfrac{dx_2^{}}{dt}={}-(\sigma +\psi(t))\;\!x_1^{}+\delta (\alpha\cos\varphi(t)+\beta\sin\varphi(t))\;\!x_3^{}+f_2^{}(t),
$
\hfill (3.21)
\\[2ex]
\mbox{}\hfill 
$
\dfrac{dx_3^{}}{dt}=\gamma (\beta\cos\varphi(t)-\alpha\sin\varphi(t))\;\!x_1^{}-
\gamma (\alpha\cos\varphi(t)+\beta\sin\varphi(t))\;\!x_2^{}+f_3^{}(t),
\hfill
$
\\[2.5ex]
where
\vspace{0.5ex} 
$\psi\colon J\to{\mathbb R},\ f_j^{}\colon J\to{\mathbb R},\, j=1,2,3,$ 
are continuous functions on an interval $J\subset {\mathbb R},$ the function 
$
\varphi\colon t\to\int\limits_{t_0^{}}^{t}\psi(\zeta)\;\!d\zeta
$ 
\vspace{0.5ex} 
for all 
$
t\in J,
\ 
t_0^{}\in J,
$
and $\alpha,\, \beta,\, \gamma,\, \delta,\, \sigma\ne 0$ are real numbers.

We claim that the differential system (3.21) is a reducible system by the 
\vspace{0.35ex} 
orthogonal group of transformations $O(J,3).$ Indeed,
the orthogonal transformation 
\\[1.5ex]                             
\mbox{}\hfill                   
$
y_1^{}\!=\cos\varphi(t)\;\!x_1^{}-\,\sin\varphi(t)\;\!x_2^{},
\quad
y_2^{}=\sin\varphi(t)\;\!x_1^{}+\cos\varphi(t)\;\!x_2^{},
\quad
y_3^{}=x_3^{}
\hfill
$
\\[1.5ex]
reduce the system (3.21) to the system with constant coefficients
\\[2ex]                             
\mbox{}\hfill                   
$
\dfrac{dy_1^{}}{dt}=\sigma\;\! y_2^{}-\beta\delta\;\! y_3^{}+
\cos\varphi(t)\;\! f_1^{}(t)-\;\!\sin\varphi(t)\;\! f_2^{}(t),
\hfill
$
\\[2.75ex]
\mbox{}\hfill                   
$
\dfrac{dy_2^{}}{dt}={}-\sigma\;\! y_1^{}+\alpha\delta\;\! y_3^{}+
\sin\varphi(t)\;\! f_1^{}(t)+\cos\varphi(t)\;\! f_2^{}(t),
\qquad
\dfrac{dy_3^{}}{dt}=\beta\gamma\;\! y_1^{}-\alpha\gamma\;\! y_2^{}+ f_3^{}(t).
\hfill
$
\\[-2.75ex]

\newpage

The matrix 
$
C=
\left\|\!\!
\begin{array}{ccc}
0& {}-\sigma\omega & \beta\gamma
\\
\sigma &  0 & {}-\alpha\gamma
\\
{}-\beta\delta & \alpha\delta & 0
\end{array}
\!\!\right\|
$
has the characteristic equation
\\[2.5ex]
\mbox{}\hfill
$
\det(C-\lambda E)=0
\iff 
\lambda\bigl(\lambda^2+(\alpha^2+\beta^2)\gamma\delta+\sigma^2\bigr)=0.
\hfill
$
\\[1.75ex]
\indent
Consider the real number $D=(\alpha^2+\beta^2)\gamma\delta+\sigma^2.$
\vspace{0.35ex}
We have three possible cases for building first integrals of the 
reducible linear differential system (3.21). 
\vspace{0.5ex}

Let $D=0.$ The matrix $C$ has the eigenvalue $\lambda_1^{}=0$
\vspace{0.5ex}
with the elementary divisor $\lambda^3$ 
corresponding to the real eigenvector $\nu^{10}=(\alpha\gamma,\beta\gamma,\sigma)$
\vspace{0.75ex}
and to the generalized eigenvectors
$
\nu^{11}=\Bigl((\alpha+\beta)\,\dfrac{\gamma}{\sigma}\,,
(\beta-\alpha)\,\dfrac{\gamma}{\sigma}\,, 1\Bigr),\
\nu^{12}=\Bigl(\dfrac{2\gamma}{\sigma^2}\,\beta,
{}-\dfrac{2\gamma}{\sigma^2}\,\alpha, \dfrac{2}{\sigma}\Bigr).$
\vspace{1ex}
Using Theorems 3.5 and 3.6, 
we can construct the basis of first integrals for 
the reducible differential system (3.21)
\\[2.5ex]
\mbox{}\hfill
$
F_1^{}\colon (t,x)\to 
\gamma \bigl(\alpha\cos\varphi(t)+\beta\sin\varphi(t)\bigr)\;\! x_1^{} +
\gamma \bigl(\beta\cos\varphi(t)-\alpha\sin\varphi(t)\bigr)\;\! x_2^{} +
\sigma\;\!x_3^{} - C_0^{}(t),
$
\hfill (3.22)
\\[2.5ex]
\mbox{}\hfill
$
F_2^{}\colon (t,x)\to \
\dfrac{\gamma}{\sigma}\, 
\bigl((\alpha+\beta)\cos\varphi(t)+(\beta-\alpha)\sin\varphi(t)\bigr)\;\! x_1^{} \ +
\hfill
$
\\[2ex]
\mbox{}\hfill
$
\displaystyle
+\ \dfrac{\gamma}{\sigma}\, 
\bigl((\beta-\alpha)\cos\varphi(t)-(\alpha+\beta)\sin\varphi(t)\bigr)\;\! x_2^{} +
x_3^{} - t\;\! F_1^{}(t,x) - C_1^{}(t)
$
\ for all 
$
(t,x)\in J\times {\mathbb R}^3, 
\hfill
$
\\[2.75ex]
\mbox{}\hfill
$
F_3^{}\colon (t,x)\to \
\dfrac{2\gamma}{\sigma^2}\, 
\bigl(\beta\cos\varphi(t)-\alpha\sin\varphi(t)\bigr)\;\! x_1^{}  -
\dfrac{2\gamma}{\sigma^2}\, 
\bigl(\alpha\cos\varphi(t)+\beta\sin\varphi(t)\bigr)\;\! x_2^{}  \ +
\hfill
$
\\[2.25ex]
\mbox{}\hfill
$
\displaystyle
+\ \dfrac{2}{\sigma}\, x_3^{} - 
t^2\;\! F_1^{}(t,x)-2t\;\! F_2^{}(t,x) - C_2^{}(t)
$
\ for all 
$
(t,x)\in J\times {\mathbb R}^3, 
\hfill
$
\\[2.25ex]
where the scalar functions
\\[2ex]
\mbox{}\hfill
$
\displaystyle
C_0^{}\colon t\to 
\int\limits_{t_0^{}}^{t}
\Bigl(
\gamma \bigl(\alpha\cos\varphi(\zeta)+\beta\sin\varphi(\zeta)\bigr)\;\! f_1^{}(\zeta) +
\gamma \bigl(\beta\cos\varphi(\zeta)-\alpha\sin\varphi(\zeta)\bigr)\;\! f_2^{}(\zeta) +
\sigma\;\!f_3^{}(\zeta)\Bigr)\;\!d\zeta,
\hfill
$
\\[2ex]
\mbox{}\hfill
$
\displaystyle
C_1^{}\colon t\to \
\int\limits_{t_0^{}}^{t}
\Bigl(
\dfrac{\gamma}{\sigma}\, 
\bigl((\alpha+\beta)\cos\varphi(\zeta)+(\beta-\alpha)\sin\varphi(\zeta)\bigr)\;\! f_1^{}(\zeta) \ +
\hfill
$
\\[1.75ex]
\mbox{}\hfill
$
\displaystyle
+\ \dfrac{\gamma}{\sigma}\, 
\bigl((\beta-\alpha)\cos\varphi(\zeta)-(\alpha+\beta)\sin\varphi(\zeta)\bigr)\;\! f_2^{}(\zeta) +
f_3^{}(\zeta) +C_0^{}(\zeta)\Bigr)\;\!d\zeta
$
\ for all 
$
t\in J,
\hfill
$
\\[2.75ex]
\mbox{}\hfill
$
\displaystyle
C_2^{}\colon t\to \
\int\limits_{t_0^{}}^{t}
\Bigl(
\dfrac{2\gamma}{\sigma^2}\, 
\bigl(\beta\cos\varphi(\zeta)-\alpha\sin\varphi(\zeta)\bigr)\;\! f_1^{}(\zeta)  -
\dfrac{2\gamma}{\sigma^2}\, 
\bigl(\alpha\cos\varphi(\zeta)+\beta\sin\varphi(\zeta)\bigr)\;\! f_2^{}(\zeta)  \ +
\hfill
$
\\[1.75ex]
\mbox{}\hfill
$
\displaystyle
+\ \dfrac{2}{\sigma}\, f_3^{}(\zeta) + 2C_1^{}(\zeta)\Bigr)\;\!d\zeta
$
\ for all 
$
t\in J,
\quad
t_0^{}\in J.
\hfill
$
\\[2.25ex]
\indent
Let $D<0.$ The matrix $C$ has the linearly independent eigenvectors
\vspace{1ex}
$\nu^{10}=(\alpha\gamma, \beta\gamma, \sigma),$
$\nu^{2}=\bigl(\sigma^2+\gamma\delta\beta^2, {}-\alpha\beta\gamma\delta-\sigma\sqrt{{}-D}\,,
\vspace{1.25ex}
\delta\bigl({}-\alpha\sigma+\beta\sqrt{{}-D}\,\bigr)\bigr),\ 
\nu^{3}=\bigl(\sigma^2+\gamma\delta\beta^2, {}-\alpha\beta\gamma\delta+\sigma\sqrt{{}-D}\,,$
${}-\delta\bigl(\alpha\sigma+\beta\sqrt{{}-D}\,\bigr)\bigr)$ 
\vspace{1ex}
corresponding to the eigenvalues 
$\lambda_1^{}\!=0,\, \lambda_2^{}\!=\sqrt{{}-D},$  
$\lambda_3^{}\!=-\,\sqrt{-\,D}.$ 
An integral basis of system (3.21) is the function (3.22) and the functions (by Theorem 3.5) 
\\[2ex]
\mbox{}\hfill
$
F_k^{}\colon (t,x)\to 
\Bigl(
\bigl((\sigma^2+\gamma\delta\beta^2)\cos\varphi(t)
-(\alpha\beta\gamma\delta+\sigma\lambda_k^{})\sin\varphi(t)\bigr)\;\!x_1^{} \ -
\hfill
$
\\[1.75ex]
\mbox{}\hfill
$
-\ \bigl((\alpha\beta\gamma\delta+\sigma\lambda_k^{})\cos\varphi(t)+
(\sigma^2+\gamma\delta\beta^2)\sin\varphi(t)\bigr)\;\!x_2^{} +
\delta (\beta\lambda_k^{}-\alpha\sigma)\;\!x_3^{}\Bigr)e^{{}-\lambda_k^{} t} \ -
\hfill
$
\\[1.75ex]
\mbox{}\hfill
$
\displaystyle
-\ \int\limits_{t_0^{}}^{t}
\Bigl(
\bigl((\sigma^2+\gamma\delta\beta^2)\cos\varphi(\zeta)
-(\alpha\beta\gamma\delta+\sigma\lambda_k^{})\sin\varphi(\zeta)\bigr)\;\!f_1^{}(\zeta)  -
\bigl((\alpha\beta\gamma\delta+\sigma\lambda_k^{})\cos\varphi(\zeta) \ +
\hfill
$
\\[1.75ex]
\mbox{}\hfill
$
 +\ 
(\sigma^2+\gamma\delta\beta^2)\sin\varphi(\zeta)\bigr)\;\!f_2^{}(\zeta) +
\delta (\beta\lambda_k^{}-\alpha\sigma)\bigr)\;\!f_3^{}(\zeta)\Bigr)e^{{}-\lambda_k^{} \zeta}\;\!d\zeta
\quad
\forall (t,x)\in J\times {\mathbb R}^3,
\quad
k=2,3.
\hfill
$
\\[2ex]
\indent
Let $D>0.$ The matrix $C$ has the linearly independent eigenvectors
\vspace{1ex}
$\nu^{10}=(\alpha\gamma,\beta\gamma,\sigma),$  
$\nu^{2}=\bigl({}-\sigma^2-\gamma\delta\beta^2,\, \alpha\beta\gamma\delta+\sigma\sqrt{D}\,i,\,
\vspace{1ex}
\delta\bigl(\alpha\sigma-\beta\sqrt{D}\,i\bigr)\bigr),\ 
\nu^{3}=\bigl({}-\sigma^2-\gamma\delta\beta^2,\, \alpha\beta\gamma\delta-\sigma\sqrt{D}\,i,$
$
\delta\bigl(\alpha\sigma+\beta\sqrt{D}\,i\bigr)\bigr)\!
$ 
corresponding to the eigenvalues 
\vspace{0.75ex}
$
\lambda_1^{}\!=0,\ \lambda_2^{}\!=\sqrt{D}\,i,$ and
$\lambda_3^{}\!=-\;\!\sqrt{D}\,i.$ 
An integral basis of system (3.21) is the function (3.22) and the functions (by Corollary 3.1) 
\\[1.5ex]
\mbox{}\hfill
$
\displaystyle
F_{k}^{}\colon (t,x_1^{},x_2^{},x_3^{})\to \
\alpha_{{}_{\scriptstyle k}}(t, x_1^{},x_2^{},x_3^{})-
\int\limits_{t_0^{}}^{t} \alpha_{{}_{\scriptstyle k}}(\zeta, f_1^{}(\zeta), f_2^{}(\zeta),f_3^{}(\zeta))\,d\zeta,
\quad 
k=2,3,
\hfill
$
\\[1.5ex]
where $t_0^{}$ is a fixed point from the interval $J,$ the scalar functions
\\[1.75ex]
\mbox{}\hfill
$
\displaystyle
\alpha_{2}^{}\colon (t,x)\to 
\bigl(\!{}-(\sigma^2+\gamma\delta\beta^2)\cos\sqrt{D}\,t\cos\varphi(t)+
\bigl(\alpha\beta\gamma\delta\cos\sqrt{D}\,t+ \sigma\sqrt{D}\,\sin\sqrt{D}\,t\bigr)\sin\varphi(t)\!\bigr)\;\!x_1^{}+
\hfill
$
\\[2ex]
\mbox{}\hfill
$
+\ \bigl(\bigl(\alpha\beta\gamma\delta\cos\sqrt{D}\,t+ \sigma\sqrt{D}\,\sin\sqrt{D}\,t\bigr)\cos\varphi(t)+
(\sigma^2+\gamma\delta\beta^2)\cos\sqrt{D}\,t\sin\varphi(t)\bigr)\;\!x_2^{}\ +
\hfill
$
\\[2ex]
\mbox{}\hfill
$
+\ \delta\bigl(\alpha\sigma\cos\sqrt{D}\,t-\beta\sqrt{D}\,\sin\sqrt{D}\,t\bigr)\;\!x_3^{}
$
\ for all 
$
(t,x)\in J\times{\mathbb R}^3,
\hfill
$
\\[2.25ex]
\mbox{}\hfill
$
\displaystyle
\alpha_{3}^{}\colon (t,x)\to 
\bigl( (\sigma^2+\gamma\delta\beta^2)\sin\sqrt{D}\,t\cos\varphi(t)+
\bigl(\sigma\sqrt{D}\,\cos\sqrt{D}\,t-\alpha\beta\gamma\delta\sin\sqrt{D}\,t\bigr)\sin\varphi(t)\!\bigr)\;\!x_1^{}+
\hfill
$
\\[2ex]
\mbox{}\hfill
$
+\ \bigl(\bigl(\sigma\sqrt{D}\,\cos\sqrt{D}\,t-\alpha\beta\gamma\delta\sin\sqrt{D}\,t\bigr)\cos\varphi(t)-
(\sigma^2+\gamma\delta\beta^2)\sin\sqrt{D}\,t\sin\varphi(t)\bigr)\;\!x_2^{}\ -
\hfill
$
\\[2ex]
\mbox{}\hfill
$
-\ \delta\bigl(\beta\sqrt{D}\,\cos\sqrt{D}\,t+\alpha\sigma\sin\sqrt{D}\,t\bigr)\;\!x_3^{}
$
\ for all 
$
(t,x)\in J\times{\mathbb R}^3.
\hfill
$
\\[2.25ex]
\indent
{\bf Corollary 3.2.}
{\it
Let the differential system {\rm (3.14)} be reducible to the system {\rm (3.15)} 
by a tran\-s\-for\-ma\-ti\-on matrix ${\rm g}\in G,$ and let 
\vspace{0.35ex}
$\lambda={\stackrel{*}{\lambda}}+\widetilde{\lambda}\,i
\ ({\rm Re}\,\lambda={\stackrel{*}{\lambda}},\ {\rm Im}\,\lambda=\widetilde{\lambda}\ne 0)$
be the complex eigenvalue of the matrix $C$ with elementary divisor of multiplicity $m\geq 2$ 
\vspace{0.35ex}
corresponding to an complex eigenvector
\vspace{0.35ex}
$\nu^0={\stackrel{*}{\nu}}\;\!{}^{0}+\widetilde{\nu}\;\!{}^{0}\,i
\ ({\rm Re}\,\nu^0={\stackrel{*}{\nu}}\;\!{}^{0},\ {\rm Im}\,\nu^0=\widetilde{\nu}\;\!{}^{0})$ 
and to generalized eigenvectors  
\vspace{0.5ex}
$\nu^k={\stackrel{*}{\nu}}\;\!{}^{k}+\widetilde{\nu}\;\!{}^{k}\,i\ 
({\rm Re}\,\nu^k={\stackrel{*}{\nu}}\;\!{}^{k},\ {\rm Im}\,\nu^k=\widetilde{\nu}\;\!{}^{k}),
\ k=1,\ldots, m-1.$ 
Then first integrals of the re\-du\-cib\-le system {\rm(3.14)} are the functions
\\[1.75ex]
\mbox{}\hfill                                  % (3.23)
$
\displaystyle
F_{{}_{\scriptstyle \theta, k+1}}\colon (t,x)\to \
\alpha_{{}_{\scriptstyle \theta k}}(t,x)-\,
\sum\limits_{\tau=0}^{k-1}
{\textstyle \binom{k}{\tau}}\, t^{k-\tau}\, F_{{}_{\scriptstyle \theta, \tau+1}}(t,x) -
C_{{}_{\scriptstyle \theta k}}(t)
\hfill                              
$
\\[0.25ex]
\mbox{}\hfill      {\rm (3.23)}                            
\\[0.25ex]
\mbox{}\hfill                                  
for all 
$
(t,x)\in J\times {\mathbb R}^n,
\quad
k=1,\ldots, m-1,
\quad
\theta=1,2,
\hfill                              
$
\\[2ex]
where the first integrals {\rm(}by Corollary {\rm 3.1)}
\\[1.75ex]
\mbox{}\hfill                                     
$
F_{{}_{\scriptstyle \theta\;\!1}}^{}\colon (t,x)\to \,
\alpha_{{}_{\scriptstyle \theta\;\! 0}}(t, x)-
C_{{}_{\scriptstyle \theta\;\! 0}}(t)
$
\ for all 
$
(t,x)\in J\!\times {\mathbb R}^n,
\quad 
\theta=1,2,
\hfill
$
\\[1.75ex]
the scalar functions
\\[1.75ex]
\mbox{}\hfill                                     
$
\displaystyle
\alpha_{{}_{\scriptstyle 1k}}\colon (t, x)\to
\bigl(\,{\stackrel{*}{\nu}}\;\!{}^k{\rm g}(t)x\,\cos\widetilde{\lambda}\,t +
\widetilde{\nu}\;\!{}^k{\rm g}(t)x\,\sin\widetilde{\lambda}\,t \bigr)
\exp\bigl({}-{\stackrel{*}{\lambda}}\,t\bigr)
$
\ for all 
$
(t,x)\in J\times {\mathbb R}^{n},
\hfill
$
\\[2ex]
\mbox{}\hfill                                     
$
\displaystyle
\alpha_{{}_{\scriptstyle 2k}}\colon (t,x)\to
\bigl(\,\widetilde{\nu}\;\!{}^k{\rm g}(t)x\,\cos\widetilde{\lambda}\,t -
{\stackrel{*}{\nu}}\;\!{}^k{\rm g}(t)x\,\sin\widetilde{\lambda}\,t\bigr)
\exp\bigl({}-{\stackrel{*}{\lambda}}\,t\bigr)
\quad
\forall (t,x)\in J\times {\mathbb R}^{n},
\hfill                              
$
\\[2ex]
\mbox{}\hfill                                  
$
\displaystyle
C_{{}_{\scriptstyle \theta k}}\colon t\to
\int\limits_{t_0^{}}^{t}\!\bigl(\alpha_{{}_{\scriptstyle \theta k}}(\zeta, f(\zeta)) +
k\;\!C_{{}_{\scriptstyle \theta, k-1}}(\zeta)\bigr)\;\! d\zeta\!
$
for all 
$
t\in J,
\ \,
\theta=1,2,
\
k=0,\ldots, m-1, 
\ t_0^{}\in J.
\hfill                              
$
\\[1.5ex]
%а 
%\vspace{0.5ex}
%$t_0^{}$ --- произвольная фиксированная точки из числового промежутка $J.$
}
\indent
{\sl Proof}.  
Since the functions (3.19) are complex-valued first integrals (by Theorem 3.6) 
of the reducible system (3.14), we see that 
the real and imaginary parts of the functions (3.19) are the real 
first integrals (3.23) of the reducible system (3.14). \k
\vspace{0.75ex}

{\bf Example 3.8.}
The linear nonhomogeneous Euler differential system 
\\[2ex]
\mbox{}\hfill                               % (3.24)
$
\dfrac{dx_1^{}}{dt} =   x_2^{} +a,
\qquad
\dfrac{dx_2^{}}{dt} =  x_3^{} + \dfrac{b}{t^2}\,\ln t,
\qquad
\dfrac{dx_3^{}}{dt} =  x_4^{} + \dfrac{c}{t^3}\,\sin\ln t,
\hfill
$
\\[0.25ex]
\mbox{}\hfill  (3.24)
\\[0.25ex]
\mbox{}\hfill
$
\dfrac{dx_4^{}}{dt} =  {}-\dfrac{1}{t^4}\,x_1^{} -\dfrac{3}{t^3}\,x_2^{} 
-\dfrac{9}{t^2}\,x_3^{}-\dfrac{6}{t}\,x_4^{}+\dfrac{d}{t^4\cos\ln t}
\qquad
(a,b,c,d\in{\mathbb R})
\hfill
$
\\[2ex]
is reducible to the linear differential system with constant coefficients
\\[2ex]
\mbox{}\hfill
$
\dfrac{dy_1^{}}{d\tau}=y_2^{}+a\;\!e^{\tau},
\qquad
\dfrac{dy_2^{}}{d\tau}=y_2^{}+y_3^{}+b\;\!\tau,
\qquad
\dfrac{dy_3^{}}{d\tau}=2y_3^{}+y_4^{}+c\sin\tau,
\hfill
$
\\[2.5ex]
\mbox{}\hfill
$
\dfrac{dy_4^{}}{d\tau}={}-y_1^{}-3y_2^{}-9y_3^{}-3y_4^{}+\dfrac{d}{\cos\tau}
\hfill
$
\\[2.25ex]
by the exponential group
${\rm Exp}(4)$ with the exponential transformation 
\\[1.5ex]
\mbox{}\hfill
$
y_1^{}=x_1^{},
\quad
y_2^{}=tx_2^{},
\quad
y_3^{}=t^2x_3^{},
\quad
y_4^{}=t^3x_4^{},
\quad
t=e^{\tau},
\hfill
$
\\[2.5ex]
\indent
The matrix 
\vspace{0.5ex}
$
C = \left\|\!\!
\begin{array}{cccc}
0& 0&  0\!\! & {}-1
\\
1 &  1& 0\!\! & {}- 3
\\
0 &  1& 2\!\! & {}- 9
\\
0 &  0& 1\!\! & {}- 3
\end{array}
\!\! \right\|
$
has the complex eigenvalues $\lambda_{1}^{}={}-i$ and $\lambda_{2}^{}=i$
with the multiple elementary divisors $(\lambda+i)^2$ and $(\lambda-i)^2,$ respectively.
\vspace{0.5ex}

The eigenvalue $\lambda_{1}^{}={}-i$ of the matrix $C$ is
\vspace{0.5ex}
corresponding to the complex eigenvector $\nu^{0}\!=(1, 1+2i, 1+3i, i)$ and to the 
generalized eigenvector of the 1-st order $\nu^{1}\!=({}-i, 2+i, i, 0).$  

By Corollaries 3.1 and 3.2, using the numbers 
${\stackrel{*}{\lambda}}_1^{}=0,\ \widetilde{\lambda}_1^{}={}-1,$  
the linearly independent vectors
${\stackrel{*}{\nu}}{}^{0}=(1, 1, 1, 0),\ \widetilde{\nu}{}^{\,0}=(0, 2, 3, 1),\ 
{\stackrel{*}{\nu}}{}^{1}=(0, 2, 0, 0), \ \widetilde{\nu}{}^{\,1}=({}-1, 1, 1, 0),$
the functions
\\[2ex]
\mbox{}\hfill
$
\alpha_{{}_{\scriptstyle 10}}\colon (t,x)\to\,
\cos\ln t\, (x_1^{}+tx_2^{}+t^2x_3^{})-t\sin\ln t\, (2x_2^{}+3tx_3^{}+t^2x_4^{})
$
for all 
$
(t,x)\in J_l^{}\times {\mathbb R}^4,
\hfill
$
\\[2.25ex]
\mbox{}\hfill
$
\alpha_{{}_{\scriptstyle 20}}\colon (t,x)\to\,
t\cos\ln t\, (2x_2^{}+3tx_3^{}+t^2x_4^{})+
\sin\ln t\,(x_1^{}+tx_2^{}+t^2x_3^{}) 
$
for all 
$
(t,x)\in J_l^{}\times {\mathbb R}^4,
\hfill
$
\\[2.25ex]
\mbox{}\hfill
$
\alpha_{{}_{\scriptstyle 11}}\!\colon (t,x)\to\,
2t\cos\ln t\, x_2^{}+
\sin\ln t\, (x_1^{}-tx_2^{}-t^2x_3^{})
$
\ for all 
$
(t,x)\in J_l^{}\times {\mathbb R}^4,
\hfill
$
\\[2.25ex]
\mbox{}\hfill
$
\alpha_{{}_{\scriptstyle 21}}\colon (t,x)\to\,
\cos\ln t\, ({}-x_1^{}+tx_2^{}+t^2x_3^{})+2t\sin\ln t\, x_2^{}
$ 
\ for all 
$
(t,x)\in J_l^{}\times {\mathbb R}^4,
\hfill
$
\\[1ex]
and
\\[1.25ex]
\mbox{}\hfill
$
\displaystyle
C_{{}_{\scriptstyle 10}}(t)=
\int \!
\biggl(\Bigl(a+\dfrac{b}{t}\,\ln t+\dfrac{c}{t}\,\sin\ln t\Bigr)\cos\ln t -
\Bigl(\dfrac{2b}{t^2}\,\ln t+\dfrac{3c}{t^2}\,\sin\ln t+\dfrac{d}{t^2\cos\ln t}\Bigr)t\sin\ln t\biggr) \;\!dt=
\hfill
$
\\[2.25ex]
\mbox{}\hfill
$
=\dfrac{a}{2}\,(\cos\ln t+\sin\ln t)\;\!t +
b\bigl(\cos\ln t-2\sin\ln t+(2\cos\ln t+\sin\ln t)\ln t\bigr) \ -
\hfill
$
\\[2.25ex]
\mbox{}\hfill
$
-\ \dfrac{c}{4}\,(6\ln t-2\sin^2\ln t-3\sin 2\ln t) +
d\ln |\cos\ln t|
$ 
\ for all 
$
t\in J_l^{},
\hfill
$
\\[2.5ex]
\mbox{}\hfill
$
\displaystyle
C_{{}_{\scriptstyle 20}}(t)=
\int \!
\biggl(
\Bigl(\dfrac{2b}{t^2}\,\ln t+\dfrac{3c}{t^2}\,\sin\ln t+\dfrac{d}{t^2\cos\ln t}\Bigr)t\cos\ln t +
\Bigl(a+\dfrac{b}{t}\,\ln t+\dfrac{c}{t}\,\sin\ln t\Bigr)\sin\ln t\biggr) \;\!dt=
\hfill
$
\\[2.25ex]
\mbox{}\hfill
$
=\dfrac{a}{2}\,(\sin\ln t-\cos\ln t)\;\!t +
b\bigl(2\cos\ln t+\sin\ln t+(2\sin\ln t-\cos\ln t)\ln t\bigr) \ +
\hfill
$
\\[2.25ex]
\mbox{}\hfill
$
+\ \dfrac{c}{4}\,(2\ln t+6\sin^2\ln t-\sin 2\ln t) +
d\ln t
$ 
\ for all 
$
t\in J_l^{},
\hfill
$
\\[2.5ex]
\mbox{}\hfill
$
\displaystyle
C_{{}_{\scriptstyle 11}}(t)=
\int
\biggl(\dfrac{2b}{t}\,\ln t\cos\ln t+
\Bigl(a-\dfrac{b}{t}\,\ln t -\dfrac{c}{t}\,\sin\ln t\Bigr)\sin\ln t+
C_{{}_{\scriptstyle 10}}(t)\biggr)\;\!dt=
\hfill
$
\\[2.25ex]
\mbox{}\hfill
$
=\dfrac{a}{2}\,(2\sin\ln t-\cos\ln t)\;\!t +
b\bigl(6\cos\ln t+\sin\ln t+4\ln t\sin\ln t\bigr) \ -
\hfill
$
\\[2.25ex]
\mbox{}\hfill
$
-\ \dfrac{c}{8}\,(2\ln t+6\ln^2 t+3\cos 2\ln t-\sin 2\ln t) +
d\bigl(\ln t\;\!\ln |\cos\ln t|+\beta(\ln t)\bigr)
$ 
\ for all 
$
t\in J_l^{},
\hfill
$
\\[2.5ex]
\mbox{}\hfill
$
\displaystyle
C_{{}_{\scriptstyle 21}}(t)=
\int
\biggl(
\Bigl({}-a+\dfrac{b}{t}\,\ln t +\dfrac{c}{t}\,\sin\ln t\Bigr)\cos\ln t+
\dfrac{2b}{t}\,\ln t\;\!\sin\ln t+C_{{}_{\scriptstyle 20}}(t)\biggr)\;\!dt=
\hfill
$
\\[2.25ex]
\mbox{}\hfill
$
={}-\dfrac{a}{2}\,(2\cos\ln t+\sin\ln t)\;\!t -
b\bigl(\cos\ln t-6\sin\ln t+4\ln t\cos\ln t\bigr) \ +
\hfill
$
\\[2.25ex]
\mbox{}\hfill
$
+\ \dfrac{c}{8}\,(6\ln t+2\ln^2 t+\cos 2\ln t-3\sin 2\ln t+4\sin^2\ln t) +
\dfrac{d}{2}\, \ln^2 t
$ 
\ for all 
$
t\in J_l^{},
\hfill
$
\\[2.5ex]
where $\beta(t)=\int t\tan t\;\!dt,$
we can build the first integrals of the reducible system (3.24)
\\[2.5ex]
\mbox{}\hfill
$
F_{{}_{\scriptstyle 11}}\colon (t,x)\to\
\cos\ln t\, (x_1^{}+tx_2^{}+t^2x_3^{})-t\sin\ln t\, (2x_2^{}+3tx_3^{}+t^2x_4^{}) \ -
\hfill
$
\\[2.25ex]
\mbox{}\hfill
$
-\ \dfrac{a}{2}\,(\cos\ln t+\sin\ln t)\;\!t -
b\bigl(\cos\ln t-2\sin\ln t+(2\cos\ln t+\sin\ln t)\ln t\bigr) \ +
\hfill
$
\\[2.25ex]
\mbox{}\hfill
$
+\ \dfrac{c}{4}\,(6\ln t-2\sin^2\ln t-3\sin 2\ln t) -
d\ln |\cos\ln t|
$ 
\ for all 
$
(t,x)\in J_l^{}\times {\mathbb R}^4,
\hfill
$
\\[2.75ex]
\mbox{}\hfill
$
F_{{}_{\scriptstyle 21}}\colon (t,x)\to\
t\cos\ln t\, (2x_2^{}+3tx_3^{}+t^2x_4^{})+
\sin\ln t\,(x_1^{}+tx_2^{}+t^2x_3^{}) \ -
\hfill
$
\\[2.25ex]
\mbox{}\hfill
$
-\ \dfrac{a}{2}\,(\sin\ln t-\cos\ln t)\;\!t -
b\bigl(2\cos\ln t+\sin\ln t+(2\sin\ln t-\cos\ln t)\ln t\bigr) \ -
\hfill
$
\\[2.25ex]
\mbox{}\hfill
$
-\ \dfrac{c}{4}\,(2\ln t+6\sin^2\ln t-\sin 2\ln t) - d\ln t
$ 
\ for all 
$
(t,x)\in J_l^{}\times {\mathbb R}^4,
\hfill
$
\\[2.75ex]
\mbox{}\hfill
$
F_{{}_{\scriptstyle 12}}\colon (t,x)\to\
2t\cos\ln t\, x_2^{}+\sin\ln t\, (x_1^{}-tx_2^{}-t^2x_3^{})
-\ln t\, F_{{}_{\scriptstyle 11}}(t,x) \ -
\hfill
$
\\[2.25ex]
\mbox{}\hfill
$
-\ \dfrac{a}{2}\,(2\sin\ln t-\cos\ln t)\;\!t -
b\bigl(6\cos\ln t+\sin\ln t+4\ln t\sin\ln t\bigr) \ +
\hfill
$
\\[2.25ex]
\mbox{}\hfill
$
+\;\!\dfrac{c}{8}\,(2\ln t+6\ln^2 t+3\cos 2\ln t-\sin 2\ln t)\! -
d\bigl(\ln t\ln |\cos\ln t|+\beta(\ln t)\bigr)\!
$ 
for all 
$\!
(t,x)\!\in\! J_l^{}\!\times\! {\mathbb R}^4,
\hfill
$
\\[2.75ex]
\mbox{}\hfill
$
F_{{}_{\scriptstyle 22}}\colon (t,x)\to\
\cos\ln t\, ({}-x_1^{}+tx_2^{}+t^2x_3^{})+2t\sin\ln t\, x_2^{}
-\ln t\, F_{{}_{\scriptstyle 21}}(t,x)\ +
\hfill
$
\\[2.25ex]
\mbox{}\hfill
$
+\ \dfrac{a}{2}\,(2\cos\ln t+\sin\ln t)\;\!t +
b\bigl(\cos\ln t-6\sin\ln t+4\ln t\cos\ln t\bigr) \ -
\hfill
$
\\[2.25ex]
\mbox{}\hfill
$
-\ \dfrac{c}{8}\,(6\ln t+2\ln^2 t+\cos 2\ln t-3\sin 2\ln t+4\sin^2\ln t) -
\dfrac{d}{2}\, \ln^2 t
$ 
\ for all 
$
(t,x)\in J_l^{}\times {\mathbb R}^4
\hfill
$
\\[2.25ex]
on any domain $J_l^{}\times {\mathbb R}^4\subset {\mathbb R}^5,\ l=0,1,2,\ldots,$ 
where the intervals
\\[1.75ex]
\mbox{}\hfill
$
J_{_0}=\Bigl(0; e^{{}^{\tfrac{\pi}{2}}}\Bigr),
\quad
J_{s}=\Bigl(e^{{}^{\scriptstyle \tfrac{\pi}{2}+\pi(s-1)}};\, 
e^{{}^{\scriptstyle\tfrac{\pi}{2}+\pi s}}\Bigr),
\ \ s=1,2,\ldots\,.
\hfill
$
\\[2ex]
\indent
The functionally independent first integrals 
\vspace{0.75ex}
$F_{11}^{},\, F_{21}^{},\, F_{12}^{},$ and $F_{22}^{}$ are 
an integral basis of the Euler differential system (3.24) on any domain $J_l^{}\times {\mathbb R}^4.$

\newpage

\mbox{}
\\[-3.5ex]

}
\end{document}